\documentclass[12pt,a4paper,french,english]{amsart}

\usepackage[utf8]{inputenc}
\usepackage[T1]{fontenc}
\usepackage{babel}
\usepackage[a4paper,top=2cm,bottom=2cm,left=3cm,right=3cm,marginparwidth=1.75cm]{geometry}

\usepackage{amssymb}
\usepackage{mathrsfs} 
\usepackage{amsmath}
\usepackage{amsthm}
\usepackage{amscd}
\usepackage{pgf,tikz}
\usetikzlibrary{arrows.meta}
\usepackage[alwaysadjust]{paralist} 
\setdefaultenum{\upshape(i)}{}{}{}
\setdefaultleftmargin{10em}{1em}{1em}{1em}{1em}{1em} %
\usepackage[pdftex,colorlinks,linkcolor=blue,filecolor=red,citecolor=blue,urlcolor=blue]{hyperref} 

\newtheoremstyle{mythmstyle}%
{1.5\baselineskip}
{\baselineskip}
{\itshape}
{}
{\bf}
{}
{0pt}
{} 

\newtheoremstyle{mydefstyle}%
{1.5\baselineskip}
{\baselineskip}
{}
{}
{\bf}
{}
{0pt}
{} 

\newtheoremstyle{mypreuvestyle}%
{\baselineskip}
{\baselineskip}
{}
{}
{\em}
{}
{0pt}
{} 

\newif\ifmynonumberenvi\mynonumberenvitrue
\theoremstyle{mythmstyle}
\newtheorem{proclaimmythm}[equation]{} 
\newtheorem*{proclaimmythm*}{}
\newenvironment{proclaim}[2][*]{\ifx*#1\mynonumberenvitrue\begin{proclaimmythm*}{\bf#2.} \ignorespaces\else\mynonumberenvifalse\begin{proclaimmythm}{.\kern0.5em\bf#2.}\label{#1} \ignorespaces\fi}{\ifmynonumberenvi\end{proclaimmythm*}\else\end{proclaimmythm}\fi}
\renewenvironment{proclaim}[2][*]{\ifx*#1\begin{proclaimmythm}{.\kern0.5em\bf#2.} \ignorespaces\else\begin{proclaimmythm}{.\kern0.5em\bf#2.}\label{#1} \ignorespaces\fi}{\end{proclaimmythm}}

\theoremstyle{mydefstyle}
\newtheorem{proclaimmydef}[equation]{} 
\newtheorem*{proclaimmydef*}{}
\newenvironment{definition}[2][*]{\ifx*#1\mynonumberenvitrue\begin{proclaimmydef*}{\bf#2.}\else\mynonumberenvifalse\begin{proclaimmydef}{.\kern0.5em\bf#2.}\label{#1} \ignorespaces\fi}{\ifmynonumberenvi\end{proclaimmydef*}\else\end{proclaimmydef}\fi}
\renewenvironment{definition}[2][*]{\ifx*#1\begin{proclaimmydef}{.\kern0.5em\bf#2.} \ignorespaces\else\begin{proclaimmydef}{.\kern0.5em\bf#2.}\label{#1} \ignorespaces\fi}{\end{proclaimmydef}}

\newenvironment{definition*}[1]{\begin{proclaimmydef*}{\bf#1.} \ignorespaces}{\end{proclaimmydef*}}

\newif\ifSymbolFinRem\SymbolFinRemtrue
\def\FinRemSymbol{\hbox{$\blacktriangleleft\!\blacktriangleright$}}
\def\FinRem{\hskip0.001em plus 40pt\null{}\kern5pt\null\nobreak\hfill
   \kern3pt{\FinRemSymbol
   }}
\newif\ifFinRemdejaplace\FinRemdejaplacefalse
\newenvironment{remarque}[2][*]%
{\ifx*#1\mynonumberenvitrue\begin{proclaimmydef*}{\bf#2.}\else\mynonumberenvifalse\begin{proclaimmydef}{\ \bf#2.}\label{#1} \ignorespaces\fi{}}%
{\ifSymbolFinRem\ifFinRemdejaplace\global\FinRemdejaplacefalse\else\FinRem\fi\fi
\ifmynonumberenvi\end{proclaimmydef*}\else\end{proclaimmydef}\fi}

\def\QEDbox{\hbox{\lower2.3pt\vbox{\hrule\hbox
   {\vrule\kern1pt\vbox{\kern1.7pt\hbox{$\scriptstyle
   QED$}\kern.6pt}\kern1pt\vrule}\hrule}}}
\def\QED{\hskip0.01em plus 40pt\null{} \null\nobreak\hfill
   \kern3pt\QEDbox} 

\newif\ifQEDdejaplace\QEDdejaplacefalse

\theoremstyle{mypreuvestyle}
\newtheorem*{proclaimmypreuve}{}
\newenvironment{preuve}[1][*]{\begin{proclaimmypreuve}{\ifx*#1{\em Proof.}\else{\em#1.}\fi} \ignorespaces}{\ifQEDdejaplace\global\QEDdejaplacefalse\else\QED\fi\end{proclaimmypreuve}}

\newif\ifmynonumberequation\mynonumberequationtrue

\makeatletter
\numberwithin{equation}{section}
\newenvironment{moneq}[1][*]{\ifx*#1\mynonumberequationtrue\begin{equation*}\else\mynonumberequationfalse\begin{equation}\label{#1}\fi}{\ifmynonumberequation\end{equation*}\@ignoretrue\else\end{equation}\@ignoretrue\fi\ignorespaces}
\makeatother

\newcommand{\recalf}[1]{(\ref{#1})}
\newcommand{\recals}[1]{{\S\ref{#1}}}
\newcommand{\recalt}[1]{{[\ref{#1}]}}

\allowdisplaybreaks

\def\oversetalign#1\to#2{\mathbin{\smash{\overset{{\text{\rlap{\hss#1}}}}{#2}}}}
\def\oversettext#1\to#2{\mathbin{\smash{\overset{{\text{{#1}}}}{#2}}}}

\def\dcap_#1{\mathchoice{%
          {\textstyle\bigcap\limits_{#1}}}%
          {\underset{#1}\cap}%
          {\underset{#1}\cap}%
          {\underset{#1}\cap}}
\def\dcup_#1{\mathchoice{%
          {\textstyle\bigcup\limits_{#1}}}%
          {\underset{#1}\cup}%
          {\underset{#1}\cup}%
          {\underset{#1}\cup}}
\def\ddcap_#1^#2{\mathchoice{%
          {\textstyle\bigcap\limits_{#1}^{#2}}}%
          {\overset{#2}{\underset{#1}\cap}}%
          {\overset{#2}{\underset{#1}\cap}}%
          {\overset{#2}{\underset{#1}\cap}}}
\def\ddcup_#1^#2{\mathchoice{%
          {\textstyle\bigcup\limits_{#1}^{#2}}}%
          {\overset{#2}{\underset{#1}\cup}}%
          {\overset{#2}{\underset{#1}\cup}}%
          {\overset{#2}{\underset{#1}\cup}}}

\newskip\myICS\myICS=0.0pt

\newcommand\intervoo[3][1]{\ifcase#1
   \mskip2.0mu\left(\mskip1mu 
   #2,#3
   \mskip1mu\right)\mskip2.0mu
\or 
   \mathopen{\mskip2mu(\mskip1mu}
   #2,#3
   \mathclose{\mskip1mu)\mskip2.0mu}
\or 
   \mathopen{\mskip2mu\bigl(\mskip1mu}
   #2\kern#1\myICS\kern-\myICS,\kern#1\myICS\kern-\myICS#3
   \mathclose{\mskip1mu\bigr)\mskip2.0mu}
\or 
   \mathopen{\mskip2mu\Bigl(\mskip1mu} 
   #2\kern#1\myICS\kern-\myICS,\kern#1\myICS\kern-\myICS#3
   \mathclose{\mskip1mu\Bigr)\mskip2.0mu}
\or 
   \mathopen{\mskip2mu\biggl(\mskip1mu} 
   #2\kern#1\myICS\kern-\myICS,\kern#1\myICS\kern-\myICS#3
   \mathclose{\mskip1mu\biggr)\mskip2.0mu}
\or 
   \mathopen{\mskip2mu\Biggl(\mskip1mu} 
   #2\kern#1\myICS\kern-\myICS,\kern#1\myICS\kern-\myICS#3
   \mathclose{\mskip1mu\Biggr)\mskip2.0mu}
\else 
   \mskip2.0mu\left(\mskip1mu 
   #2\kern#1\myICS\kern-\myICS,\kern#1\myICS\kern-\myICS#3
   \mskip1mu\right)\mskip2.0mu
\fi}

\newcommand\intervof[3][1]{\ifcase#1
   \mskip2.0mu\left(\mskip1mu 
   #2,#3
   \mskip2.5mu\right]\mskip0.5mu
\or 
   \mathopen{\mskip2.0mu(\mskip1mu}
   #2,#3
   \mathclose{\mskip2.5mu]\mskip0.5mu}
\or 
   \mathopen{\mskip2.0mu\bigl(\mskip1mu}
   #2\kern#1\myICS\kern-\myICS,\kern#1\myICS\kern-\myICS#3
   \mathclose{\mskip2.5mu\bigr]\mskip0.5mu}
\or 
   \mathopen{\mskip2mu\Bigl(\mskip1mu} 
   #2\kern#1\myICS\kern-\myICS,\kern#1\myICS\kern-\myICS#3
   \mathclose{\mskip2.5mu\Bigr]\mskip0.5mu}
\or 
   \mathopen{\mskip2.0mu\biggl(\mskip1mu} 
   #2\kern#1\myICS\kern-\myICS,\kern#1\myICS\kern-\myICS#3
   \mathclose{\mskip2.5mu\biggr]\mskip0.5mu}
\or 
   \mathopen{\mskip2.0mu\Biggl(\mskip1mu} 
   #2\kern#1\myICS\kern-\myICS,\kern#1\myICS\kern-\myICS#3
   \mathclose{\mskip2.5mu\Biggr]\mskip0.5mu}
\else 
   \mskip2.0mu\left(\mskip1mu 
   #2\kern#1\myICS\kern-\myICS,\kern#1\myICS\kern-\myICS#3
   \mskip2.5mu\right]\mskip0.5mu
\fi}

\newcommand\intervfo[3][1]{\ifcase#1
   \mskip0.5mu\left[\mskip2mu 
   #2,#3
   \mskip1mu\right)\mskip2.0mu
\or 
   \mathopen{\mskip0.5mu[\mskip2mu}
   #2,#3
   \mathclose{\mskip1mu)\mskip2.0mu}
\or 
   \mathopen{\mskip0.5mu\bigl[\mskip2mu}
   #2\kern#1\myICS\kern-\myICS,\kern#1\myICS\kern-\myICS#3
   \mathclose{\mskip1mu\bigr)\mskip2.0mu}
\or 
   \mathopen{\mskip0.5mu\Bigl[\mskip2mu}
   #2\kern#1\myICS\kern-\myICS,\kern#1\myICS\kern-\myICS#3
   \mathclose{\mskip1mu\Bigr)\mskip2.0mu}
\or 
   \mathopen{\mskip0.5mu\biggl[\mskip2mu}
   #2\kern#1\myICS\kern-\myICS,\kern#1\myICS\kern-\myICS#3
   \mathclose{\mskip1mu\biggr)\mskip2.0mu}
\or 
   \mathopen{\mskip0.5mu\Biggl[\mskip2mu}
   #2\kern#1\myICS\kern-\myICS,\kern#1\myICS\kern-\myICS#3
   \mathclose{\mskip1mu\Biggr)\mskip2.0mu}
\else 
   \mskip0.5mu\left[\mskip2mu 
   #2\kern#1\myICS\kern-\myICS,\kern#1\myICS\kern-\myICS#3
   \mskip1mu\right)\mskip2.0mu
\fi}

\newcommand\intervff[3][1]{\ifcase#1
   \mskip0.5mu\left[\mskip2mu 
   #2,#3
   \mskip2.5mu\right]\mskip0.5mu
\or 
   \mathopen{\mskip0.5mu[\mskip2mu}
   #2,#3
   \mathclose{\mskip2.5mu]\mskip0.5mu}
\or 
   \mathopen{\mskip0.5mu\bigl[\mskip2mu}
   #2\kern#1\myICS\kern-\myICS,\kern#1\myICS\kern-\myICS#3
   \mathclose{\mskip2.5mu\bigr]\mskip0.5mu}
\or 
   \mathopen{\mskip0.5mu\Bigl[\mskip2mu}
   #2\kern#1\myICS\kern-\myICS,\kern#1\myICS\kern-\myICS#3
   \mathclose{\mskip2.5mu\Bigr]\mskip0.5mu}
\or 
   \mathopen{\mskip0.5mu\biggl[\mskip2mu}
   #2\kern#1\myICS\kern-\myICS,\kern#1\myICS\kern-\myICS#3
   \mathclose{\mskip2.5mu\biggr]\mskip0.5mu}
\or 
   \mathopen{\mskip0.5mu\Biggl[\mskip2mu}
   #2\kern#1\myICS\kern-\myICS,\kern#1\myICS\kern-\myICS#3
   \mathclose{\mskip2.5mu\Biggr]\mskip0.5mu}
\else 
   \mskip0.5mu\left[\mskip2mu 
   #2\kern#1\myICS\kern-\myICS,\kern#1\myICS\kern-\myICS#3
   \mskip2.5mu\right]\mskip0.5mu
\fi}

\newcommand\norme[2][1]{\ifcase#1
   \left\Vert#2\right\Vert 
   \or\Vert#2\Vert 
   \or\bigl\Vert#2\bigr\Vert 
   \or\Bigl\Vert#2\Bigr\Vert 
   \or\biggl\Vert#2\biggr\Vert 
   \or\Biggl\Vert#2\Biggr\Vert 
   \else\Vert#2\Vert\fi} 

\newcommand\inprod[3][1]{\ifcase#1
   \left\langle\mskip1.5mu #2,#3\mskip1.5mu \right\rangle 
   \or\langle\mskip1.5mu #2,#3\mskip1.5mu \rangle 
   \or\bigl\langle\mskip1.5mu #2,#3\mskip1.5mu \bigr\rangle 
   \or\Bigl\langle\mskip1.5mu #2,#3\mskip1.5mu \Bigr\rangle 
   \or\biggl\langle\mskip1.5mu #2,#3\mskip1.5mu \biggr\rangle 
   \or\Biggl\langle\mskip1.5mu #2,#3\mskip1.5mu \Biggr\rangle 
   \else\langle\mskip1.5mu #2,#3\mskip1.5mu \rangle\fi} 
   
\newcommand\inprodsym{\inprod\cdot\cdot}

\newcommand\inprodd[3][1]{\ifcase#1
   \left\langle\mkern-3mu\left\langle#2,#3
        \right\rangle\mkern-3mu\right\rangle 
   \or\langle\mkern-3mu\langle#2,#3
        \rangle\mkern-3mu\rangle 
   \or\bigl\langle\mkern-4mu\bigl\langle#2,#3
        \bigr\rangle\mkern-4mu\bigr\rangle 
   \or\Bigl\langle\mkern-6mu\Bigl\langle#2,#3
        \Bigr\rangle\mkern-6mu\Bigr\rangle 
   \or\biggl\langle\mkern-8mu\biggl\langle#2,#3
        \biggr\rangle\mkern-8mu\biggr\rangle 
   \or\Biggl\langle\mkern-10mu\Biggl\langle#2,#3
        \Biggr\rangle\mkern-10mu\Biggr\rangle 
   \else\langle\mkern-3mu\langle#2,#3
        \rangle\mkern-3mu\rangle\fi} 

\newcommand\supinprsym{\mathscr{S}}

\newcommand\superinprod[3][1]{\ifcase#1
   \supinprsym\left(#2,#3\right) 
   \or\supinprsym(#2,#3) 
   \or\supinprsym\bigl(#2,#3\bigr) 
   \or\supinprsym\Bigl(#2,#3\Bigr) 
   \or\supinprsym\biggl(#2,#3\biggr) 
   \or\supinprsym\Biggl(#2,#3\Biggr) 
   \else\supinprsym(#2,#3)\fi} 

\newcommand\homsuperinprod[4][1]{\ifcase#1
   \supinprsym{#4}\left(#2,#3\right) 
   \or\supinprsym#4(#2,#3) 
   \or\supinprsym#4\bigl(#2,#3\bigr) 
   \or\supinprsym#4\Bigl(#2,#3\Bigr) 
   \or\supinprsym#4\biggl(#2,#3\biggr) 
   \or\supinprsym#4\Biggl(#2,#3\Biggr) 
   \else\supinprsym#4(#2,#3)\fi} 

\newcommand\contrfoper{{\iota}}
\newcommand\contrs[2]{{\langle\mskip1.5mu}
    #1\mskip1.5mu |\mskip1.5mu #2
    {\mskip1.5mu\rangle}}

\newcommand\contrf[3][1]{\ifcase#1
   \contrfoper\left(#2\right)#3 
   \or\contrfoper(#2)#3 
   \or\contrfoper\bigl(#2\bigr)#3 
   \or\contrfoper\Bigl(#2\Bigr)#3 
   \or\contrfoper\biggl(#2\biggr)#3 
   \or\contrfoper\Biggl(#2\Biggr)#3 
   \else\contrfoper(#2)#3\fi} 

\newcount\mysubsectioncount\mysubsectioncount=0
\newcommand\firstofmysubsection{\mysubsectioncount=0}
\newcommand\mysubsection[2]{\vskip2\baselineskip\advance\mysubsectioncount by1\noindent{\textbf{\ref{#1}.\the\mysubsectioncount{} #2}}\medskip}

\newif\ifKortereTekst\KortereTekstfalse
\newcommand\KorteOfLangeTekst[2]{\ifKortereTekst#1\else#2\fi}
\long\def\dontprintcomp#1{#1}

\newcommand\bigrestricted{{\kern1pt\vrule height3.3ex depth1.7ex width0.6pt\kern1pt}}
\newcommand\caprestricted{{\kern1pt\vrule height1.7ex depth0.5ex width0.6pt\kern1pt}}
\newcommand\restricted[1][1.3]{{\kern1pt\vrule height#1ex depth0.5ex width0.6pt\kern1pt}}

\newcommand\Ad{{\operatorname{Ad}}}
\def\aglgrp/{an \glgrp/}
\def\agmfd/{an \gmfd/}
\def\agvs/{an \gvs/}
\newcommand\ah{{\blauw{\widehat a}}}
\newcommand\alphah{\blauw{\widehat \alpha}}
\newcommand\adh{\rood{
a}}
\newcommand\alphadh{\rood{
\alpha}}
\newcommand\Aut{\operatorname{Aut}}
\newcommand\betah{\blauw{\widehat\beta}}
\newcommand\bh{\blauw{\widehat b}}
\newcommand\betadh{\rood{
\beta}}
\newcommand\bdh{\rood{
b}}
\newcommand\blauw[1]{#1}
\newcommand\body{{\mathbf B}}
\newcommand\CA{{\mathcal A}}
\newcommand\CC{{\mathbf C}}
\newcommand\ch{\blauw{\widehat c}}
\newcommand\cdh{\rood{
c}}

\newcommand\Coad{{\operatorname{Coad}}}
\newcommand\conjugate{\mathfrak{C}}
\newcommand\datum{\the\day/\the\month/\the\year}
\newcommand\Dense{\mathcal{D}}

\newcommand\Denset{\widetilde{\Dense}}
\newcommand\eexp{{\operatorname{e}}}
\newcommand\End{\operatorname{End}}
\newcommand\extder{{\operatorname{d}}}
\newcommand\extdersep{\kern0.75em}
\newcommand\etab{{\bar\eta}}
\newcommand\FgroupHilbert{\mathbf{F}}

\newcommand\Fourier{\mathcal{F}}
\newcommand\Fourierodd{\smash{\overset{\kern2pt\raise-2pt\hbox{$\scriptscriptstyle \mathrm o$}}{\Fourier}}}
\newcommand\Fourieroddp{\Fourierodd^{\,'}}
\newcommand\Fourieroddpp{\Fourierodd^{\,''}}
\newcommand\fracp[2]{\frac{\partial #1}{\partial #2}}
\newcommand\gammah{\blauw{\widehat\gamma}}
\newcommand\gammadh{\rood{
\gamma}}
\newcommand\gh{{\blauw{\hat g}}}
\newcommand\GL{\mathrm{GL}}
\def\glalg/{$\CA$-Lie algebra}
\def\glgrp/{$\CA$-Lie group}
\def\gmfd/{$\CA$-mani\-fold}
\newcommand\grel{{(\ah,\bh,\ch)}}
\def\gvs/{$\CA$-vector space}
\newcommand\hBigl[1]{\raise0.6pt\hbox{\footnotesize$\Bigl#1$}}
\newcommand\hBigr[1]{\raise0.6pt\hbox{\footnotesize$\Bigr#1$}}
\newcommand\Hilbert{\mathcal{H}}
\newcommand\id{\textrm{id}}
\newcommand\Identif{\mathcal{I}}
\newcommand\ie{i.e.}
\newcommand\itemizedef[1][œ]{\ifx#1*\else\par\medskip\noindent\fi{\large$\bullet$} \ignorespaces}
\newcommand\Itemizedef[1][œ]{\ifx#1*\else\par\bigskip\noindent\fi$\vrule width1.2ex height1.2ex depth0pt$\kern0.5em {}\ignorespaces}
\newcommand\itemizepreuve[2][œ]{\ifx#1*\else\par\medskip\noindent\fi$\bullet$ (\ref{#2}). \ignorespaces}
\newcommand\KK{\mathbf{K}}
\newcommand\KKA{\mathcal{K}}
\newcommand\labelFamilySUR{pqFSUR\kern0.1em\relax}
\newcommand\labelIsomSHS{IsoSH\kern0.1em\relax}
\newcommand\labelMainSUR{SUR\kern0.1em\relax}
\newcommand\labelMainSUReq{SUREQ\kern0.1em\relax}
\newcommand\LD[2][0]{\ifcase#1 {}^*#2\or \kern1pt{}^*\kern-1pt#2\or\kern2pt{}^*\kern-2pt#2\else{}^*#2\fi}
\renewcommand\LD[2][0]{#2^*}
\newcommand\Leb{\mathcal{L}\kern-1pt\mathit{eb}}
\newcommand\Liealg[1]{\mathfrak{#1}}
\newcommand\Lieder{{\mathcal L}}

\newcommand\mapob{\kern.3em} 
\newcommand\mfdmetric{\mathbf{g}}
\newcommand\mo{^{-1}}
\newcommand\myquote[1]{``#1''}
\newcommand\nilpotent{\mathcal{N}}
\newcommand\NN{{\mathbf N}}
\newcommand\numeral[1]{$^\textrm{#1}$}
\newcommand\Orbit{{\mathcal O}}
\newcommand\parity[1]{\vert#1\vert}
\newcommand\Phih{{\widehat\Phi}}
\newcommand\psih{\hat\psi}
\newcommand\psur{proto super unitary representation}

\newcommand\refmetnaam[2]{#1\ref{#2}}
\def\repcomm#1,#2{\bigl[[#1,(#2)]\bigr]}
\newcommand\repprod[1]{(\!(#1)\!)}
\newcommand\rhob{{\bar\rho}}
\newcommand\rhoh{{\hat\rho}}
\newcommand\rhot{{\tilde\rho}}
\newcommand\rood[1]{#1}
\newcommand\RR{{\mathbf R}}
\newcommand\Schwartz[2]{%
\def\EmPtY{}\def\NrTwoEmpty{#2}\ifx\NrTwoEmpty\EmPtY
                   \mathcal{S}(\RR^{#1})\else\mathcal{S}(\RR^{#1\vert#2})\fi}
\newcommand\scirc{\,{\raise 0.8pt\hbox{$\scriptstyle\circ$}}\,}
\newcommand\shifttag[1]{\kern#1&\kern-#1}
\newcommand\smooths{\mathbf{F}}
\newcommand\Sphere{{\mathbf S}}
\newcommand\stresd[1]{{\emph{#1}}}
\newcommand\stress[1]{\textbf{\textit{#1}}}
\newcommand\subhead[1]{\vskip\baselineskip\noindent\textbf{#1}.}
\newcommand\SuperH{\mathit{S}\kern-0.25em\mathit{H}}
\newcommand\supscal{super scalar }
\newcommand\SwitchHats{%
\renewcommand\ah{a}\renewcommand\adh{{\hat a}}
\renewcommand\bh{b}\renewcommand\bdh{{\hat b}}
\renewcommand\ch{c}\renewcommand\cdh{{\hat c}}
\renewcommand\alphah{\alpha}\renewcommand\alphadh{{\hat\alpha}}
\renewcommand\betah{\beta}\renewcommand\betadh{{\hat\beta}}
\renewcommand\gammah{\gamma}\renewcommand\gammadh{{\hat\gamma}}
\renewcommand\gh{g}
}
\newcommand\BackSwitchHats{%
\renewcommand\ah{{\hat a}}\renewcommand\adh{a}
\renewcommand\bh{{\hat b}}\renewcommand\bdh{b}
\renewcommand\ch{{\hat c}}\renewcommand\cdh{c}
\renewcommand\alphah{{\hat\alpha}}\renewcommand\alphadh{\alpha}
\renewcommand\betah{{\hat\beta}}\renewcommand\betadh{\beta}
\renewcommand\gammah{{\hat\gamma}}\renewcommand\gammadh{\gamma}
\renewcommand\gh{{\hat g}}
}

\newcommand\xb{{\bar x}}
\newcommand\xh{{\hat x}}
\newcommand\xib{{\bar\xi}}
\newcommand\xih{{\hat \xi}}
\newcommand\xnul{x^o}
\newcommand\xnulb{\xb^o}
\newcommand\yb{{\bar y}}
\newcommand\yh{{\hat y}}
\newcommand\zb{{\bar z}}

\newcommand\ZZ{{\mathbf Z}}

\begin{document}
\renewcommand\labelitemi{\textbullet}
\setcounter{MaxMatrixCols}{20}

\title[Decomposing super unitary representations]{Decomposing super unitary representations using Berezin integration}
\address{Laboratoire Paul Painlevé, U.M.R. CNRS 8524 et Département de Mathématiques, Faculté des Sciences et Technologies, Université de Lille, 59655 Villeneuve d'Ascq Cedex, France}
\email{FirstName[dot]LastName[at]univ-lille[dot]fr}
\author{Gijs M. Tuynman 
}
\thanks{Publishing history: an earlier version of this paper \cite{Tuynman..2009} was never published and \cite{Tuynman..2017} was (for reasons of length) only partially published as \cite{Tuynman..2018}. This new version \myquote{thus} is a rewrite of \cite{Tuynman..2009}, incorporating the notion of a super unitary representation as introduced in \cite{Tuynman..2018} and replacing\slash improving the discussion of the Berezin-Fourier transform with material from \cite{Tuynman..2017} that was excluded in \cite{Tuynman..2018}.}

\begin{abstract}  
I introduce the notion of a family of super unitary representations depending upon odd parameters as well as the notion of a direct Berezin integral of such representations. 
I use this to decompose the (left-)regular representation of the super Heisenberg group and of a super Heisenberg-like group into smaller (irreducible) subrepresentations. 
I compare the obtained subrepresentations to the representations obtained by applying (super) geometric quantization of the coadjoint orbits of the corresponding groups. 
It turns out that for the super Heisenberg group we need coadjoint orbits with an odd symplectic form in order to describe all subrepresentations found in the decomposition and for the Heisenberg-like super group we even need orbits with a non-homogeneous symplectic form. 
\end{abstract}

\maketitle

\tableofcontents

%
%

\section{Introduction\slash overview}

In this paper I present some kind of \myquote{experimental} mathematics in the sense that I perform experiments\slash computations and that I try to find a theory that best fits these experimental data. 
More precisely, I do two different experiments: decomposing the left-regular representation of a (super) Lie group into irreducible components and computing the representations obtained by (super) geometric quantization of coadjoint orbits. 
These experiments start with the standard Heisenberg group $H$ as template, simply because this is a well understood case\slash experiment\slash computation. 
I then look at the super Heisenberg group $\SuperH$ by declaring that two of the three coordinates of $H$ are odd and I simply imitate\slash copy the computations done for the standard Heisenberg group. 
This can be done without any problem, except that in the ungraded case we use direct integrals of representations in the decomposition of $H$ (Fourier decomposition), which in the super case become Berezin-Fourier integrals. 
And I finish with a Heisenberg-like super Lie group $G$ where I imitate the kind of computations done for $\SuperH$ (the decomposition of the left-regular representations as well as for a comparison with representations obtained by geometric quantization of coadjoint orbits). 

In imitating the decomposition of the Heisenberg group using Fourier transforms in the super context, these Fourier integrals become (partially) integrals over odd coordinates. 
This has two consequences for a theoretical framework: we have to adapt the underlying Hilbert spaces to odd coordinates and we have to show that, with the chosen definition of these underlying super Hilbert spaces, our Berezin-Fourier transforms are super unitary equivalences. 
In \recals{SECTInitialDefinitionsAndMore}-\ref{SECTpVertqFamiliesOfSURs} I provide a suitable theoretical framework that fits all computations presented in this paper. 

For the second set of \myquote{experiments} --- the computation of the representations associated to coadjoint orbits and the link with the irreducible representations found in the decomposition of the left-regular representation---the situation is less straightforward. 
Not all representations obtained in the decomposition of the regular representation can be linked to coadjoint orbits with an even symplectic form. 
To obtain a satisfactory link between coadjoint orbits and the irreducible representations in the decomposition of the regular representation, it turns out that for the super Heisenberg group $\SuperH$ we need coadjoint orbits with an odd symplectic form and for the Heisenberg-like super group $G$ we even need coadjoint orbits with a non-homogeneous symplectic form as described in \cite{Tuynman..2010}.\footnote{This was and is the main reason for adding this group to my experiments.} 

When performing the computations for geometric quantization of the coadjoint orbits, I need (a generalisation of) the notion of a polarization. 
And here appear the first differences with the non-super case: even though the total dimension of a polarization is constant (at least in these examples), its graded dimension is not. 
Moreover, the representations associated to a same orbit but with polarizations of different graded dimensions are not equivalent. 
Now this may seem disconcerting, but it fits nicely with computations for the super Heisenberg group $\SuperH$. 
When decomposing the classical Heisenberg group $H$ using the Fourier transform, there are two natural choices, choices that are decribed by physicists as the position representation and the momentum representation respectively. In terms of polarizations this corresponds to the vertical and horizontal polarization. 
And of course, in the classical case these two representations are equivalent. 
But for the super Heisenberg group the horizontal and vertical polarizations have different graded dimensions, yielding (obviously) non-equivalent representations. 
This fits nicely with the two natural choices in the decomposition: the two natural Fourier transforms for $H$ (which yield the position and momentum representations) become different for $\SuperH$, as one will be done with an even coordinate, whereas the other has to be done with an odd coordinate. 
We thus have two essentially different ways to decompose $\SuperH$, but these different ways correspond to different choices for a polarization on the same coadjoint orbit. 

Contrary to the situation for the decomposition, I only have a partial theoretical framework (given in \recals{GQgeneral}) for the computation of a representation associated to a super coadjoint orbit via super geometric quantization. 
I have a framework for super prequantization, the notion of a super polarization and a way to associate function spaces to these notions on which the (super) group acts. 
But I do not have a way to generalise the notion of half-densities or half-forms that would allow me to define a (natural) structure of a super Hilbert space on these function spaces. 
Moreover, the geometric quantization procedure suggests (strongly) to introduce an odd Berezin-Fourier transform and this Berezin-Fourier transform introduces an odd (Fourier) parameter $\kappa$. 
And although it nicely gives us representations that correspond to representations obtained in the decomposition of the regular representation, it also poses a problem when comparing to $p\vert q$-parameter families of representations that appear in such a decomposition: all odd parameters (all $q$ of them) have to be real multiples of this same $\kappa$. 

Let me finish this introduction with the remark that I will work with the geometric $H^\infty$ version of DeWitt supermanifolds,
which is equivalent to the theory of graded manifolds of Leites and Kostant
(see \cite{DW84}, \cite{Ko77}, \cite{Le80}, \cite{Ro07}, \cite{Tu04}). Any reader using a
(slightly) different version of supermanifolds should be able to translate the
results to her\slash his version of supermanifolds. 
Readers unfamiliar with my approach will find an extremely succinct overview of $\CA$-manifold theory in the appendix.

\section{\label{SECTInitialDefinitionsAndMore}Super Hilbert spaces}

Let me start very naively: a representation of a Lie group $G$ on a vector space $E$ is a homomorphism $\rho:G\to \Aut(E)$ and two representations of $G$ on $E_1$ and $E_2$ are equivalent if there exists an isomorphism (linear bijection) $\Phi:E_1\to E_2$ that intertwines the two representations. 
If the vector space $E$ is equipped with an inner product, we can define the notion of a unitary representation as a representation that \myquote{preserves} the inner product. 
And then two unitary representations are equivalent if there exists a unitary isomorphism $\phi:E_1\to E_2$ that intertwines the representations. 

When we transcribe this to the super setting (more or less putting everywhere the adjective super\slash graded), we arrive at the description that a (super) representation of \aglgrp/ $G$ on a graded vector space $E$ is a homomorphism $\rho:G\to \Aut(E)$. 
For a super unitary representation we then will need a \supscal product that should be preserved by the representation, and two super unitary representations are equivalent when there exists a linear bijection $\phi:E_1\to E_2$ that intertwines the representations and is compatible with the \supscal products. 

When we remain in the realm of finite dimensional (graded) vector spaces, this transcription does not pose any problems. 
However, when going to infinite dimensional (graded) vector spaces, the situation becomes less clear. 
For (ordinary) unitary representations, one requires that the vector space be a Hilbert space, \ie, complete with respect to the topology induced by the inner product. 
Unfortunately, a \supscal product is not positive definite and there is no natural way to define a topology associated to a \supscal product, hence no way to define a notion of completeness. 
So we have to ask ourselves the question: is this a problem?
Or, changing our point of view, why does one want a Hilbert space in the ungraded case? 
One (possible) answer is quite obvious: to restrict the number of inequivalent unitary representations. 
When we allow non-complete vector spaces, any invariant dense subspace of a representation space generates a new unitary representation. 
And on the other hand, any unitary representation on a non-complete vector space induces a unique unitary representation on its completion. 
In short: restricting attention to Hilbert spaces does not make an essential difference and reduces the number of inequivalent unitary representations. 

The preceding discussion (strongly?) suggests that we should introduce a topology and\slash or an inner product on a representation space for a super unitary representation in order to achieve a similar reduction of super unitary representations. 
Once we admit that there should be two structures on a graded vector space in order to be a representation space for a super unitary representation: a \supscal product $\supinprsym$ for the notion of super-unitary and an inner product $\inprodsym$ to define a topology and the notion of completeness, we should specify if and how these two structures are related. 
The definition usually adopted requires a direct link between the \supscal product $\supinprsym$ and the inner product $\inprodsym$. 
In \cite{Tuynman..2018} I argued that in this approach for most \glgrp/s the (left) regular representation will not be super unitary. 
To circumvent this restriction I weakened this condition by requiring only that $\supinprsym$ be continuous with respect to the topology induced by $\inprodsym$. 
In \recalt{ExampleWithDifferentTopologiesInducedByChangeOfCoordinates} I will give another reason to weaken the link between the \supscal product and the inner product by showing that when there is a fixed (algebraic) relation between the \supscal product $\supinprsym$ and the inner product $\inprodsym$, then a simple change of coordinates on an underlying super manifold will yield inequivalent super unitary representations. 
Having defined a super inner product space as a graded vector space equipped with two structures $\supinprsym$ and $\inprodsym$, we can define a super Hilbert space as a super inner product space that is complete with respect to the inner product $\inprodsym$. 
This gives us the following definitions.

\begin{definition}[DefinitionOfSuperInnerProductSpaces]{Definitions: super inner product spaces}
Let $\KKA$ denote any of the four rings\slash fields $\RR$, $\CC$, $\CA$ or $\CA^\CC$ and let $V$ be a graded vector space over $\KKA$. 

\itemizedef
I will reserve the name \stresd{inner product} for the ungraded notion of a positive definite symmetric sesquilinear form on a real or complex vector space (with or without grading).\footnote{And I will follow the physicists convention that it is linear in the second variable and anti-linear in the first.} 

\itemizedef
A \stresd{\supscal product $\supinprsym$ on a graded vector space $V$} (over $\KKA$) is a smooth graded symmetric non-degenerate graded right-sesquilinear form on $V$, \ie, a map $\supinprsym$ defined on $V\times V$ and with values in $\KKA$ satisfying the conditions
\begin{enumerate}
\item
$\superinprod{v_1+v_2\lambda}w = \superinprod{v_1}w + \superinprod{v_2}{\overline\lambda w}$ (right-anti-linear in the first entry);

\item
$\superinprod v{w_1+w_2 \lambda} = \superinprod v{w_1} + \superinprod v{w_2}\cdot \lambda$ (right-linear in the second entry);

\item
$\superinprod w v = (-1)^{\parity w \parity v}\,\overline{\superinprod vw}$ (graded symmetry);

\item
$\forall v\neq0\ \exists w:\superinprod vw\neq0$ (non-degeneracy);

\item
$\superinprod{\body V}{\body V}\subset \body\KKA$ (smoothness);

\end{enumerate}
Note that we explicitly do not require $\supinprsym$ to be even nor homogeneous (its \myquote{parity} is \myquote{concentrated} in the symbol $\supinprsym$). 

\itemizedef
A right-linear map $\Phi:V\to V$ is said to be \stresd{graded skew-symmetric} with respect to a \supscal product $\supinprsym$ on $V$ if we have (see \recalf{DefinitionOfConjugationOperator} for the definition of $\conjugate$)
\begin{moneq}
\superinprod{\Phi v}w + \superinprod v{\Phi_0w} + \superinprod[2]{\conjugate(v)}{\Phi_1w}
=
0
\mapob,
\end{moneq}
where $\Phi=\Phi_0+\Phi_1$ is the decomposition of $\Phi$ into its homogeneous parts.

\itemizedef
A \stresd{super inner product space} is a triple $(E,\inprodsym, \supinprsym)$ in which $E$ is a (complex) graded vector space, $\inprodsym$ is an inner product on $E$ (\ie, $(E,\inprodsym)$ is a hermitian space) and $\supinprsym$ is a \supscal product on $E$ such that 
\begin{enumerate}

\item
$\inprod{E_0}{E_1}=0$;

\item
$\supinprsym$ is continuous with respect to the topology of $(E,\inprodsym)$.

\end{enumerate}
The triple $(E,\inprodsym, \supinprsym)$ will be called a \stresd{super Hilbert space} when $(E,\inprodsym)$ is an ordinary Hilbert space, \ie, complete with respect to the inner product $\inprodsym$. 

\itemizedef
Let $(E, \inprodsym_E, \supinprsym_E)$ and $(F, \inprodsym_F, \supinprsym_F)$ be two super inner product spaces and let $k$ be a non-zero real constant. 
We will say that these two spaces are \stresd{$k$-isomorphic}\footnote{Motivations for this definition of equivalence with $k=1$ can be found in \cite[Rem. 3.4]{Tuynman..2018}. The reason to add the non-zero constant $k$ is explained in \recalt{BerezinIntegrationAndMotivationkEquivalence}.} if there exists a linear bijection $\Phi:E\to F$ satisfying the following conditions (where $\Phi=\Phi_0+\Phi_1$ and $x=x_0+x_1$ are the decompositions of $\Phi$ respectively $x$ into homogeneous parts).
\begin{enumerate}[\kern0.5em\labelIsomSHS1.]
\item\label{conditiondirectsuminequivalenceSHS}
The homogeneous parts of $\Phi$ do not interfere in the sense that $\ker(\Phi_0) + \ker(\Phi_1) = E$ (and thus $E=\ker(\Phi_0) \oplus \ker(\Phi_1)\,$).

\item\label{conditionequivalenceofSHSisunitary}
$\Phi$ intertwines the inner products: $\forall x,y\in E$: $\inprod{\Phi x}{\Phi y}_F = \inprod xy_E$.

\item\label{conditionoddisoforequivalenceofSHS}
$\Phi$ intertwines the \supscal products up to the constant $k$ in the sense that for all $x,y\in E$ we have
\begin{moneq}
k\cdot \homsuperinprod{\Phi x}{\Phi y}{_F}
=
\begin{cases}
\homsuperinprod{x}{y}{_E}
&\qquad x \in E\ ,\ y\in \ker(\Phi_1)
\\
\homsuperinprod{\conjugate x}{y}{_E}
&\qquad x\in \ker(\Phi_1) \ , \ y\in \ker(\Phi_0)
\\
i\,\homsuperinprod{\conjugate x}{y}{_E}
&\qquad x,y\in \ker(\Phi_0) 
\mapob.
\end{cases}
\end{moneq}

\end{enumerate}
Alternatively we will say that a map $\Phi$ satisfying these conditions is a \stresd{$k$-isomorphism} between the two super inner product spaces. 
In the case $k=1$ we will simply say that $\Phi$ is an isomorphism instead of a $1$-isomorphism and that $E$ and $F$ are isomorphic instead of $1$-isomorphic. 

\end{definition}

\begin{remarque}{Remark}
As (by definition) any graded vector space $V$ over $\CA$ respectively $\CA^\CC$ is of the form $V=E\otimes\CA$ (resp. $E\otimes \CA^\CC$) for a graded vector space $E$ over $\RR$ (resp. $\CC$), the smoothness condition for a super scalar product implies that restriction to $E\cong E\otimes \{1\}\subset V$ establishes a bijection between super scalar products on $V$ and super scalar products on $E$. 
\end{remarque}

\begin{proclaim}[ExtensionSuperScalarProductToCompletion]{Lemma}
Let $(E, \inprodsym_E, \supinprsym_E)$ be a super inner product space. 
\begin{enumerate}
\item
Let $(\Hilbert,\inprodsym_\Hilbert)$ be the completion of $(E,\inprodsym_E)$. 
If the extension $\supinprsym_\Hilbert$ of $\supinprsym_E$ to $\Hilbert$ by continuity is non-degenerate, then $(\Hilbert,\inprodsym_\Hilbert, \supinprsym_\Hilbert)$ is a super Hilbert space. 

\item
If $F\subset E$ is a dense subspace, then $(F, \inprodsym_E\caprestricted_F, \supinprsym_E\caprestricted_F)$ is a super inner product space. 

\end{enumerate}

\end{proclaim}

\begin{remarque}{Remark}
Let $(E,\inprodsym, \supinprsym)$ be a super inner product space. 
In view of \recalt{ExtensionSuperScalarProductToCompletion}, it might be better to replace the condition that $\supinprsym$ should be continuous with respect to the topology of $(E,\inprodsym)$ by the condition that there exist two constants $c,C>0$ such that
\begin{moneq}
\forall x\in E\setminus\{0\}
:
c\le \sup_{y\in E\setminus\{0\}}\ \frac{\vert \superinprod xy\vert}{\norme x\cdot\norme y} \le C
\mapob.
\end{moneq}
The second inequality is the condition that $\supinprsym$ is continuous and the first inequality is slightly stronger than non-degeneracy (together they say that the map $x\mapsto \superinprod x\cdot$ is a bicontinuous isomorphism from $E$ to its image in $E^*$, the topological dual of $E$). 
It has the advantage that it passes on to a completion, thus assuring that the extension of $\supinprsym$ to a completion is again (slightly more than) non-degenerate. 
\end{remarque}

\begin{definition}[BerezinIntegrationAndMotivationkEquivalence]{Berezin integration and orientation}
If $\psi(x,\xi)$ is a (smooth) function of $p$ even and $q$ odd variables with compact support (in the even variables), then Berezin integration of $\psi$ is defined as follows:
we write (see \recalt{DefinitionRepeatedProducts} and \recalt{structureofsmoothfunctions} for notation)
\begin{moneq}[StructureOfSmoothFunctionInPowersXiAndFunctionsX]
\psi(x,\xi) = \sum_{I\subset \{1, \dots, q\}} \repprod\xi^I\cdot \psi_I(x)
\end{moneq}
and we define
\begin{moneq}
\int_{\CA_1^q} \extder \xi_q\cdots\extder\xi_1\extdersep  \psi(x,\xi)
=
\psi_{\{1, \dots, q\}}(x)
\mapob.
\end{moneq}
Before discussing orientation, let me make some observations concerning this definition. 
In the first place, writing the \myquote{super measure} $\extder\xi_q\cdots\extder\xi_1$ on the left of the function to be integrated implies that we see this as a right-$\CA$-linear map defined on functions with values in $\CA^\CC$ (or $\CA$). 
We thus have in particular for homogeneous $\lambda\in \CA$:
\begin{align*}
\int\extder \xi_q\cdots\extder\xi_1\extdersep  \lambda\cdot \psi(x,\xi)
&
=
(-1)^{q\parity\lambda}\lambda \cdot \int\extder \xi_q\cdots\extder\xi_1\extdersep   \psi(x,\xi)
\\
\int\extder \xi_q\cdots\extder\xi_1\extdersep  \bigl(\psi(x,\xi)\cdot\lambda\bigr)
&
=
\biggl(\,\int\extder \xi_q\cdots\extder\xi_1\extdersep   \psi(x,\xi)\biggr)\cdot\lambda
\mapob.
\end{align*}
This is in complete agreement with the observation that Berezin integration is in fact a derivation: the partial derivative with respect to all odd coordinates evaluated at $\xi=0$:
\begin{moneq}
\int_{\CA_1^q} \extder \xi_q\cdots\extder\xi_1\extdersep  \psi(x,\xi)
=
(\partial_{\xi_q}\cdots\partial_{\xi_1}\psi)(x,\xi=0)
\mapob,
\end{moneq}
an operation that is indeed right-linear. 
This also explains the reason to write the $\extder\xi_j$ in the given order, as it is exactly the order of the partial derivatives. 

Coming back to the general definition of the Berezin integral, we note that changing the order of the odd variables by an odd permutation changes the sign of this integral. 
That is why, in the context of integration on \agmfd/, one has to require that it is oriented in the odd directions (the associated Batchelor bundle should be oriented). 
This is exactly the same phenomenon as for integration of differential (top) forms on ordinary manifolds, where one needs an orientation of the tangent bundle. 
And just as one hardly ever provides such an orientation when integrating a differential form on an ordinary manifold, for \gmfd/s one hardly ever provides this orientation of the odd directions. 
Unfortunately for us, we will encounter orientation reversing applications, so we have to take care of that issue. 
But fortunately, we can avoid a formal discussion of orientation (and in particular we can avoid having to invent a notation for it), as we can, in our context, \myquote{hide} the aspect of orientation reversing maps in the \supscal product. 

To understand what happens, we first look at the non-super situation of a unitary representation of a group on a Hilbert space $(\Hilbert, \inprodsym)$. 
Replacing the inner product $\inprodsym$ by any multiple $\lambda\cdot \inprodsym$ will not change the unitary aspect. 
However, we have to require $\lambda>0$, else we will loose the positive-definite property of an inner product. 
Apart from that, there is no restriction on $\lambda$. 
For instance, using the Lebesgue measure on $\RR$ or any positive multiple of it will not change the unitary aspect of a representation on $L^2(\RR)$. 
And even though for compact groups one can \myquote{normalise} the inner product (so that the associated volume of the group becomes $1$), we do not have that luxury for non-compact groups. 
Similarly, replacing a \supscal product $\supinprsym$ by any multiple $\lambda\cdot \supinprsym$ will not change the super-unitary aspect of a representation. 
But contrary to the case of an ordinary inner product, we do not have a positive-definite condition for \supscal products (only a non-degeneracy condition), so we can use any $\lambda\in\RR\setminus\{0\}$ (we cannot use complex values, as these interfere with the symmetry condition of a \supscal product). 

Now suppose we have two ordinary inner product spaces that are equivalent, meaning that we have a unitary isomorphism. 
When we change the inner product on one of these spaces by a multiple $\lambda>0$, the unitary isomorphisme will no longer be unitary, as it maps the inner product on the source to a multiple of the inner product on the target space. 
However, changing the isomorphism by a factor $\sqrt\lambda$ will give us a unitary isomorphism between the new inner products, meaning that we can maintain the definition of equivalence as the existence of a \stress{unitary} isomorphism. 
On the other hand, for a super inner product space we have two structures: the \supscal product and the inner product, both of which can be changed by a factor. 
Changing an isomorphism by a factor will be able to absorb the factor involving the inner product, but it can never absorb both. 
And it cannot absorb a negative factor for the \supscal product, because multiplying an isomorphism $A$ by a factor $\lambda\in \CC\setminus\{0\}$ will change the expression $\superinprod{Ax}{Ay}$ by a factor $\vert\lambda\vert^2>0$. 
The above arguments motivate the definition of a $k$-isomorphism given in \recalt{DefinitionOfSuperInnerProductSpaces}. 
It also solves our problem of orientation reversing maps, as those can be seen as $(-1)$-isomorphisms, changing the sign of the \supscal product. 
In this way we will not have to discuss orientation problems, we only have to indicate how \supscal products are related.  

\end{definition}

\begin{definition}[TheSuperManifoldsRpVertq]{The \gmfd/s $\RR^{p\vert q}$ and associated function spaces}
Even though I strongly dislike this notation, I admit that it is very suggestive and short, so \dots{} 
let us denote by $\RR^{p\vert q}$ the \gmfd/ $\CA_0^p\times \CA_1^q$, (isomorphic to) the even part of a $p\vert q$-dimensional \gvs/, which we occasionally will interpret as an abelian \glgrp/. 
We now consider a (complex) super Hilbert space $(\Hilbert,\inprodsym_\Hilbert, \supinprsym_\Hilbert)$ and a dense subspace $E\subset \Hilbert$ equipped with the restrictions of $\inprodsym_\Hilbert$ and $\supinprsym_\Hilbert$ (a super inner product space \recalt{ExtensionSuperScalarProductToCompletion}). 
Next we consider the space 
\begin{moneq}
F=C^\infty_c(\RR^{p\vert q}; E\otimes \CA^\CC)
\end{moneq}
of compactly supported smooth functions on $\RR^{p\vert q}$ with values in $E \otimes \CA^\CC$. 
Any such function $\psi$ is determined by $2^q$ compactly supported smooth functions $\psi_I:\RR^p\to E$, $I\subset \{1,\dots,q\}$ according to \recalf{StructureOfSmoothFunctionInPowersXiAndFunctionsX}, 
providing a bijection between $F$ and $C^\infty_c(\RR^p; E)^{2^q}$ (see \recalt{structureofsmoothfunctions}). 
For completeness sake we remind the reader that the parity of the function $(x,\xi)\mapsto \repprod\xi^I\cdot \psi_I(x)\in E$ is the parity of the function $\psi_I:\RR^p\to E$ (provided it is homogeneous) plus the cardinal of $I$ modulo $2$. 

On $C^\infty_c(\RR^p;E)$ we define the (ordinary) inner product $\inprodd\cdot\cdot$ by
\begin{moneq}[InnerProductWithLebesgueIntegral]
\inprodd\chi\psi= \int_{\RR^p}\extder\Leb^{(p)}(x)\extdersep  \inprod{\chi(x)}{\psi(x)}_\Hilbert
\mapob,
\end{moneq}
making it a dense subspace of the (ordinary) Hilbert space $\bigl(L^2(\RR^p, \Hilbert), \inprodd\cdot\cdot\bigr)$. 
We next define on $F\times F$ the maps $\inprodsym_{p\vert q}^\Hilbert$ and $\supinprsym_{p\vert q}^\Hilbert$ by
\begin{align}
\inprod\chi\psi_{p\vert q}^\Hilbert
&
=
\sum_{I\subset \{1,\dots,q\}}\int_{\RR^p}\extder\Leb^{(p)}(x)\extdersep  \inprod{\chi_I(x)}{\psi_I(x)}_\Hilbert
\oversettext{\ \recalf{InnerProductWithLebesgueIntegral}\ }\to\equiv
\sum_{I\subset \{1,\dots, q\}} \inprodd{\chi_I}{\psi_I}
\label{ScalarProductInFunctionHilbertSpace}
\\
\homsuperinprod\chi\psi{_{p\vert q}^\Hilbert}
&
=
\int_{\RR^p}\extder \Leb^{(p)}(x)\ \int_{\CA_1^q} \extder \xi_q\dots\extder\xi_2\,\extder \xi_1\extdersep  \homsuperinprod[2]{\chi(x, \xi)}{ \psi(x,\xi)}{_\Hilbert}
\notag
\\&
=
\sum_{I\subset \{1, \dots, q\}} 
\int_{\RR^p} \extder\Leb^{(p)}(x)\extdersep  (-1)^{\varepsilon(I,I^c)}\, \homsuperinprod[2]{\conjugate^{\vert I\vert}\chi_I(x)}{\conjugate^q \psi_{I^c}(x)}{_\Hilbert}
\mapob,
\label{SuperScalarProductInFunctionHilbertSpace}
\end{align}
where, for any two disjoint subsets $I,J\subset \mathcal I$ of a finite set of integers $\mathcal I$, the integer $\varepsilon(I,J)$ denotes the number of transpositions (modulo $2$) needed to put the sequence
\begin{moneq}
i_1, \dots, i_p, j_1, \dots, j_q
\end{moneq}
in increasing order, when we have
\begin{moneq}
\begin{aligned}
I &= \{i_1, \dots, i_p\}
\text{ with }i_1<i_2<\cdots<i_p
\\
J &= \{j_1, \dots, j_q\}
\text{ with }j_1<j_2<\cdots<j_q
\mapob.
\end{aligned}
\end{moneq}
If either $I$ or $J$ is empty, we put $\varepsilon(I,J)=0$. 
In the special case $I\cup J=\mathcal I$, we have $J=\mathcal I\setminus I$ and we note $J=I^c$. 

The construction of the triple $(F,\inprodsym_{p\vert q}^\Hilbert,\supinprsym_{p\vert q}^\Hilbert)$ is a special case of the direct (Berezin) integral of the super inner product space $E$ over the \gmfd/ $M=\RR^{p\vert q}$ equipped with the standard (invariant) super metric as described in \cite[5.19]{Tuynman..2018}, where it is shown that it is a super inner product space, and in particular that $\supinprsym_{p\vert q}^\Hilbert$ is continuous with respect to (the topology defined by) $\inprodsym_{p\vert q}^\Hilbert$. 
When we forget about the \supscal product $\supinprsym_{p\vert q}^\Hilbert$, it is immediate that the completion of $F$ with respect to its inner product $\inprodsym_{p\vert q}^\Hilbert$, a space we will denote as $L^2(\RR^{p\vert q};\Hilbert)$, is (isomorphic to) $2^q$ copies of $L^2(\RR^p;\Hilbert)$, equipped with the inner product $\inprodsym_{p\vert q}^\Hilbert$ \recalf{ScalarProductInFunctionHilbertSpace}. 
However, after performing this completion, we can no longer say that $L^2(\RR^{p\vert q};\Hilbert)$ consists of (smooth?) functions of $q$ odd variables, as that aspect gets lost in the completion procedure. 
On the other hand, officially the aspect of a function on $\RR^p$ also gets lost in the passage from $C^\infty_c(\RR^p;E)$ to $L^2(\RR^p;\Hilbert)$. 
The best one can do (when $E=\Hilbert$, if not, it gets even more complicated) is to say that it will be isomorphic to an equivalence class of functions that are almost everywhere the same. 
Nevertheless, nobody has any problem to \myquote{identify} elements of $L^2(\RR^p;\Hilbert)$ with a function (on $\RR^p$ with values in $\Hilbert$). 

In our super context we will do exactly the same: we will identify any element of $L^2(\RR^{p\vert q}; \Hilbert)$, \ie, a family of functions $\psi_I\in L^2(\RR^p;\Hilbert)$, $I\subset \{1, \dots, q\}$, with a function $\psi:\RR^{p\vert q}\to \Hilbert\otimes \CA^\CC$ according to \recalf{StructureOfSmoothFunctionInPowersXiAndFunctionsX}. 
And under that identification, we can say\slash pretend that the extension of the \supscal product $\supinprsym_{p\vert q}^\Hilbert$ to $L^2(\RR^{p\vert q};\Hilbert)$, the completion of $F$, is still given by Lebesgue-Berezin integration as in \recalf{SuperScalarProductInFunctionHilbertSpace}. 
Standard arguments involving bump functions on $\RR^p$ then allow us to affirm that $\supinprsym^\Hilbert_{p\vert q}$ is non-degenerate on $L^2(\RR^{p\vert q}; \Hilbert)$. 

To conclude: the space $L^2(\RR^{p\vert q};\Hilbert)$, the completion of $F$, together with the inner product $\inprodsym_{p\vert q}^\Hilbert$ \recalf{ScalarProductInFunctionHilbertSpace} and the \supscal product $\supinprsym_{p\vert q}^\Hilbert$ \recalf{SuperScalarProductInFunctionHilbertSpace} is a super Hilbert space, one we will denote as $L^2(\RR^{p\vert q};\Hilbert;\inprodsym_{p\vert q}^\Hilbert; \supinprsym_{p\vert q}^\Hilbert)$. 
And \dots{} as is usual, when no confusion is possible, we will leave out the inner product and\slash or the \supscal product and denote this space simply as $L^2(\RR^{p\vert q};\Hilbert)$. 
In the particular case $\Hilbert=\CC$ with $\inprod xy_\Hilbert=\homsuperinprod xy{_\Hilbert} = \overline x\,y$ we will omit the superscript $^\CC \equiv\, ^\Hilbert$ on the inner product and \supscal product and we will write $\inprodsym_{p\vert q}$ and $\supinprsym_{p\vert q}$. 
It yields the super Hilbert space $L^2(\RR^{p\vert q};\CC)$ of the left-regular representation of $\RR^{p\vert q}$ (seen as an abelian \glgrp/). 
The fact that this is a super unitary representation reflects the fact that Lebesgue integration and Berezin integration are invariant under translations. 
And, as before if no confusion is possible, we will denote this space as $L^2(\RR^{p\vert q})$. 
When we identify an element $\psi$ of $L^2(\RR^{p\vert q})$ with $2^q$ ordinary functions $\psi_I:\RR^p\to \CC$ (see \recalf{StructureOfSmoothFunctionInPowersXiAndFunctionsX}) the inner product $\inprod{\chi}{\psi}_{p\vert q}$ \recalf{ScalarProductInFunctionHilbertSpace} simplifies as 
\begin{align}
\inprod\chi\psi_{p\vert q}
&
=
\sum_{I\subset \{1, \dots, q\}} \inprodd{\chi_I}{\psi_I}
=
\sum_{I\subset \{1, \dots, q\}} \int_{\RR^p}  \extder\Leb^{(p)}(x)\extdersep  \overline{\chi_I(x)}\cdot \psi_I(x)
\label{SimplifiedInnerProductOnL2Rpq}
\\
\noalign{\noindent and the \supscal product $\supinprsym_{p\vert q}$  \recalf{SuperScalarProductInFunctionHilbertSpace} simplifies as}
\homsuperinprod{\chi}{\psi}{_{p\vert q}}
&=
\sum_{I\subset \{1, \dots, q\}}(-1)^{\varepsilon(I,I^c)}\cdot \inprodd{\chi_I}{\psi_{I^c}}
\mapob.
\label{SimplifiedSuperScalarProductOnL2Rpq}
\end{align}

For future reference we also introduce the Schwartz space $\Schwartz pq \subset  L^2(\RR^{p\vert q})$ as follows. 
Just as any element of $L^2(\RR^{p\vert q})$ \myquote{is} a function on $\RR^{p\vert q}$ of the form
\begin{moneq}
\psi(x,\xi) = \sum_{I\subset \{1, \dots, q\}} \repprod\xi^I \cdot \psi_I(x)
\end{moneq}
with $\psi_I\in L^2(\RR^p)$, in the same way a function $\psi$ on $\RR^{p\vert q}$ belongs to $\Schwartz pq$ if all \myquote{components} $\psi_I$ belong to the Schwartz space $\Schwartz p{} \subset L^2(\RR^p)$ of rapidly decreasing smooth functions on $\RR^p$ defined as
\begin{multline*}
\quad
\Schwartz p{}
=
\bigl\{\, \psi\in C^\infty(\RR^p)\,\bigm\vert\,
\forall i,j\in\NN^p
\quad:\quad
\\
x_1^{i_1}\cdots x_p^{i_p}\cdot (\partial_1^{j_1}\dots\partial_p^{j_p}\psi)(x)
\in L^2(\RR^p)
\,\bigr\}
\mapob.
\quad
\end{multline*}
We thus have the inclusions $C^\infty_c(\RR^{p\vert q}) \subset \Schwartz pq \subset C^\infty(\RR^{p\vert q})\cap L^2(\RR^{p\vert q})$ as dense graded subspaces.

\end{definition}

\begin{remarque}{Remark}
It should be noted that the definition of the \supscal product $\supinprsym_{p\vert q}$ is more \myquote{intrinsic} than the definition of the inner product  $\inprodsym_{p\vert q}$ on $L^2(\RR^{p\vert q})$. 
The definition of $\supinprsym_{p\vert q}$ uses (Berezin) integration over $\RR^{p\vert q}$ and as such behaves in a \myquote{nice} way with respect to a change of coordinates on $\RR^{p\vert q}$. 
On the other hand, $\inprodsym_{p\vert q}$ depends upon the identification $C^\infty(\RR^{p\vert q})\cong \bigl(C^\infty(\RR^p)\bigr)^q$, an identification which does not behave nicely with respect to a change of coordinates (see also \recalt{ExampleWithDifferentTopologiesInducedByChangeOfCoordinates}). 
\end{remarque}

\begin{proclaim}[OddPermutationsChangeSignOfSuperScalarProductsNEW]{Lemma}
Let $\sigma$ be a permutation of $\{1, \dots, q\}$ and let $\varepsilon(\sigma)$ denote its parity. 
Then the map $\Identif_\sigma:L^2(\RR^{p\vert q};\Hilbert)
\to  
L^2(\RR^{p\vert q};\Hilbert)$ defined as
\begin{moneq}
\bigl(\Identif_\sigma(\psi)\bigr)(x_1, \dots, x_p,\xi_1, \dots, \xi_q) 
= \psi(x_1, \dots, x_p, \xi_{\sigma(1)}, \dots, \xi_{\sigma(q)})
\end{moneq}
is an even $\varepsilon(\sigma)$-isomorphism of super Hilbert spaces. 

\end{proclaim}

\begin{proclaim}[DirectIntegralOfHilbertSpacesPotentiallyChangeSignSuperScalarProductsNEW]{Lemma}
The map $\Identif:L^2(\RR^{p+r\vert q+s};\Hilbert)
\to  
L^2\bigl(\RR^{p\vert q}; L^2(\RR^{r\vert s};\Hilbert)\bigr)$ defined by
\begin{multline*}
\bigl(\Identif(\psi)\bigr)(x_1, \dots, x_p,\xi_1, \dots, \xi_q) : (y_1, \dots, y_r, \eta_1,\dots, \eta_s) \\
\mapsto \psi(x_1, \dots, x_p,y_1, \dots, y_r, \xi_1, \dots, \xi_q,\eta_1, \dots, \eta_s)
\end{multline*}
is an even $(-1)^{qs}$-isomorphism of super Hilbert spaces. 

\end{proclaim}

\begin{preuve}
Forgetting about the \supscal product and concentrating on the standard inner product, we have isomorphisms\slash identifications
\begin{moneq}[IdentificationL2SpacesOnProductInProofEqReedSimon]
L^2(\RR^{p+r\vert q+s};\Hilbert) \cong L^2(\RR^p; \Hilbert)^{2^{q+s}}
\end{moneq}
and
\begin{moneq}[IdentificationRepeatedL2SpacesInProofEqReedSimon]
L^2\bigl(\RR^{p\vert q}; L^2(\RR^{r\vert s};\Hilbert)\bigr)
\cong
L^2\bigl(\RR^p; L^2(\RR^r;\Hilbert)^{2^s}\bigr)^{2^q}
\cong
L^2\bigl(\RR^p; L^2(\RR^q; \Hilbert)\bigr)^{2^{q+s}}
\mapob.
\end{moneq}
More precisely, according to \recalf{StructureOfSmoothFunctionInPowersXiAndFunctionsX} slightly rewritten, any element $\psi\in L^2(\RR^{p+r\vert q+s};\Hilbert)$ is determined by $2^{q+s}$ elements $\psi_{I,J}\in L^2(\RR^{p+r}; \Hilbert)$ (with $I\subset \{1, \dots, q\}$ and $J\subset\{1, \dots, s\}$) via the formula 
\begin{moneq}
\psi(x,y,\xi_1, \dots, \xi_{q},\eta_1, \dots, \eta_s)
=
\sum_{I\subset \{1, \dots, q\}} \sum_{J\subset\{1, \dots, s\}} \repprod\xi^I\,\repprod\eta^J\,\psi_{I,J}(x,y)
\mapob.
\end{moneq}
Similarly, any element $\chi\in L^2\bigl(\RR^{p\vert q}; L^2(\RR^{r\vert s};\Hilbert)\bigr)$ is determined by $2^q$ elements $\chi_I\in L^2\bigl(\RR^p; L^2(\RR^{r\vert s};\Hilbert)\bigr)$ (with $I\subset \{1, \dots, q\}$) via the formula
\begin{moneq}
\chi(x,\xi) = \sum_{I\subset \{1, \dots, q\}} \repprod\xi^I\,\chi_I(x)
\mapob.
\end{moneq}
and any element $\chi_I(x)\in L^2(\RR^{r\vert s}; \Hilbert)$ is determined by $2^s$ elements $\chi_{I}^J(x)\in L^2(\RR^r; \Hilbert)$ via the formula
\begin{moneq}
\bigl(\chi_I(x)\bigr)(y,\eta) = \sum_{J\subset \{1, \dots, s\}} \repprod\eta^J \bigl(\chi_{I}^J(x)\bigr)(y)
\mapob.
\end{moneq}
According to the definition of $\Identif$ it follows immediately that we have
\begin{moneq}[FunctionsOfFunctionsDecompositionForFunctionHilbertSpacesNEW]
\Identif(\psi)_{I}(x):(y,\eta)
\mapsto
\sum_{J\subset \{1, \dots, s\}} \repprod\eta^J\,\psi_{I,J}(x,y)
\mapob,
\end{moneq}
\ie, $\bigl(\Identif(\psi)_{I}^J(x)\bigr)(y)=\psi_{I,J}(x,y)$. 
We now recall (see \cite[Thm II.10]{ReedSimon1..1980}) that the map $\mathcal{J}:L^2(\RR^{p+r}; \Hilbert)\to L^2\bigl(\RR^p; L^2(\RR^r; \Hilbert)\bigr)$ defined by
\begin{moneq}
\forall \psi\in L^2(\RR^{p+r}; \Hilbert)\ \forall x\in \RR^p
\ :\ 
(\mathcal{J}\psi)(x):y\mapsto \psi(x,y)
\end{moneq}
is an isomorphism of Hilbert spaces.
Applied to our case it tells us that we have the equality $\Identif(\psi)_I^J = \mathcal{J}(\psi_{I,J})$. 
In other words, in terms of the identifications \recalf{IdentificationL2SpacesOnProductInProofEqReedSimon} and \recalf{IdentificationRepeatedL2SpacesInProofEqReedSimon} the map $\Identif$ is given by
\begin{moneq}
\Identif\hBigl((\psi_{I,J})_{I\subset\{1, \dots, q\}, J\subset\{1, \dots, s\}}\hBigr)
=
\hBigl(\mathcal{J}(\psi_{I,J})\hBigr)_{I\subset\{1, \dots, q\}, J\subset\{1, \dots, s\}}
\mapob.
\end{moneq}
As $\mathcal{J}$ is an isomorphism, $\Identif$ is too. 

To prove that $\Identif$ is even, suppose $\psi$ has parity $\varepsilon$. 
It then follows that the separate functions $\psi_{I,J}$ have parity $(-1)^{\parity I + \parity J}\cdot \varepsilon$. 
As the map $\mathcal{J}$ does not change parity, the maps $\mathcal{J}(\psi_{I,J})$ also have parity $(-1)^{\parity I + \parity J}\cdot \varepsilon$. 
For fixed $x\in \RR^p$ the map $(y,\eta)\mapsto \sum_{J}\repprod\eta^J\,\bigl(\mathcal{J}(\psi_{I,J})(x)\bigr)(y)$ then has parity $(-1)^{\parity I}\cdot \varepsilon$, and thus the map $\Identif(\psi)$ will have parity $\varepsilon$, showing that $\Identif$ is even as claimed. 

Remains to show that $\Identif$ preserves the inner products and the \supscal products in the sense
\begin{moneq}
\inprod{\Identif(\chi)}{\Identif(\psi)}_{p\vert q}
=
\inprod \chi \psi_{p+r\vert q+s}
\quad\text{and}\quad
(-1)^{qs}\,\homsuperinprod[2]{\Identif(\chi)}{\Identif(\psi)}{_{p\vert q}}
=
\homsuperinprod \chi \psi{_{p+r\vert q+s}}
\mapob.
\end{moneq}
For the (ordinary) inner product this is already implicit in the previous argument, but we can compute directly:
\begin{align*}
\inprod{\Identif(\chi)}{\Identif(\psi)}_{p\vert q}
\kern1em
&
\oversetalign{\recalf{ScalarProductInFunctionHilbertSpace}}\to=
\kern0.5em
\sum_{I\subset \{1, \dots, q,\}} \inprod{\Identif(\chi)_I}{\Identif(\psi)_I}_{r\vert s}
\\&
=
\kern0.5em
\sum_{I\subset \{1, \dots, q,\}} \ \sum_{J\subset \{1, \dots, s\}} \inprod{\Identif(\chi)_I^J}{\Identif(\psi)_I^J}_{\Hilbert}
\\&
\oversetalign{\recalf{FunctionsOfFunctionsDecompositionForFunctionHilbertSpacesNEW}}\to=
\kern0.5em
\sum_{I\subset \{1, \dots, q,\}} \ \sum_{J\subset \{1, \dots, s\}} \inprod{\chi_{I,J}}{\psi_{I,J}}_{\Hilbert}
\oversettext{\recalf{ScalarProductInFunctionHilbertSpace}}\to=
\inprod \chi \psi_{p+r\vert q+s}
\mapob.
\end{align*}
And for the \supscal product we compute:
\begin{align*}
\homsuperinprod[2]{\Identif(\chi)}{\Identif(\psi)}{_{p\vert q}}
&
=
\int_{\RR^p}\extder\Leb^{(p)}(x)\kern0.5em
\int_{\CA_1^q} \extder\xi_q\cdots\extder\xi_1\extdersep 
\homsuperinprod[2]{\Identif(\chi)(x,\xi)}{\Identif(\psi)(x,\xi)}{_{r\vert s}}
\\&
=
\int_{\RR^p}\extder\Leb^{(p)}(x)\extdersep 
\int_{\CA_1^q} \extder\xi_q\cdots\extder\xi_1\extdersep 
\int_{\RR^r}\extder\Leb^{(r)}(y)\kern0.5em
\\&
\kern6em
\int_{\CA_1^s}\extder\eta_s\cdots\extder\eta_1\extdersep 
\homsuperinprod[2]{\chi(x,y,\eta,\xi)}{\psi(x,y,\eta,\xi)}{_\Hilbert}
\\&
=
\int_{\RR^{p+r}}\extder\Leb^{(p+r)}(x,y)\kern0.5em
\int_{\CA_1^{q+s}} \extder\xi_q\cdots\extder\xi_1\,\extder\eta_s\cdots\extder\eta_1\extdersep 
\\&
\kern6em
\homsuperinprod[2]{\chi(x,y,\eta,\xi)}{\psi(x,y,\eta,\xi)}{_\Hilbert}
\\&
=
(-1)^{qs}
\int_{\RR^{p+r}}\extder\Leb^{(p+r)}(x,y)\kern0.5em
\int_{\CA_1^{q+s}} \extder\eta_s\cdots\extder\eta_1\,\extder\xi_q\cdots\extder\xi_1\extdersep 
\\&
\kern6em
\homsuperinprod[2]{\chi(x,y,\eta,\xi)}{\psi(x,y,\eta,\xi)}{_\Hilbert}
\\&
=
(-1)^{qs}
\homsuperinprod \chi \psi{_{p+r\vert q+s}}
\mapob.
\end{align*}

\vskip-1.5\baselineskip
\end{preuve}

\begin{proclaim}[PartialBerezeinFourierTransformIsIsomorphismSHS]{Proposition}
For $r\le p$ and $s\le q$ the \myquote{partial} Berezin-Fourier transform $\Fourierodd_{r\vert s}:L^2(\RR^{p\vert q})\to L^2(\RR^{p\vert q})$ defined by
\begin{multline*}
\quad
(\Fourierodd_{r\vert s}\psi)(k_1, \dots, k_r,x_{r+1}, \dots,x_p, \lambda_1, \dots, \lambda_s, \xi_{s+1}, \dots, \xi_q)
=
\\[2\jot]
\frac1{(2\pi)^{r/2}}\int_{\RR^r}\extder\Leb^{r}(x)\kern0.5em \int_{\CA_1^s}\extder \xi_s\cdots\extder\xi_1\extdersep 
\eexp^{-i\sum_{j=1}^rk_jx_j - i\sum_{j=1}^s\xi_j\lambda_j}\,\psi(x,\xi)
\quad
\end{multline*}
is an isomorphism of super Hilbert spaces of the same parity as $s$ with inverse given by
\begin{multline*}
\quad
(\Fourierodd_{r\vert s}^{\,\,-1} \psi)(x,\xi)
=
\frac{i^{s^2}}{(2\pi)^{r/2}}\cdot \int_{\RR^r}\extder\Leb^{r}(k)\kern0.5em \int_{\CA_1^s}\extder \lambda_s\cdots\extder\lambda_1\extdersep 
\\[2\jot]
\eexp^{i\sum_{j=1}^rk_jx_j + i\sum_{j=1}^s\xi_j\lambda_j}\,\psi(k_1, \dots, k_r,x_{r+1}, \dots,x_p, \lambda_1, \dots, \lambda_s, \xi_{s+1}, \dots, \xi_q)
\mapob.
\quad
\end{multline*}
It maps $\Schwartz pq \subset L^2(\RR^{p\vert q})$ to itself. 

\end{proclaim}

\begin{preuve}
We first recall that the partial Fourier transform $\Fourier_r:L^2(\RR^p)\to L^2(\RR^p)$ defined as
\begin{moneq}
(\Fourier_r\psi)(k_1, \dots, k_r, x_{r+1}, \dots, x_p)
=
(2\pi)^{-r/2}\cdot \int_{\RR^r} \extder\Leb^{(r)}(x)\kern0.5em \eexp^{-i\sum_{j=1}^rk_jx_j
}\,\psi(x_1, \dots, x_p)
\end{moneq}
is a unitary isomorphism. 
Next we recall that any function $\psi\in L^2(\RR^{p\vert q})$ is determined by $2^q$ functions $\psi_I:\RR^p\to \CC$ according to \recalt{StructureOfSmoothFunctionInPowersXiAndFunctionsX}, a formula that we can rewrite, separating the first $s$ odd coordinates from the remaining ones, as
\begin{moneq}
\psi(x,\xi) = \sum_{J'\subset \{1, \dots, s\}}\sum_{J\subset \{s+1, \dots, q\}} \repprod\xi^{J'}\,\repprod\xi^{J}\, \psi_{J',J}(x)
\mapob.
\end{moneq}
For such a function we thus have
\begin{align*}
\shifttag{4em}
(\Fourierodd_{r\vert s}\psi)(k_1, \dots, k_r,x_{k+1}, \dots,x_p, \lambda_1, \dots, \lambda_s, \xi_{s+1}, \dots, \xi_q)
\\[3\jot]
&
=
\frac1{(2\pi)^{r/2}}\int_{\RR^r}\extder\Leb^{r}(x)\kern0.5em \int_{\CA_1^s}\extder \xi_s\cdots\extder\xi_1\extdersep  \sum_{J'\subset \{1, \dots, s\}}\sum_{J\subset \{s+1, \dots, q\}}
\\[1\jot]&
\kern9em
\eexp^{- i\sum_{j=1}^s\xi_j\lambda_j}\,\repprod\xi^{J'}\,\repprod\xi^{J}\,\eexp^{-i\sum_{j=1}^rk_jx_j}\,\psi_{J',J}(x)
\\[3\jot]&
=
\int_{\CA_1^s}\extder \xi_s\cdots\extder\xi_1\kern0.5em \sum_{I\subset \{1, \dots, s\}} \sum_{J'\subset \{1, \dots, s\}}\sum_{J\subset \{s+1, \dots, q\}} 
\\[1\jot]&
\kern9em
\repprod{-i\xi\lambda}^I\,\repprod\xi^{J'}\,\repprod\xi^{J}\,(\Fourier_r \psi_{J',J})(k,x)
\\[3\jot]&
=
\int_{\CA_1^s}\extder \xi_s\cdots\extder\xi_1\kern0.5em \sum_{I\subset \{1, \dots, s\}} \sum_{J'\subset \{1, \dots, s\}}\sum_{J\subset \{s+1, \dots, q\}} 
\\[1\jot]&
\kern5em
i^{\parity I}\,(-1)^{\parity I(\parity I-1)/2}\,\repprod{\lambda}^I\,\repprod{\xi}^I\,\repprod\xi^{J'}\,\repprod\xi^{J}\,(\Fourier_r \psi_{J',J})(k,x)
\\[3\jot]&
=
\sum_{I\subset \{1, \dots, s\}}\sum_{J\subset \{s+1, \dots, q\}}\ i^{\parity I}\,(-1)^{s\parity I+\parity I(\parity I-1)/2}\,\repprod{\lambda}^I\cdot
\\[1\jot]&
\kern5em
\int_{\CA_1^s}\extder \xi_s\cdots\extder\xi_1\extdersep 
\repprod{\xi}^I\,\repprod\xi^{I^c}\,\repprod\xi^{J}\,(\Fourier_r \psi_{I^c,J})(k,x)
\\[3\jot]&
=
\sum_{I\subset \{1, \dots, s\}}\sum_{J\subset \{s+1, \dots, q\}}\ i^{\parity I}\,(-1)^{s\parity I+\parity I(\parity I-1)/2+\varepsilon(I,I^c)}\cdot
\\&
\kern15em
\repprod{\lambda}^I\,\repprod\xi^{J}\,(\Fourier_r \psi_{I^c,J})(k,x)
\mapob.
\end{align*}
In terms of the decomposition of $\Fourierodd_{r\vert s}\psi$ with respect to its odd variables
\begin{moneq}
(\Fourierodd_{r\vert s}\psi)(k,x, \lambda, \xi)
=
\sum_{I\subset\{1, \dots, s\}}\sum_{J\subset \{s+1, \dots, q\}} 
\repprod\lambda^I\,\repprod\xi^J\,
(\Fourierodd_{r\vert s}\psi)_{I,J}(k,x)
\end{moneq}
we thus have
\begin{moneq}[PartialBerezinFourierInComponents]
(\Fourierodd_{r\vert s}\psi)_{I,J}
=
i^{\parity I}\,(-1)^{s\parity I+\parity I(\parity I-1)/2+\varepsilon(I,I^c)}\cdot
(\Fourier_r \psi_{I^c,J})
\mapob.
\end{moneq}
For the formula given for $\Fourierodd_{r\vert s}^{\,\,-1}$ we find in a similar way the result
\dontprintcomp{
\begin{align*}
(\Fourierodd_{r\vert s}^{\,\,-1} \psi)(x,\xi)
&
=
\frac{i^{s^2}}{(2\pi)^{r/2}}\cdot \int_{\RR^r}\extder\Leb^{r}(k)\kern0.5em \int_{\CA_1^s}\extder \lambda_s\cdots\extder\lambda_1\extdersep 
\\&
\kern3em
\eexp^{i\sum_{j=1}^rk_jx_j + i\sum_{j=1}^s\xi_j\lambda_j}\,\psi(k,x, \lambda, \xi)
\\&
=
\frac{i^{s^2}}{(2\pi)^{r/2}}\cdot \int_{\RR^r}\extder\Leb^{r}(k)\kern0.5em \int_{\CA_1^s}\extder \lambda_s\cdots\extder\lambda_1\kern0.5em \sum_{J'\subset \{1, \dots, s\}} \sum_{J\subset \{s+1, \dots, q\}} \sum_{I\subset \{1, \dots, s\}}
\\&
\kern3em
\repprod{i\xi\lambda}^I\,\repprod\lambda^{J'}\,\repprod\xi^J\,
\eexp^{i\sum_{j=1}^rk_jx_j}\,\psi_{J',J}(k,x)
\\&
=
i^{s^2}\cdot \int_{\CA_1^s}\extder \lambda_s\cdots\extder\lambda_1\kern0.5em \sum_{J\subset \{s+1, \dots, q\}} \sum_{I\subset \{1, \dots, s\}}
\\&
\kern3em
\repprod\lambda^{I^c}\,\repprod{i\xi\lambda}^I\,\repprod\xi^J\,
(\Fourier_r^{\,\,-1}\psi_{I^c,J})(x)
\\&
=
i^{s^2}\cdot  \sum_{J\subset \{s+1, \dots, q\}} \sum_{I\subset \{1, \dots, s\}}
\\&
\kern3em
i^{\parity I}\,(-1)^{\parity I(\parity I+1)/2+\varepsilon(I^c,I)}
\repprod{\xi}^I\,\repprod\xi^J\,
(\Fourier_r^{\,\,-1}\psi_{I^c,J})(x)
\mapob.
\end{align*}
And thus:
}
\begin{moneq}
(\Fourierodd_{r\vert s}^{\,\,-1} \psi)_{I,J}
=
i^{s^2} \, i^{\parity I}\,(-1)^{\parity I(\parity I+1)/2+\varepsilon(I^c,I)}\,
(\Fourier_r^{\,\,-1}\psi_{I^c,J})
\mapob.
\end{moneq}
Combining this with \recalf{PartialBerezinFourierInComponents} we find
\begin{moneq}
\bigl(\Fourierodd_{r\vert s}^{\,\,-1}(\Fourierodd_{r\vert s} \psi)\bigr)_{I,J}
=
\psi_{I,J}
\mapob,
\end{moneq}
proving that these two operations are indeed each other's inverse. 

As the partial Fourier transform $\Fourier_r$ is a unitary isomorphism of $L^2(\RR^p)$, it follows immediately that we have
\begin{moneq}
\inprod[2]{\Fourierodd_{r\vert s}\chi}{\Fourierodd_{r\vert s}\psi}_{p\vert q}
=
\inprod \chi \psi_{p\vert q}
\mapob.
\end{moneq}
The parity claim is an immediate consequence of the fact that we perform a Berezin integration over $s$ odd variables, but it can also be deduced from \recalf{PartialBerezinFourierInComponents}: the function $\repprod\xi^{I^c}\cdot\repprod\xi^J\cdot \psi_{I^c,J}$ of parity $\parity{I^c}+\parity J$ is mapped (up to a constant) to $\repprod\lambda^I\cdot\repprod\xi^J\cdot \Fourier_r\psi_{I^c,J}$ of parity $\parity{I}+\parity J$. 
The parity change thus is $(-1)^{\parity{I^c}+\parity I} = (-1)^s$ as claimed. 

It thus remains to show the equivalence property with respect to the \supscal product, for which we have to be careful. 
As the Berezin-Fourier transform is homogeneous of parity $(-1)^s$, we can summarise condition \refmetnaam{\labelIsomSHS}{conditionoddisoforequivalenceofSHS} for an equivalence as the condition
\begin{moneq}
\homsuperinprod{\Fourierodd_{r\vert s}f}{\Fourierodd_{r\vert s}\psi}{_{p\vert q}}
=
i^{s^2}\,\homsuperinprod{\conjugate^s f}\psi{_{p\vert q}}
\mapob.
\end{moneq}
We thus compute:
\begin{align*}
\shifttag{4em}
\homsuperinprod{\Fourierodd_{r\vert s}\chi}{\Fourierodd_{r\vert s}\psi}{_{p\vert q}}
\oversettext{\ \recalf{SuperScalarProductInFunctionHilbertSpace}\ }\to=
\sum_{I\subset\{1,\dots, s\}}\sum_{J\subset \{s+1, \dots, q\}} 
\\&
\kern3em
(-1)^{\varepsilon(I\cup J,(I\cup J)^c)}\, \inprodd{\conjugate^{\parity I+\parity J} (\Fourierodd_{r\vert s}\chi)_{I,J}}{\conjugate^q (\Fourierodd_{r\vert s}\psi)_{I^c,J^c}}
\\[3\jot]
&
=
\sum_{I\subset\{1,\dots, s\}}\sum_{J\subset \{s+1, \dots, q\}} 
(-1)^{\parity{I^c}\parity J+\varepsilon(I,I^c)+\varepsilon(J,J^c)}\, \inprodd{ (\Fourierodd_{r\vert s}\chi)_{I,J}}{(\Fourierodd_{r\vert s}\psi)_{I^c,J^c}}
\\[3\jot]
&
\oversetalign{\recalf{PartialBerezinFourierInComponents}}\to=
\kern0.5em
\sum_{I\subset\{1,\dots, s\}}\sum_{J\subset \{s+1, \dots, q\}} 
(-1)^{\parity{I^c}\parity J+\varepsilon(I,I^c)+\varepsilon(J,J^c)}\, 
\\[1\jot]&\kern6em
(-i)^{\parity I}\,(-1)^{s\parity I+\parity I(\parity I-1)/2+\varepsilon(I,I^c)}\,
\\[1\jot]&\kern9em
i^{\parity{I^c}}\,(-1)^{s\parity{I^c}+\parity{I^c}(\parity{I^c}-1)/2+\varepsilon(I^c,I)}\,
\inprodd{ \Fourier_r\chi_{I^c,J}}{\Fourier_r\psi_{I,J^c}}
\\[3\jot]&
=
i^{3s^2}\,\sum_{I\subset\{1,\dots, s\}}\sum_{J\subset \{s+1, \dots, q\}} 
(-1)^{\parity{I^c}\parity J-s\parity I+\varepsilon(I^c,I)+\varepsilon(J,J^c)}\, 
\inprodd{ \chi_{I^c,J}}{\psi_{I,J^c}}
\\[3\jot]&
=
i^{3s^2}\,\sum_{I\subset\{1,\dots, s\}}\sum_{J\subset \{s+1, \dots, q\}} 
(-1)^{\parity{I}\parity J-s\parity{I^c}+\varepsilon(I,I^c)+\varepsilon(J,J^c)}\, 
\inprodd{ \chi_{I,J}}{\psi_{I^c,J^c}}
\\[3\jot]&
=
i^{s^2}\sum_{I\subset\{1,\dots, s\}}\sum_{J\subset \{s+1, \dots, q\}} 
(-1)^{\parity{I^c}\parity J+\varepsilon(I,I^c)+\varepsilon(J,J^c)+s(\parity I+\parity J)}\, 
\,\inprodd{ \chi_{I,J}}{\psi_{I^c,J^c}}
\\&
=
i^{s^2}\,\homsuperinprod{\conjugate^{s}\chi}{\psi}{_{p\vert q}}
\mapob.
\end{align*}
We thus may conclude that $\Fourierodd_{r\vert s}$ is indeed an equivalence of super Hilbert spaces as claimed. 
As $\Fourier_r:L^2(\RR^p)\to L^2(\RR^p)$ maps $\Schwartz p{}$ to itself, it is immediate that $\Fourierodd_{r\vert s}$ maps $\Schwartz pq$ to itself. 
\end{preuve}

\section{\label{SECTSuperUnitaryRepresentations}Super unitary representations}

Once we have decided what a super inner product space should be, comes the question how a super unitary representation should behave with respect to the inner product $\inprodsym$. 
A natural idea that comes to mind is that we should stay as close as possible to the ungraded notion of a unitary representation, meaning that we should require that the action not only preserves the \supscal product (we have \myquote{superized} the whole notion of unitary representations), but also the additional structure of an inner product. 
However, this is far too restrictive, as already for the left-regular representation of the simplest super group, the $0\vert1$-dimensional odd line $G=\CA_1$, this cannot be fulfilled, as this is the representation on $\GL(1\vert1)$ given by $\rho(\xi) = \bigl(\begin{smallmatrix} 1 & \xi \\ 0 & 1 \end{smallmatrix}\bigr)$. 
And there is no (positive definite) inner product on $\GL(1\vert1)$ that is preserved by this action. 
The next best requirement is that at least the action of the underlying ungraded (ordinary) Lie group $\body G$ preserves the inner product (making it an honest unitary representation of $\body G$). 
This is the usual approach and the one I adopted in \cite{Tuynman..2018}. 

Unfortunately, even if an ordinary unitary representation extends in a natural way to a completion, the same is not true for super unitary representations: in general there will (to my knowledge) be no way to extend the (infinitesimal versions of) the odd generators to the completion. 
We thus are lead to introduce a dense subspace on which the odd part of the \glgrp/ acts, either in the way of definition \recalt{DefinitionsSuperScalarProductSuperUnitaryRep} below, or in terms of the infinitesimal version linked to the Harish-Chandra pair \recalt{equivalentDefSuperUnitaryRepNEW}.

\begin{definition}[DefinitionsSuperScalarProductSuperUnitaryRep]{Definitions}
\itemizedef[*]
A \stresd{\psur{}} of \aglgrp/ $G$ on a (complex) super Hilbert space $(\Hilbert, \inprodsym, \supinprsym)$ is a couple $(\Dense, \rho)$ with the following properties: 
$\Dense\subset \Hilbert$ is a dense graded subspace and $\rho:G\to \Aut(\Dense\otimes\CA^\CC)$ a group homomorphism satisfying the following conditions.
\begin{enumerate}[\kern0.5em\labelMainSUR1.]
\item\label{MainSURlabelrhopreservessupinpr}
For all $g\in  G$, $\rho(g)$ preserves $\supinprsym$ (or more accurately: preserves the restriction to $\Dense \otimes \CA^\CC$ of the extension of $\supinprsym$ to $\Hilbert \otimes \CA^\CC$), \ie, 
\begin{moneq}
\forall g\in G\ 
\forall \chi,\psi\in \Dense: 
\superinprod{\rho(g)\chi}{\rho(g)\psi} = \superinprod\chi\psi
\mapob.
\end{moneq}

\item\label{MainSURlabelrhopsiissmooth}
For all $\psi\in \Dense$ the map $\FgroupHilbert_\psi:G\to \Dense \otimes \CA^\CC$, $\FgroupHilbert_\psi(g) = \rho(g)\psi$ is smooth.

\item\label{MainSURlabelrhoonbodyisunitary}
For all $g\in \body G$ the map $\rho(g)$ preserves $\inprodsym$.

\end{enumerate}

\itemizedef
If $(\Dense,\rho)$ and $(\Dense',\rho')$ are two \psur{}s of $G$ on the super Hilbert space $(\Hilbert,\inprodsym, \supinprsym)$, then $(\Dense',\rho')$ \stresd{extends} $(\Dense,\rho)$ if $\Dense\subset \Dense'$ and for all $g\in G$ we have $\rho'(g)\restricted_{\Dense} = \rho(g)$. We will say that a \psur{} $(\Dense, \rho)$ is a \stresd{super unitary representation of $G$ on $\Hilbert$} if it does not admit a non-trivial extension. 

\itemizedef
A \psur{} $(\Dense,\rho)$ is said to be \stresd{irreducible} is there exist no non-trivial closed graded subspaces of $\Dense$ that are invariant, \ie, if $F\subset \Dense$ is a closed graded subspace such that $\rho(g)(F)\subset F$ for all $g\in G$, then either $F=\{0\}$ or $F=\Dense$. 

\end{definition}

\begin{proclaim}[ProtoSuperUnitaryRepExtendsToUniqueSuperUnitaryRep]{Lemma}
Let $(\Dense, \inprodsym_\Dense, \supinprsym_\Dense)$ be a super inner product space, let $(\Hilbert, \inprodsym, \supinprsym)$ be its completion (which we suppose to be a super Hilbert space, see \recalt{ExtensionSuperScalarProductToCompletion}), let $G$ be \aglgrp/ and let $\rho:G\to \Aut(\Dense\otimes \CA^\CC)$ be a group homomorphism satisfying the conditions \refmetnaam{\labelMainSUR}{MainSURlabelrhopreservessupinpr}--\ref{MainSURlabelrhoonbodyisunitary}, \ie, $(\Dense, \rho)$ is a \psur. 
Then there exists a unique super unitary representation $(\Dense',\rho')$ of $\,G$ on $(\Hilbert, \inprodsym, \supinprsym)$ that extends $(\Dense, \rho)$, \ie, any \psur{} determines a unique super unitary representation. 

\end{proclaim}

\begin{preuve}
The notion of an extension is such that it is immediate that we can apply Zorn's lemma to prove the existence of a maximal extension, which will be a super unitary representation. 
In order to prove uniqueness, we take two extensions $(\Dense_1, \rho_1)$ and $(\Dense_2, \rho_2)$ of $(\Dense, \rho)$ and we prove that they coincide on $\Dense_1\cap\Dense_2$. 
It then follows that there exists an extension defined on $\Dense_1+\Dense_2$. 
Applying this to two maximal extensions shows that they must coincide, proving uniqueness of the maximal extension. 

So let $(\Dense_1, \rho_1)$ and $(\Dense_2, \rho_2)$ be two extensions of $(\Dense, \rho)$, let $\psi\in \Dense_1\cap \Dense_2$, let $g\in G$ and define $\chi=\rho_1(g)\psi - \rho_2(g)\psi\in \Hilbert\otimes\CA^\CC$.
It follows that there exist $\chi_1, \dots, \chi_n\in\Hilbert$ and $\lambda_1, \dots, \lambda_n\in \CA^\CC\setminus\{0\}$ such that 
\begin{moneq}
\chi = \sum_{j=1}^n \chi_j\otimes \lambda_j
\mapob.
\end{moneq}
When we assume that $n$ is minimal, it follows that $\lambda_1, \dots, \lambda_n$ are independent over $\CC$. 

We now start our computations with the observation that $\rho_j$ is a group homomorphism with values in $\Aut(\Dense_j\otimes \CA^\CC)$ preserving $\supinprsym$ and that they extend $\rho$, hence for any $\phi\in \Dense$ we have:
\begin{moneq}
\superinprod\phi{\rho_1(g)\psi}
=
\superinprod{\rho_1(g\mo)\phi}{\psi}
=
\superinprod{\rho_2(g\mo)\phi}{\psi}
=
\superinprod\phi{\rho_2(g)\psi}
\mapob.
\end{moneq}
It follows that we have the property
\begin{moneq}
\forall \phi\in \Dense: 0 = \superinprod\phi\chi=\sum_j \superinprod\phi{\chi_j}\cdot \lambda_j
\mapob.
\end{moneq}
As the $\lambda_j$ are independent over $\CC$, it follows that we have
\begin{moneq}
\forall \phi\in \Dense\ \forall j: \superinprod\phi{\chi_j}=0
\mapob.
\end{moneq}
By continuity of $\supinprsym$ (with respect to the topology defined by $\inprodsym$) and the fact that $\Dense$ is dense in $\Hilbert$ it follows that we have
\begin{moneq}
\forall j\ \forall \phi\in \Hilbert: \superinprod\phi{\chi_j}=0
\mapob.
\end{moneq}
By non-degeneracy of $\supinprsym$ it then follows that $\chi_j=0$ for all $j$ and thus $\chi=0$, proving $\rho_1(g)\psi=\rho_2(g)\psi$ as wanted. 
\end{preuve}

\begin{remarque}{Remarks}
\itemizedef[*]
In \cite[Thm 7.22]{Tuynman..2018} I showed that the (maximal) dense domain $\Dense$ of a super unitary representation of $G$ on a super Hilbert space $\Hilbert$ is the super analogue of the space of smooth vectors of an ordinary unitary representation. 

\itemizedef
In \cite[Thm 3.6]{Tuynman..2018} repeated below (see also \cite[Thm 4.18]{deGoursacMichel..2015}) I showed that the definition of a super unitary representation as given in \recalt{DefinitionsSuperScalarProductSuperUnitaryRep} is equivalent to a caracterisation closer to the commonly accepted definition in terms of Harish-Chandra pairs for \aglgrp/ as can be found for instance in \cite{CCTV..2006}.
\end{remarque}

\begin{proclaim}[equivalentDefSuperUnitaryRepNEW]{Theorem \cite[Thm 3.6]{Tuynman..2018}}
Let $(\Hilbert, \inprodsym, \supinprsym)$ be a (complex) super Hilbert space, let $G$ be a connected $\CA$-Lie group and let $(\Dense, \rho)$ be a \psur{} of $G$ on $\Hilbert$. 
If we denote by $\Liealg g \cong T_eG$ the graded Lie algebra of $G$, then the map $\tau : \body \Liealg g \to \End(\Dense)$ defined by (see \refmetnaam{\labelMainSUR}{MainSURlabelrhopsiissmooth})
\begin{moneq}[defofthemaptauforSUR]
\tau(X)\psi = X\FgroupHilbert_\psi
\end{moneq}
is an even graded Lie algebra morphism with the following properties. 
\begin{enumerate}
\item\label{alternateSUR1}
For each $X\in \body \Liealg g_0$ the map $\tau(X)$ is the restriction of the infinitesimal generator of $\rho\bigl( \exp(tX)\bigr)$ to $\Dense$. 

\item\label{alternateSUR2}
For all $X\in \body \Liealg g$ the map $\tau(X)$ is graded skew-symmetric with respect to $\supinprsym$.

\item\label{alternateSUR3}
For all $g\in \body G$ and all $X\in \body \Liealg g_1$ we have
\begin{moneq}
\tau\bigl(\Ad(g)X\bigr) = \rho_o(g)\scirc \tau(X) \scirc \rho_o(g\mo)
\mapob.
\end{moneq}

\end{enumerate}
The triple $(\rho_o=\rho\restricted_{\body G}, \Dense, \tau)$ is called \emph{the infinitesimal form of $(\Dense, \rho)$}.

\smallskip

Conversely, if a triple $(\rho_o, \Dense,\tau)$ satisfies the conditions (\ref{alternateSUR1})--(\ref{alternateSUR3}) with $\rho_o$ an even unitary representation of $\body G$ on $\Hilbert$ (meaning that $\rho_o(g)$ preserves $\Hilbert_\alpha$), $\Dense\subset C^\infty(\rho_o)$ a dense graded subspace of $\Hilbert$ and invariant under the maps $\rho_o(g)$, $g\in \body G$ and $\tau:\body \Liealg g\to \End( \Dense )$ an even graded Lie algebra morphism, then there exists a unique group homomorphism $\rho:G\to \Aut(\Dense \otimes \CA^\KK)$ satisfying the conditions \refmetnaam{\labelMainSUR}{MainSURlabelrhopreservessupinpr}--\ref{MainSURlabelrhoonbodyisunitary} such that the map $\tau$ is given by \recalf{defofthemaptauforSUR} and such that the restriction of $\rho$ to $\body G$ is $\rho_o$. 

\end{proclaim}

\begin{remarque}{Nota Bene}
In the \myquote{experiments} of \recals{SECTSuperHeisenbergGroup} and \recals{SECTHeisenbergLikeSuperGroup} I will never comment on the condition \refmetnaam{\labelMainSUR}{MainSURlabelrhopsiissmooth} of a \psur. 
In all cases it will be a (more or less direct) consequence of the well known fact that, for $\psi\in \Schwartz{}{}$, the map
\begin{moneq}
\FgroupHilbert_\psi:\RR\to \Schwartz{}{}\ ,\qquad
\FgroupHilbert_\psi(t) \ :\ x\mapsto \psi(x+t)
\end{moneq}
is smooth (with $(\FgroupHilbert_\psi)' = \FgroupHilbert_{\psi'}$). 
\end{remarque}

Once we have decided on this definition of a super unitary representation, we have to address the question of equivalence. 
The most obvious answer would be to say that two super unitary representations $(\rho_i,\Dense_i)$ of $G$ on super Hilbert spaces $(\Hilbert_i,\inprodsym_i,\supinprsym_i)$ are equivalent if there exists a $k$-isomorphism of super Hilbert spaces $\Phi : \Hilbert_1\to \Hilbert_2$ that maps $\Dense_1$ onto $\Dense_2$ and that intertwines the two representations. 
However, there is a serious problem with this naive idea that will become obvious when we look at the following example.

\begin{remarque}[ExampleWithDifferentTopologiesInducedByChangeOfCoordinates]{Example}
Consider the abelian \glgrp/ $G=\RR^{1\vert2}$ of (graded) dimension $1\vert2$ with one global even coordinate $a$ and two odd coordinates $\beta,\gamma$. 
In terms of a Harish-Chandra pair, this means that $\body G= \RR$ is the abelian real line and that the graded Lie algebra $\Liealg g$ is the abelian super Lie algebra of (graded) dimension $1\vert2$ with (homogeneous) basis $e_1,f_2,f_3$ with $e_1$ even and $f_2,f_3$ odd. 
We are now going to define three super unitary representations $\rho_\varepsilon$, $\varepsilon\in \{0,\pm1\}$ of $G$ on the super Hilbert space $L^2(\RR^{1\vert2})$ \recalt{TheSuperManifoldsRpVertq}, for which we recall that we have an identification $L^2(\RR^{1\vert2})\cong L^2(\RR)^4$ in terms of which the inner product $\inprodsym_{1\vert2}$ and the \supscal product $\supinprsym_{1\vert2}$ are given by \recalf{SimplifiedInnerProductOnL2Rpq} and \recalf{SimplifiedSuperScalarProductOnL2Rpq}
\begin{align*}
\inprod\chi\psi_{1\vert2}
&
=
\inprodd{\chi_0}{\psi_0} + \inprodd{\chi_1}{\psi_1} + \inprodd{\chi_2}{\psi_2} + \inprodd{\chi_{12}}{\psi_{12}} 
\\
\homsuperinprod\chi\psi{_{1\vert2}}
&
=
\inprodd{\chi_0}{\psi_{12}} + \inprodd{\chi_1}{\psi_2} - \inprodd{\chi_2}{\psi_1} + \inprodd{\chi_{12}}{\psi_{0}} 
\mapob,
\end{align*}
where $\inprodd\cdot\cdot$ denotes the standard inner product on $L^2(\RR)$. 
All three representations are defined on the same (maximal) dense subspace $\Dense\subset L^2(\RR^{1\vert2})\cong \bigl(L^2(\RR)\bigr)^4$ given by
\begin{moneq}
\Dense = \{\,(\psi_0,\psi_1,\psi_2,\psi_{12})\in C^\infty(\RR)^4\mid \forall I\subset \{1,2\}\ \forall n\in\NN: \psi_I^{(n)}\in L^2(\RR)\,\}
\mapob.
\end{moneq}
With these preparations, the three maps $\rho_\varepsilon:G\to \Aut(\Dense\otimes \CA^\CC)$ are given by the formula
\begin{align}
\rho_\varepsilon(a,\beta,\gamma)(\psi_0,\psi_1,\psi_2,\psi_{12})
&
=
\hBigl(
T_a\psi_0 
-\beta T_a\psi_1 
-\gamma T_a\psi_2 
+ \beta\gamma T_a(\psi_{12}
-\varepsilon\psi_0')
\ ,
\notag
\\&
\kern2.25em
T_a\psi_1 
-\gamma T_a(\psi_{12} - \varepsilon\psi_0') \ ,\ 
\notag
\\&
\kern2.625em
T_a\psi_2 
+\beta T_a(\psi_{12}-\varepsilon\psi_0')  \ ,\ 
\label{TripleRepresentationAbelianSuperGroupe1vert2}
\\&
\kern3em
T_a\psi_{12}
+ \varepsilon\beta\gamma T_a(\psi_{12}-\varepsilon\psi_0')' 
-  \varepsilon\gamma T_a\psi_2' 
-  \varepsilon\beta T_a\psi_1'
\hBigr)
\mapob,
\notag
\end{align}
where $T_a$ denotes the translation operator $(T_af)(x) = f(x-a)$. 
It is immediate that the restriction of $\rho_\varepsilon$ to $\body G=\RR$ is the standard unitary representation $\rho_o$ of $\RR$ on $L^2(\RR)^4$ (or better: four copies of the standard representation) given by
\begin{moneq}
\rho_o(a)(\psi_0,\psi_1,\psi_2,\psi_{12}) = (T_a\psi_0, T_a\psi_1,T_a\psi_2,T_a\psi_{12})
\mapob,
\end{moneq}
and it is a standard result that the set of smooth functions for this representation is $\Dense$. 
It follows that in terms of Harish-Chandra pairs the infinitesimal form of these representations (see \recalt{equivalentDefSuperUnitaryRepNEW}\,) is given by this standard unitary representation $\rho_o$ of $\RR$ on $\Dense$ and that the (even) graded Lie algebra morphisms $\tau_\varepsilon:\body \Liealg g\to \End(\Dense)$ are given by
\begin{moneq}
\begin{matrix}
\tau_\varepsilon(e_1)(\psi_0,\psi_1,\psi_2,\psi_{12}) & = & (\ -\psi_0' &,& -\psi_1' &,& -\psi_2' &,& -\psi_{12}'\ )
\\[1.5\jot]
\tau_\varepsilon(f_2)(\psi_0,\psi_1,\psi_2,\psi_{12}) & = & (\ -\psi_1 &,&  0 &,&  -\psi_{12} + \varepsilon\psi_0' &,&  -\varepsilon\psi_1'\ )
\\[1.5\jot]
\tau_\varepsilon(f_3)(\psi_0,\psi_1,\psi_2,\psi_{12}) & = & (\ -\psi_2 &,&  \psi_{12} - \varepsilon\psi_0' &,&  0 &,&  -\varepsilon\psi_2'\ )
\rlap{\mapob.}
\end{matrix}
\end{moneq}
I leave it to the reader to verify that the maps $\rho_\varepsilon$ and $\tau_\varepsilon$ preserve the inner product and the \supscal product in the appropriate ways for a super unitary representation.

A direct computation shows that the map $\Phi:\Dense \to \Dense$ defined as
\begin{moneq}
\Phi(\psi_0,\psi_1,\psi_2,\psi_{12}) = (\psi_0,\psi_1,\psi_2,\psi_{12}+\psi_0')
\end{moneq}
is an even (right-linear) bijection intertwining the \supscal product $\supinprsym_{1\vert2}$ on source and target and intertwining the representations $\rho_\varepsilon$ and $\rho_{\varepsilon+1}$ ($\varepsilon=-1,0$):
\begin{align*}
\superinprod[2]{\Phi(\chi)}{\Phi(\psi)}
&
=
\inprodd{\Phi(\chi)_0}{\Phi(\psi)_{12}} + \inprodd{\Phi(\chi)_1}{\Phi(\psi)_{2}} 
\\&
\kern6em
- \inprodd{\Phi(\chi)_2}{\Phi(\psi)_{1}} + \inprodd{\Phi(\chi)_{12}}{\Phi(\psi)_{0}}
\\&
=
\inprodd{\chi_{0}}{\psi_{12}+\psi_0'} + \inprodd{\chi_1}{\psi_{2}} - \inprodd{\chi_2}{\psi_{1}} + \inprodd{\chi_{12}+\chi_0'}{\psi_{0}}
\\&
=
\superinprod\chi\psi + \inprodd{\chi_0}{\psi_0'} + \inprodd{\chi_0'}{\psi_0} 
\\&
=
\superinprod\chi\psi + 
\int_\RR\extder x\kern0.5em\hBigl(\ \overline{\chi_0(x)}\cdot \psi_0'(x) + \overline{\chi_0'(c)}\cdot \psi_0(x) \hBigr)
\\&
=
\superinprod\chi\psi + 
\int_\RR\extder x\kern0.5em\bigl(\ \overline{\chi_0}\cdot \psi_0\bigr)'(x)
=
\superinprod\chi\psi
\mapob.
\end{align*}
And
\begin{align*}
\Phi\bigl(\rho_\varepsilon(a,\beta,\gamma)\psi\bigr)
\kern0.5em
&
\oversetalign{\recalf{TripleRepresentationAbelianSuperGroupe1vert2}}\to=
\kern0.5em
\Phi\hBigl(
T_a\psi_0 
-\beta T_a\psi_1 
-\gamma T_a\psi_2 
+ \beta\gamma T_a(\psi_{12}
-\varepsilon\psi_0')
\ ,
\notag
\\&
\kern2.25em
T_a\psi_1 
-\gamma T_a(\psi_{12} - \varepsilon\psi_0') \ ,\ 
\notag
\\&
\kern2.625em
T_a\psi_2 
+\beta T_a(\psi_{12}-\varepsilon\psi_0')  \ ,\ 
\notag
\\&
\kern3em
T_a\psi_{12}
+ \varepsilon\beta\gamma T_a(\psi_{12}-\varepsilon\psi_0')' 
-  \varepsilon\gamma T_a\psi_2' 
-  \varepsilon\beta T_a\psi_1'
\hBigr)
\mapob,
\notag
\\&
=
\Bigl(
T_a\psi_0 
-\beta T_a\psi_1 
-\gamma T_a\psi_2 
+ \beta\gamma T_a(\psi_{12}
-\varepsilon\psi_0')
\ ,
\notag
\\&
\kern2.25em
T_a\psi_1 
-\gamma T_a(\psi_{12} - \varepsilon\psi_0') \ ,\ 
\notag
\\&
\kern2.625em
T_a\psi_2 
+\beta T_a(\psi_{12}-\varepsilon\psi_0')  \ ,\ 
\notag
\\&
\kern3em
T_a\psi_{12}
+ \varepsilon\beta\gamma T_a(\psi_{12}-\varepsilon\psi_0')' 
-  \varepsilon\gamma T_a\psi_2' 
-  \varepsilon\beta T_a\psi_1'
\\&
\kern3.375em
+\hBigl(
T_a\psi_0 
-\beta T_a\psi_1 
-\gamma T_a\psi_2 
+ \beta\gamma T_a(\psi_{12}
-\varepsilon\psi_0')
\hBigr)'\ 
\Bigr)
\\&
=
\Bigl(
T_a\psi_0 
-\beta T_a\psi_1 
-\gamma T_a\psi_2 
+ \beta\gamma T_a\bigl((\psi_{12}+\psi_0')
-(\varepsilon+1)\psi_0'\bigr)
\ ,
\notag
\\&
\kern2.25em
T_a\psi_1 
-\gamma T_a\bigl((\psi_{12}+\psi_0')-(\varepsilon+1)\psi_0'\bigr) \ ,\ 
\notag
\\&
\kern2.625em
T_a\psi_2 
+\beta T_a\bigl((\psi_{12}+\psi_0')-(\varepsilon+1)\psi_0'\bigr)  \ ,\ 
\notag
\\&
\kern3em
T_a(\psi_{12}+\psi_0')
+ (\varepsilon+1)\beta\gamma T_a\bigl((\psi_{12}+\psi_0')-(\varepsilon+1)\psi_0'\bigr)' 
\\&
\kern7.375em
-  (\varepsilon+1)\gamma T_a\psi_2' 
-  (\varepsilon+1)\beta T_a\psi_1'
\,\Bigr)
\\&
=
\rho_{\varepsilon+1}(a,\beta,\gamma)\Phi(\psi)
\mapob.
\notag
\end{align*}
However, $\Phi$ certainly does not preserve the inner product, as it is not continuous due to the presence of the derivative in its expression. 
It thus cannot be extended to a $k$-isomorphism of super Hilbert spaces and we are led to conclude that these three representations are not equivalent. 

On the other hand, $\Phi$ is the bijection induced on $\Dense$ by the change of coordinates 
\begin{moneq}
(\xh, \hat \eta,\hat \zeta) \equiv \phi(x,\eta,\zeta)=(x-\eta\zeta, \eta,\zeta) 
\end{moneq}
in the manifold $\RR^{1\vert2}$. 
The expression $\psih$ of $\psi$ in terms of these new coordinates is given by $\psih(\xh, \hat \eta,\hat \zeta) = \psi(x,\eta,\zeta)$, \ie,
\begin{align*}
\psih(\xh,\hat\eta,\hat\zeta)
&
=
\psih_0(\xh) + \hat\eta\,\psih_1(\xh) + \hat\zeta\,\psih_2(\xh) + \hat\eta\hat\zeta\,\psih_{12}(\xh)
\\&
=
\psih_0(x-\eta\zeta) + \eta\,\psih_1(x-\eta\zeta) + \zeta\,\psih_2(x-\eta\zeta) + \eta\zeta\,\psih_{12}(x-\eta\zeta)
\\&
=
\psih_0(x) + \eta\,\psih_1(x) + \zeta\,\psih_2(x) + \eta\zeta\,\bigl(\psih_{12}(x)-\psih_0'(x)\bigr)
\\
\text{\footnotesize but also\kern4em}
&
=
\psi(x, \eta,\zeta)
\\&
=
\psi_0(x) + \eta\, \psi_1(x) + \zeta \,\psi_2(x) + \eta\zeta\,\psi_{12}(x)
\mapob,
\end{align*}
giving us the identification (the bijection $\Phi$)
\begin{moneq}
(\psih_0,\psih_1,\psih_2,\psih_{12}) \cong (\psi_0,\psi_1,\psi_2,\psi_{12}+\psi_0')
\mapob.
\end{moneq}
In other words, the three representations $\rho_\varepsilon$ are manifestations of the same representation of $G=\RR^{1\vert2}$ on $L^2(\RR^{1\vert2})$ in terms of different coordinate systems. 
More precisely, $\rho_0$ is the standard left-regular representation given by
\begin{moneq}
\bigl(\rho_0(a,\beta,\gamma)\psi\bigr)(x,\eta,\zeta) = \psi(x-a, \eta-\beta,\zeta-\gamma)
\end{moneq}
and $\rho_{\pm1}$ are obtained by changing the coordinates $(x,\eta,\zeta)$ by $\phi^{\pm1}$. 

A completely different outlook of this phenomenon (three different representations of $G$ on the same super Hilbert space and the same dense subspace linked by a change of coordinates) is obtained when we pass through an \myquote{obvious} intermediate. 
As said, the map $\Phi:\Dense\to \Dense$ intertwines the \supscal product $\supinprsym_{1\vert2}$ on source and target, but not the inner product $\inprodsym_{1\vert2}$. 
But when we use $\Phi$ to transport the inner product $\inprodsym_{1\vert2}$ on the source\slash target to the target\slash source, we get the inner products
\begin{moneq}
\inprod\chi\psi_\varepsilon
=
\inprodd{\chi_0}{\psi_0} + \inprodd{\chi_1}{\psi_1} + \inprodd{\chi_2}{\psi_2} + \inprodd{\chi_{12}-\varepsilon\chi_0'}{\psi_{12}-\varepsilon\psi_0'}
\mapob,
\end{moneq}
and then the map
\begin{moneq}
\Phi = (\Dense, \inprodsym_\varepsilon, \supinprsym_{1\vert2}) \to (\Dense, \inprodsym_{\varepsilon+1}, \supinprsym_{1\vert2})
\end{moneq}
will be an isomorphism of super inner product spaces for $\varepsilon= -1,0$, as we have
\begin{moneq}
\inprod{\Phi(\chi)}{\Phi(\psi)}_{\varepsilon+1} = \inprod\chi\psi_\varepsilon
\mapob.
\end{moneq}
The completion of $(\Dense, \inprodsym_\varepsilon, \supinprsym_{1\vert2})$ will still be (isomorphic to) the super Hilbert space $\bigl(L^2(\RR^{1\vert2}), \inprodsym_{1\vert2}, \supinprsym_{1\vert2}\bigr)$, but the embedding of $\Dense$ for $\varepsilon=\pm1$ as a dense subspace is not the \myquote{natural} one, but uses $\Phi^{\pm1}$. 
We now can \myquote{disguise} this change by defining three different \psur s $\rhoh_\varepsilon$ of $G$ using \recalt{ProtoSuperUnitaryRepExtendsToUniqueSuperUnitaryRep} as follows. 
We define $\rhoh_\varepsilon$ on $(\Dense_\varepsilon, \inprodsym_\varepsilon, \supinprsym_{1\vert2})$ by 
\begin{moneq}
\forall \varepsilon\in \{-1,0,1\} 
\quad:\quad 
\Dense_\varepsilon=\Dense
\quad\text{and}\quad
\rhoh_\varepsilon=\rho_0
\mapob.
\end{moneq}
We thus have a single graded vector space equipped with a \supscal product $(\Dense, \supinprsym_{1\vert2})$ and a single homomorphism $\rho_0:G\to \Aut(\Dense\otimes\CA^\CC)$ preserving the \supscal product. 
But we have three different non-equivalent inner products that all three make it into a \psur. 
\end{remarque}

Given that I find it highly undesirable that a simple change of coordinates produces an inequivalent representation, I conclude that we have to abandon the naïve idea that an equivalence of super unitary representations is given by a $k$-isomorphism of super Hilbert spaces that intertwines the representations. 
The fact that we can disguise a change of coordinates as different inner products (topologies) on a super inner product space is another argument in favour of my choice not to impose an (obvious) link between the \supscal product and the inner product used in the definition of a super inner product space (see \recals{SECTInitialDefinitionsAndMore}).

\begin{definition}[DefinitionEquivalenceSuperUnitaryRep]{Definition}
Let $(\rho_1,\Dense_1)$ and $(\rho_2,\Dense_2)$ be two super unitary representations of \aglgrp/ $G$ on the super Hilbert spaces $(\Hilbert_1,\inprodsym_1, \supinprsym_1)$ and $(\Hilbert_2,\inprodsym_2, \supinprsym_2)$ respectively. 
We will say that these two super unitary representations are \stresd{equivalent} if there exist dense graded subspaces $\Dense_1'\subset \Dense_1$, $\Dense_2'\subset \Dense_2$ and a right-linear bijection $\Phi:\Dense_1'\to\Dense_2'$ satisfying the two conditions
\begin{enumerate}[\kern0.5em\labelMainSUReq1.]
\item\label{MainSUREQlabelPhiIntertwinesRepresentations}
$\Phi$ intertwines the two representations: $\forall g\in G:\rho_2(g)\scirc \Phi = \Phi \scirc \rho_1(g)$;

\item\label{MainSUREQlabelPhiIntertwinesSuperScalarProducts}
$\Phi$ is partially a $k$-isomorphism of super inner product spaces in the sense that it should satisfy conditions \refmetnaam{\labelIsomSHS}{conditiondirectsuminequivalenceSHS} and \refmetnaam{\labelIsomSHS}{conditionoddisoforequivalenceofSHS}, but not (necessarily) \refmetnaam{\labelIsomSHS}{conditionequivalenceofSHSisunitary}. 

\end{enumerate}

\end{definition}

\begin{remarque}{Remark}
Once we have admitted that an equivalence between two super unitary representations should not necessarily be (bi-)continuous\slash unitary, the reader might wonder why I introduced the notion of equivalence using dense subspaces. 
A more natural choice would be that two super unitary representations $(\Dense_i,\rho_i)$ of \aglgrp/ $G$ on $\Hilbert_i$ are equivalent if there exists a linear bijection $\Phi:\Dense_1\to \Dense_2$ that intertwines the representations $\rho_i$ and the \supscal products $\supinprsym_i$. 
In all the examples I know of, the equivalence bijection $\Phi=\Dense_1'\to \Dense_2'$ (with $\Dense_i'\subset \Dense_i$ dense subspaces) extends to an equivalence bijection $\Phi=\Dense_1\to \Dense_2$, but I have no formal proof that this will be the case in general.
On the other hand, according to \recalt{ProtoSuperUnitaryRepExtendsToUniqueSuperUnitaryRep}, any \psur{} extends in a unique way to a super unitary representation. 
It follows that in order to discuss super unitary representations, we can restrict attention to \psur{}s (a similar idea is expressed in \cite[\S2.3, Prop.2]{CCTV..2006}). 
And my definition of equivalence simply says that we must have a bijection between the representations spaces of \psur{}s. 
\end{remarque}

\section{\label{SECTpVertqFamiliesOfSURs} \texorpdfstring{$p\vert q$}{pq}-parameter families of super unitary representations}

In order to discuss our decompositions of super unitary representations, we need the notion of a direct Berezin integral of a family of (proto super unitary) representations and the associated notion of a $p\vert q$-parameter family of \psur s (in particular the dependence on odd parameters). 
Now we \stress{could} copy the non-super approach, defining a $p\vert q$-parameter family of super unitary representations as a super unitary representation $\rho^m$ for every value $m\in \RR^{p\vert q}$, all defined on a common dense subspace $\Dense$ of a fixed super Hilbert space $\Hilbert$. 
And then to define the direct Berezin-integral of such a family as the (proto super unitary) representation $\rhot$ defined on $C^\infty_c(\RR^{p\vert q}, \Dense)\subset L^2(\RR^{p\vert q}, \Hilbert)$ \recalt{TheSuperManifoldsRpVertq} by
\begin{moneq}
\bigl(\rhot(g)\psi\bigr)(m) = \bigl(\rho^m(g)\bigr)\bigl(\psi(m)\bigr)
\mapob.
\end{moneq}
However, as we will see later in naturally occurring families (e.g., the $1\vert 1$-parameter family $\rho^{k,\mu}$ of $\SuperH$ \recals{SubSECTDecompSuperHeisenberg}), there is no guarantee that for $m\notin \body \RR^{p\vert q}\equiv \RR^p$ (and in particular with non-zero values for the odd parameters), the maps $\rho^m(g)$ will preserve the inner product of the Hilbert space $\Hilbert$. 
As a consequence, we have no way to control the behaviour of $\rhot$ on $C^\infty_c(\RR^{p\vert q}, \Dense)$ with respect to the inner product $\inprodsym_{p\vert q}^\Hilbert$, meaning that we have no way to guarantee that $\rhot$ will be a \psur. 
On the other hand, we have seen in \recalt{ExampleWithDifferentTopologiesInducedByChangeOfCoordinates} that we have a certain leeway in the choice of an inner product when defining a super unitary representation. 
The solution I adopt below (all examples I know of satisfy this definition) is to require that at least for $m\in \RR^p = \body\RR^{p\vert q}$ the map $\rho^m$ is a \psur. 
And to define the direct integral representation $\rhot$ on $C^\infty_c(\RR^{p\vert q}, \Dense)$ only if there exists an inner product on this space making $\rhot$ a \psur.

\begin{definition}[DefOfpqFamilyOfPSURs]{Definitions}
\itemizedef[*]
A \stresd{$p\vert q$-parameter family of \psur{}s} of \aglgrp/ $G$ on a (complex) super Hilbert space $(\Hilbert, \inprodsym_\Hilbert, \supinprsym_\Hilbert)$ is a couple $(\Dense, \rho)$ with the following properties: 
$\Dense\subset \Hilbert$ is a dense graded subspace and $\rho:\RR^{p\vert q}\times G\to \Aut(\Dense\otimes\CA^\CC)$ is a map satisfying the following conditions.
\begin{enumerate}[\kern0.5em\labelFamilySUR1.]
\item\label{FamilySUR1}
For each fixed value $m\in \RR^{p\vert q}$ the map $\rho^m:G\to \Aut(\Dense\otimes\CA^\CC)$ defined as $\rho^m(g) = \rho(m,g)$ is a group homomorphism.

\item\label{FamilySUR2}
For all $(m,g)\in  \RR^{p\vert q}\times G$, $\rho^m(g)$ preserves $\supinprsym_\Hilbert$, \ie, 
\begin{moneq}
\forall (m,g)\in \RR^{p\vert q}\times G\ 
\forall \chi,\psi\in \Dense: 
\homsuperinprod[2]{\rho^m(g)\chi}{\rho^m(g)\psi}{_\Hilbert} = \homsuperinprod\chi\psi{_\Hilbert}
\mapob.
\end{moneq}

\item\label{FamilySUR3}
For all $\psi\in \Dense$ the map $\FgroupHilbert_\psi:\RR^{p\vert q}\times G\to \Dense \otimes \CA^\CC$, $\FgroupHilbert_\psi(m,g) 
\equiv \rho(m,g)\psi$ is smooth.

\item\label{FamilySUR3.5}
For all $(m,g)\in \body(\RR^{p\vert q}\times G)\cong \RR^p\times \body G$, $\rho^m(g)$ preserves $\inprodsym_\Hilbert$.

\item\label{FamilySUR4}
There exists an inner product $\inprodsym_{\Denset}$ on $\Denset = C^\infty_c(\RR^{p\vert q};\Dense)$ such that the couple $(\Denset,\rhot)$ is a \psur{} on $(\Denset,\inprodsym_{\Denset}, \supinprsym_{p\vert q}^\Dense)$, where $\rhot:G\to \Aut(\Denset\otimes \CA^\CC)$ is given by
\begin{moneq}
\forall \psi\in \Denset\ 
\forall m\in \RR^{p\vert q}
\quad:\quad
\bigl(\rhot(g)\psi\bigr)(m) = \rho^m(g)\psi(m)
\mapob.
\end{moneq}

\end{enumerate}
Note that by definition of a $p\vert q$-parameter family, the map $\rhot$ automatically is a homomorphism preserving the \supscal product $\supinprsym_{p\vert q}^\Hilbert$. 
However, there is no guarantee that the standard inner product $\inprodsym_{p\vert q}^\Hilbert$ on $\Denset$ will make $(\Denset,\rhot)$ a \psur{} (see below for an example where there is none). 
We will call $\rhot$ the \stresd{direct integral representation of the family $\rho$} on $(\Denset, \inprodsym_{\Denset}, \supinprsym_{p\vert q}^\Hilbert)$ (or in a sloppy notation we will talk about the family $\rho^m$). 
It is immediate that two different inner products on $\Denset$ yield equivalent \psur s via the identity isomorphism, simply because with our definition \recalt{DefinitionEquivalenceSuperUnitaryRep} of equivalence we do not look at the inner product. 
This justifies that we talk about \stress{the} direct integral representation.

\itemizedef
Let $\rho:\RR^{p\vert q}\times \Dense\to \Aut(\Dense\otimes \CA^\CC)$ be a $p\vert q$-parameter family of \psur s and let $m\in \RR^{p\vert q}$ be fixed. 
We will say that the \myquote{representation} $\rho^m$ of $G$ is \stresd{irreducible} if there exists no invariant closed graded subspace $F\subset \Dense$, \ie, if $F\subset \Dense$ is a closed graded subspace such that $\rho^m(g)(F)\subset F\otimes \CA^\CC$ for all $g\in G$, then $F=\{0\}$ or $F=\Dense$.

\end{definition}

\begin{remarque}{Example\slash Remark}
\itemizedef[*]
The map $\bigl(\rho^\mu(t)\psi\bigr)(\xi) = \psi(\xi+t\mu)$ satisfies the conditions \refmetnaam{\labelFamilySUR}{FamilySUR1}-\ref{FamilySUR3.5} of a $0\vert1$-parameter family of \psur s of $G=\CA_0$ (essentially the real line) on $\Dense\equiv L^2(\RR^{0\vert 1})\cong\CC^2$. 
On the other hand, there exists no inner product on $\Denset\equiv C_c^\infty(\RR^{0\vert1};\Dense)\cong \CC^4$ making $\rhot$ a \psur{} on $\Denset$, \ie, the family $\rho^\mu$, $\mu\in \CA_1$ is not a $0\vert1$-parameter family of \psur s in the sense of \recalt{DefOfpqFamilyOfPSURs}. 

\itemizedef
I have hesitated a long time whether or not to include the condition \refmetnaam{\labelFamilySUR}{FamilySUR4} in the definition of a $p\vert q$-parameter family of \psur s. 
I finally decided to include it for two reasons: it avoids having two different kind of families, satisfying \refmetnaam{\labelFamilySUR}{FamilySUR4} on not, and because I see no real interest in families for which the map $\rhot$ (which could be called the algebraic direct integral representation) is not a \psur.  
\end{remarque}

\section{\label{GQgeneral}Geometric quantization of coadjoint orbits}

\renewcommand\LD[2][0]{#2^*}

The final ingredient we will need for our \myquote{experiments} is the geometric quantization procedure of coadjoint orbits (the orbit method) to link coadjoint orbits of a (super) Lie group $G$ to (irreducible) unitary representations of $G$. 
But before we go into the details, first a remark on notation.

\begin{remarque}{Remarks on notation}
\itemizedef[*]
Contraction of a $k$-form $\alpha$ with a vector field $X$ will be denoted as $\contrf{X}\alpha$ and, for $k\ge2$, contraction with two vector fields $X$ and $Y$ will be denoted a $\contrf{X,Y}\alpha$ and is defined as repeated contraction:
\begin{moneq}
\contrf{X,Y}\alpha = \contrf{X}{\bigl(\contrf Y\alpha\bigr)}
\qquad\text{or}\qquad
\contrf{X,Y}{} = \contrf{X}{} \scirc \contrf{Y}{}
\mapob.
\end{moneq}

\itemizedef
The Lie derivative of a $k$-form $\alpha$ in the direction of a vector field $X$ is denoted as $\Lieder(X)\alpha$. It is defined by the formula
\begin{moneq}
\Lieder(X) = d\scirc \contrf X{} + \contrf{X}{}\scirc d
\mapob.
\end{moneq}

\itemizedef
As is usual, we identify the (super) Lie algebra $\Liealg g$ of a (super) Lie group $G$ with the set of left-invariant vector fields on $G$ and as a consequence its dual $\Liealg g^*$ as the set of left-invariant $1$-forms. 
Here the qualifier \myquote{left} refers to left-multiplication in $G$, not to left-linearity with respect to $\CA$. 
But we have defined the cotangent bundle $T^*G$ as the left-dual of the tangent bundle $TG$ (so we have coherence with the standard notation of contraction on the left of a $k$-form with a vector field). 
This implies that our dual Lie algebra $\Liealg g^*$ will be the left-dual of $\Liealg g$, where here the qualifier \myquote{left} refers to left-linearity over $\CA$. 
We thus will denote evaluation of an element $\mu\in \Liealg g^*$ on an element $X\in \Liealg g$ as $\contrs{X}{\mu}$ with $X$ to the left of $\mu$ (by analogy to the contraction of a vector field with a $1$-form). 
This means that we have the equality
\begin{moneq}
\contrs{\lambda\cdot X}\mu= \lambda\cdot\contrs{X}{\mu}
\mapob,
\end{moneq}
for any $\lambda\in \CA$. 
Had we written $\mu(X)$ for the evaluation, we would have had a conflict with standard conventions about signs: we would (still) have $\mu(\lambda X) = \lambda\,\mu(X)$ because $\mu$ is left-linear. 
However, standard conventions would tell us that we should add a sign according to the parity of $\lambda$ and $\mu$ as $\mu(\lambda  X) = (-1)^{\vert \lambda\vert\cdot\vert\mu\vert}\, \lambda\,\mu(X)$, simply because we interchange $\lambda$ and $\mu$.

\itemizedef
If $\Phi : G \times M \to M$ denotes the (left) action of a (super) Lie group $G$ on a (super) manifold $M$, then for all $v \in \Liealg g = T_eG$ the associated
\stress{fundamental vector field} $v^M$ on $M$ is defined as $v^M\restricted_m =
- T_{(e,m)}\Phi (v,0)$. The minus sign is conventional and ensures that the map
from $\Liealg g$ to vector fields on $M$ is a homomorphism of super Lie algebras.
\end{remarque}

We start the description of the geometric quantization procedure with a half-baked summary of the classical case. 
Given a coadjoint orbit $\Orbit_{\mu_o}\subset \LD[1]{\Liealg g}$ through $\mu_o\in \LD[1]{\Liealg g}$ of a Lie group $G$ (with Lie algebra $\Liealg g$), one constructs a principal $\Sphere^1\equiv U(1)$ bundle $Y\to \Orbit_{\mu_o}$ with a connection $1$-form $\Gamma$ whose curvature is the (Kostant-Kirillov-Souriau) symplectic form $\omega_o$ of the coadjoint orbit, which is determined by the conditions
\begin{moneq}[DefiningFormulaForKKSForm]
\forall v,w\in \Liealg g
\quad:\quad
-\contrf{v^{\LD[1]{\Liealg g}},w^{\LD[1]{\Liealg g}}}{\omega_{\mu}} = 
\contrs{[v,w]}{\mu}
\mapob,
\end{moneq}
where $v^{\LD[1]{\Liealg g}}$ denotes the fundamental vector field on $\LD[1]{\Liealg g}$ associated to $v\in \Liealg g$ for the Coadjoint action. 
The structure group of $Y$ is $\RR/d\ZZ\cong \Sphere^1$ with $d\in (0,\infty)$ (see \cite{TW87}) on which we will use the coordinate $t\in \RR$ (modulo $d$). 
As in our examples the symplectic form $\omega_o$ is exact and the coadjoint orbits admit a single chart, it follows that $Y=\Orbit_{\mu_o}\times (\RR/d\ZZ)$ and that the connection $1$-form $\Gamma$ is of the form
\begin{moneq}
\Gamma = \extder t + \Theta
\mapob,
\end{moneq}
with $\Theta$ a potential of $\omega_o$. 
The presence of the (canonical) momentum map on a coadjoint orbit implies that the infinitesimal action of $G$ on the orbit can be lifted to an  infinitesimal action on the prequantum bundle preserving the connection $1$-form. 
More precisely, for $X \in \Liealg g$ the infinitesimal generator $X^{Y}$ on $Y$ is determined by the equations
\begin{moneq}[LiftingEquationsOfGroupActionOnPrequantumBundle]
\contrf{X^{Y}\restricted_{(\mu,t)}}\Gamma = \contrs{X}{\mu}
\qquad\text{and}\qquad
\Lieder(X^{Y})\Gamma = 0
\mapob.
\end{moneq}
This action integrates to an action $\phi$ of $G$ on $Y$ preserving the connection $1$-form. 
(Remark: in general the lifting of the coadjoint action on an orbit to a prequantum bundle is equivalent to the existence of a character of the isotropy subgroup of the base point whose infinitesimal form is the given base point, which in turn corresponds to an integrality condition on the symplectic form.) 
With this lifted action $\phi$ of $G$ on $Y$, we obtain the natural representation space of functions $\psi$ on $Y$ with the standard action $\Phi(g)$ of $g\in G$:
\begin{moneq}[InducedActionOnPrequantumFunctions]
\bigl(\Phi(g)\psi\bigr)(y) = \psi\bigl(\phi(g)\mo(y)\bigr)
\mapob.
\end{moneq}
However, this representation space will certainly not be irreducible. 
First one observes that a Fourier transform with respect to the circle coordinate (the coordinate of the structure group) allows a decomposition of the space of functions on $Y$ into various Fourier modes, all of which are invariant under the lifted action of $G$. 
A (not completely arbitrary) choice of a non zero Fourier mode still is not irreducible, so the final step in geometric quantization is the introduction of an invariant polarization, \ie, a foliation that is maximal isotropic with respect to the symplectic form on the orbit and which is invariant under the group action; such objects are in bijection with subalgebras $\Liealg h\subset \Liealg g$ containing the stabilizer subalgebra $\Liealg g_{\mu_o}$ and which are maximal with respect to the condition $\contrs{[\Liealg h, \Liealg h]}{\mu_o} = 0$. 
One then shows that the subspace of functions on $Y$ that satisfy the two conditions
\begin{itemize}
\item
belonging to the non-zero Fourier mode and

\item
constant in the direction of the horizontally lifted (invariant) polarization
\end{itemize}
is an invariant subspace under the natural action $\Phi$ given above. 
This is the representation space constructed by geometric quantization of the orbit. 
Except that we need some technical details (such as half-densities or half-forms) to obtain an inner product on this representation space, details that I will gloss over. 

\medskip

For geometric quantization of a symplectic \gmfd/ in the context of the coadjoint orbits of \aglgrp/ $G$, I just imitate the above approach. 
As explained in \cite{Tuynman..2010}, since the symplectic form on a (super)coadjoint orbit can be non-homogeneous, the structure group of the prequantum bundle $Y$ will be $(\CA_0/d\ZZ)\times \CA_1$ of dimension $1\vert1$ with $d\in (0,\infty)$. 
On $(\CA_0/d\ZZ)\times \CA_1$ we will use coordinates $(t,\tau)$ with $t$ the canonical even coordinate on $\CA_0$ modulo $d$ and $\tau$ the canonical coordinate on $\CA_1$. 
In our case the orbits all have a single global chart of the form $(\CA_0)^p \times (\CA_1)^q\equiv \RR^{p\vert q}$, hence are simply connected and thus the prequantum bundle is unique and of the form $Y=\Orbit_{\mu_o} \times (\CA_0/d\ZZ)\times \CA_1$. Moreover, there exists a global potential $\Theta$ for the symplectic form $\omega$, \ie, $\extder\Theta = \omega$. 
The connection $1$-form $\Gamma$ on $Y$, an even Lie algebra (of the structure group $(\CA_0/d\ZZ)\times \CA_1$) valued $1$-form on $Y$, then takes the form
\begin{moneq}
\Gamma = (\Theta_0+\extder t) \otimes \partial_t + (\Theta_1+\extder \tau) \otimes \partial_\tau
\mapob,
\end{moneq}
with $\Theta=\Theta_0+\Theta_1$ the splitting of $\Theta$ into its homogeneous parts. 
It is useful to identify $\Gamma$ with the mixed $1$-form $(\Theta_0+\extder t) + (\Theta_1+\extder \tau) = \Theta + \extder (t+\tau)$ on $Y$ and in the sequel we will often pretend that we have this equality.

As in the non-super case, we can lift the (coadjoint) $G$-action on the orbit to a connection preserving action $\phi$ on $Y$ whose infinitesimal form is determined by exactly the same equations \recalf{LiftingEquationsOfGroupActionOnPrequantumBundle}. 
And as in the non-super case, we get a natural representation $\Phi$ on the space of $\CA^\CC$-valued functions on $Y$, a representation that certainly will not be irreducible. 
Just as in the ungraded case we can perform a Berezin-Fourier transform with respect to the coordinates of the structure group $\Sphere^1\times \CA_1$ giving a family of invariant subspaces. 
And then we introduce the notion of a polarisation, \ie, an invariant isotropic foliation of maximal rank. 
And just as in the non-super case, such objects are in bijection with polarization subalgebras, \ie, (graded) subalgebras $\Liealg h\subset \Liealg g$ containing the stabilizer subalgebra $\Liealg g_{\mu_o}$ and which are maximal with respect to the condition $\contrs{[\Liealg h, \Liealg h]}{\mu_o} = 0$. 
The idea then is to look at functions that satisfy the same conditions as in the ungraded case: beloning to a non-zero Berezin-Fourier mode and constant in the directions of the polarization. 
The second condition translates easily to the super case as
\begin{moneq}[secondGQcondition]
\forall X\in \Liealg h
\quad:\quad
X^h\psi=0
\mapob,
\end{moneq}
where $\Liealg h$ is the polarization subalgebra and where $X^h$ is the horizontal lift to $Y$ of the infinitesimal action $X^{\LD[1]{\Liealg g}}$ (the fundamental vector field on the orbit $\Orbit_{\mu_o}$ associated to $X\in \Liealg g$), a vector field of the form 
\begin{moneq}
X^h = X^{\LD[1]{\Liealg g}} + p_0 \cdot \partial_t + p_1\cdot \partial_\tau
\mapob,
\end{moneq}
where the coefficients $p_i$ are determined by the equations of being horizontal:
\begin{moneq}[DefinitionOfHorizontalLiftFundamentalVectorFields]
0 = \contrf{X^{h}}{\Gamma_0} \equiv \contrf{X^{h}}{(\Theta_0+\extder t)}
\quad\text{and}\quad
0 = \contrf{X^{h}}{\Gamma_1} \equiv \contrf{X^{h}}{(\Theta_1+\extder \tau)}
\mapob.
\end{moneq}

For the first condition we have a small problem. 
The Berezin-Fourier modes $(k,\kappa)$ are described as functions on $Y$ satisfying the condition
\begin{moneq}
\psi\bigl(y\cdot (t\,\textrm{mod}\,d,\tau)\bigr) = \eexp^{2\pi i tk/d+i\tau\kappa}\cdot \psi(y)
\mapob.
\end{moneq}
Now in the non-super case (without the coordinate $\tau$), we can choose any Fourier mode $k\neq0$ and decide to write $\hbar=d/(2k\pi )$ (this notation is inherited from physics). 
However, we can not choose a fixed value for the odd parameter $\kappa$, as then we will obtain results that are no longer in the category of \gmfd/s (see \cite{Tu04} for details); we thus are obliged to keep it as a free parameter. 
The upshot of this discussion is that I will interpret the choice of a Berezin-Fourier mode as looking at functions on $Y$ satisfying the condition
\begin{moneq}[firstGQcondition]
\psi\bigl(y\cdot (t\,\textrm{mod}\,d,\tau)\bigr) = \eexp^{i\tau\kappa+it/\hbar}\cdot \psi(y)
\mapob,
\end{moneq}
where $\kappa\in \CA_1$ is a free parameter.

\begin{remarque}{Remark}
In our context, with $Y=\Orbit_{\mu_o} \times (\CA_0/d\ZZ)\times \CA_1$, we can, without loss of generality, put $\hbar=1$ (meaning we see the even coordinate $t$ really as an angle), but we will not do so in order to maintain a minimum of equality between the even and odd coordinates on the structure group. 
\end{remarque}

\section{\label{SECT/HeisenbergTemplate}The classical Heisenberg group as template}

\SwitchHats

As our purpose is to mimic the standard decomposition of the classical Heisenberg group, we provide here a short version of what one does. 
The Heisenberg group $H$ is (isomorphic to) the subgroup of $\GL(3,\RR)$ consisting of matrices of the form
\begin{moneq}
\begin{pmatrix}
1 & a & c \\ 0 & 1 & b \\ 0 & 0 & 1 \end{pmatrix}
\mapob.
\end{moneq}
In terms of the three coordinates $(a,b,c)$ the group law is given by
\begin{moneq}
(\ah,\bh, \ch)\cdot(\adh,\bdh,\cdh) = (\ah+\adh, \bh+\bdh, \ch+\cdh+\ah \bdh)
\mapob.
\end{moneq}

\medskip
\noindent\textbf{A short preview of results.}
The left-regular representation of the Heisenberg group decomposes as a direct integral of a $2$-parameter family of (mostly) irreducible representations. This can be done in an infinite number of ways, of which we present the two \myquote{natural} ones (by physicists these are called the position and momentum representations). When one of the parameters in the family is zero, the corresponding representation decomposes as a direct integral of a $1$-parameter family of $1$-dimensional representations. Comparison with coadjoint orbits tells us that the $2$-parameter family corresponds to $2$-dimensional orbits and the $1$-dimensional representations to point orbits. Moreover, the two natural ways to decompose the regular representations correspond to two different polarizations (the vertical and horizontal) of the coadjoint orbit.

\subsection{Decomposition of the left-regular representation}

\vrule width0pt
\vskip0.5\baselineskip

The (left) regular representation $\rho$ of $H$ is defined, for any $\gh=(\ah,\bh,\ch)\in H$ on the Hilbert space $L^2(\RR^3)$ by
\begin{align*}
\forall \psi\in L^2(\RR^3)
\quad:\quad
(\rho_\grel \psi)(\adh,\bdh,\cdh) 
&
= 
\psi\bigl((\ah,\bh,\ch)\mo\cdot(\adh,\bdh,\cdh)\bigr)
\\&
=
\psi(\adh-\ah,\bdh-\bh,\cdh-\ch+\ah\bh-\ah \bdh)
\mapob.
\end{align*}
This is a reducible unitary representation.  
It is not hard to see (with complete disregard for the square-integrability condition) that functions of the form
\begin{moneq}
\psi(\adh,\bdh,\cdh) = \chi(\bdh)\,\eexp^{i(\adh k+\cdh m)}
\end{moneq}
form an invariant subspace as for such functions we have
\begin{align*}
(\rho_\grel \psi)(\adh,\bdh,\cdh) 
&
= 
\psi(\adh-\ah,\bdh-\bh,\cdh-\ch+\ah\bh-\ah \bdh)
\dontprintcomp{
\\&
=
\chi(\bdh-\bh)\,\eexp^{i(\adh-\ah)k + i(\cdh-\ch+\ah\bh-\ah \bdh)m}
}%
\\&
=
\chi(\bdh-\bh)\,\eexp^{- i\ah \bdh m}\,\eexp^{-i(\ah k +(\ch-\ah\bh)m)}\,\eexp^{i(\adh k + \cdh m)}
\mapob.
\end{align*}
We \myquote{thus} obtain a $2$-parameter family of unitary representations $\rho^{k,m}$ of $H$ on $L^2(\RR)$ given by
\begin{moneq}
(\rho^{k,m}_\grel \chi)(\bdh) = \eexp^{- i\ah \bdh m}\,\eexp^{-i(\ah k +(\ch-\ah\bh)m)}\,\chi(\bdh-\bh)
\mapob.
\end{moneq}
And indeed, the partial Fourier transform $\Fourier_2:L^2(\RR^3)\to L^2\bigl(\RR^2, L^2(\RR)\bigr)\cong L^2(\RR^3)$ given by
\begin{moneq}
\bigl((\Fourier_2\psi)(k,m)\bigr)(\bdh)
=
\frac1{2\pi}\int\extder \adh\,\extder\cdh\extdersep 
\eexp^{-i(\adh k+\cdh m)}\psi(\adh,\bdh,\cdh)
\end{moneq}
intertwines the left-regular representation $\rho$ with the direct integral over $(k,m)$ of the $2$-parameter family $\rho^{k,m}$ of representations:
\KorteOfLangeTekst{
\begin{align*}
\biggl(\Bigl(\Fourier_2\bigl(\rho_\grel \psi\bigr)\Bigr)(k,m)\biggr)(\bdh)
=
\Bigl(\rho^{k,m}_\grel \bigl((\Fourier_2\psi)(k,m)\bigr)\Bigr)(\bdh)
\mapob.
\end{align*}
}{
\begin{align*}
\shifttag{6em}
\biggl(\Bigl(\Fourier_2\bigl(\rho_\grel \psi\bigr)\Bigr)(k,m)\biggr)(\bdh)
=
\frac1{2\pi}\, \int \extder \adh\,\extder \cdh\extdersep   \eexp^{-i(\adh k+\cdh m)}\,(\rho_\grel \psi)(\adh,\bdh,\cdh)
\\&
=
\frac1{2\pi}\, \int \extder \adh\,\extder \cdh\extdersep   \eexp^{-i(\adh k+\cdh m)}\,\psi(\adh-\ah,\bdh-\bh,\cdh-\ch+\ah\bh-\ah \bdh)
\\
\intertext{\footnotesize change of variables $\bar a = \adh-\ah$ and $\bar c= \cdh-\ch+\ah\bh-\ah \bdh$}
&
=
\frac1{2\pi}\, \int \extder \bar a\,\extder \bar c\extdersep   
\eexp^{-i((\bar a+\ah)k+(\bar c+\ch-\ah\bh+\ah \bdh)m)}
\,\psi(\bar a,\bdh-\bh,\bar c)
\\&
=
\eexp^{-i(\ah k+(\ch-\ah\bh+\ah \bdh)m)}
\,\bigl((\Fourier_2\psi)(k,m)\bigr)(\bdh-\bh)
\\&
=
\Bigl(\rho^{k,m}_\grel \bigl((\Fourier_2\psi)(k,m)\bigr)\Bigr)(\bdh)
\mapob.
\end{align*}}

But is also elementary to see that functions of the form
\begin{moneq}
\psi(\adh,\bdh,\cdh) = \chi(\adh)\,\eexp^{i\bdh\ell +i(\cdh-\adh \bdh)m}
\end{moneq}
form an invariant subspace as for such functions we have
\begin{align*}
(\rho_\grel \psi)(\adh,\bdh,\cdh) 
&
= 
\psi(\adh-\ah,\bdh-\bh,\cdh-\ch+\ah\bh-\ah \bdh)
\dontprintcomp{
\\&
=
\chi(\adh-\ah)\,\eexp^{i(\bdh-\bh)\ell + i(\cdh-\ch+\ah\bh-\ah \bdh-(\adh-\ah)(\bdh-\bh))m}
}%
\\&
=
\Bigl(\chi(\adh-\ah)\,\eexp^{-i\bh\ell + i(\adh\bh-\ch)m}\Bigr)\,\eexp^{i\bdh\ell+i(\cdh-\adh\bdh)m}
\mapob.
\end{align*}
This gives us another $2$-parameter family of unitary representations $\rhoh^{\ell,m}$ of $H$ on the space $L^2(\RR)$ given by
\begin{moneq}
(\rhoh^{\ell,m}_\grel \chi)(\adh) = \eexp^{-i\bh\ell + i(\adh\bh-\ch)m}\,\chi(\adh-\ah)
\mapob.
\end{moneq}
When we perform the (slightly modified) partial Fourier transform $\Fourier'_2:L^2(\RR^3)\to L^2\bigl(\RR^2;L^2(\RR)\bigr)\cong L^2(\RR^3)$ given by
\begin{moneq}
\bigl((\Fourier'_2\psi)(\ell,m)\bigr)(\adh) = \frac1{2\pi}\, \int_{\RR^2} \extder \bdh\,\extder \cdh\extdersep   
\eexp^{-i\bdh\ell-i(\cdh-\adh \bdh)m}\,\psi(\adh,\bdh,\cdh)
\mapob,
\end{moneq}
we find
\KorteOfLangeTekst{
\begin{align*}
\Bigl((\Fourier'_2\bigl(\rho_\grel \psi)\bigr)(\ell,m)\Bigr)(\adh)
=
\Bigl(\rhoh^{\ell,m}_\grel \bigl((\Fourier'_2\psi)(\ell,m)\bigr)\Bigr)(\adh)
\mapob,
\end{align*}
}{
\begin{align*}
\shifttag{6em}
\Bigl((\Fourier'_2\bigl(\rho_\grel \psi)\bigr)(\ell,m)\Bigr)(\adh)
=
\frac1{2\pi}\, \int \extder \bdh\,\extder \cdh\extdersep   \eexp^{-i\bdh\ell-i(\cdh-\adh \bdh)m}\,(\rho_\grel \psi)(\adh,\bdh,\cdh)
\\&
=
\frac1{2\pi}\, \int \extder \bdh\,\extder \cdh\extdersep   \eexp^{-i\bdh\ell-i(\cdh-\adh \bdh)m}\,\psi(\adh-\ah,\bdh-\bh,\cdh-\ch+\ah\bh-\ah \bdh)
\\
\intertext{\footnotesize change of variables $\bar b = \bdh-\bh$ and $\bar c= \cdh-\ch+\ah\bh-\ah \bdh$}
&
=
\frac1{2\pi}\, \int \extder \bar b\,\extder \bar c\extdersep   \eexp^{-i\bar b\ell-i(\bar c-(\adh-\ah)\bar b)m}\,\eexp^{-i\bh\ell-i(\ch-\adh\bh)m}\,\psi(\adh-\ah,\bar b,\bar c)
\\&
=
\eexp^{-i\bh\ell-i(\ch-\adh\bh)m}\,\bigl((\Fourier'_2\psi)(\ell,m)\bigr)(\adh-\ah)
\\&
=
\Bigl(\rhoh^{\ell,m}_\grel \bigl((\Fourier'_2\psi)(\ell,m)\bigr)\Bigr)(\adh)
\mapob,
\end{align*}}%
showing that $\Fourier'_2$ intertwines the left-regular representation of $H$ on $L^2(\RR^3)$ with the direct integral over $\ell,m$ of the $2$-parameter family of representations $\rhoh^{\ell,m}$ of $H$ on $L^2(\RR)$. 

Now for fixed $m\neq0$ all representations $\rhoh^{\ell,m}$ are irreducible (and unitarily equivalent to the standard irreducible representation $\rhoh^{0,m}$ of $H$ by the unitary equivalence $T_{\ell/m}:L^2(\RR)\to L^2(\RR)$ defined as $(T_{\ell/m}\psi)(\adh) = \psi(\adh-\ell/m)$). 
However, for $m=0$ the representations $\rhoh^{\ell,0}$ are reducible: one easily verifies that functions of the form
\begin{moneq}
\chi(\adh) = C\,\eexp^{i\adh k}
\ ,\qquad C\in \CC
\end{moneq}
form an invariant subspace as for such functions we have
\begin{align*}
(\rhoh^{\ell,0}_\grel \chi)(\adh) 
&
=
\eexp^{-i\bh\ell}\,\chi(\adh-\ah)
=
\eexp^{-i\bh\ell}\,C\,\eexp^{i(\adh-\ah)k}
=
\eexp^{- i\ah k -i\bh\ell }\,\chi(\adh)
\mapob.
\end{align*}
This gives us, for a fixed value of $\ell$, a $1$-parameter family of unitary (irreducible) representations $\rhob^{k,(\ell,0)}$ of $H$ on $\CC$ given  by
\begin{moneq}
\rhob^{k,(\ell,0)}_\grel C = \eexp^{-ik\ah-ik\bh}\,C
\mapob.
\end{moneq}
When we perform the Fourier transform $\Fourier_1:L^2(\RR)\to L^2(\RR)$ given by
\begin{moneq}
(\Fourier_1\chi)(k) = (2\pi)^{-1/2}\, \int_{\RR} \extder \adh\extdersep   \eexp^{-i\adh k}\,\chi(\adh)
\mapob,
\end{moneq}
we find
\KorteOfLangeTekst{
\begin{align*}
\bigl(\Fourier_1(\rhoh^{\ell,0}_\grel \chi)\bigr)(k)
=
\rhob^{k,(\ell,0)}_\grel \bigl((\Fourier_1\chi)(k)\bigr)
\mapob,
\end{align*}
}{
\begin{align*}
\bigl(\Fourier_1(\rhoh^{\ell,0}_\grel \chi)\bigr)(k)
&
=
\int \extder \adh\extdersep   \eexp^{-i\adh k}\,\rhoh^{\ell,0}_\grel \chi(\adh)
=
\int \extder a\extdersep   \eexp^{-i\adh k}\,\eexp^{-i\bh\ell}\,\chi(\adh-\ah)
\\
\intertext{\footnotesize change of variables $\bar a = \adh-\ah$}
&
=
\int \extder \bar a\extdersep   \eexp^{-i\bh\ell-i\ah k}\,\eexp^{-i\adh k}\,\chi(\adh)
\\&
=
\eexp^{-i\ah k-i\bh\ell}\,(\Fourier_1\chi)(k)
=
\rhob^{k,(\ell,0)}_\grel \bigl((\Fourier_1\chi)(k)\bigr)
\mapob,
\end{align*}}%
showing that $\Fourier_1$ intertwines the representation $\rhoh^{\ell,0}$ on $L^2(\RR)$ with the direct integral over $k$ of the $1$-parameter family of representations $\rhob^{k,(\ell,0)}$ on $\CC$.

\subsection{Representations associated to coadjoint orbits}

\BackSwitchHats 
\renewcommand\grel{{(a,b,c)}}

\vrule width0pt
\vskip0.5\baselineskip

Let us now compute the representations obtained by geometric quantization of coadjoint orbits. 
Using the standard basis $e_1=\partial_a, e_2=\partial_b, e_3=\partial_c$ of the Lie algebra $\Liealg g = T_eH$ with its dual basis $\LD[1]e_1, \LD[1]e_2, \LD[1]e_3$ of $\LD[1]{\Liealg g}$, we obtain the Coadjoint action in terms of the coordinates $(x,y,z)$ with respect to this dual basis as
\begin{align*}
\Coad(a,b,c)(x,y,z)
=
(x+zb, y-za,  z)
\mapob.
\end{align*}
The fundamental vector fields $e_i^{\LD[1]{\Liealg g}}$ associated to the Coadjoint action are given by
\begin{moneq}
e_1^{\LD[1]{\Liealg g}} =
z \partial_y
\quad,\quad
e_2^{\LD[1]{\Liealg g}} =
-z\partial_x
\quad,\quad
e_3^{\LD[1]{\Liealg g}} = 0
\mapob.
\end{moneq}
When looking at the orbit through a point $(x_o,y_o,0)$, we obtain the point itself, and when starting at a point $(x_o,y_o,z_o)$ with $z_o\neq0$, we obtain the set $\RR^2\times \{z_o\}$ as orbit. 

\itemizedef
In the first case we have $\Gamma = \extder t$ and the infinitesimal version $e_i^Y=c_i\partial_t$ of the lifted action to the bundle $Y=\{(x_o,y_o,0)\}\times (\RR/d\ZZ)\to \{(x_o,y_o,0)\}$ is given by the condition \recalf{LiftingEquationsOfGroupActionOnPrequantumBundle}
\begin{moneq}
\contrf{c_i\partial_t\restricted_{(x_o,y_o,0,t)}}\extder t = \contrs{e_i}{(x_o,y_o,0)}
\qquad\text{and}\qquad
\Lieder(c_i\partial_t)\extder t = 0
\mapob,
\end{moneq}
whose solutions are given by
\begin{moneq}
e_1^Y = x_o\partial_t
\ ,\qquad
e_2^Y = y_o\partial_t
\ ,\qquad
e_3^Y = 0
\mapob.
\end{moneq}
This integrates to the action $\phi$ of $H$ on $Y$ as
\begin{moneq}
\phi_\grel:(x_o,y_o,0,t)\mapsto (x_o, y_o, 0, t-ax_o-by_o)
\mapob.
\end{moneq}
Functions on $Y$ that are equivariant with respect to the action of the structure group are of the form $\psi(x_o,y_o,0,t) = C\cdot \eexp^{it/\hbar}$ for $C\in \CC$ and thus the representation $\Phi^{(x_o,y_o)}$ associated to such a point orbit (there is only the trivial polarization) is given by
\begin{align*}
\bigl(\Phi^{(x_o,y_o)}_\grel\psi\bigr)(x_o,y_o,0,t) 
&
= 
\psi\bigl(\phi_\grel\mo(x_o,y_o,0,t)\bigr)
\\&
=
\psi(x_o,y_o,0,t+ax_o+by_o)
\\&
=
C\cdot \eexp^{i(t+ax_o+by_o)/\hbar}
=
\eexp^{i(ax_o+by_o)/\hbar}\cdot \psi(x_o,y_o,0,t)
\mapob.
\end{align*}
Comparing this with the family $\rhob^{k,\ell,0}$ of representations of $H$, we see that they are \myquote{the same,} provided we make the identifications $-x_o/\hbar = k$ and $-y_o/\hbar = \ell$. 

\itemizedef
For the orbit $\RR^2\times \{z_o\}$ with $z_o\neq 0$ the Kirillov-Kostant-Souriau symplectic form $\omega$ is defined by \recalf{DefiningFormulaForKKSForm}, which reduces to 
\begin{moneq}
-\contrf{e_1^{\LD[1]{\Liealg g}},e_2^{\LD[1]{\Liealg g}}}{\omega_{(x,y,z)}} = 
\contrf{[e_1,e_2]}{(x,y,z)}
\mapob.
\end{moneq}
This gives $\omega = \extder x\wedge \extder y/z_o$ when using the coordinates $(x,y)$ for the point $(x,y,z_o)$ of the orbit. 
We now choose (any other choice will yield an equivalent description) the symplectic potential $\Theta = x\,\extder y/z_o$ and thus the connection $1$-form $\Gamma = \extder t + x\,\extder y/z_o$. 
In this case the infinitesimal version of the lifted action of $H$ on $Y$ is determined by $e_i^Y = e_i^{\LD[1]{\Liealg g}} + d_i\,\partial_t$ satisfying \recalf{LiftingEquationsOfGroupActionOnPrequantumBundle}, giving the solutions
\begin{moneq}
e_1^Y = z_0\partial_y 
\ ,\qquad
e_2^Y = -z_0\partial_x + y\partial_t
\ ,\qquad
e_3^Y = z_0\partial_t
\mapob.
\end{moneq}
This integrates to the action $\phi$ of $H$ on $Y$ given by
\begin{moneq}
\phi_\grel:(x,y,z_o,t)\mapsto (x+z_ob, y-az_o, z_o,t-yb - z_oc +z_o ab)
\mapob.
\end{moneq}
For this orbit the stabilizer subalgebra is generated by $e_3$, and in this (completely even) case, all $2$-dimensional subspaces containing $e_3$ will be polarizations. 
According to theory all will give equivalent representations, but we will concentrate on two: $\Liealg h = \langle e_3,e_1\rangle$ and $\Liealg h = \langle e_3,e_2\rangle$. 
The horizontally lifted vector fields $e_i^h$ are determined by the conditions \recalf{DefinitionOfHorizontalLiftFundamentalVectorFields} and are given by
\begin{moneq}
e_1^h = z_0\partial_y - x\partial_t
\ ,\qquad
e_2^h = -z_0\partial_x
\ ,\qquad
e_3^h = 0
\mapob.
\end{moneq}
---
The representation space associated to the polarization $\langle e_1,e_3\rangle$ consists of functions $\psi$ on $Y$ of the form $\psi(x,y,z_o,t)=\eexp^{it/\hbar}\cdot \psi(x,y,z_o,0)$ satisfying $b^h\psi=0$ for $b=e_1,e_3$. 
This yields the space of functions of the form
\begin{moneq}
\psi(x,y,z_o,t) = \eexp^{it/\hbar}\cdot \eexp^{ix y/(\hbar z_0)}\cdot \chi(x)
\mapob,
\end{moneq}
with $\chi$ a function of a single variable. 
The actual representation $\Phi^{z_o}$ then is given by
\begin{align*}
\bigl(\Phi^{z_o}_\grel\psi\bigr)(x,y,z_o,t)
&
=
\psi\bigl(\phi_\grel\mo(x,y,z_o,t)\bigr)
\\&
=
\psi(x-z_ob,y+z_oa,z_o,t+yb+z_oc)
\\&
=
\eexp^{it/\hbar}\cdot \eexp^{ixy/(\hbar z_0)} \cdot \eexp^{ ixa/\hbar}\cdot \eexp^{  iz_0(c-ab)/\hbar}\cdot  \chi(x-z_0b)
\mapob,
\end{align*}
which yields in terms of the function $\chi$ of a single variable the representation
\begin{moneq}
\bigl(\Phi^{z_o}_\grel \chi\bigr)(x) = \eexp^{ ixa/\hbar}\cdot \eexp^{  iz_0(c-ab)/\hbar}\cdot  \chi(x-z_0b)
\mapob.
\end{moneq}
When we perform a simple change of the coordinate $x = x_o+z_o \xh $ and writing $\hat \chi(\xh) = \chi(x)$, we find in terms of the functions $\hat \chi$ the representation
\begin{align*}
\bigl(\Phi^{z_o}_\grel\hat \chi\bigr)(\xh) 
&
= 
\eexp^{ i(x_o+z_o\xh)a/\hbar}\, \eexp^{  iz_0(c-ab)/\hbar}\,  \chi(x_o+ z_o\xh-z_0b)
\\&
=
\eexp^{ i\xh az_o/\hbar}\cdot\eexp^{i x_oa/\hbar}\, \eexp^{  iz_0(c-ab)/\hbar}\,  \hat \chi(\xh - b)
\mapob.
\end{align*}
This is exactly the representation $\rho^{k,m}$ above, provided we make the identifications $-z_o/\hbar = m$ and $-x_o/\hbar = k$. 

\noindent---
We next turn our attention to the representation space associated to the polarization $\langle e_2,e_3\rangle$ consisting of functions $\psi$ on $Y$ of the form $\psi(x,y,z_o,t)=\eexp^{it/\hbar}\cdot \psi(x,y,z_o,0)$ satisfying $b^h\psi=0$ with $b=e_2,e_3$. 
This yields the space of functions of the form
\begin{moneq}
\psi(x,y,z_o,t) = \eexp^{it/\hbar}\, f(y)
\mapob,
\end{moneq}
with $\chi$ a function of a single variable. 
The actual representation $\Phih^{z_o}$ then is given by
\begin{align*}
\bigl(\Phih^{z_o}_\grel\psi\bigr)(x,y,z_o,t)
&
=
\psi\bigl(\phi(a,b,c)\mo(x,y,z_o,t)\bigr)
\\&
=
\psi(x-z_ob,y+z_oa,z_o,t+yb+z_oc)
\\&
=
\eexp^{it/\hbar}\, \eexp^{ iyb/\hbar}\, \eexp^{  iz_0c/\hbar}\cdot  \chi(y+z_0a)
\mapob,
\end{align*}
which yields in terms of the function $\chi$ of a single variable the representation
\begin{moneq}
\bigl(\Phih^{z_o}_\grel \chi\bigr)(x) = \eexp^{ iyb/\hbar}\, \eexp^{  iz_0c/\hbar}\,  \chi(y+z_0a)
\mapob.
\end{moneq}
When we perform a simple change of the coordinate $y =y_o- z_o \yh$ and writing $\hat \chi(\yh) = \chi(y)$, we find in terms of the functions $\hat \chi$ the representation
\begin{align*}
\bigl(\Phih^{z_o}_\grel\hat \chi\bigr)(\yh) 
&
= 
\eexp^{ i(y_o-z_o\yh)b/\hbar}\, \eexp^{  iz_0c/\hbar}\,  \chi(y_o-z_o\yh+z_0a)
\\&
=
\eexp^{ i\yh bz_o/\hbar}\cdot\eexp^{i by_o/\hbar}\, \eexp^{  iz_0c/\hbar}\,  \hat \chi(\yh - a)
\mapob.
\end{align*}
And this is exactly the representation $\rhoh^{\ell,m}$ above, provided we make the identifications $-z_o/\hbar = m$ (as before) and $-y_0/\hbar = \ell$.

\begin{remarque}{Concluding remarks}
\itemizedef[*]
We see that the $1$-dimensional representations obtained in the decomposition of the regular representations are all linked to point orbits and that the two initial ways to decompose it, yielding the representations $\rho^{k,m}$ and $\rhoh^{\ell,m}$, are related to $2$-dimensional orbits using two different polarizations.

\itemizedef
It is known that different polarizations (at least in simple cases) always yield equivalent representations, which is corroborated by the fact that the two representations $\rho^{k,m}$ and $\rhoh^{\ell,m}$ are indeed equivalent (via a Fourier transform). 
The link with the representation associated to the $2$-dimensional orbit also shows (a well known fact) that the representations $\rho^{k,m}$ for various $k$ are all equivalent, as are the representations $\rhoh^{\ell,m}$ for various $\ell$. 
\end{remarque}

\section{\label{SECTSuperHeisenbergGroup}A super version of the Heisenberg group}

\SwitchHats
\renewcommand\grel{{(\ah,\betah,\gammah)}}

In this section I consider a super version of the classical Heisenberg group by taking the variables $b$ and $c$ to be odd instead of even. 
More precisely, I consider the \glgrp/ $\SuperH$ of graded dimension $1\vert2$, a subgroup of $\GL(2\vert1,\CA)$ consisting of matrices of the form
\begin{moneq}
\begin{pmatrix}
1 & a\rlap{$\smash{\kern0.4em\vrule height1.8ex depth 6ex}$} & \gamma \\ 0 & 1 & \beta \\ 
\hline 0 & 0 & 1 \end{pmatrix}
\ ,\qquad a\in \CA_0\ \&\ \beta,\gamma\in \CA_1
\mapob.
\end{moneq}
As \agmfd/ we have $\SuperH = \RR^{1\vert2}$ and in terms of the three coordinates $(a,\beta,\gamma)$ the group law is given by
\begin{moneq}
(\ah,\betah, \gammah)\cdot(\adh,\betadh,\gammadh) = (\ah+\adh, \betah+\betadh, \gammah+\gammadh+\ah \betadh)
\mapob.
\end{moneq}
In terms of Harish-Chandra pairs this is the pair $(\RR,\Liealg h)$ in which the (underlying real) Lie group is the additive real line $\RR$ and $\Liealg h$ the super Lie algebra of dimension $1\vert2$ with homogeneous basis $e_1,f_2,f_3$ ($e_1$ even and $f_i$ odd) with (only non-zero) commutator $[e_1,f_2]=f_3$.

\bigskip
\noindent\textbf{A short preview of results.}
The left-regular representation of the super Heisenberg group decomposes as a direct integral of a $2$-parameter family of (mostly) irreducible representations. This can be done in two ways (we no longer have an infinity of possibilities due to grading). When one of the parameters in the family is zero, the corresponding representation decomposes as a direct integral of a $1$-parameter family of $1$-dimensional representations. Comparison with coadjoint orbits tells us that the $2$-parameter family corresponds to $2$-dimensional orbits and the $1$-dimensional representations to point orbits. Moreover, the two ways to decompose the regular representations correspond to two different polarizations of the coadjoint orbit. So far it seems that the results are (nearly) completely analogous to the non-super case. However, a closer look reveals some important differences. Due to grading, the two different ways to decompose the regular representation are not equivalent, despite being linked to the same coadjoint orbit: the two different polarizations have a different graded dimension.

\subsection{Decomposition of the left-regular representation}
\label{SubSECTDecompSuperHeisenberg}

\vrule width0pt
\vskip0.5\baselineskip

My purpose now is to copy the computations done for the classical Heisenberg group, simply replacing the even variables\slash coordinates $b$ and $c$ by the odd ones $\beta$ and $\gamma$. 
Of course there will be small changes, in particular that I will define the (left) regular representation $\rho$ of $\SuperH$ on the space $\Schwartz12 \subset C^\infty(\RR^{1\vert2})\cap L^2(\RR^{1\vert2})$. 
This choice is inspired by the fact that the Schwartz space $\Schwartz{}{}$ is mapped by the (ordinary) Fourier transform on $L^2(\RR)$ to itself and that it is closed under differentiation and multiplication by a polynomial. 
On this space the left-regular representation $\rho$ of $\SuperH$ is defined by
\begin{align*}
\forall \psi\in \Schwartz12
\ :\ 
(\rho_\grel \psi)(\adh,\betadh,\gammadh) 
&
= 
\psi\bigl((\ah,\betah,\gammah)\mo\cdot(\adh,\betadh,\gammadh)\bigr)
\\&
=
\psi(\adh-\ah,\betadh-\betah,\gammadh-\gammah+\ah\betah-\ah \betadh)
\mapob.
\end{align*}
It is not hard to see that functions of the form
\begin{moneq}
\psi(\adh,\betadh,\gammadh) = \chi(\betadh)\,\eexp^{i(\adh k+\gammadh\mu)}
\end{moneq}
with $k\in \CA_0$ even, $\mu\in \CA_1$ odd and $\chi$ a function of a single odd variable, form an invariant subspace as for such functions we have
\begin{align*}
(\rho_\grel \psi)(\adh,\betadh,\gammadh) 
&
= 
\psi(\adh-\ah,\betadh-\betah,\gammadh-\gammah+\ah\betah-\ah \betadh)
\dontprintcomp{
\\&
=
\chi(\betadh-\betah)\,\eexp^{i(\adh-\ah)k + i(\gammadh-\gammah+\ah\betah-\ah \beta)\mu}
}
\\&
=
\bigl(\,\chi(\betadh-\betah)\,\eexp^{- i\ah \betadh\mu}\,\eexp^{-i(\ah k +(\gammah-\ah\betah)\mu)}\,\bigr)\,\eexp^{i(\adh k + \gammadh\mu)}
\mapob.
\end{align*}
In this way we obtain a $1\vert1$-parameter family of representations $\rho^{k,\mu}$ of $\SuperH$ on the super Hilbert space $L^2(\RR^{0\vert1})\equiv C^\infty(\CA_1)\cong \CC^2$ of smooth functions of a single odd variable given by
\begin{moneq}
\bigl(\rho^{k,\mu}_\grel \chi\bigr)(\betadh) = \eexp^{- i\ah \betadh\mu}\,\eexp^{-i(\ah k +(\gammah-\ah\betah)\mu)}\,\chi(\betadh-\betah)
\mapob.
\end{moneq}
To show that $\rho$ decomposes as the direct integral of this family $\rho^{k,\mu}$, we consider the modified partial Berezin-Fourier transform $\Fourieroddp_{1\vert1}:\Schwartz12\to \Schwartz12$ given by $\Fourieroddp_{1\vert1} = \Fourierodd_{1\vert1}\scirc \Identif_\sigma$ with $\Fourierodd_{1\vert1}$ as in \recalt{PartialBerezeinFourierTransformIsIsomorphismSHS} and $\Identif_\sigma$ as in \recalt{OddPermutationsChangeSignOfSuperScalarProductsNEW} with $\sigma = \bigl( \begin{smallmatrix} 1&2 \\ 2&1 \end{smallmatrix}\bigr)$ the (unique) transposition of $\{1,2\}$, \ie, given by the formula 
\begin{moneq}
\bigl((\Fourieroddp_{1\vert1}\psi)(k,\mu)\bigr)(\betadh)
=
(2\pi)^{-1/2}\,\int\extder \adh \,\extder\gammadh\extdersep 
\eexp^{-i(\adh k+\gammadh\mu)}\psi(\adh,\betadh,\gammadh)
\mapob.
\end{moneq}
And indeed, $\Fourieroddp_{1\vert1}$ intertwines the left-regular representation $\rho$ with the direct integral over $(k,\mu)$ of the $1\vert1$-parameter family $\rho^{k,\mu}$ of representations:
\KorteOfLangeTekst{
\begin{align*}
\Bigl(\hBigl(\Fourieroddp_{1\vert1}\bigl(\rho_\grel \psi\bigr)\hBigr)(k,\mu)\Bigr)(\betadh)
=
\hBigl(\rho^{k,\mu}_\grel \bigl((\Fourieroddp_{1\vert1}\psi)(k,\mu)\bigr)\hBigr)(\betadh)
\mapob.
\end{align*}
}{
\begin{align*}
\shifttag{6em}
\Bigl(\hBigl(\Fourieroddp_{1\vert1}\bigl(\rho_\grel \psi\bigr)\hBigr)(k,\mu)\Bigr)(\betadh)
=
(2\pi)^{-1/2}\, \int \extder \adh\,\extder \gammadh\extdersep   \eexp^{-i(\adh k+\gammadh\mu)}\,(\rho_\grel \psi)(\adh,\betadh,\gammadh)
\\&
=
(2\pi)^{-1/2}\, \int \extder \adh\,\extder \gammadh\extdersep   \eexp^{-i(\adh k+\gammadh\mu)}\,\psi(\adh-\ah,\betadh-\betah,\gammadh-\gammah+\ah\betah-\ah \betadh)
\\
\intertext{\footnotesize change of variables $\bar a = \adh-\ah$ and $\bar \gamma= \gammadh-\gammah+\ah\betah-\ah \betadh$}
&
=
(2\pi)^{-1/2}\, \int \extder \bar a\,\extder \bar \gamma\extdersep   
\eexp^{-i((\bar a+\ah)k+(\bar \gamma+\gammah-\ah\betah+\ah \betadh)\mu)}
\,\psi(\bar a,\betadh-\betah,\bar \gamma)
\\&
=
\eexp^{-i(\ah k+(\gammah-\ah\betah+\ah \betadh)\mu)}
\,\bigl((\Fourier_2\psi)(k,\mu)\bigr)(\betadh-\betah)
\\&
=
\hBigl(\rho^{k,\mu}_\grel \bigl((\Fourieroddp_{1\vert1}\psi)(k,\mu)\bigr)\hBigr)(\betadh)
\mapob.
\end{align*}}

But is also elementary to see that functions of the form
\begin{moneq}
\psi(\adh,\betadh,\gammadh) = \chi(\adh)\,\eexp^{i\betadh\lambda +i(\gammadh-\adh\betadh)\mu}
\end{moneq}
with $\lambda$ and $\mu$ odd and $\chi$ a function of a single even variable, form an invariant subspace as for such functions we have
\begin{align*}
(\rho_\grel \psi)(\adh,\betadh,\gammadh) 
&
= 
\psi(\adh-\ah,\betadh-\betah,\gammadh-\gammah+\ah\betah-\ah \betadh)
\dontprintcomp{
\\&
=
\chi(\adh-\ah)\,\eexp^{i(\betadh-\betah)\lambda + i(\gammadh-\gammah+\ah\betah-\ah \betadh-(\adh-\ah)(\betadh-\betah))\mu}
}%
\\&
=
\Bigl(\chi(\adh-\ah)\,\eexp^{-i\betah\lambda + i(\adh\betah-\gammah)\mu}\Bigr)\,\eexp^{i\betadh\lambda+i(\gammadh-\adh\betadh)\mu}
\mapob.
\end{align*}
This gives us a $0\vert2$-parameter family of representations $\rhoh^{\lambda,\mu}$ of $\SuperH$ on the space $\Schwartz{}{}$ given by
\begin{moneq}
(\rhoh^{\lambda,\mu}_\grel \chi)(\adh) = \eexp^{-i\betah\lambda + i(\adh\betah-\gammah)\mu}\,\chi(\adh-\ah)
\mapob.
\end{moneq}
When we perform the modified partial Berezin-Fourier transform $\Fourieroddp_{0\vert2}:\Schwartz12 \to \Schwartz12$ given by (see \recalt{PartialBerezeinFourierTransformIsIsomorphismSHS}, but with a change of coordinates replacing $\gammadh$ by $\gammadh-a\betadh$ in the exponent)
\begin{moneq}
\bigl((\Fourieroddp_{0\vert2}\psi)(\lambda,\mu)\bigr)(\adh) = \int \extder \betadh\,\extder \gammadh\extdersep   
\eexp^{-i\beta\lambda-i(\gamma-\adh\betadh)\mu}\,\psi(\adh,\betadh,\gammadh)
\mapob,
\end{moneq}
we find
\KorteOfLangeTekst{
\begin{align*}
\hBigl((\Fourieroddp_{0\vert2}\bigl(\rho_\grel \psi)\bigr)(\lambda,\mu)\hBigr)(\adh)
=
\hBigl(\rhoh^{\lambda,\mu}_\grel \bigl((\Fourieroddp_{0\vert2}\psi)(\lambda,\mu)\bigr)\hBigr)(\adh)
\mapob,
\end{align*}
}{
\begin{align*}
\shifttag{6em}
\hBigl((\Fourieroddp_{0\vert2}\bigl(\rho_\grel \psi)\bigr)(\lambda,\mu)\hBigr)(\adh)
=
\int \extder \betadh\,\extder \gammadh\extdersep   \eexp^{-i\betadh\lambda-i(\gamma-\adh\betadh)\mu}\,(\rho_\grel \psi)(\adh,\betadh,\gammadh)
\\&
=
\int \extder \betadh\,\extder \gammadh\extdersep   \eexp^{-i\betadh\lambda-i(\gammadh-\adh\betadh)\mu}\,\psi(\adh-\ah,\betadh-\betah,\gammadh-\gammah+\ah\betah-\ah \betadh)
\\
\intertext{\footnotesize change of variables $\bar \beta = \betadh-\betah$ and $\bar \gamma= \gammadh-\gammah+\ah\betah-\ah \betadh$}
&
=
\int \extder \bar \beta\,\extder \bar \gamma\extdersep   \eexp^{-i\bar \beta\lambda-i(\bar \gamma-(\adh-\ah)\bar \beta)\mu}\,\eexp^{-i\betah\lambda-i(\gammah-\adh\betah)\mu}\,\psi(\adh-\ah,\bar \beta,\bar \gamma)
\\&
=
\eexp^{-i\betah\lambda-i(\gammah-\adh\betah)\mu}\,\bigl((\Fourieroddp_{0\vert2}\psi)(\lambda,\mu)\bigr)(\adh-\ah)
\\&
=
\hBigl(\rhoh^{\lambda,\mu}_\grel \bigl((\Fourieroddp_{0\vert2}\psi)(\lambda,\mu)\bigr)\hBigr)(\adh)
\mapob,
\end{align*}}%
showing that $\Fourieroddp_{0\vert2}$ intertwines the left-regular representation of $\SuperH$ on $\Schwartz12$ with the direct integral over $\lambda,\mu$ of the $0\vert2$-parameter family of representations $\rhoh^{\lambda,\mu}$ of $\SuperH$ on $\Schwartz{}{}$. 

It is easy to see that for $\mu=0$ the representations $\rhoh^{\lambda,0}$ are reducible: functions of the form
\begin{moneq}
\chi(\adh) = C\,\eexp^{i\adh k}
\ ,\qquad C\in \CC
\end{moneq}
form an invariant subspace as for such functions we have
\begin{align*}
(\rhoh^{\lambda,0}_\grel \chi)(\adh) 
&
=
\eexp^{-i\betah\lambda}\,\chi(\adh-\ah)
=
\eexp^{-i\betah\lambda}\,C\,\eexp^{i(\adh-\ah)k}
=
\eexp^{- i\ah k -i\betah\lambda }\,\chi(\adh)
\mapob.
\end{align*}
This gives us, for a fixed value of $\lambda$, a $1\vert0$-parameter family of unitary representations $\rhob^{k,(\lambda,0)}$ of $\SuperH$ on $\CC$ given  by
\begin{moneq}
\rhob^{k,(\lambda,0)}{(\ah,\betah,\gammah)}\,C = \eexp^{-i\ah k-i\betah\lambda}\,C
\mapob.
\end{moneq}
When we perform the Fourier transform $\Fourier_1:\Schwartz{}{}\to \Schwartz{}{}$ given by
\begin{moneq}
(\Fourier_1\chi)(k) = (2\pi)^{-1/2}\, \int \extder \adh\extdersep   \eexp^{-i\adh k}\,\chi(\adh)
\mapob,
\end{moneq}
we find
\KorteOfLangeTekst{
\begin{align*}
\bigl(\Fourier_1(\rhoh^{\lambda,0}_\grel \chi)\bigr)(k)
=
\rhob^{k,(\lambda,0)}_\grel \bigl((\Fourier_1\chi)(k)\bigr)
\mapob,
\end{align*}
}{
\begin{align*}
\bigl(\Fourier_1(\rhoh^{\lambda,0}_\grel \chi)\bigr)(k)
&
=
\int \extder \adh\extdersep   \eexp^{-i\adh k}\,\rhoh^{\lambda,0}_\grel \chi(\adh)
=
\int \extder \adh\extdersep   \eexp^{-i\adh k}\,\eexp^{-i\betah\lambda}\,\chi(\adh-\ah)
\\
\intertext{\footnotesize change of variables $\bar a = \adh-\ah$}
&
=
\int \extder \bar a\extdersep   \eexp^{-i\betah\lambda-i\ah k}\,\eexp^{-i\bar a k}\,\chi(\bar a)
\\&
=
\eexp^{-i\ah k-i\betah\lambda}\,(\Fourier_1\chi)(k)
=
\rhob^{k,(\lambda,0)}_\grel \bigl((\Fourier_1\chi)(k)\bigr)
\mapob,
\end{align*}}%
showing that $\Fourier_1$ intertwines the representation $\rhoh^{\lambda,0}$ on $\Schwartz{}{}$ with the direct integral over $k$ of the $1$-parameter family of representations $\rhob^{k,(\lambda,0)}$ on $\CC$.

\begin{remarque}{Additional details concerning these computations}
\itemizedef[*]
It is elementary to check by a direct computation that the representation $\rho$ preserves the \supscal product $\supinprsym_{1\vert2}$ on $\Schwartz12$, but it can also be deduced from the fact that the Berezinian of the coordinate transformation
\begin{moneq}
(\adh,\betadh,\gammadh) \mapsto
(\adh-\ah,\betadh-\betah,\gammadh-\gammah+\ah\betah-\ah \betadh)
\end{moneq}
is $1$ and thus that the integral over $\RR^{1\vert2}$ is invariant. 
On the other hand, the standard inner product $\inprodsym_{1\vert2}$ is not invariant under the action of $\rho(\body\SuperH)$: the restriction of $\rho$ to $\body \SuperH\cong\RR$ is given by
\begin{moneq}
\bigl(\rho_{(\ah,0,0)}\psi\bigr)(\adh,\betadh,\gammadh) = \psi(\adh-\ah, \betadh,\gammadh - \ah\betadh)
\mapob,
\end{moneq}
which does not preserve the standard inner product $\inprodsym_{1\vert2}$ \recalf{ScalarProductInFunctionHilbertSpace}. 
However, the action of $\rho(\body\SuperH)$ preserves the inner product $\inprodsym_\mfdmetric$ defined as
\begin{moneq}
\inprod{\chi}{\psi}_\mfdmetric
=
\inprodd{\chi_0}{\psi_0} + \inprodd{\chi_1+\phi\chi_2}{\psi_1+\phi\psi_2} + \inprodd{\chi_2}{\psi_2} + \inprodd{\chi_{12}}{\psi_{12}}
\mapob,
\end{moneq}
where $\inprodd\cdot\cdot$ is the standard inner product on $\Schwartz{}{} \subset L^2(\RR)$, and where $\phi:\RR\to \RR$ denotes the identity function $\phi(x) = x$. 
This modified inner product is the one obtained via a (left-invariant) metric on $\RR^{1\vert2}=\SuperH$ according to the procedure described in \cite[4.19, 4.22]{Tuynman..2018}. 
It is directly linked to the change of coordinates $\gammadh\mapsto \gammah = \gamma-a\beta$ in $\SuperH$ (this is the change of coordinates $\Phi(a,\beta,\gammah) = (a,0,0)\cdot \exp(\beta,\gammah)$ of \cite[(3.12)]{Tuynman..2018}). 
Note that we use here that $\Schwartz{}{}$ is closed under multiplication by the (degree $1$) polynomial $\phi$. 
It follows that $\rho$ is a \psur{} on the super inner product space $\bigl(\Schwartz12, \inprodsym_\mfdmetric, \supinprsym_{1\vert2}\bigr)$. 
Its completion is isomorphic to the super Hilbert space $L^2(\RR^{1\vert2})$, but the embedding\slash realisation of $\Schwartz12$ as a dense subspace is given by the map
\begin{multline*}
\qquad
(\psi_0, \psi_1,\psi_2,\psi_{12})\in \Schwartz{}{}^4\cong \bigl(\Schwartz12, \inprodsym_\mfdmetric, \supinprsym_{1\vert2})
\\
\mapsto
(\psi_0, \psi_1+\phi\psi_2, \psi_2,\psi_{12})\in L^2(\RR)^4\cong \bigl(L^2(\RR^{1\vert2}), \inprodsym_{1\vert2}, \supinprsym_{1\vert2}\bigr)
\mapob.
\qquad
\end{multline*}

\itemizedef
One verifies easily that the representations $\rho^{k,\mu}$ of $\SuperH$ on $L^2(\RR^{0\vert1})\equiv C^\infty_c(\RR^{0\vert1})\equiv \CC^2$ satisfy the conditions \refmetnaam{\labelFamilySUR}{FamilySUR1}-\ref{FamilySUR3.5} when we use the standard inner product $\inprodsym_{0\vert1}$ and the standard \supscal product $\supinprsym_{0\vert1}$ on $\Dense=L^2(\RR^{0\vert1})\cong\CC^2$. 
In order to show that this family also satisfies \refmetnaam{\labelFamilySUR}{FamilySUR4}, we proceed as follows. 
Any $\psi\in \Denset\equiv C^\infty_c(\RR^{1\vert1}; \Dense)$ is determined by two smooth compactly supported functions $\psi_i:\RR\to \Dense\cong\CC^2$ according to
\begin{moneq}
\forall (k,\mu)\in \RR^{1\vert1} 
:
\psi(k,\mu) = \psi_0(k)+\mu\,\psi_1(k)
\mapob.
\end{moneq}
And for any $k\in\RR$ an element $\psi_i(k)\in \Dense$ is determined by two (complex) numbers $\psi_{ij}(k)$ according to
\begin{moneq}
\bigl(\psi_i(k)\bigr)(\betadh) = \psi_{i0}(k) + \betadh\,\psi_{i1}(k)
\mapob.
\end{moneq}
It then is elementary to show that $\psi_i:\RR\to \Dense$ is smooth if and only if the functions $\psi_{ij}:\RR\to \CC$ are smooth. 
Any element $\psi\in \Denset$ thus is determined by four functions $\psi_{ij}\in C^\infty_c(\RR)$ according to
\begin{moneq}
\bigl(\psi(k,\mu)\bigr)(\betadh) = \psi_{00}(k) + \betadh\,\psi_{01}(k) + \mu\,\psi_{10}(k) + \mu\betadh\,\psi_{11}(k)
\mapob.
\end{moneq}
What we have done here is essentially proving that the isomorphism $L^2(\RR^{1\vert2})\to L^2\bigl(\RR^{1\vert1}; L^2(\RR^{0\vert1})\bigr)$ of \recalt{DirectIntegralOfHilbertSpacesPotentiallyChangeSignSuperScalarProductsNEW} restricted to $C^\infty_c(\RR^{1\vert2})$ has $\Denset\equiv C^\infty_c(\RR^{1\vert1}; \Dense)$ as image. 
Writing $\psi(k,\mu,\betadh)$ for $\bigl(\psi(k,\mu)\bigr)(\betadh)$, the action of $\rhot$ on $\Denset$ is given by
\begin{align*}
\bigl(\rhot_\grel \psi\bigr)(k,\mu,\betadh)
&
\equiv\rho^{k,\mu}_\grel \psi(k,\mu,\betadh)
\\&
=
\eexp^{- i\ah \betadh\mu}\,\eexp^{-i(\ah k +(\gammah-\ah\betah)\mu)}\,\psi(k,\mu,\betadh-\betah)
\mapob.
\end{align*}
It then is straightforward to show that the inner product $\inprodsym_{\Denset}$ on $\Denset$ given by
\begin{moneq}[InnerProductOnDensetSuperHeisenberg]
\inprod{\psi}{\chi}_{\Denset}
=
\inprodd{\psi_{00}}{\chi_{00}} + \inprodd{\psi_{01}}{\chi_{01}} + \inprodd{\psi_{10}}{\chi_{10}} + \inprodd{\psi_{11}-\psi_{00}'}{\chi_{11}-\chi'_{00}}
\end{moneq}
is preserved by each $\rhot(\ah,0,0)$ with $\ah\in \RR$, as their action in terms of the four functions $\psi_{ij}$ is given by
\begin{moneq}
\rhot_{(\ah,0,0)}(\psi_{00},\psi_{01}, \psi_{10}, \psi_{11}) 
=
\eexp^{-i\ah k}\cdot(\psi_{00}, \psi_{01}, \psi_{10}, \psi_{11} - i\ah\,\psi_{00})
\mapob,
\end{moneq}
where the expression $\eexp^{-i\ah k}$ is abuse of notation for the function $k\mapsto \eexp^{-i\ah k}$. 
As we already know that $\supinprsym_{1\vert2}$ is preserved by $\rhot$ (it is given by Berezin integration over $\RR^{1\vert2}$), it follows easily that $(\Denset,\rhot)$ is a \psur{} as required for \refmetnaam{\labelFamilySUR}{FamilySUR4}.

\itemizedef
To prove that the modified Berezin-Fourier transform $\Fourieroddp_{1\vert1}$ is an equivalence of super unitary representations, we first note that the formulæ for $\rhot$ and $\inprodsym_{\Denset}$ in the identification $C^\infty_c(\RR^{1\vert1}; \Dense)\cong C^\infty_c(\RR^{1\vert2})$ extend naturally to the space $\Schwartz12$. 
It then follows from the given computations above that $\Fourieroddp_{1\vert1}$ intertwines the representations $\rho$ on the source with $\rhot$ (the direct integral of the family $\rho^{k,\mu}$) on the target. 
But according to \recalt{PartialBerezeinFourierTransformIsIsomorphismSHS} and \recalt{OddPermutationsChangeSignOfSuperScalarProductsNEW} it also intertwines the \supscal product $\supinprsym_{1\vert2}$ on source and target, proving that $\Fourieroddp_{1\vert1}$ is an equivalence of super unitary representations in the sense of \recalt{DefinitionEquivalenceSuperUnitaryRep}. 
But we can say even more: not only does $\smash{\Fourieroddp_{1\vert1}}$ intertwine the standard inner product $\inprodsym_{1\vert2}$ on source and target (as stated in \recalt{PartialBerezeinFourierTransformIsIsomorphismSHS} and \recalt{OddPermutationsChangeSignOfSuperScalarProductsNEW}, but something we are not interested in), it also intertwines the inner products $\inprodsym_\mfdmetric$ on the source and $\inprodsym_{\Denset}$ on the target, making it an odd $(-1)$-isomorphism of super inner product spaces (something not required of an equivalence of super unitary representations).

\itemizedef
One verifies easily that the representations $\rhoh^{\lambda,\mu}$ of $\SuperH$ on $\Schwartz{}{}$ satisfy the conditions \refmetnaam{\labelFamilySUR}{FamilySUR1}-\ref{FamilySUR3.5} when we use the standard inner product and \supscal product $\inprodsym_{1\vert0}=\supinprsym_{1\vert0}\equiv\inprodd\cdot\cdot$ on $\Dense=\Schwartz{}{}$. 
In order to show that this family also satisfies \refmetnaam{\labelFamilySUR}{FamilySUR4}, we note that, with arguments similar to those used for the family $\rho^{k,\mu}$, the isomorphism $L^2(\RR^{2\vert1})\to L^2\bigl(\RR^{0\vert2}; L^2(\RR^{1\vert0})\bigr)$ of \recalt{DirectIntegralOfHilbertSpacesPotentiallyChangeSignSuperScalarProductsNEW} maps $\Schwartz12$ to $\Denset \equiv C^\infty_c(\RR^{0\vert2}; \Dense) \equiv C^\infty_c\bigl(\RR^{0\vert2}; \Schwartz{}{}\bigr)$. 
Writing $\psi(\adh,\lambda,\mu)$ for $\bigl(\psi(\lambda,\mu)\bigr)(\adh)$, the action of $\rhot$ on $\Denset$ is given by
\begin{align*}
\bigl(\rhot_\grel \psi\bigr)(\adh,\lambda,\mu)
&
\equiv\rhoh^{\lambda,\mu}(\ah,\betah,\gammah)\psi(\adh,\lambda,\mu)
\\&
=
\eexp^{-i\betah\lambda + i(\adh\betah-\gammah)\mu}\,\psi(\adh-\ah, \lambda,\mu)
\mapob.
\end{align*}
It then is straightforward to show that the standard inner product $\inprodsym_{1\vert2}$ and the standard \supscal product $\supinprsym_{1\vert2}$ on $L^2(\RR^{1\vert2})$ restricted to $\Denset$ (in the identification $\Denset\cong\Schwartz12$) make $(\Denset,\rhot)$ a \psur{} as required for \refmetnaam{\labelFamilySUR}{FamilySUR4}.

\itemizedef
The modified partial Berezin-Fourier transform $\Fourierodd'_{0\vert2}:\Schwartz12 \to \Schwartz12$ is an even equivalence between the super unitary representation $\rho$ and the direct integral over $\lambda,\mu$ of the $0\vert2$-parameter family of representations $\rhoh^{\lambda,\mu}$. 
It is also an equivalence of super Hilbert spaces, as it does intertwine the inner products $\inprodsym_\mfdmetric$ (the inner product on the representation space of $\rho$) and $\inprodsym_{1\vert2}$ (the inner product on the direct integral of the representation spaces of $\rhoh^{\lambda,\mu}$), due to the modification in the Berezin-Fourier transform with the factor $\eexp^{-i(\gammadh-\adh\betadh)\mu}$ instead of $\eexp^{-i\gammadh\mu}$.

\itemizedef
As for the other families, it is easy to show that the family $\rhob^{k,(\lambda,0)}$ on $\CC$ satisfies the conditions \refmetnaam{\labelFamilySUR}{FamilySUR1}-\ref{FamilySUR3.5} when we equip $\CC$ with the standard inner product and \supscal product. 
In order to show that it satisfies \refmetnaam{\labelFamilySUR}{FamilySUR4}, we recall that $\rhot$ on $\Denset=C^\infty_c(\RR^{1\vert1}; \CC)$ is given by
\begin{moneq}
\bigl(\rhot_\grel \psi\bigr)(k,\lambda) = \eexp^{-i\ah k - i\betah\lambda}\,\psi(k,\lambda)
\mapob.
\end{moneq}
It then is (again) immediate that $\rhot$ preserves the standard inner product $\inprodsym_{1\vert1}$ and the standard \supscal product $\supinprsym_{1\vert1}$ on $\Denset$, making it into a \psur. 

As the expression for $\rhot$ extends naturally to $\Schwartz11$, it follows easily that $\Fourier_{1\vert0}$ is an even isomorphism of super Hilbert spaces that intertwines $\rhot$ with the direct integral (over $\lambda$) of the representations $\rhoh^{\lambda,0}$. 
\end{remarque}

\subsection{Irreducibility}
\label{SUBSECTIrreducibilityOfSuperH}

In order to prove that the \psur s $\rho^{k,\mu}$ and $\rhoh^{\lambda,\mu}$, $\mu\neq0$ of $\SuperH$ are irreducible, we start with the observation that invariance of a graded subspace $F\subset \Dense$ under the maps $\rho^{k,\mu}(g)$ implies invariance under the infinitesimal generators $\tau^{k,\mu}(X)\in \End(\Dense\otimes \CA^\CC)$, $X\in \Liealg h\cong T_e\SuperH$ defined as (see \recalf{MainSURlabelrhopsiissmooth})
\begin{moneq}
\tau^{k,\mu}(X)\psi = X\FgroupHilbert_\psi = \frac{\extder}{\extder t}\bigrestricted_{t=0}\ \rho^{k,\mu}(\eexp^{tX})\psi
\mapob.
\end{moneq}
We thus investigate invariance under these infinitesimal operators, where we can ignore those that are multiples of the identity, as these do not impose any restriction on the (supposedly invariant) subspace.

\itemizedef
The $1\vert1$-parameter family of \psur{} $\rho^{k,\mu}$ of $\SuperH$ on $L^2(\RR^{0\vert1})\equiv C^\infty(\CA_1)\cong \CC^2$ (with $\Dense=L^2(\RR^{0\vert1})\cong \CC^2$) is given by
\begin{moneq}
\bigl(\rho^{k,\mu}_\grel \psi\bigr)(\betadh) = \eexp^{- i\ah \betadh\mu}\,\eexp^{-i(\ah k +(\gammah-\ah\betah)\mu)}\,\psi(\betadh-\betah)
\mapob.
\end{moneq}
The infinitesimal actions that do not reduce to a multiple of the identity are
\begin{align*}
\tau^{k,\mu}(\partial_{\ah})\psi
&=
-i\betadh\mu\psi - ik\psi
&
\tau^{k,\mu}(\partial_{\betah})\psi
&=
-\partial_{\betadh}\psi
\end{align*}
If $F\subset \Dense\cong \CC^2$ is a closed invariant subspace that is invariant under $\rho^{k,\mu}(g)$ and if $\psi\in F$, then necessarily the above images should belong to $F\otimes \CA^\CC$, \ie, the elements
\begin{moneq}
\partial_{\betadh}\chi
\qquad\text{and}\qquad
\mu\betadh\chi
\end{moneq}
should belong to $F\otimes \CA^\CC$. 
When $\mu\neq0$, this implies that
\begin{moneq}
\partial_{\betadh}\chi
\qquad\text{and}\qquad
\betadh\chi
\end{moneq}
should belong to $F$. 
Now the only non-trivial graded subspaces of $\Dense\cong \CC\oplus\CC$ are $\CC\oplus\{0\}$ and $\{0\}\oplus\CC$. 
But neither can be invariant, as $\partial_\betadh(\{0\}\oplus\CC) = \CC\oplus\{0\}$ and $\betadh(\CC\oplus\{0\}) = \{0\}\oplus \CC$. 
It follows that the representation $\rho^{k,\mu}$ is irreducible for $\mu\neq0$.

\itemizedef
The $0\vert2$-parameter family of \psur{} $\rhoh^{\lambda,\mu}$ of $\SuperH$ on the space $\Schwartz{}{}$ is given by
\begin{moneq}
(\rhoh^{\lambda,\mu}_\grel \chi)(\adh) = \eexp^{-i\betah\lambda + i(\adh\betah-\gammah)\mu}\,\chi(\adh-\ah)
\mapob.
\end{moneq}
The infinitesimal actions that do not reduce to a multiple of the identity are
\begin{align*}
\tau^{\lambda,\mu}(\partial_{\ah})\psi
&=
-\partial_{\adh}\psi
&
\tau^{\lambda,\mu}(\partial_{\betah})\psi
&=
-i\lambda\psi + i\mu \adh\psi
\end{align*}
If $\psi$ belongs to an invariant closed subspace $F\subset \Schwartz{}{}$, it follows that the elements
\begin{moneq}
\partial_{\adh}\chi
\qquad\text{and}\qquad
\mu\adh\chi
\end{moneq}
should also belong to this subspace. 
When $\mu\neq0$, this implies that
\begin{moneq}
\partial_{\adh}\chi
\qquad\text{and}\qquad
\adh\chi
\end{moneq}
should belong to $F$. 
But these two operations are (two out of three of) the generators of the infinitesimal version of the standard irreducible representation of the Heisenberg group on $L^2(\RR)$, implying that $F$ must be a trivial subspace. 
And thus $\rhoh^{\lambda,\mu}$ is irreducible when $\mu\neq0$.

\subsection{Representations associated to coadjoint orbits}

\BackSwitchHats
\renewcommand\grel{{(a,\beta,\gamma)}}
\vrule width0pt
\vskip0.5\baselineskip

We now turn our attention to geometric quantization of the Coadjoint orbits of $\SuperH$, and we \myquote{copy} the arguments for the standard Heisenberg group. 
Using the standard basis $e_1=\partial_a, f_2=\partial_\beta, f_3=\partial_\gamma$ of the Lie algebra $\Liealg g = T_e\SuperH$ with its left-dual basis $\LD[1]e_1, \LD[2]f_2, \LD[2]f_3$ of $\LD[1]{\Liealg g}$, we obtain the Coadjoint action in terms of the left-coordinates $(x,y,z)\in \CA^3$ with respect to this dual basis as
\begin{align*}
\Coad(a,\beta,\gamma)(x\LD[1] e+y\LD[2]f_1+z\LD[2]f_2)
\cong
\Coad(a,\beta,\gamma)(x,y,z)
=
(x+z\beta, y-za,  z)
\mapob.
\end{align*}
When we split the coordinates $x,y,z\in \CA$ into their homogeneous parts $x=x_0+x_1$ etc., we find
\begin{moneq}
\Coad(a,\beta,\gamma)(x_0, x_1, y_0, y_1, z_0, z_1)
=
(x_0+z_1\beta, x_1+z_0\beta, y_0-az_0, y_1-az_1, z_0, z_1)
\mapob.
\end{moneq}
In terms of these six (three even and three odd) coordinates, the fundamental vector fields associated to the Coadjoint actions are given by
\begin{moneq}
e_1^{\LD[1]{\Liealg g}} =
z_0 \partial_{y_0} + z_1\partial_{y_1} 
\quad,\quad
f_2^{\LD[1]{\Liealg g}} =
z_1\partial_{x_0}-z_0\partial_{x_1}
\quad,\quad
f_3^{\LD[1]{\Liealg g}} = 0
\mapob.
\end{moneq}
Using the six homogeneous coordinates and starting at a point with real coordinates, \ie, $x_1=y_1=z_1=0$, we find the orbits described by
\begin{moneq}
(x_0,0, y_0,0, z_0, 0) \mapsto (x_0, \beta z_0, y_0-a z_0, 0, z_0,0)
\mapob,
\end{moneq}
which shows that for $z_0=0$ we have point orbits, and for $z_0\neq0$ we have $1\vert1$-dimensional orbits of the form $\{x_0\}\times \CA_1\times \CA_0\times \{(0,z_0,0)\}$.

\medskip\noindent\textbf{\textit{Point orbits.}}
For a point orbit $\mu_o=(x_0,0,y_0,0,0,0)$ we have $Y=\{\mu_o\}\times (\CA_0/d\ZZ)\times \CA_1$ and $\Gamma = \extder t + \extder \tau$, with which the infinitesimal lifted action of $\SuperH$ on $Y$ is given by the vector fields $b^Y$, $b=e_1,f_2,f_3$ satisfying \recalf{LiftingEquationsOfGroupActionOnPrequantumBundle}, giving
\begin{moneq}
e_1^Y = x_0\partial_t
\ ,\qquad
f_2^Y = y_o\partial_\tau
\ ,\qquad
f_3^Y = 0
\mapob.
\end{moneq}
This integrates to the action $\phi$ of $\SuperH$ on $Y$ as
\begin{moneq}
\phi_\grel:(\mu_o,t,\tau)\mapsto (\mu_o, t-ax_0, \tau-\beta y_0)
\mapob.
\end{moneq}
Functions on $Y$ that are equivariant with respect to the action of the structure group are of the form $\psi(\mu_o,t,\tau) = C\cdot \eexp^{it/\hbar+i\tau\kappa}$ for $C\in \CC$ and thus the representation $\Phi^{(x_0,y_0)}$ associated to such a point orbit (there is only the trivial polarization) is given by
\begin{align*}
\bigl(\Phi^{x_0,y_0}_\grel \psi\bigr)(\mu_o,t,\tau) 
&
= 
\psi\bigl(\phi_\grel\mo(\mu_o,t,\tau)\bigr)
\dontprintcomp{
=
\psi(\mu_o,t+ax_0, \tau+\beta y_0)
}%
\\&
=
C\cdot \eexp^{i(t+ax_0)/\hbar + i(\tau+\beta y_0)\kappa}
=
\eexp^{iax_0/\hbar+i\beta y_0\kappa}\cdot \psi(\mu_o,t,\tau)
\mapob.
\end{align*}
Comparing this with the family $\rhob^{k,(\lambda,0)}$ of representations of $\SuperH$ found above, we see that they are \myquote{the same,} provided we make the identifications $-x_0/\hbar = k$ and $-y_0\kappa = \lambda$.

\medskip\noindent\textbf{\textit{Orbits of dimension $1\vert1$.}}
For the $1\vert1$-dimensional orbit $\Orbit_{\mu_o}=\{x_0\}\times \CA_1\times \CA_0\times \{(0,z_0,0)\}$ through the point $\mu_o = (x_0,0,y_0,0,z_0,0)$ with $z_0\neq 0$ we use the coordinates $(\xi,y)\in \CA_1\times \CA_0$ to describe the point $(x_0,\xi,y,0,z_0,0,0)\in \Orbit_{\mu_o}$. 
In terms of these coordinates the Coadjoint action is given as $(\xi,y)\mapsto (\xi+\beta z_0, y-az_0)$ and the Kirillov-Kostant-Souriau symplectic form $\omega$ (still defined by \recalf{DefiningFormulaForKKSForm}) then is given as
\begin{moneq}
\omega = \frac1{z_0}\,\extder \xi\wedge \extder y
\mapob.
\end{moneq}
We now choose (any other choice will yield an equivalent description) the symplectic potential $\Theta = \xi\,\extder y/z_0$ and thus the connection $1$-form $\Gamma = \extder t + \extder\tau + \xi\,\extder y/z_0$. 
In this case the infinitesimal version of the lifted action of $\SuperH$ on $Y$ is determined by $b^Y = b^{\LD[1]{\Liealg g}} + d\,\partial_t+\delta\,\partial_\tau$, $b=e_1,f_2,f_3$ satisfying \recalf{LiftingEquationsOfGroupActionOnPrequantumBundle}, giving the solutions
\begin{moneq}
e_1^Y = z_0\partial_y +x_0\partial_t
\ ,\qquad
f_2^Y = -z_0\partial_\xi + y\partial_\tau
\ ,\qquad
f_3^Y = z_0\partial_\tau
\mapob.
\end{moneq}
This integrates to the action $\phi$ of $\SuperH$ on $Y$ given by
\begin{moneq}
\phi_\grel:(\xi,y,t,\tau)\mapsto (\xi+z_0\beta, y-az_0 ,t-ax_0, \tau -y\beta - z_0\gamma +z_0 a\beta)
\mapob.
\end{moneq}
For this orbit the stabilizer subalgebra is generated by $f_3$, but, contrary to the standard Heisenberg group, now, due to grading, there are only two polarization subalgebras: $\langle e_1,f_3\rangle$ of dimension $1\vert1$ and $\langle f_2,f_3\rangle$ of dimension $0\vert 2$. 
The horizontally lifted vector fields $b^h$, $b=e_1,f_2,f_3$ determined by the conditions \recalf{DefinitionOfHorizontalLiftFundamentalVectorFields} are given by
\begin{moneq}
e_1^h = z_0\partial_y - \xi\partial_\tau
\ ,\qquad
f_2^h = -z_0\partial_\xi
\ ,\qquad
f_3^h = 0
\mapob.
\end{moneq}

\itemizedef
The representation space associated to the polarization $\langle e_1,f_3\rangle$ consists of functions $\psi$ on $Y$ of the form $\psi(\xi,y,t,\tau)=\eexp^{it/\hbar+i\tau\kappa}\cdot \psi(\xi,y,0,0)$ satisfying $b^h\psi=0$ for $b=e_1,f_3$. 
This yields the space of functions of the form
\begin{moneq}
\psi(\xi,y,t,\tau) = \eexp^{it/\hbar+i\tau\kappa}\cdot \eexp^{i\xi y\kappa/z_0}\cdot \chi(\xi)
\mapob,
\end{moneq}
with $\chi$ a function of a single odd variable. 
The actual representation $\Phi^{x_0,z_0}$ then is given by
\begin{align*}
\bigl(\Phi^{x_0,z_0}_\grel \psi\bigr)(\xi,y,t,\tau)
&
=
\psi\bigl(\phi_\grel\mo(\xi,y,t,\tau)\bigr)
\dontprintcomp{
\\&
=
\psi(\xi-z_0\beta,y+z_0a,t+ax_0, \tau+y\beta+z_0\gamma)
}%
\\&
=
\eexp^{it/\hbar+i\tau\kappa}\, \eexp^{i\xi y\kappa/ z_0} \, \eexp^{iz_0(\gamma-a\beta)\kappa + iax_0/\hbar} \cdot \eexp^{i\xi\kappa a}\,  \chi(\xi-z_0\beta)
\mapob,
\end{align*}
which yields in terms of the function $\chi$ of a single odd variable the representation
\begin{moneq}
\bigl(\Phi^{x_0,z_0}_\grel \chi\bigr)(\xi) = \eexp^{ i\xi a\kappa}\, \eexp^{  iz_0(\gamma-a\beta)\kappa}\,\eexp^{ ix_0a/\hbar}\cdot  \chi(\xi-z_0\beta)
\mapob.
\end{moneq}
When we perform a simple change of the coordinate $\xi = z_0 \xih$ and writing $\hat \chi(\xih) = \chi(\xi)$, we find in terms of the functions $\hat \chi$ the representation
\begin{align*}
\bigl(\Phi^{x_0,z_0}_\grel \hat \chi\bigr)(\xih) 
&
= 
\eexp^{ iz_0\xih a\kappa}\, \eexp^{  iz_0(\gamma-a\beta)\kappa}\,\eexp^{ ix_0a/\hbar}\cdot  \chi(z_0\xih-z_0\beta)
\\&
=
\eexp^{ iz_0\xih a\kappa}\, \eexp^{  iz_0(\gamma-a\beta)\kappa}\,\eexp^{ ix_0a/\hbar}\cdot  \hat \chi(\xih-\beta)
\mapob.
\end{align*}
This is exactly the representation $\rho^{k,\mu}$ above, provided we make the identifications $-x_0/\hbar=k$ and $-z_0\kappa = \mu$. 

\itemizedef
We next turn our attention to the representation space associated to the polarization $\langle f_2,f_3\rangle$ consisting of functions of the form $\psi(\xi,y,t,\tau)=\eexp^{it/\hbar+i\tau\kappa}\, \psi(\xi,y,0,0)$ satisfying $b^h\psi=0$ with $b=f_2,f_3$. 
This yields the space of functions of the form
\begin{moneq}
\psi(\xi,y,t,\tau) = \eexp^{it/\hbar+i\tau\kappa}\,\eexp^{-ix_0y/(z_0\hbar)}\, \chi(y)
\mapob,
\end{moneq}
with $\chi$ a function of a single (even, real) variable (the reason to separate a factor $\eexp^{-ix_0y/(z_0\hbar)}$ will become clear soon). 
The actual representation $\Phih^{x_0,z_0}$ then is given by
\begin{align*}
\bigl(\Phih^{x_0,z_0}_\grel \psi\bigr)(\xi,y,t,\tau)
&
=
\psi\bigl(\phi_\grel\mo(\xi,y,t,\tau)\bigr)
\dontprintcomp{
\\&
=
\psi(\xi-z_0\beta,y+z_0a,t+x_0a, \tau+y\beta+z_0\gamma)
}%
\\&
=
\eexp^{i\tau\kappa+it/\hbar}\,\eexp^{-ix_0y/(z_0\hbar)}\, \eexp^{i(y\beta + z_0\gamma)\kappa}\,\chi(y+az_0)
\mapob,
\end{align*}
which yields in terms of the function $\chi$ of a single variable the representation (independent of $x_0$, the reason we introduced the factor $\eexp^{-ix_0y/(z_0\hbar)}$ above)
\begin{moneq}
\bigl(\Phih^{x_0,z_0}_\grel \chi\bigr)(y) = \eexp^{iy\beta\kappa + i\gamma z_0\kappa }\,  \chi(y+z_0a)
\mapob.
\end{moneq}
When we perform the simple change of the coordinate $y =y_0- z_0 \yh$ and write $\hat \chi(\yh) = \chi(y)$, we find in terms of the functions $\hat \chi$ the representation
\begin{align*}
\bigl(\Phih^{x_0,z_0}_\grel\hat \chi\bigr)(\yh) 
&
= 
\eexp^{i(y_0- z_0 \yh)\beta\kappa + i\gamma z_0\kappa }\,  \chi(y_0-z_0\yh+z_0a)
\\&
=
\eexp^{-i z_0 \yh\beta\kappa + i y_0\beta\kappa + i\gamma z_0\kappa }\,  \hat \chi(\yh - a)
\mapob.
\end{align*}
And this is exactly the representation $\rhoh^{\lambda,\mu}$ above, provided we make the identifications $-z_0\kappa = \mu$ (as before) and $-y_0\kappa = \lambda$.

\begin{remarque}{Concluding remarks}
\itemizedef[*]
We see that the $1$-dimensional representations obtained in the decomposition of the regular representations are all linked to point orbits and that the two initial ways to decompose it, yielding the two families of representations $\rho^{k,\mu}$ and $\rhoh^{\lambda,\mu}$, are related to the $1\vert1$-dimensional orbit using two different polarizations. 
However, contrary to the standard Heisenberg group, here these two different ways to decompose the regular representation do \stress{not} yield equivalent (sub)representations, even though the link with coadjoint orbits is the same. 
In this graded case the representations spaces are obviously different (one being finite dimensional, the other infinite dimensional). 
Whereas in the ungraded case (under good circumstances) the representation associated to a Coadjoint orbit is independent of the choice of a polarization, in the graded case the best one can hope is that representations using two polarizations \stress{of the same graded dimension} are equivalent.

\itemizedef
In order to link representations associated to Coadjoint orbits with representations found in our decomposition of the regular representation, we had to make identifications between the parameters found in the decomposition and the parameters of the orbit, as well for the standard Heisenberg group $H$ as for the super Heisenberg group $\SuperH$. 
Now in all except one case this poses no problems at all, but in the last case we have a singularity: in the decomposition we found a $0\vert2$-parameter family of representations with odd parameters $\lambda$ and $\mu$. 
In order to make a \myquote{perfect} link with the orbits, we had to identify $\mu =-z_0\kappa$ and $\lambda = -y_0\kappa$. 
But this means that the two odd parameters are a real multiple of the same odd Fourier parameter $\kappa$ from our decomposition of the natural representation of $\SuperH$ on the prequantum bundle, implying that we cannot choose both freely. 
We will encounter the same phenomenon in the next example, but I have no explanation\slash solution. 
\end{remarque}

\section{\label{SECTHeisenbergLikeSuperGroup}A Heisenberg-like supergroup}

\SwitchHats
\renewcommand\grel{{(\ah,\ch,\alphah,\gammah)}}

In this section I want to imitate what we did for the super Heisenberg group in \recals{SECTSuperHeisenbergGroup} to a particular Heisenberg-like \glgrp/ $G$ (Heisenberg-like because, just as the classical Heisenberg group, it is a central extension of an abelian \glgrp/). 
As \agmfd/ it will be $G=\RR^{4\vert4}$ of (graded) dimension $4\vert4$ with  coordinates $(a_1, a_2, a_3, c, \alpha_4, \alpha_5, \alpha_6, \gamma)$, of which $(a_1, a_2, a_3, c)$ are even and $(\alpha_4, \alpha_5, \alpha_6, \gamma)$ odd, equipped with the group law given by
\begin{moneq}
\begin{pmatrix} \ah_1\\ \ah_2\\ \ah_3\\ \ch\\ \alphah_4\\ \alphah_5\\ \alphah_6\\ \gammah
\end{pmatrix}
\cdot
\begin{pmatrix} \adh_1\\ \adh_2\\ \adh_3\\ \cdh\\ \alphadh_4\\ \alphadh_5\\ \alphadh_6\\ \gammadh
\end{pmatrix}
=
\begin{pmatrix} \ah_1+\adh_1\\ \ah_2+\adh_2\\ \ah_3+\adh_3
\\ 
\ch+\cdh + \adh_1\ah_2 + \tfrac12(\alphah_6\alphadh_6 - \alphah_5 \alphadh_5)
\\ 
\alphah_4+\alphadh_4\\ \alphah_5+\alphadh_5\\ \alphah_6+\alphadh_6
\\ 
\gammah+\gammadh + \alphah_4\adh_1 - \ah_3\alphadh_5
\end{pmatrix}
\mapob.
\end{moneq}
The underlying (ordinary) real Lie group $\body G$ is the direct product $H\times \RR$ of the Heisenberg groupe $H$ of \recals{SECT/HeisenbergTemplate} (coordinates $a_2,a_1,c$) with $\RR$ (coordinate $a_3$). 
In terms of a Harish-Chandra pair, it is the pair $(H\times \RR, \Liealg g)$ with $\Liealg g$ a $4\vert4$-dimensional graded Lie algebra, which will be described in \recals{SECTDeterminationOfCoadjointOrbitsHlike}. 
The neutral element has all coordinates $0$ and the inverse is given by
\begin{moneq}
(a,c,\alpha,\gamma)\mo = (-a, -c+a_1a_2, -\alpha, -\gamma+a_1\alpha_4-a_3\alpha_5)
\mapob.
\end{moneq}
Hence we have
\begin{align*}
(\ah,\ch,\alphah,\gammah)\mo \cdot (\adh,\cdh,\alphadh,\gammadh)
&
=
(\adh-\ah, \cdh -\ch+\ah_1\ah_2 - \adh_1\ah_2 + \tfrac12(\alphah_5\alphadh_5 - \alphah_6 \alphadh_6), 
\\&
\kern4em
\alphadh-\alphah , \gammadh -\gammah+\ah_1\alphah_4-\ah_3\alphah_5 -\alphah_4\adh_1 + \ah_3\alphadh_5 )
\mapob.
\end{align*}

\subsection{Decomposition of the left-regular representation}
\label{SUBSECTDecompLeftRegularRepHeisenbergLikeSuperGroup}

\vrule width0pt
\vskip0.5\baselineskip

For the same reasons as given in \recals{SECTSuperHeisenbergGroup}, I will define the left-regular representation $\rho$ of $G$ on the Schwartz space $\Schwartz44$; it is given by the standard formula
\begin{moneq}
(\rho_\grel\psi)(\adh,\cdh,\alphadh,\gammadh)
=
\psi\bigl((\ah,\ch,\alphah,\gammah)\mo \cdot (\adh,\cdh,\alphadh,\gammadh)\bigr)
\mapob.
\end{moneq}
This will be a super unitary representation when we use the standard \supscal product $\supinprsym_{4\vert4}$ and a \myquote{modified} inner product $\inprodsym_\mfdmetric$ associated to an invariant super metric $\mfdmetric$ (see \cite[4.19, 4.22]{Tuynman..2018} for details; it is related to the change of coordinates $\eta=\gammadh+ \adh_3\alphadh_5$) given by 
\begin{align*}
\inprod\chi\psi_\mfdmetric
&
=
\inprodd{\chi_0}{\psi_0}
+
\inprodd{\chi_4}{\psi_4} + \inprodd{{\chi_5-\phi_3 \chi_7}}{\psi_5-\phi_3 \psi_7} + \inprodd{\chi_6}{\psi_6} + \inprodd{\chi_7}{\psi_7} 
\\&
\kern2em
+ \inprodd{{\chi_{45}-\phi_3 \chi_{47}}}{\psi_{45}-\phi_3 \psi_{47}} 
+ \inprodd{\chi_{46}}{\psi_{46}} 
+ \inprodd{\chi_{47}}{\psi_{47}} 
\\&
\kern2em
+ \inprodd{{\chi_{56}+\phi_3 \chi_{67}}}{\psi_{56}+\phi_3 \psi_{67}}
+ \inprodd{\chi_{57}}{\psi_{57}} 
+ \inprodd{\chi_{67}}{\psi_{67}} 
\\&
\kern3.5em
+\inprodd{{\chi_{456}+\phi_3 \chi_{467}}}{\psi_{456}+\phi_3 \psi_{467}} + \inprodd{\chi_{467}}{\psi_{467}} 
\\&\kern3.5em
+ \inprodd{\chi_{457}}{\psi_{457}} + \inprodd{\chi_{567}}{\psi_{567}}
+ \inprodd{\chi_{4567}}{\psi_{4567}}
\mapob.
\end{align*}
In this formula $\phi_3:\RR^4\to \RR$ stands for the function $\phi_3(\adh,\cdh) = \adh_3$ and I used the notation (assimilating the odd coordinate $\gammadh$ with $\alphadh_7$)
\begin{moneq}[FunctionPsiInacalphagammaWithGammaSeven]
\psi(\adh,\cdh,\alphadh,\gammadh) = \sum_{I\subset \{4,\dots,7\}} \repprod\alphadh^I\cdot\psi_I(\adh,\cdh)
\mapob,
\end{moneq}
with $\psi_I\in \Schwartz4{}$ (and $\psi_0\equiv\psi_\emptyset$).

As for the classical (super) Heisenberg group it is not hard to see that functions on $G$ of the form
\begin{moneq}
\psi(\adh,\cdh, \alphadh,\gammadh)
=
\chi(\adh_1, \alphadh_5, \alphadh_6) \cdot
\eexp^{i\ell_2\adh_2+i\ell_3\adh_3+i\alphadh_4\lambda_4+i\ell_0\cdh+i\gammadh\lambda_7}
\mapob,
\end{moneq}
with $\chi$ a function of one even and two odd coordinates form an invariant subspace, as for such functions we have
\begin{align*}
\shifttag{3.5em}
(\rho_\grel \psi)(\adh, \cdh, \alphadh, \gammadh)
=
\psi\bigl(\gh\mo\cdot(\adh, \cdh, \alphadh, \gammadh)\bigr)
\\&
=
\Bigl(\,\chi(\adh_1 - \ah_1, \alphadh_5 - \alphah_5, \alphadh_6 - \alphah_6) \cdot\eexp^{-i\ell_0( \ah_2\adh_1 - \frac12\alphah_5\alphadh_5 +\frac12 \alphah_6 \alphadh_6)}\cdot
\eexp^{-i\ell_0( \ch-\ah_1\ah_2 )}
\\&
\kern3em
\cdot
\eexp^{-i(\alphah_4\adh_1 - \ah_3\alphadh_5)\lambda_7}
\cdot
\eexp^{-i(\gammah-\ah_1\alphah_4+\ah_3\alphah_5)\lambda_7}
\cdot
\eexp^{-i\ell_2\ah_2-i\ell_3\ah_3-i\alphah_4\lambda_4}
\,\Bigr)
\\&
\kern7em
{}\cdot
\eexp^{i\ell_2\adh_2+i\ell_3\adh_3+i\alphadh_4\lambda_4+i\ell_0\cdh+i\gammadh\lambda_7}
\mapob.
\end{align*}
We {thus} obtain a $3\vert2$-parameter family of representations $\rhoh^{(\ell_0,\ell_2,\ell_3,\lambda_4, \lambda_7)}$ of $G$ on $\Schwartz12$ given by
\begin{align}
\shifttag{5em}
(\rhoh^{(\ell_0,\ell_2,\ell_3,\lambda_4, \lambda_7)}_\grel \chi)(\adh_1,\alphadh_5,\alphadh_6)
\notag
\\
&
=
\chi(\adh_1 - \ah_1, \alphadh_5 - \alphah_5, \alphadh_6 - \alphah_6) \cdot\eexp^{-i\ell_0( \ah_2\adh_1 - \frac12\alphah_5\alphadh_5 +\frac12 \alphah_6 \alphadh_6)}\cdot
\eexp^{-i\ell_0( \ch-\ah_1\ah_2 )}
\notag
\\&
\kern3em
\cdot
\eexp^{-i(\alphah_4\adh_1 - \ah_3\alphadh_5)\lambda_7}
\cdot
\eexp^{-i(\gammah-\ah_1\alphah_4+\ah_3\alphah_5)\lambda_7}
\cdot
\eexp^{-i\ell_2\ah_2-i\ell_3\ah_3-i\alphah_4\lambda_4}
\mapob.
\label{FirstDecompositionOfLeftRegularRepOfHLikeSuperGroup}
\end{align}
We introduce the modified partial Berezin-Fourier transform $\Fourieroddp_{3\vert2} = \Fourierodd_{3\vert2}\scirc \Identif_\sigma$ with $\sigma$ the permutation $\sigma=\bigl(\begin{smallmatrix} 4&5&6&7 \\ 4&7&5&6 \end{smallmatrix}\bigr)$ \recalt{OddPermutationsChangeSignOfSuperScalarProductsNEW} \recalt{PartialBerezeinFourierTransformIsIsomorphismSHS}. 
It is given by
\begin{multline*}
\quad
\bigl((\Fourieroddp_{3\vert2}\psi)(\ell_0,\ell_2,\ell_3,\lambda_4, \lambda_7)\bigr)(\adh_1,\alphadh_5,\alphadh_6)
\\
=
(2\pi)^{-3/2}\cdot
\int \extder \adh_2\,\extder \adh_3\,\extder\cdh\,\extder\gammadh\,\extder \alphadh_4 \extdersep 
\eexp^{-i\ell_2\adh_2-i\ell_3\adh_3-i\alphadh_4\lambda_4-i\ell_0\cdh-i\gammadh\lambda_7}\cdot\psi(\adh,\cdh,\alphadh,\gammadh)
\mapob.
\quad
\end{multline*}
Applied to $\rho_\grel\psi$ we find
\KorteOfLangeTekst{
\begin{align*}
\shifttag{4em}
\Bigl(\bigl(\Fourieroddp_{3\vert2}(\rho_\grel\psi)\bigr)(\ell_0,\ell_2,\ell_3,\lambda_4, \lambda_7)\Bigr)(\adh_1,\alphadh_5,\alphadh_6)
\\&
=
\Bigl(\rhoh^{(\ell_0,\ell_2,\ell_3,\lambda_4, \lambda_7)}_\grel \bigl((\Fourieroddp_{3\vert2}\psi)(\ell_0,\ell_2,\ell_3,\lambda_4, \lambda_7)\bigr)\Bigr)(\adh_1,\alphadh_5,\alphadh_6)
\mapob.
\end{align*}
}{
\begin{align*}
\shifttag{4em}
\Bigl(\bigl(\Fourieroddp_{3\vert2}(\rho_\grel\psi)\bigr)(\ell_0,\ell_2,\ell_3,\lambda_4, \lambda_7)\Bigr)(\adh_1,\alphadh_5,\alphadh_6)
\\&
=
(2\pi)^{-3/2}\cdot
\int \extder \adh_2\,\extder \adh_3\,\extder\cdh\,\extder\gammadh\,\extder \alphadh_4 \extdersep 
\\&
\kern7em
\eexp^{-i\ell_2\adh_2-i\ell_3\adh_3-i\alphadh_4\lambda_4-i\ell_0\cdh-i\gammadh\lambda_7}\cdot(\rho_\grel\psi)(\adh,\cdh,\alphadh,\gammadh)
\\&
=
(2\pi)^{-3/2}\cdot
\int \extder\adh_2\,\extder\adh_3\,\extder\cdh\,\extder\gammadh\,\extder \alphadh_4 \extdersep 
\eexp^{-i\ell_2\adh_2-i\ell_3\adh_3-i\alphadh_4\lambda_4-i\ell_0\cdh-i\gammadh\lambda_7}
\\&
\kern7em
{}\cdot\psi(\adh-\ah, \cdh -\ch+\ah_1\ah_2 - \adh_1\ah_2 + \tfrac12(\alphah_5\alphadh_5 - \alphah_6 \alphadh_6), 
\\&
\kern10em
\alpha-\alphah , \gammadh -\gammah+\ah_1\alphah_4-\ah_3\alphah_5 -\alphah_4\adh_1 + \ah_3\alphadh_5 )
\\
\intertext{\footnotesize change of variables $\bar a_{2/3} = \adh_{2/3}-\ah_{2/3}$, $\bar\alpha_4=\alphadh_4-\alphah_4$, $\bar c= \cdh -\ch+\ah_1\ah_2 - a_1\ah_2 + \tfrac12(\alphah_5\alphadh_5 - \alphah_6 \alpha_6)$ and $\bar\gamma = \gammadh -\gammah+\ah_1\alphah_4-\ah_3\alphah_5 -\alphah_4\adh_1 + \ah_3\alphadh_5$ (of Berezinian equal $1$)}
&
=
(2\pi)^{-3/2}\cdot
\int \extder\bar a_2\,\extder\bar a_3\,\extder\bar c\,\extder\bar \gammadh\,\extder\bar \alphadh_4 \extdersep 
\eexp^{-i\ell_2\bar a_2-i\ell_3\bar a_3-i\bar\alphadh_4\lambda_4-i\ell_0\bar c-i\bar\gammadh\lambda_7}
\\&
\kern7em
{}\cdot\eexp^{-i\ell_2\ah_2-i\ell_3\ah_3-i\alphah_4\lambda_4-i\ell_0(\ch-\ah_1\ah_2 + \adh_1\ah_2 - \frac12(\alphah_5\alphadh_5 - \alphah_6 \alphadh_6))}
\\&
\kern8.5em
{}\cdot\eexp^{-i(\gammah-\ah_1\alphah_4+\ah_3\alphah_5 +\alphah_4\adh_1 - \ah_3\alphadh_5)\lambda_7}
\\&
\kern10em
{}\cdot\psi(\adh_1-\ah_1, \bar a_2, \bar a_3, \bar c,\bar\alpha_4, \alphadh_5-\alphah_5, \alphadh_6-\alphah_6, \bar\gamma)
\\&
=
\eexp^{-i\ell_2\ah_2-i\ell_3\ah_3-i\alphah_4\lambda_4-i\ell_0(\ch-\ah_1\ah_2 + a_1\ah_2 - \frac12(\alphah_5\alpha_5 - \alphah_6 \alpha_6))}
\\&
\kern2em
{}\cdot\eexp^{-i(\gammah-\ah_1\alphah_4+\ah_3\alphah_5 +\alphah_4\adh_1 - \ah_3\alphadh_5)\lambda_7}
\\&
\kern3.5em
{}\cdot
\bigl(\Fourieroddp_{3\vert2}\psi(\ell_0,\ell_2,\ell_3,\lambda_4, \lambda_7)\bigr)(\adh_1-\ah_1,\alphadh_5-\alphah_5,\alphadh_6-\alphah_6)
\\&
=
\Bigl(\rhoh^{(\ell_0,\ell_2,\ell_3,\lambda_4, \lambda_7)}_\grel \bigl((\Fourieroddp_{3\vert2}\psi)(\ell_0,\ell_2,\ell_3,\lambda_4, \lambda_7)\bigr)\Bigr)(\adh_1,\alphadh_5,\alphadh_6)
\mapob.
\end{align*}}%
This shows that $\Fourieroddp_{3\vert2}$  intertwines the left-regular representation of $G$ on $\Schwartz44$ with the direct integral (over three even and two odd variables) of the $3\vert2$-parameter family of representations $\rhoh^{(\ell_0,\ell_2,\ell_3,\lambda_4, \lambda_7)}$ of $G$ on $\Schwartz12$.

\begin{remarque}{Some additional remarks concerning representation spaces}
\itemizedef[*]
To show that the family $\rhoh^{(\ell_0,\ell_2,\ell_3,\lambda_4, \lambda_7)}$ is a $3\vert2$-parameter family of \psur s on $(\Schwartz12, \inprodsym_{1\vert2}, \supinprsym_{1\vert2})$, we first note that its expression \recalf{FirstDecompositionOfLeftRegularRepOfHLikeSuperGroup} obviously maps $\Schwartz12$ to itself and that it preserves both $\inprodsym_{1\vert2}$ (for $\gh\in \body G$) and $\supinprsym_{1\vert2}$.
The conditions \refmetnaam{\labelFamilySUR}{FamilySUR1}-\ref{FamilySUR3.5} then follow immediately. 
To show that it also satisfies \refmetnaam{\labelFamilySUR}{FamilySUR4}, we first note that $\rhot$ on $\Denset\equiv C^\infty_c\bigl(\RR^{3\vert2}; \Schwartz12\bigr)$ is given by
\begin{align*}
\shifttag{5em}
\bigl((\rhot_\grel\psi)(\ell_0,\ell_2,\ell_3,\lambda_4, \lambda_7)\bigr)( \adh_1,\alphadh_5,\alphadh_6)
\\&
=
\bigl(\psi(\ell_0,\ell_2,\ell_3,\lambda_4, \lambda_7)\bigr)( \adh_1 - \ah_1, \alphadh_5 - \alphah_5, \alphadh_6 - \alphah_6) 
\\&
\kern1.5em
\cdot\eexp^{-i\ell_0( \ah_2\adh_1 - \frac12\alphah_5\alphadh_5 +\frac12 \alphah_6 \alphadh_6)}\cdot
\eexp^{-i\ell_0( \ch-\ah_1\ah_2 )}
\notag
\\&
\kern3em
\cdot
\eexp^{-i(\alphah_4\adh_1 - \ah_3\alphadh_5)\lambda_7}
\cdot
\eexp^{-i(\gammah-\ah_1\alphah_4+\ah_3\alphah_5)\lambda_7}
\cdot
\eexp^{-i\ell_2\ah_2-i\ell_3\ah_3-i\alphah_4\lambda_4}
\mapob.
\end{align*}
As in the proof of \recalt{DirectIntegralOfHilbertSpacesPotentiallyChangeSignSuperScalarProductsNEW}, any element $\psi\in \Denset$ can be identified with a family of functions $\psi_{I,J}\in C^\infty_c\bigl(\RR^3; \Schwartz{}{}\bigr)$ with $I\subset \{4,7\}$ and $J\subset \{5,6\}$ according to
\begin{moneq}
\bigl(\psi(\ell_0,\ell_2,\ell_3,\lambda_4, \lambda_7)\bigr)( \adh_1 , \alphadh_5, \alphadh_6)
=
\sum_{I,J}\repprod{\lambda}^I\cdot\repprod\alphadh^J\cdot \bigl(\psi_{I,J}(\ell_0,\ell_2,\ell_3)\bigr)( \adh_1)
\mapob.
\end{moneq}
In terms of this identification, we define the inner product $\inprodsym_{\Denset}$ on $\Denset$ as
\begin{align}
\inprod\chi\psi_{\Denset}
&
=
\inprodd{\chi_0}{\psi_0}
+
\inprodd{\chi_4}{\psi_4} + \inprodd{{\chi_5}}{\psi_5} + \inprodd{\chi_6}{\psi_6} + \inprodd{\chi_7}{\psi_7} 
\notag
\\&
\kern2em
+ \inprodd{{\chi_{45}}}{\psi_{45}} 
+ \inprodd{\chi_{46}}{\psi_{46}} 
+ \inprodd{\chi_{47}}{\psi_{47}} 
+ \inprodd{{\chi_{56}}}{\psi_{56}}
\notag
\\&
\kern2em
+ \inprodd{\chi_{57}+\partial_{\ell_3}\chi_0}{\psi_{57}+\partial_{\ell_3}\psi_0} 
+ \inprodd{\chi_{67}}{\psi_{67}} 
\notag
\\&
\kern3.0em
+\inprodd{{\chi_{456}}}{\psi_{456}} + \inprodd{\chi_{467}}{\psi_{467}} 
+ \inprodd{\chi_{457}+\partial_{\ell_3}\chi_4}{\psi_{457}+\partial_{\ell_3}\psi_4} 
\label{InnerProductOnDensetFor44HeisenbergLikeGroupIn32Decomp}
\\&\kern3.0em
+ \inprodd{\chi_{567}-\partial_{\ell_3}\chi_6}{\psi_{567}-\partial_{\ell_3}\psi_6}
+ \inprodd{\chi_{4567}-\partial_{\ell_3}\chi_{46}}{\psi_{4567}-\partial_{\ell_3}\psi_{46}}
\mapob,
\notag
\end{align}
with the \myquote{obvious} identification on the indices and where $\inprodd\cdot\cdot$ stands for (see also \recalf{InnerProductWithLebesgueIntegral})
\begin{align*}
\inprodd{\chi}{\psi}
&
=
\int_{\RR^3} \  \inprod{\chi(\ell_0,\ell_2,\ell_3)}{\psi(\ell_0,\ell_2,\ell_3)}_{\Schwartz{}{}}\extdersep\extder\ell_0\,\extder\ell_2\,\extder\ell_3
\\&
=
\int_{\RR^3} \ \biggl(\ \int_{\RR} \overline{\bigl(\chi(\ell_0,\ell_2,\ell_3)\bigr)( \adh_1 )}\cdot \bigl(\psi(\ell_0,\ell_2,\ell_3)\bigr)( \adh_1)\extdersep\extder \adh_1\,\biggr)\extdersep\extder\ell_0\,\extder\ell_2\,\extder\ell_3
\\&
\cong
\int_{\RR^4} \  \overline{\bigl(\chi(\ell_0,\ell_2,\ell_3)\bigr)( \adh_1 )}\cdot \bigl(\psi(\ell_0,\ell_2,\ell_3)\bigr)( \adh_1)\extdersep\extder \adh_1\,\extder\ell_0\,\extder\ell_2\,\extder\ell_3
\mapob.
\end{align*}
(Nota Bene: the partial derivative $\partial_{\ell_3}$ denotes the partial derivative of a function $f:\RR^3\to \Schwartz{}{}$, a function of the $3$ variables $\ell_0, \ell_2,\ell_3$, as in \recalf{InnerProductOnDensetFor44HeisenbergLikeGroupIn32Decomp} we deal with elements $f\in C^\infty_c\bigl(\RR^3; \Schwartz{}{}\bigr)$.) 
It then is elementary to show that the maps $\rhot_\gh$ with $\gh\in \body G$ (\ie, $\alphah=0=\gammah$) preserve this inner product. 
As $\rhot$ also preserves the \myquote{standard} \supscal product $\supinprsym_{3\vert2}^{\Schwartz12}$ on $C^\infty_c\bigl(\RR^{3\vert2}; \Schwartz12\bigr)$, it follows that $\rhot$ is indeed a \psur{} as required by \refmetnaam{\labelFamilySUR}{FamilySUR4}. 

\itemizedef
In order to prove that $\Fourieroddp_{3\vert2}$ provides an equivalence between $\rho$ and $\rhot$, the direct integral of the $3\vert2$-parameter family $\rhoh^{(\ell_0,\ell_2,\ell_3,\lambda_4, \lambda_7)}$ of representations, we first note that we have the inclusions
\begin{moneq}
C^\infty_c(\RR^{4\vert4}) \subset C^\infty_c\bigl(\RR^{3\vert2}; C^\infty_c(\RR^{1\vert2})\bigr) \subset C^\infty_c\bigl(\RR^{3\vert2}, \Schwartz12\bigr)\equiv \Denset
\end{moneq}
and that these inclusions do not depend upon the inner product on $\Denset$, but only on the inner products $\inprodsym_{1\vert2}$ on $\Schwartz12$ and the standard inner product on $\CC$ (in order to be able to define the notion of smooth functions). 
We next note that the inner product $\inprodsym_{\Denset}$ on $\Denset$, restricted to $C^\infty_c(\RR^{4\vert4})$ can be extended by continuity to $\Schwartz44$ (by the same formula). 
Now the modified Berezin-Fourier transform $\Fourieroddp_{3\vert2}$ is a bijection from $\Schwartz44$ to $\Schwartz44$ and according to \recalt{OddPermutationsChangeSignOfSuperScalarProductsNEW} and \recalt{PartialBerezeinFourierTransformIsIsomorphismSHS} it intertwines the (standard) inner products and \supscal products $\inprodsym_{4\vert4}$ and $\supinprsym_{4\vert4}$. 
An elementary calculation shows that it also intertwines the modified inner product $\inprodsym_\mfdmetric$ on the source and the inner product $\inprodsym_{\Denset}$ on the target space. 
It follows that $\Fourieroddp_{3\vert2}$ establishes a bijection between a subspace $\Dense'\subset \Schwartz44$ and $C^\infty_c(\RR^{4\vert4})\subset \Schwartz44$ and that it is an equivalence of super inner product spaces between
\begin{moneq}
(\Dense', \inprodsym_\mfdmetric, \supinprsym_{4\vert4}) \subset \bigl(\Schwartz44, \inprodsym_\mfdmetric, \supinprsym_{4\vert4}\bigr)
\quad\text{the representation space of $\rho$}
\end{moneq}
and
\begin{moneq}
\bigl(C^\infty_c(\RR^{4\vert4}), \inprodsym_{\Denset}, \supinprsym_{4\vert4}\bigr)
\subset (\Denset, \inprodsym_{\Denset}, \supinprsym_{2\vert3}^{\Schwartz12})
\quad\text{the representation space of $\rhot$.}
\end{moneq}

There thus remains a single detail to check before we may conclude that it is an equivalence between $\rho$ and $\rhot$: we have to show that the given subspaces are dense. 
This seems obvious, but we are not dealing with the standard inner product (for which it is indeed obvious), but with modified inner products: $\inprodsym_\mfdmetric$ on the source and $\inprodsym_{\Denset}$ on the target. 
In order to show this last detail, we start with the definition of the map $\Phi:\Denset\to\Denset$ given by $\Phi\bigl((\psi_{I,J})_{I,J}\bigr) = \bigl((\psih_{I,J})_{I,J}\bigr)$ with $\psih_{I,J} = \psi_{I,J}$, except for four cases:
\begin{align*}
\psih_{57} &= \psi_{57} - \partial_{\ell_3}\psi_0
&
\psih_{567} &= \psi_{567} - \partial_{\ell_3}\psi_6
\\
\psih_{457} &= \psi_{457} - \partial_{\ell_3}\psi_4
& 
\psih_{4567} &= \psi_{4567} - \partial_{\ell_3}\psi_{46}
\mapob.
\end{align*}
It then follows immediately that $\Phi$ is a linear bijection and that it maps the \myquote{standard} inner product $\inprodsym_{4\vert4}$ (a slight abuse of notation) on $\Denset$ given by
\begin{moneq}
\inprod\chi\psi_{4\vert4} = \sum_{I,J}\ \inprodd{\chi_{I,J}}{\psi_{I,J}}
\end{moneq}
to the inner product $\inprodsym_{\Denset}$, meaning that
\begin{moneq}
\forall \chi,\psi\in \Denset:
\inprod{\Phi(\chi)}{\Phi(\psi)}_{\Denset} = \inprod\chi\psi_{4\vert4}
\mapob.
\end{moneq}
Said differently, $\Phi:(\Denset, \inprodsym_{4\vert4})\to (\Denset, \inprodsym_{\Denset})$ is an isometry of (non-super) inner product spaces. 
Moreover, it maps $C^\infty_c(\RR^{4\vert4})$ to itself. 
As $\bigl(C^\infty_c(\RR^{4\vert4}), \inprodsym_{4\vert4}\bigr)$ is dense in $L^2(\RR^{4\vert4})$, it is dense in $(\Denset, \inprodsym_{4\vert4})$, and thus $\bigl(C^\infty_c(\RR^{4\vert4}), \inprodsym_{\Denset}\bigr)$ is dense in $(\Denset, \inprodsym_{\Denset})$ as the image under an isometry of a dense subspace. 

Analogous to \recalt{ExampleWithDifferentTopologiesInducedByChangeOfCoordinates}, the isomorphism $\Phi$ is related to the (simple) change of coordinates on $\RR^{4\vert4}$ (elements of $\Denset$ can be seen as functions on $\RR^{4\vert4}$) given by 
\begin{moneq}
\phi(\ell_0,\ell_2,\ell_3, \lambda_4,\lambda_7,\adh_1,\alphadh_5,\alphadh_6)
=
(\ell_0,\ell_2,\ell_3+\alphadh_5\lambda_7, \lambda_4,\lambda_7,\adh_1,\alphadh_5,\alphadh_6)
\mapob.
\end{moneq}
To finish our last detail (in a hurry), we note that an analogous argument involving the change of coordinates
\begin{moneq}
\phi(\adh_1,\adh_2,\adh_3,\cdh,\alphadh_4,\alphadh_5,\alphadh_6, \gammadh)
=
(\adh_1,\adh_2,\adh_3,\cdh,\alphadh_4,\alphadh_5,\alphadh_6, \gammadh-\adh_3\alphadh_5)
\end{moneq}
shows that $(\Dense',\inprodsym_\mfdmetric)$, the source space of $\Fourieroddp_{3\vert2}$, is dense in $\bigl(\Schwartz44, \inprodsym_\mfdmetric\bigr)$, finishing the proof that $\rho$ and the direct integral $\rhot$ of the $3\vert2$-parameter family $\rhoh^{(\ell_0,\ell_2,\ell_3,\lambda_4, \lambda_7)}$ of representations are equivalent. 
\end{remarque}

Just as in the case of the classical Heisenberg group, we can decompose some representations among this $3\vert2$-parameter family of representations even further. 
For the classical Heisenberg group we had to distinguish just two cases: $m=0$ and $m\neq0$, but here we have to distinguish four cases according to whether $\ell_0$ and\slash or $\lambda_7$ is zero. 

\itemizedef
We start with the case of $\lambda_7=0$ and $\ell_0\neq0$, in which case it is not hard to see that functions of the form
\begin{moneq}[lambda7zeroInvariantSubspaces]
\chi(\adh_1, \alphadh_5, \alphadh_6)
=
h(\adh_1, \alphadh_5+\epsilon\alphadh_6) 
\cdot \eexp^{-\frac i2\ell_0\epsilon\alphadh_5\alphadh_6}
\end{moneq}
with $\epsilon=\pm1$ and $h$ a function of one even and one odd variable, form an invariant subspace, denoted by $\Schwartz12^\epsilon$, under the representation $\rhoh^{(\ell_0,\ell_2,\ell_3, \lambda_4,0)}$, as for such functions we have
\begin{align*}
(\rhoh^{(\ell_0,\ell_2,\ell_3, \lambda_4,0)}_\grel \chi)(\adh_1, \alphadh_5, \alphadh_6)
&
=
\chi(\adh_1 - \ah_1, \alphadh_5 - \alphah_5, \alphadh_6 - \alphah_6) 
{}\cdot
\eexp^{-i\ell_0( \ah_2\adh_1 - \frac12\alphah_5\alphadh_5 +\frac12 \alphah_6 \alphadh_6)}
\\&
\kern3.5em
{}\cdot
\eexp^{-i\ell_0( \ch-\ah_1\ah_2 )}
\cdot
\eexp^{-i\ell_2\ah_2}\cdot
\eexp^{-i\ell_3\ah_3}\cdot
\eexp^{-i\alphah_4\lambda_4}
\dontprintcomp{
\\&
=
h\bigl(\adh_1 - \ah_1, (\alphadh_5+\epsilon\alphadh_6) - (\alphah_5+\epsilon\alphah_6)\bigr) 
\\&
\kern2em
{}\cdot \eexp^{-\frac i2\ell_0\epsilon(\alphadh_5-\alphah_5)(\alphadh_6-\alphah_6)}
{}\cdot
\eexp^{-i\ell_0( \ah_2\adh_1 - \frac12\alphah_5\alphadh_5 +\frac12 \alphah_6 \alphadh_6)}
}
\\&
\kern3.5em
{}\cdot
\eexp^{-i\ell_0( \ch-\ah_1\ah_2 )}
\cdot
\eexp^{-i\ell_2\ah_2}\cdot
\eexp^{-i\ell_3\ah_3}\cdot
\eexp^{-i\alphah_4\lambda_4}
\\&
=
h\bigl(\adh_1 - \ah_1, (\alphadh_5+\epsilon\alphadh_6) - (\alphah_5+\epsilon\alphah_6)\bigr) \cdot \eexp^{-\frac i2\ell_0\epsilon\alphadh_5\alphadh_6}
\\&
\kern2em
{}\cdot
\eexp^{-i\ell_0( \ah_2\adh_1 - \frac12(\alphah_5-\epsilon\alphah_6)(\alphadh_5+\epsilon\alphadh_6) )}
\\&
\kern3.5em
{}\cdot
\eexp^{-i\ell_0( \ch-\ah_1\ah_2 +\frac12\epsilon\alphah_5\alphah_6)}
\cdot
\eexp^{-i\ell_2\ah_2}\cdot
\eexp^{-i\ell_3\ah_3}\cdot
\eexp^{-i\alphah_4\lambda_4}
\mapob.
\end{align*}
We {thus} obtain, for fixed values of $(\ell_0,\ell_2,\ell_3, \lambda_4)$, two representations $\rhoh^{(\ell_0,\ell_2,\ell_3, \lambda_4,0),\pm}$ of $G$ on $\Schwartz11$ given by
\begin{align}
(\rhoh^{(\ell_0,\ell_2,\ell_3, \lambda_4,0),\epsilon}_\grel h)(\adh_1, \xi)
&
=
h\bigl(\adh_1 - \ah_1, \xi - (\alphah_5+\epsilon\alphah_6)\bigr) 
{}\cdot
\eexp^{-i\ell_0( \ah_2\adh_1 - \frac12(\alphah_5-\epsilon\alphah_6)\xi )}
\notag
\\&
\kern3.5em
{}\cdot
\eexp^{-i\ell_0( \ch-\ah_1\ah_2 +\frac12\epsilon\alphah_5\alphah_6)}
\cdot
\eexp^{-i\ell_2\ah_2-i\ell_3\ah_3-i\alphah_4\lambda_4}
\mapob.
\label{RhohDecompositionPlusMoinsAtLambda7Zero}
\end{align}

As it is easy to show that $\Schwartz12$ is the direct sum of the two subspaces $\Schwartz12^\pm$, it is tempting to think that the super unitary representation $\rhoh^{(\ell_0,\ell_2,\ell_3, \lambda_4,0)}$ is the direct sum of the two super unitary representations $\rhoh^{(\ell_0,\ell_2,\ell_3, \lambda_4,0),\pm}$. 
But there is a snag: the restriction of the \supscal product $\supinprsym_{1\vert2}$ to each of the subspaces $\Schwartz12^\pm$ is identically zero. 
And thus $\rhoh^{(\ell_0,\ell_2,\ell_3, \lambda_4,0)}$ is \stress{not} the direct sum of the two representations $\rhoh^{(\ell_0,\ell_2,\ell_3, \lambda_4,0),\pm}$. 

What we will show is that $\Schwartz12$ decomposes as a direct integral of the (constant, independent of the odd parameter of the family) $0\vert1$-parameter family of super unitary representations $\rhoh^{(\ell_0,\ell_2,\ell_3, \lambda_4,0),+}$. 
More precisely, we consider the $0\vert1$-parameter family of representations $\rhob^\kappa$ on $\bigl(\Schwartz11,\inprodsym_{1\vert1}, \supinprsym_{1\vert1}\bigr)$ defined by (see \recalf{RhohDecompositionPlusMoinsAtLambda7Zero})
\begin{moneq}
(\rhob^\kappa_\grel h)(\adh_1,\xi)
=
(\rhoh^{(\ell_0,\ell_2,\ell_3, \lambda_4,0),+}_\grel h)(\adh_1,\xi)
\mapob,
\end{moneq}
independent of the odd parameter $\kappa$. 
It is elementary to show that this family satisfies the conditions \refmetnaam{\labelFamilySUR}{FamilySUR1}-\ref{FamilySUR3.5} and that $\rhot$ is a \psur{} when we use  the standard inner product $\inprodsym_{1\vert2}$ and the standard \supscal product $\supinprsym_{1\vert2}$ on $C^\infty_c\bigl(\RR^{0\vert1}; \Schwartz11\bigr) \cong\Schwartz12$ as required by \refmetnaam{\labelFamilySUR}{FamilySUR4}, proving that it indeed is a $0\vert1$-parameter family of \psur s. 

In order to show that $\rhoh^{(\ell_0,\ell_2,\ell_3, \lambda_4,0)}$ is equivalent to the direct integral of the $\rhob^\kappa$, we first define the isomorphisms $I_\epsilon:\Schwartz11 \to \Schwartz12^\epsilon$ by (for $h(\adh_1,\xi) = h_0(\adh_1)+\xi\,h_1(\adh_1)$\,)
\begin{align*}
\bigl(I_\epsilon(h)\bigr)(\adh_1,\alphadh_5,\alphadh_6) 
&
= 
h(\adh_1,\alphadh_5+\epsilon\alphadh_6)\cdot \eexp^{-\frac i2\ell_0\alphadh_5\alphadh_6}
\\&
=
h_{0}(\adh_1)+\alphadh_5\,h_{1}(\adh_1)+\epsilon\alphadh_6\, h_{1}(\adh_1) -\tfrac i2\ell_0\epsilon\alphadh_5\alphadh_6\,h_{0}(\adh_1)
\mapob.
\end{align*}
It follows (essentially by definition) that we have the intertwining property
\begin{moneq}
\rhoh^{(\ell_0,\ell_2,\ell_3, \lambda_4,0)}_\grel \scirc I_\epsilon  = I_\epsilon \scirc \rhoh^{(\ell_0,\ell_2,\ell_3, \lambda_4,0),\epsilon}_\grel 
\mapob.
\end{moneq}
Next we observe that the (slightly modified) odd Berezin-Fourier transform $\Fourieroddp_{0\vert1}:\Schwartz11 \to \Schwartz11$ given by
\begin{align*}
(\Fourieroddp_{0\vert1}h)(\adh_1,\eta)
&
=
\int\extder\xi \extdersep \eexp^{\frac12 i\ell_0\xi\eta}\cdot h(\adh_1,\xi)
=
h_1(\adh_1) +\tfrac12i\ell_0\eta\,h_0(\adh_1)
\mapob,
\end{align*}
intertwines the representations $\rhoh^{(\ell_0,\ell_2,\ell_3, \lambda_4,0),+}$ and $\rhoh^{(\ell_0,\ell_2,\ell_3, \lambda_4,0),-}$\,:
\begin{moneq}
\rhoh^{(\ell_0,\ell_2,\ell_3, \lambda_4,0),-\epsilon}_\grel \scirc \Fourieroddp_{0\vert1}
=
\Fourieroddp_{0\vert1} \scirc \rhoh^{(\ell_0,\ell_2,\ell_3, \lambda_4,0),\epsilon}_\grel
\mapob.
\end{moneq}
Now recall that any (smooth) function $h:\RR^{0\vert1}\to \Schwartz11$ is determined by two elements $h_0,h_1\in \Schwartz11$ according to $h(\kappa) = h_0+\kappa\,h_1$. 
Using this identification we define the odd bijection $\Phi:C^\infty_c\bigl(\RR^{0\vert1}; \Schwartz11\bigr)\cong\Schwartz12 \to \Schwartz12$ by
\begin{align*}
\Phi(h)
&
=
I_+(h_1) + I_-\bigl(\Fourieroddp_{0\vert1} (h_0)\bigr)
\mapob.
\end{align*}
That it indeed is a bijection follows from the facts that $I_\epsilon\bigl(\Schwartz11\bigr)=\Schwartz12^\epsilon$ and $\Schwartz12=\Schwartz12^+\oplus\Schwartz12^-$. 
To show that it intertwines the direct integral of the representations $\rhob^\kappa$ on $C^\infty_c\bigl(\RR^{0\vert1}; \Schwartz11\bigr)$ and the representation $\rhoh^{(\ell_0,\ell_2,\ell_3, \lambda_4,0)}$ on $\Schwartz12$, we compute:
\begin{align*}
\rhoh^{(\ell_0,\ell_2,\ell_3, \lambda_4,0)}_\grel\bigl(\Phi(h)\bigr)
&
=
\rhoh^{(\ell_0,\ell_2,\ell_3, \lambda_4,0)}_\grel\Bigl( I_+(h_1) + I_-\bigl(\Fourieroddp_{0\vert1}( h_0)\bigr) \Bigr)
\\&
=
I_+\bigl(\rhoh^{(\ell_0,\ell_2,\ell_3, \lambda_4,0),+}_\grel (h_1)\bigr) + I_- \hBigl(\rhoh^{(\ell_0,\ell_2,\ell_3, \lambda_4,0),-}_\grel \bigl(\Fourieroddp_{0\vert1} (h_0)\bigr)\hBigr)
\\&
=
I_+\bigl(\rhoh^{(\ell_0,\ell_2,\ell_3, \lambda_4,0),+}_\grel (h_1)\bigr) + I_- \hBigl(\Fourieroddp_{0\vert1} \bigl( \rhoh^{(\ell_0,\ell_2,\ell_3, \lambda_4,0),+}_\grel (h_0)\bigr)\hBigr)
\\&
=
\Phi\bigl(\rhob_\grel^{\kappa}( h)\bigr)
\mapob.
\end{align*}
And we finish by showing that $\Phi$ is partially an odd $1/\ell_0$-isomorphism of super inner product spaces in the sense of \refmetnaam{\labelMainSUReq}{MainSUREQlabelPhiIntertwinesSuperScalarProducts}. 
As $\Phi$ is odd, the condition \refmetnaam{\labelIsomSHS}{conditiondirectsuminequivalenceSHS} is automatically satisfied. As for \refmetnaam{\labelIsomSHS}{conditionoddisoforequivalenceofSHS}, we start with the decomposition of an element $h\in C^\infty_c\bigl(\RR^{0\vert1};\Schwartz11\bigr)\cong\Schwartz12$:
\begin{moneq}
h(\kappa) = h_0+\kappa h_1\text{ with } h_i\in \Schwartz11
\ \Rightarrow\ 
h_i(x,\xi) = h_i^0(x) + \xi h_i^1(x)
\text{ with } h_i^j\in \Schwartz{}{}
\mapob.
\end{moneq}
For $h,k\in C^\infty_c\bigl(\RR^{0\vert1}; \Schwartz11\bigr)$ we thus have
\begin{align*}
\homsuperinprod hk{_{0\vert1}^{\Schwartz11}}
&
=
\int\extder\lambda\extdersep 
\homsuperinprod{h_0+\lambda h_1}{k_0+\lambda k_1}{_{1\vert1}}
\\&
=
\homsuperinprod{h_0}{\conjugate k_1}{_{1\vert1}} + \homsuperinprod{\conjugate h_1}{\conjugate k_0}{_{1\vert1}}
\\&
=
-\inprodd{h^0_0}{k^1_1} + \inprodd{h_0^1}{k_1^0} - \inprodd{h_1^0}{k_0^1} - \inprodd{h^1_1}{k^0_0}
\\&
=
-\hBigl(\inprodd{h_{0}^0}{k_{1}^1} - \inprodd{h_{0}^1}{k_{1}^0} + \inprodd{h_{1}^0}{k_{0}^1} + \inprodd{h_{1}^1}{k_{0}^0}
\hBigr)
\mapob.
\end{align*}
But also
\begin{align*}
\Phi(h)(\adh_1,\alphadh_5,\alphadh_6)
&
=
\bigl(I_+(h_1)\bigr)(\adh_1,\alphadh_5,\alphadh_6) + \hBigl(I_-\bigl(\Fourieroddp (h_0)\bigr)\hBigr)(\adh_1,\alphadh_5,\alphadh_6)
\\&
=
\bigl(h_1^0(\adh_1) + \alphadh_5\,h_1^1(\adh_1) + \alphadh_6\,h_1^1(\adh_1) -\tfrac i2 \ell_0 \alphadh_5\alphadh_6\,h_1^0(\adh_1)\bigr)
\\&
\kern2.5em
+
\bigl( h_0^1(\adh_1) + \alphadh_5\,\tfrac i2\ell_0 h_0^0(\adh_1) - \alphadh_6\,\tfrac i2\ell_0 h^0_0(\adh_1)   + \tfrac i2\ell_0\alphadh_5\alphadh_6\,h_0^1(\adh_1)
\bigr)
\\&
=
\bigl(h_1^0(\adh_1) + h_0^1(\adh_1)\bigr) + \alphadh_5\bigl(h^1_1(\adh_1) + \tfrac i2\ell_0 h^0_0(\adh_1)\bigr) 
\\&
\kern2.5em 
+ \alphadh_6\bigl(h^1_1(\adh_1) - \tfrac i2\ell_0 h^0_0(\adh_1)\bigr)
+ \tfrac i2\ell_0 \alphadh_5\alphadh_6\bigl( h_0^1(\adh_1) - h_1^0(\adh_1)\bigr)
\end{align*}
and thus
\begin{align*}
\homsuperinprod{\Phi(h)}{\Phi(k)}{_{1\vert2}}
&
=
\inprodd{h^0_1 + h^1_0}{\tfrac i2\ell_0( k_0^1 - k_1^0)} + \inprodd{\tfrac i2\ell_0 (h_0^1 - h_1^0)}{k^0_1 + k^1_0}
\\&
\kern3em
+ \inprodd{h^1_1 + \tfrac i2\ell_0 h^0_0}{k^1_1 - \tfrac i2\ell_0 k^0_0} - \inprodd{h^1_1 - \tfrac i2\ell_0 h^0_0}{k^1_1 + \tfrac i2\ell_0 k^0_0}
\\&
=
\tfrac i2\ell_0 \hBigl(\inprodd{h^0_1 + h^1_0}{k_0^1 - k_1^0} - \inprodd{h_0^1 - h_1^0}{k^0_1 + k^1_0}\hBigr)
\\&
\kern3em
+ \inprodd{h^1_1 }{ - i\ell_0 k^0_0} + \inprodd{\tfrac i2\ell_0 h^0_0}{2k^1_1 }
\\&
=
\tfrac i2\ell_0 \hBigl(\inprodd{h^0_1 }{ 2k^1_0} + \inprodd{h^1_0}{-2k^0_1}\hBigr)
- i\ell_0\hBigl( \inprodd{h^1_1 }{  k^0_0} + \inprodd{h^0_0}{k^1_1 } \hBigr)
\\&
=
-i\ell_0\hBigl( \inprodd{h^0_0}{k^1_1 } - \inprodd{h_1^0}{k_0^1} + \inprodd{h_0^1 }{ k_1^0} + \inprodd{h^1_1 }{  k^0_0} 
\hBigr)
\\&
=
i\ell_0 \,\homsuperinprod{\conjugate h}k{_{0\vert1}^{\Schwartz11}}
\mapob,
\end{align*}
proving that $\Phi$ satisfies \refmetnaam{\labelIsomSHS}{conditionoddisoforequivalenceofSHS} with $k=1/\ell_0$ as claimed.

\itemizedef
Next we consider the case $\ell_0=0$ and $\lambda_7\neq0$, in which case it is not hard to see that functions of the form
\begin{moneq}
\chi(\adh_1, \alphadh_5, \alphadh_6)
=
h(\adh_1, \alphadh_5)\cdot \eexp^{i\alphadh_6\lambda_6}
\mapob,
\end{moneq}
with $h$ a function of one even and one odd variable, form an invariant subspace under the representation $\rhoh^{(\ell_0=0,\ell_2,\ell_3,\lambda_4, \lambda_7)}$, as for such functions we have
\begin{align*}
(\rhoh^{(0,\ell_2,\ell_3,\lambda_4, \lambda_7)}_\grel \chi)(\adh_1, \alphadh_5, \alphadh_6)
&
=
\chi(\adh_1 - \ah_1, \alphadh_5 - \alphah_5, \alphadh_6 - \alphah_6) 
\cdot
\eexp^{-i(\alphah_4\adh_1 - \ah_3\alphadh_5)\lambda_7}
\\&
\kern5em
{}\cdot
\eexp^{-i(\gammah-\ah_1\alphah_4+\ah_3\alphah_5)\lambda_7}
\cdot
\eexp^{-i\ell_2\ah_2-i\ell_3\ah_3-i\alphah_4\lambda_4}
\\&
=
\Bigl(h(\adh_1-\ah_1, \alphadh_5-\alphah_5) 
\cdot
\eexp^{-i(\alphah_4\adh_1 - \ah_3\alphadh_5)\lambda_7 
-i(\gammah-\ah_1\alphah_4+\ah_3\alphah_5)\lambda_7}
\\&
\kern5em
{}
\cdot
\eexp^{-i\ell_2\ah_2-i\ell_3\ah_3-i\alphah_4\lambda_4-i\alphah_6\lambda_6}
\Bigr)\cdot \eexp^{i\alphadh_6\lambda_6}
\mapob.
\end{align*}
We {thus} obtain, for fixed values of $(\ell_2,\ell_3,\lambda_4, \lambda_7)$, a $0\vert1$-parameter family of representations $\rhoh^{(0,\ell_2,\ell_3,\lambda_4, \lambda_7),\lambda_6}$ of $G$ on $\Schwartz11$ given by
\begin{align}
(\rhoh^{(0,\ell_2,\ell_3,\lambda_4, \lambda_7),(\lambda_6)}_\grel h)(\adh_1,\alphadh_5)
&
=
h(\adh_1-\ah_1, \alphadh_5-\alphah_5) 
\cdot
\eexp^{-i(\alphah_4\adh_1 - \ah_3\alphadh_5)\lambda_7}
\notag
\\&
\kern2em
{}\cdot
\eexp^{-i(\gammah-\ah_1\alphah_4+\ah_3\alphah_5)\lambda_7}
\cdot
\eexp^{-i\ell_2\ah_2-i\ell_3\ah_3-i\alphah_4\lambda_4-i\alphah_6\lambda_6}
\mapob.
\label{ell0zerobutnotlambda0action}
\end{align}
To show that $\rhoh^{(0,\ell_2,\ell_3,\lambda_4, \lambda_7)}$ decomposes as the direct integral of this family, we consider the the modified partial Berezin-Fourier transform $\Fourieroddp_{0\vert1}:\Schwartz12\to \Schwartz12$ given by $\Fourieroddp_{0\vert1} = \Fourierodd_{0\vert1}\scirc \Identif_\sigma$ with $\Fourierodd_{0\vert1}$ as in \recalt{PartialBerezeinFourierTransformIsIsomorphismSHS} and $\Identif_\sigma$ as in \recalt{OddPermutationsChangeSignOfSuperScalarProductsNEW} with $\sigma = \bigl( \begin{smallmatrix} 1&2 \\ 2&1 \end{smallmatrix}\bigr)$ the (unique) transposition of $\{1,2\}$, \ie, given by the formula 
\begin{moneq}
\bigl((\Fourieroddp_{0\vert1}\chi)(\lambda_6)\bigr)(\adh_1,\alphadh_5)
=
\int \extder \alphadh_6 \extdersep 
\eexp^{-i\alphadh_6\lambda_6}\cdot \chi(\adh_1,\alphadh_5,\alphadh_6)
\mapob,
\end{moneq}
and 
\KorteOfLangeTekst{%
a direct computation gives
\begin{align*}
\shifttag{4em}
\Bigl(\bigl(\Fourieroddp_{0\vert1}(\rhoh^{(0,\ell_2,\ell_3,\lambda_4, \lambda_7)}_\grel \chi)\bigr)(\lambda_6)\Bigr)(\adh_1,\alphadh_5)
=
\Bigl(\rhoh^{(0,\ell_2,\ell_3,\lambda_4, \lambda_7),(\lambda_6)}_\grel \bigl((\Fourieroddp_{0\vert1}\chi)(\lambda_6)\bigr)\Bigr)(\adh_1,\alphadh_5)
\mapob.
\end{align*}
}{%
we compute
\begin{align*}
\shifttag{4em}
\Bigl(\bigl(\Fourieroddp_{0\vert1}(\rhoh^{(0,\ell_2,\ell_3,\lambda_4, \lambda_7)}_\grel \chi)\bigr)(\lambda_6)\Bigr)(\adh_1,\alphadh_5)
\\&
=
\int \extder \alphadh_6 \extdersep  \eexp^{-i\alphadh_6\lambda_6}
\cdot(\rhoh^{(0,\ell_2,\ell_3,\lambda_4, \lambda_7)}_\grel \chi)(\adh_1,\alphadh_5,\alphadh_6)
\\&
=
\int \extder \alphadh_6 \extdersep  \eexp^{-i\alpha_6\lambda_6}
{}\cdot 
\chi(\adh_1 - \ah_1, \alphadh_5 - \alphah_5, \alphadh_6 - \alphah_6) 
{}\cdot
\eexp^{-i(\alphah_4\adh_1 - \ah_3\alphadh_5)\lambda_7}
\\&\kern3.5em
{}\cdot
\eexp^{-i(\gammah-\ah_1\alphah_4+\ah_3\alphah_5)\lambda_7}
\cdot
\eexp^{-i\ell_2\ah_2}\cdot
\eexp^{-i\ell_3\ah_3}\cdot
\eexp^{-i\alphah_4\lambda_4}
\\
\intertext{\footnotesize change of variables $\bar\alpha_6=\alphadh_6-\alphah_6$}
&
=
\int \extder \bar\alpha_6 \extdersep  \eexp^{-i\bar\alpha_6\lambda_6}
{}\cdot 
\chi(\adh_1 - \ah_1, \alphadh_5 - \alphah_5, \bar\alpha_6) 
{}\cdot
\eexp^{-i(\alphah_4\adh_1 - \ah_3\alphadh_5)\lambda_7}
\\&\kern3.5em
{}\cdot
\eexp^{-i(\gammah-\ah_1\alphah_4+\ah_3\alphah_5)\lambda_7}
\cdot
\eexp^{-i\ell_2\ah_2-i\ell_3\ah_3-i\alphah_4\lambda_4-i\alphah_6\lambda_6}
\\&
=
\eexp^{-i(\alphah_4\adh_1 - \ah_3\alphadh_5)\lambda_7}
\cdot
\eexp^{-i(\gammah-\ah_1\alphah_4+\ah_3\alphah_5)\lambda_7}
\cdot
\eexp^{-i\ell_2\ah_2-i\ell_3\ah_3-i\alphah_4\lambda_4-i\alphah_6\lambda_6}
\\&
\kern3.5em
{}\cdot
\bigl((\Fourieroddp_{0\vert1} \chi)(\lambda_6)\bigr)(\adh_1-\ah_1,\alphadh_5-\alphah_5)
\\&
=
\Bigl(\rhoh^{(0,\ell_2,\ell_3,\lambda_4, \lambda_7),(\lambda_6)}_\grel \bigl((\Fourieroddp_{0\vert1}\chi)(\lambda_6)\bigr)\Bigr)(\adh_1,\alphadh_5)
\mapob.
\end{align*}}%
This shows that $\Fourieroddp_{0\vert1}$ intertwines the representation $\rhoh^{(0,\ell_2,\ell_3,\lambda_4, \lambda_7)}$ of $G$ on $\Schwartz12$ with the direct integral representation of the $0\vert1$-parameter family $\rhoh^{(0,\ell_2,\ell_3,\lambda_4, \lambda_7),(\lambda_6)}$ of $G$ on $\Schwartz11$.

\begin{remarque}{Some additional remarks concerning representation spaces}
In the discussion of the case $\lambda_7=0$ above I have been a bit cavalier on the spaces on which the family are defined. 
So let me be a bit more precise in this case (which is slightly more \myquote{complicated} than the case $\lambda_7=0$, where each $\rhoh^{(\ell_0,\ell_2,\ell_3,\lambda_4,0),\epsilon}_{(\ah,\ch,0,0)}$ with $(\ah,\ch,0,0)\in \body G$ preserves the inner product $\inprodsym_{1\vert2}$ on $\Schwartz12$, as here the maps $\rhoh^{(0,\ell_2,\ell_3,\lambda_4, \lambda_7)}_{(\ah,\ch,0,0)}$ do not.
This means that we cannot see the $\rhoh^{(0,\ell_2,\ell_3,\lambda_4, \lambda_7)}$ separately as \psur s as in the case of $\lambda_7=0$ above, and thus that we have to keep track of the whole family. 
More precisely, we first have to show that the subfamily $\rhoh^{(0,\ell_2,\ell_3,\lambda_4, \lambda_7)}$ is indeed a $2\vert2$-parameter family of \psur s on $\Schwartz12$, that the family $\rhoh^{(0,\ell_2,\ell_3,\lambda_4, \lambda_7),(\lambda_6)}$ is a $2\vert3$-parameter family of \psur s on $\Schwartz11$ and then that the direct integral of the subfamily $\rhoh^{(0,\ell_2,\ell_3,\lambda_4, \lambda_7)}$ is equivalent to the direct integral of the family $\rhoh^{(0,\ell_2,\ell_3,\lambda_4, \lambda_7),(\lambda_6)}$. 

For the first point we note that the subfamily automatically satisfies \refmetnaam{\labelFamilySUR}{FamilySUR1}-\ref{FamilySUR3.5}, and that the inner product $\inprodsym_{\Denset}$ on $C^\infty_c\bigl(\RR^{2\vert2}; \Schwartz12\bigr)$ given by the same expression \recalf{InnerProductOnDensetFor44HeisenbergLikeGroupIn32Decomp} but now with
\begin{moneq}
\inprodd{\chi}{\psi}
=
\int_{\RR^3} \  \overline{\bigl(\chi(\ell_2,\ell_3)\bigr)( \adh_1 )}\cdot \bigl(\psi(\ell_2,\ell_3)\bigr)( \adh_1)\extdersep\extder \adh_1\,\extder\ell_2\,\extder\ell_3
\end{moneq}
makes that this family also satisfies \refmetnaam{\labelFamilySUR}{FamilySUR4}.

For the family $\rhoh^{(0,\ell_2,\ell_3,\lambda_4, \lambda_7),(\lambda_6)}$ it is elementary to show that it satisfies \refmetnaam{\labelFamilySUR}{FamilySUR1}-\ref{FamilySUR3.5} when $\Schwartz11$ is equipped with the standard inner product $\inprodsym_{1\vert1}$ and \supscal product $\supinprsym_{1\vert1}$. 
The existence of an inner product on $C^\infty_c\bigl(\RR^{2\vert3}; \Schwartz11\bigr)$ that will make it satisfy \refmetnaam{\labelFamilySUR}{FamilySUR4} will be postponed. 

For the equivalence between these two (direct integrals of) representations, we just extend the bijection $\Fourieroddp_{0\vert1}:\Schwartz12\to C^\infty_c\bigl(\RR^{0\vert1}; \Schwartz11\bigr)\cong\Schwartz12$ to the added parameters by defining $\Fourieroddpp_{0\vert1}:C^\infty_c\bigl(\RR^{2\vert2}; \Schwartz12\bigr)  \to C^\infty_c\bigl(\RR^{2\vert3}; \Schwartz11\bigr)$ by
\begin{moneq}
(\Fourieroddpp_{0\vert1}\psi)(\ell_2,\ell_3,\lambda_4,\lambda_7,\lambda_6)
=
\hBigl(\Fourieroddp_{0\vert1}\bigl(\psi(\ell_2,\ell_3,\lambda_4,\lambda_7)\bigr)\hBigr)(\lambda_6)
\mapob.
\end{moneq}
It is (again) elementary to verify that this is an odd bijection that intertwines the \supscal products. 
As we have seen above that the family $\rho^{(0,\ell_2,\ell_3,\lambda_4, \lambda_7)}$ satisfies \refmetnaam{\labelFamilySUR}{FamilySUR4}, we can use this odd bijection $\Fourieroddpp_{0\vert1}$ to transport the inner product on $C^\infty_c\bigl(\RR^{2\vert2}; \Schwartz12\bigr)$ to $C^\infty_c\bigl(\RR^{2\vert3}; \Schwartz11\bigr)$. 
It follows immediately that equipped with this transported inner product, the family $\rhoh^{(0,\ell_2,\ell_3,\lambda_4, \lambda_7),(\lambda_6)}$ satisfies \refmetnaam{\labelFamilySUR}{FamilySUR4}. 
\end{remarque}

\itemizedef
The last case to discuss is $\ell_0=\lambda_7=0$, in which case it is not hard to show that functions of the form
\begin{moneq}
\chi(\adh_1,\alphadh_5,\alphadh_6) = C\cdot \eexp^{i\ell_1 \adh_1+i \alphadh_5\lambda_5+i\alphadh_6\lambda_6}
\end{moneq}
with $C\in \CC$ form an invariant subspace (isomorphic to $\CC$) under the representation $\rhoh^{(0,\ell_2,\ell_3, \lambda_4,0)}$, as for such functions we have
\begin{align*}
(\rhoh^{(0,\ell_2,\ell_3, \lambda_4,0)}_\grel \chi)(\adh_1, \alphadh_5, \alphadh_6)
&
=
\chi(\adh_1 - \ah_1, \alphadh_5 - \alphah_5, \alphadh_6 - \alphah_6) 
\cdot
\eexp^{-i\ell_2\ah_2-i\ell_3\ah_3-i\alphah_4\lambda_4}
\\&
=
C\cdot 
\eexp^{i\ell_1 \adh_1+i \alphadh_5\lambda_5+i\alphadh_6\lambda_6} \cdot \eexp^{-i\ell_1 \ah_1-i \alphah_5\lambda_5-i\alphah_6\lambda_6}
\\&
\kern5em
{}\cdot
\eexp^{-i\ell_2\ah_2-i\ell_3\ah_3-i\alphah_4\lambda_4}
\mapob.
\end{align*}
We {thus} obtain, for fixed values of $(\ell_2,\ell_3, \lambda_4)$, a $1\vert2$-parameter family of representations $\rhoh^{(0,\ell_2,\ell_3, \lambda_4,0),(\ell_1,\lambda_5,\lambda_6)}$ of $G$ on $\CC$ given by
\begin{align}
\rhoh^{(0,\ell_2,\ell_3, \lambda_4,0),(\ell_1,\lambda_5,\lambda_6)}_\grel \,C
&
=
C\cdot \eexp^{-i\ell_1 \ah_1-i \alphah_5\lambda_5-i\alphah_6\lambda_6-i\ell_2\ah_2-i\ell_3\ah_3-i\alphah_4\lambda_4}
\mapob.
\label{ell0zeroandlambda0zeroaction}
\end{align}
When we perform the (full) Berezin-Fourier transform $\Fourierodd_{1\vert2} : \Schwartz12 \to \Schwartz12$ given by
\begin{moneq}
\quad
(\Fourierodd_{1\vert2}\chi)(\ell_1,\lambda_5,\lambda_6)\
=
\int \extder \adh_1\,\extder\alphadh_6\,\extder \alphadh_5 \extdersep 
\eexp^{-i\ell_1\adh_1-i\alphadh_5\lambda_5-i\alphadh_6\lambda_6}
\cdot \chi(\adh_1,\alphadh_5,\alphadh_6)
\mapob,
\quad
\end{moneq}
we find 
\KorteOfLangeTekst{
\begin{align*}
\shifttag{4em}
\bigl(\Fourierodd_{1\vert2}(\rhoh^{(0,\ell_2,\ell_3, \lambda_4,0)}_\grel \chi)\bigr)(\ell_1,\lambda_5,\lambda_6)
=
\rhoh^{(0,\ell_2,\ell_3, \lambda_4,0),(\ell_1,\lambda_5,\lambda_6)}_\grel \bigl((\Fourierodd_{1\vert2}\chi)(\ell_1,\lambda_5,\lambda_6)\bigr)
\mapob.
\end{align*}
}{
\begin{align*}
\shifttag{4em}
\bigl(\Fourierodd_{1\vert2}(\rhoh^{(0,\ell_2,\ell_3, \lambda_4,0)}_\grel \chi)\bigr)(\ell_1,\lambda_5,\lambda_6)
\\&
=
\int \extder \adh_1\,\extder\alphadh_6\,\extder \alphadh_5 \extdersep 
\eexp^{-i\ell_1\adh_1-i\alphadh_5\lambda_5-i\alphadh_6\lambda_6}
\cdot(\rhoh^{(0,\ell_2,\ell_3, \lambda_4,0)}_\grel \chi)(\adh_1,\alphadh_5,\alphadh_6)
\\&
=
\int \extder \adh_1\,\extder\alphadh_6\,\extder \alphadh_5 \extdersep 
\eexp^{-i\ell_1\adh_1-i\alphadh_5\lambda_5-i\alphadh_6\lambda_6}
{}\cdot 
\chi(\adh_1 - \ah_1, \alphadh_5 - \alphah_5, \alphadh_6 - \alphah_6) 
\\&\kern3.5em
{}\cdot
\eexp^{-i\ell_2\ah_2-i\ell_3\ah_3-i\alphah_4\lambda_4}
\\
\intertext{\footnotesize change of variables $\bar a_1=\adh_1-\ah_1$, $\bar\alpha_5=\alphadh_5-\alphah_5$, $\bar\alpha_6=\alphadh_6-\alphah_6$}
&
=
\int \extder \bar a_1\,\extder\bar\alpha_6\,\extder \bar\alpha_5 \extdersep 
\eexp^{-i\ell_1\bar a_1-i\bar\alpha_5\lambda_5-i\bar\alpha_6\lambda_6}
{}\cdot 
\chi(\bar a_1, \bar\alpha_5, \bar\alpha_6) 
\\&\kern3.5em
{}\cdot
\eexp^{-i\ell_1\ah_1-i\ell_2\ah_2-i\ell_3\ah_3-i\alphah_4\lambda_4-i\alphah_5\lambda_5-i\alphah_6\lambda_6}
\\&
=
\eexp^{-i\ell_1\ah_1-i\ell_2\ah_2-i\ell_3\ah_3-i\alphah_4\lambda_4-i\alphah_5\lambda_5-i\alphah_6\lambda_6}
\cdot
(\Fourierodd_{1\vert2} \chi)(\ell_1,\lambda_5,\lambda_6)
\\&
=
\rhoh^{(0,\ell_2,\ell_3, \lambda_4,0),(\ell_1,\lambda_5,\lambda_6)}_\grel \bigl((\Fourierodd_{1\vert2}\chi)(\ell_1,\lambda_5,\lambda_6)\bigr)
\mapob.
\end{align*}}%
This shows that $\Fourierodd_{1\vert2}$ intertwines the representation $\rhoh^{(0,\ell_2,\ell_3, \lambda_4,0)}$ of $G$ on $\Schwartz12$ with the direct (Berezin-)integral (over one even and two odd variables) of the $1\vert2$-parameter family of representations $\rhoh^{(0,\ell_2,\ell_3, \lambda_4,0),(\ell_1,\lambda_5,\lambda_6)}$ of $G$ on $\CC$.

\begin{remarque}{Some additional remarks concerning representation spaces}
As in the case of $\lambda_7=0$ and $\ell_0\neq0$, the individual elements of the subfamily $\rhoh^{(0,\ell_2,\ell_3, \lambda_4,0)}$ preserve the inner product $\inprodsym_{1\vert2}$ on $\Schwartz12$, so they are all (separately) \psur s. 
And then it is easy to show that our decomposition satisfies the requirements of a $1\vert2$-parameter family of \psur s. 
\end{remarque}

\medskip\noindent\stress{To summarise.}
We have decomposed the regular representation in the following invariant \myquote{subspaces}:
\begin{moneq}
\begin{tikzpicture}[line width=0.6pt]
\small
\def\TZVARecart{5}\def\TZVARfirsth{2}\def\TZVARsecondh{3}\def\TZVARthirdh{6}
\def\TZVARbatonl{\vrule width0pt height0ex depth1.75ex}
\def\TZVARbatonh{\vrule width0pt height3ex depth0ex}

\draw[->] (0,\TZVARthirdh) node[above]{$\TZVARbatonl\rho\text{ on }\Schwartz44$} --++(0,-1) node[midway,right]{$\int \scriptstyle \extder \ell_0 \,\extder \ell_2\,\extder \ell_3 \,\extder \lambda_7 \,\extder \lambda_4$} node[below]{$\TZVARbatonh\rhoh^{(\ell_0,\ell_2,\ell_3,\lambda_4,\lambda_7)}\text{ on }\Schwartz12$};

\draw[<->] (-\TZVARecart,\TZVARsecondh)--++(0,0.5)--++(2*\TZVARecart,0)--++(0,-0.5);
\draw[<-] (0,\TZVARsecondh)--++(0,1);

\draw (-\TZVARecart,\TZVARfirsth) node[above]{$\ell_0\neq0=\lambda_7$} node[below]{$\rhoh^{(\ell_0,\ell_2,\ell_3,\lambda_4,0)}\text{ on }\Schwartz12$};

\draw (0,\TZVARfirsth) node[above]{$\ell_0=0\neq\lambda_7$} node[below]{$\rhoh^{(0,\ell_2,\ell_3,\lambda_4, \lambda_7)}\text{ on }\Schwartz12$};

\draw (\TZVARecart,\TZVARfirsth) node[above]{$\ell_0=0=\lambda_7$} node[below]{$\rhoh^{(0,\ell_2,\ell_3,\lambda_4,0)}\text{ on }\Schwartz12$};

\draw[->] (-\TZVARecart,1)--(-\TZVARecart,0) node[midway,right]{$\int \scriptstyle\extder\kappa$} node[below]{$\TZVARbatonh \rhoh^{(\ell_0,\ell_2,\ell_3,\lambda_4,0),+}\text{ on }\Schwartz11$};

\draw[->] (0,1)--(0,0) node[midway,right]{$\int \scriptstyle\extder\lambda_6$} node[below]{$\TZVARbatonh \rhoh^{(0,\ell_2,\ell_3,\lambda_4, \lambda_7),(\lambda_6)}\text{ on }\Schwartz11$};

\draw[->] (\TZVARecart,1)--(\TZVARecart,0) node[midway,right]{$\int \scriptstyle \extder\ell_1\,\extder\lambda_5\, \extder\lambda_6$} node[below]{$\TZVARbatonh \rhoh^{(0,\ell_2,\ell_3,\lambda_4,0),(\ell_1, \lambda_5,\lambda_6)}\text{ on }\CC$};

\end{tikzpicture}
\end{moneq}

\subsection{Irreducibility}

In order to prove that the \psur s $\rhoh^{(\ell_0,\ell_2,\ell_3,\lambda_4, \lambda_7)}$, $\rhoh^{(\ell_0,\ell_2,\ell_3, \lambda_4,0),\pm}$ and $\rhoh^{(0,\ell_2,\ell_3,\lambda_4, \lambda_7),\lambda_6}$ are irreducible (for $\ell_0\neq0\neq\lambda_7$), we proceed as in \recals{SUBSECTIrreducibilityOfSuperH}. 

\itemizedef
The $3\vert2$-parameter family $\rhoh^{(\ell_0,\ell_2,\ell_3,\lambda_4, \lambda_7)}$ of $G$ on $\Schwartz12$ is given by
\begin{align}
\shifttag{5em}
(\rhoh^{(\ell_0,\ell_2,\ell_3,\lambda_4, \lambda_7)}_\grel \psi)(\adh_1,\alphadh_5,\alphadh_6)
\notag
\\
&
=
\psi(\adh_1 - \ah_1, \alphadh_5 - \alphah_5, \alphadh_6 - \alphah_6) \cdot\eexp^{-i\ell_0( \ah_2\adh_1 - \frac12\alphah_5\alphadh_5 +\frac12 \alphah_6 \alphadh_6)}\cdot
\eexp^{-i\ell_0( \ch-\ah_1\ah_2 )}
\notag
\\&
\kern3em
\cdot
\eexp^{-i(\alphah_4\adh_1 - \ah_3\alphadh_5)\lambda_7}
\cdot
\eexp^{-i(\gammah-\ah_1\alphah_4+\ah_3\alphah_5)\lambda_7}
\cdot
\eexp^{-i\ell_2\ah_2-i\ell_3\ah_3-i\alphah_4\lambda_4}
\mapob.
\notag
\end{align}%
The infinitesimal actions that do not reduce to a multiple of the identity are
\begin{align*}
\tau(\partial_{\ah_1})\psi
&=
-\partial_{\adh_1}\psi
&
\tau(\partial_{\alphah_4})\psi
&=
-i\lambda_7\adh_1\psi - i\lambda_4\psi
\\
\tau(\partial_{\ah_2})\psi
&=
-i\ell_0\adh_1\psi - i\ell_2\psi
&
\tau(\partial_{\alphah_5})\psi
&=
-\partial_{\alphadh_5}\psi + \tfrac12i\ell_0\alphadh_5\psi
\\
\tau(\partial_{\ah_3})\psi
&=
-i\lambda_7\alphadh_5\psi - i\ell_3\psi
&
\tau(\partial_{\alphah_6})\psi
&=
-\partial_{\alphadh_6}\psi - \tfrac12i\ell_0\alphadh_6\psi
\end{align*}
If $\psi$ belongs to an invariant closed subspace $F\subset \Schwartz12$, it follows that the elements
\begin{moneq}
\partial_{\adh_1}\psi
\ ,\quad
\ell_0\ah_1\psi
\ ,\quad
\lambda_7\alphadh_5\psi
\ ,\quad
\lambda_7\adh_1\psi
\ ,\quad
-\partial_{\alphadh_5}\psi + \tfrac12i\ell_0\alphadh_5\psi
\ ,\quad
-\partial_{\alphadh_6}\psi - \tfrac12i\ell_0\alphadh_6\psi
\end{moneq}
should belong to $F\otimes\CA^\CC$. 
When $\ell_0\neq0\neq\lambda_7$, this implies that
\begin{moneq}[InvarianceOperationsForGWithTwoNonZero]
\partial_{\adh_1}\psi
\ ,\quad
\adh_1\psi
\ ,\quad
\alphadh_5\psi
\ ,\quad
\partial_{\alphadh_5}\psi
\ ,\quad
\partial_{\alphadh_6}\psi + \tfrac12i\ell_0\alphadh_6\psi
\end{moneq}
should belong to $F$.\footnote{Here we use the standard abuse of notation by writing $\adh_1\psi$ for the function $(\adh_1,\alphadh_5,\alphadh_6)\mapsto \adh_1\cdot \psi(\adh_1,\alphadh_5,\alphadh_6)$ and similarly for $\alphadh_5\psi$ and the others.}
Now let $\psi\cong(\psi_0,\psi_5,\psi_6,\psi_{56})\in \Schwartz{}{}^4$, $\psi=\psi_0+\alphadh_5\psi_5+\alphadh_6\psi_6+\alphadh_5\alphadh_6\psi_{56}$ belong to the invariant graded subspace $F$, then (by definition of a graded subspace), the elements $(\psi_0,0,0,\psi_{56})$ and $(0,\psi_5,\psi_6,0)$ (the even and odd parts of $\psi$) also belong to $F$. 
Applying any sequence of operations given in \recalf{InvarianceOperationsForGWithTwoNonZero} thus still belongs to $F$, in particular the sequences (in this order) \myquote{5\numeral{th}, 4\numeral{th}, 5\numeral{th}, 3\numeral{d}} and \myquote{5\numeral{th}, 3\numeral{d}, 5\numeral{th}, 4\numeral{th}.} 
This tells us that
\begin{align*}
(0,-\tfrac12 i\ell_0\psi_5,0,0)
&
=
\alphadh_5\scirc (\partial_{\alphadh_6}\psi + \tfrac12i\ell_0\alphadh_6\psi)\scirc \partial_{\alphadh_5}\scirc(\partial_{\alphadh_6}\psi + \tfrac12i\ell_0\alphadh_6\psi)(0,\psi_5,\psi_6,0)
\\
(0, 0, -\tfrac12i\ell_0\psi_6, 0)
&
=
\partial_{\alphadh_5}\scirc (\partial_{\alphadh_6}\psi + \tfrac12i\ell_0\alphadh_6\psi)\scirc \alphadh_5\scirc (\partial_{\alphadh_6}\psi + \tfrac12i\ell_0\alphadh_6\psi)(0,\psi_5,\psi_6,0)
\end{align*}
belong to $F$ when $\psi=(\psi_0,\psi_5,\psi_6,\psi_{56})$ belongs. 
A similar argument shows that $(-\tfrac12i\ell_0\psi_0,0,0, 0)$ and $(0, 0, 0, -\tfrac12i\ell_0\psi_{56})$ belong to $F$ when $\psi$ does. 
In short, if $\psi = (\psi_0,\psi_5,\psi_6,\psi_{56})$ belongs to the invariant subspace $F$, so do the separate components. 

When we concentrate on the first component (the argument will be similar for the others), we note that 
\begin{moneq}
F^{(0)}=F\cap \{\,(\psi_0,0,0,0)\mid \psi_0\in \Schwartz{}{}\,\}
\end{moneq}
is a closed subspace of $(\Schwartz12, \inprodsym_{1\vert2})$. 
By the previous argument it is also invariant. 
But when we look at the first two operations in \recalf{InvarianceOperationsForGWithTwoNonZero}, we note that this implies that $F^{(0)}$ is a closed subspace of $\Schwartz{}{}$ that is invariant under the standard representation of the classical Heisenberg algebra, which is irreducible. 
Hence $F^{(0)}$ is trivial: either $F^{(0)}=\{0\}$ or $F^{(0)}=\Schwartz{}{}$. 
As the same holds for the other components, we conclude that if a component of $\psi\in F$ is non-zero, this component can take all values in $\Schwartz{}{}$. 

Now suppose $(\psi_0,0,0,0)\in F$ with $\psi_0\neq0$, then 
\begin{moneq}
\alphadh_5\scirc (\partial_{\alphadh_6}\psi + \tfrac12i\ell_0\alphadh_6\psi)(\psi_0,0,0,0) = (0,0,0,\tfrac12 i\ell_0\psi_0)
\end{moneq}
also belongs to $F$, meaning that there is an element with the $(56)$-component non-zero. 
A similar argument shows the converse, implying that if the even part of $F$ contains a non-zero element, this even part should be the whole even part of $\Schwartz12$. 
As a similar argument applies to the odd part, it follows that if the odd part contains a non-zero element, it should be the whole odd part of $\Schwartz12$. 
And then, when either the odd or the even part contains a non-zero element, we can use the terms $\alphadh_5\psi$ or $\partial_{\alphadh_5}\psi$ in \recalf{InvarianceOperationsForGWithTwoNonZero} to conclude that the odd or the even part also contains a non-zero element. 
It follows that the \psur s $\rhoh^{(\ell_0,\ell_2,\ell_3,\lambda_4, \lambda_7)}$ with $\ell_0\neq0\neq\lambda_7$ are all irreducible.

\itemizedef
The two $3\vert1$-families of \psur s $\rhoh^{(\ell_0,\ell_2,\ell_3, \lambda_4,0),\pm}$ of $G$ on $\Schwartz11$ with $\ell_0\neq0$ are given by
\begin{align}
(\rhoh^{(\ell_0,\ell_2,\ell_3, \lambda_4,0),\epsilon}_\grel \psi)(\adh_1, \xi)
&
=
\psi\bigl(\adh_1 - \ah_1, \xi - (\alphah_5+\epsilon\alphah_6)\bigr) 
{}\cdot
\eexp^{-i\ell_0( \ah_2\adh_1 - \frac12(\alphah_5-\epsilon\alphah_6)\xi )}
\notag
\\&
\kern3.5em
{}\cdot
\eexp^{-i\ell_0( \ch-\ah_1\ah_2 +\frac12\epsilon\alphah_5\alphah_6)}
\cdot
\eexp^{-i\ell_2\ah_2-i\ell_3\ah_3-i\alphah_4\lambda_4}
\mapob.
\notag
\end{align}%
The infinitesimal actions that do not reduce to a multiple of the identity are
\begin{align*}
\tau(\partial_{\ah_1})\psi
&=
-\partial_{\adh_1}\psi
&
\tau(\partial_{\alphah_5})\psi
&=
-\partial_{\xi}\psi + \tfrac12i\ell_0\xi \psi
\\
\tau(\partial_{\ah_2})\psi
&=
-i\ell_0\adh_1\psi - i\ell_2\psi
&
\tau(\partial_{\alphah_6})\psi
&=
-\epsilon\partial_{\xi}\psi - \tfrac12i\ell_0\epsilon\xi \psi
\end{align*}
If $\psi$ belongs to an invariant closed subspace $F\subset \Schwartz11$, it follows that the elements
\begin{moneq}
\partial_{\adh_1}\psi
\ ,\quad
\ell_0\ah_1\psi
\ ,\quad
-\partial_{\xi}\psi + \tfrac12i\ell_0\xi \psi
\ ,\quad
\partial_{\xi}\psi + \tfrac12i\ell_0\xi \psi
\end{moneq}
should belong to this subspace. 
When $\ell_0\neq0$, this implies that
\begin{moneq}
\partial_{\adh_1}\psi
\ ,\quad
\adh_1\psi
\ ,\quad
\xi \psi
\ ,\quad
\partial_{\xi}\psi
\end{moneq}
should belong to this subspace.
Applying the last two operators in both orders to $\psi\cong (\psi_0,\psi_1)$ tells us that both $(\psi_0,0)$ and $(0,\psi_1)$ belong to the invariant graded subspace. 
But then the first two operators, representing the irreducible representation of the Heisenberg algebra\slash group tells us that each component should belong to the whole space $\Schwartz{}{}$.
As the actions $\xi\psi$ and $\partial_\xi\psi$ interchanges the even and odd terms, it follows that these \psur s are irreducible.

\itemizedef
The $2\vert3$-family of \psur s $\rhoh^{(0,\ell_2,\ell_3,\lambda_4, \lambda_7),\lambda_6}$ of $G$ on $\Schwartz11$ is given by
\begin{align}
(\rhoh^{(0,\ell_2,\ell_3,\lambda_4, \lambda_7),(\lambda_6)}_\grel \psi)(\adh_1,\alphadh_5)
&
=
\psi(\adh_1-\ah_1, \alphadh_5-\alphah_5) 
\cdot
\eexp^{-i(\alphah_4\adh_1 - \ah_3\alphadh_5)\lambda_7}
\notag
\\&
\kern2em
{}\cdot
\eexp^{-i(\gammah-\ah_1\alphah_4+\ah_3\alphah_5)\lambda_7}
\cdot
\eexp^{-i\ell_2\ah_2-i\ell_3\ah_3-i\alphah_4\lambda_4-i\alphah_6\lambda_6}
\mapob.
\notag
\end{align}%
The infinitesimal actions that do not reduce to a multiple of the identity are
\begin{align*}
\tau(\partial_{\ah_1})\psi
&=
-\partial_{\adh_1}\psi
&
\tau(\partial_{\alphah_4})\psi
&=
-i\lambda_7\adh_1\psi - i\lambda_4\psi
\\
\tau(\partial_{\ah_3})\psi
&=
-i\lambda_7\alphadh_5\psi - i\ell_3\psi
&
\tau(\partial_{\alphah_5})\psi
&=
-\partial_{\alphadh_5}\psi
\end{align*}
If $\psi$ belongs to a closed invariant subspace $F\subset \Schwartz11$, it follows that the elements
\begin{moneq}
\partial_{\adh_1}h
\ ,\quad
\lambda_7\alphadh_5h
\ ,\quad
\lambda_7\adh_1\psi
\ ,\quad
-\partial_{\alphadh_5}h
\end{moneq}
should belong to $F\otimes \CA^\CC$. 
When $\lambda_7\neq0$, this implies that
\begin{moneq}
\partial_{\adh_1}\psi
\ ,\quad
\adh_1\psi
\ ,\quad
\alphadh_5\psi
\ ,\quad
\partial_{\alphadh_5}\psi
\end{moneq}
should belong to $F$.
As before, this shows that an invariant subspace must be trivial, proving irreducibility.

\subsection{\label{SECTDeterminationOfCoadjointOrbitsHlike}Representations associated to coadjoint orbits}

\BackSwitchHats
\renewcommand\grel{{(a,c,\alpha,\gamma)}}
\vrule width0pt
\vskip0.5\baselineskip

In this section we determine the coadjoint orbits of our \glgrp/ $G$.  It turns out that there are four families of orbits: a family of $0\vert0$ dimensional orbits with the trivial symplectic form, a family of $2\vert2$ dimensional orbits with an even symplectic form, a family of $2\vert2$ dimensional orbits with an odd symplectic form and a fourth family of $3\vert3$ dimensional orbits with a mixed symplectic form. The existence of orbits with a mixed symplectic form was already established in \cite{Tuynman..2010}.

A basis for the left-invariant vector fields (the Lie algebra $\Liealg g$ of $G$, isomorphic to $T_eG$) is given by $k_0=\partial_{c}$, $k_7=\partial_{\gamma}$ and
\begin{moneq}
\begin{aligned}
&e_1=\partial_{a_1} +  a_2\, \partial_{c} +  \alpha_4\, \partial_{\gamma}
&&,\quad
e_2=\partial_{a_2} 
&&,\quad
e_3=\partial_{a_3}
\\
&e_4=\partial_{\alpha_4}
&&,\quad
e_5=\partial_{\alpha_5} + \tfrac12 \alpha_5 \,\partial_{c} -  a_3 \,\partial_{\gamma}
&&,\quad
e_6=\partial_{\alpha_6} - \tfrac12 \alpha_6\,\partial_{c} 
\mapob.
\end{aligned}
\end{moneq}
In this list $k_0,e_1,e_2,e_3$ are even and $e_4,e_5,e_6,k_7$ are odd. 
A simple computation tells us that the only non-zero commutators among these basis vectors are given by
\begin{gather*}
[e_2,e_1] = - [e_1,e_2] =[e_5,e_5] = -[e_6,e_6] = k_0
\\
[e_4,e_1] = - [e_1,e_4] = [e_5,e_3] = - [e_3,e_5] = k_7
\rlap{\mapob.}
\end{gather*}

A direct computation shows that the triple product $g\gh g\mo$ (for $g=(a,c,\alpha,\gamma)$) is given by
\begin{multline*}
g\gh g\mo = (\ah_1, \ah_2, \ah_3, \ch + a_2 \ah_1 - a_1 \ah_2 + \alphah_5 \alpha_5 - \alphah_6 \alpha_6, 
\\
\alphah_4, \alphah_5, \alphah_6, \gammah + \ah_1 \alpha_4 - a_1 \alphah_4 + \ah_3 \alpha_5 - a_3 \alphah_5)
\mapob.
\end{multline*}
For the adjoint action we get the formulæ $\Ad(g)k_0 = k_0$, $\Ad(g)k_7 = k_7$ and
\begin{gather*}
\Ad(g)e_1 = e_1 + a_2 k_0 + \alpha_4 k_7
\quad,\quad
\Ad(g)e_2 = e_2 - a_1 k_0
\quad,\quad
\Ad(g)e_3 = e_3 + \alpha_5 k_7
\\
\Ad(g)e_4 = e_4 - a_1 k_7
\quad,\quad
\Ad(g)e_5 = e_5 + \alpha_5 k_0 - a_3 k_7
\quad,\quad
\Ad(g)e_6 = e_6 - \alpha_6 k_0
\mapob.
\end{gather*}

Seen as \agmfd/ the left-dual $\LD[1]{\Liealg g}$ has dimension $8\vert8$ with even coordinates $y_0$, $x_1$, $x_2$, $x_3$, $\xb_4$, $\xb_5$, $\xb_6$, $\yb_7$ and odd coordinates $\etab_0$, $\xib_1$, $\xib_2$, $\xib_3$, $\xi_4$, $\xi_5$, $\xi_6$, $\eta_7$. These are the even and odd parts of the left-coordinates of a point $\mu\in \LD[1]{\Liealg g}$ with respect to the (left-dual) basis $\LD[1]k_0$, $\LD[1]e_1$, $\LD[1]e_2$, $\LD[1]e_3$, $\LD[1]e_4$, $\LD[1]e_5$, $\LD[1]e_6$, $\LD[1]k_7$ of $\LD[1]{\Liealg g}$ determined by the formula
\begin{moneq}
\mu = 
(y_0+\etab_0)\cdot \LD[1]k_0 + \sum_{i=1}^3 
(x_i+\xib_i)\cdot \LD[1]e_i + 
\sum_{i=4}^6 (\xb_i+\xi_i)\cdot \LD[1]e_i + 
(\yb_7+\eta_7)\cdot \LD[1]k_7
\mapob.
\end{moneq}
The names are chosen such that the unbarred coordinates represent the even part $\mu_0$ of $\mu\in \LD[1]{\Liealg g}$ and the barred coordinates represent the odd part $\mu_1$, while latin characters represent even coordinates and greek characters represent odd coordinates. 

For the coadjoint action we obtain (from the adjoint action) that the basis vectors $\LD[1]e_i$ are left unchanged and that we have
\begin{align*}
\Coad(g)\LD[1]k_0 
&
= \LD[1]k_0 - \LD[1]e_1 a_2 + \LD[1]e_2 a_1 -\LD[1]e_5 \alpha_5 + \LD[1]e_6 \alpha_6
\\
\Coad(g)\LD[1]k_7 
&
= \LD[1]k_7 - \LD[1]e_1 \alpha_4 - \LD[1]e_3 \alpha_5 + \LD[1]e_4 a_1 + \LD[1]e_5 a_3
\mapob,
\end{align*}
which gives us in terms of coordinates (with as always $g=(a,c,\alpha,\gamma)$)
\begin{gather*}
x_1\mapsto x_1 - y_0 \cdot a_2 - \eta_7 \cdot \alpha_4
\quad ,\qquad
x_2 \mapsto x_2 + y_0 \cdot a_1
\quad ,\qquad
x_3\mapsto x_3 - \eta_7 \cdot \alpha_5
\\
\xi_4\mapsto \xi_4 + \eta_7\cdot a_1
\quad ,\qquad
\xi_5 \mapsto \xi_5 + y_0 \cdot \alpha_5 + \eta_7 \cdot a_3
\quad ,\qquad
\xi_6 \mapsto \xi_6 - y_0 \cdot \alpha_6 
\end{gather*}
and 
\begin{gather*}
\xib_1\mapsto \xib_1 - \etab_0 \cdot a_2 - \yb_7 \cdot \alpha_4
\quad ,\qquad
\xib_2 \mapsto \xib_2 + \etab_0 \cdot a_1
\quad ,\qquad
\xib_3\mapsto \xib_3 - \yb_7 \cdot \alpha_5
\\
\xb_4\mapsto \xb_4 + \yb_7\cdot a_1
\quad ,\qquad
\xb_5 \mapsto \xb_5 + \etab_0 \cdot \alpha_5 + \yb_7 \cdot a_3
\quad ,\qquad
\xb_6 \mapsto \xb_6 - \etab_0 \cdot \alpha_6 
\mapob,
\end{gather*}
while the coordinates $y_0$, $\yb_7$, $\etab_0$ and $\eta_7$ remain unchanged. \KorteOfLangeTekst{}{For the fundamental vector fields  $b^{\LD[1]{\Liealg g}}$ associated to the basis vectors $b=k_0, e_1$, \dots, $e_6, k_7\in \Liealg g$ we obtain $k_0^{\LD[1]{\Liealg g}} = k_7^{\LD[1]{\Liealg g}} = 0$ and
\begin{align*}
e_1^{\LD[1]{\Liealg g}} 
&
= 
-y_0\partial_{x_2} - \yb_7\partial_{\xb_4} - \etab_0\partial_{\xib_2} - \eta_7\partial_{\xi_4} 
&
e_2^{\LD[1]{\Liealg g}} 
&
= 
y_0\partial_{x_1} + \etab_0\partial_{\xib_1}
\\
e_3^{\LD[1]{\Liealg g}} 
&
= 
- \yb_7\partial_{\xb_5} -\eta_7\partial_{\xi_5} 
&
e_4^{\LD[1]{\Liealg g}} 
&
= 
\yb_7\partial_{\xib_1}-\eta_7\partial_{x_1}  
\\
e_5^{\LD[1]{\Liealg g}} 
&
= 
\yb_7\partial_{\xib_3} - y_0\partial_{\xi_5} -\eta_7\partial_{x_3} +\etab_0\partial_{\xb_5}
&
e_6^{\LD[1]{\Liealg g}} 
&
= 
y_0\partial_{\xi_6} - \etab_0\partial_{\xb_6} 
\mapob.
\end{align*}
}%

We now look at the orbit $\Orbit_{\mu_o}$ through a fixed point $\mu_o\in \LD[1]{\Liealg g}$ with real coordinates (which we need if we want to be sure that the orbit is a bona fide \gmfd/ embedded in $\LD[1]{\Liealg g}$). 
Since the coordinates of $\mu_o$ are real, we have in particular that the odd coordinates of $\mu_o$ are zero: $\xi_i(\mu_o) = \xib_j(\mu_o) = \etab_0(\mu_o) = \eta_7(\mu_o)=0$. 
Since the coordinates $y_0$ and $\yb_7$ are constant on the orbit, they remain the same real constant. 
However, the other even coordinates are in general not constant on $\Orbit_{\mu_o}$, so we will denote the remaining even coordinates of $\mu_o$ by $\xnul_1, \xnul_2, \xnul_3, \xnulb_4, \xnulb_5, \xnulb_6\in \RR$. 
The odd coordinates $\etab_0 = \eta_7 = 0$ are also invariant under the coadjoint action, so they remain constant equal zero on the orbit. 
Hence the formulæ for the coadjoint action reduce to:
\begin{gather*}
x_1 \mapsto x_1 - y_0 \cdot a_2
\ ,\ \ 
x_2 \mapsto x_2 + y_0 \cdot a_1
\ ,\ \ 
\xi_5 \mapsto \xi_5 + y_0 \cdot \alpha_5 
\ ,\ \ 
\xi_6 \mapsto \xi_6 - y_0 \cdot \alpha_6 
\\
\xib_1 \mapsto \xib_1 - \yb_7 \cdot \alpha_4
\ ,\ \ 
\xib_3 \mapsto \xib_3 - \yb_7 \cdot \alpha_5
\ ,\ \ 
\xb_4 \mapsto \xb_4 + \yb_7 \cdot a_1 
\ ,\ \ 
\xb_5 \mapsto \xb_5 + \yb_7 \cdot a_3
\rlap{\mapob,}
\end{gather*}
while the other $8$ coordinates remain unchanged. 
For the fundamental vector fields $b^{\LD[1]{\Liealg g}}$ associated to the basis vectors $b=k_0, e_i, k_7\in \Liealg g$ we find $k_0^{\LD[1]{\Liealg g}} = k_7^{\LD[1]{\Liealg g}} = 0$ and
\begin{align*}
e_1^{\LD[1]{\Liealg g}} &= -y_0\,\partial_{x_2} - \yb_7\,\partial_{\xb_4}
\ ,\quad&
e_2^{\LD[1]{\Liealg g}} &= y_0\,\partial_{x_1}
\ ,\quad&
e_3^{\LD[1]{\Liealg g}} &= -\yb_7\,\partial_{\xb_5}
\\
e_4^{\LD[1]{\Liealg g}} &= \yb_7\,\partial_{\xib_1} 
\ ,\quad&
e_5^{\LD[1]{\Liealg g}} &= -y_0\,\partial_{\xi_5} + \yb_7\,\partial_{\xib_3}
\ ,\quad&
e_6^{\LD[1]{\Liealg g}} &= y_0\,\partial_{\xi_6}
\end{align*}
We now distinguish four cases: (i) $y_0=\yb_7=0$, (ii) $y_0 \neq 0$ (but real!) and $\yb_7 = 0$,
(iii) $\yb_7 \neq 0$ and $y_0 = 0$, and (iv) $y_0 \cdot \yb_7 \neq 0$.

\subhead{Case (i)}
For $y_0=\yb_7=0$, the base point $\mu_o$ of the orbit $\Orbit_{\mu_o}$ is determined by the $6$ real values $\xnul_1, \xnul_2, \xnul_3, \xnulb_4, \xnulb_5, \xnulb_6$, while all other coordinates are zero. Since the coadjoint action reduces to the trivial action on such a point, $\Orbit_{\mu_o}=\{\mu_o\}$ is a point and its symplectic form is the $2$-form which is identically zero (which is indeed symplectic as there are no non-zero tangent vectors). This gives us a $6$-dimensional family of $0$-dimensional orbits labeled by $\xnul_1, \xnul_2, \xnul_3, \xnulb_4, \xnulb_5, \xnulb_6$.

\subhead{Case (ii)}
For $\yb_7=0$ and $y_0\neq0$, the base point $\mu_o$ is determined by the $7$ real values $\xnul_1$, $\xnul_2$, $\xnul_3$, $\xnulb_4$, $\xnulb_5$, $\xnulb_6$, $y_0$, while all other coordinates are zero. The orbit has dimension $2|2$ with even coordinates $x_1, x_2$ and odd coordinates $\xi_5, \xi_6$ (all other coordinates on $\LD[1]{\Liealg g}$ remain unchanged for such an orbit). Substituting $(v,w) = (e_1,e_2)$, $(e_5,e_5)$ and $(e_6,e_6)$ in \recalf{DefiningFormulaForKKSForm} \KorteOfLangeTekst{allow us to}{we find the following identities
\begin{moneq}
\contrf{\fracp{}{x_1}, \fracp{}{x_2}}{\omega} = \frac{1}{y_0}
\quad,\quad
\contrf{\fracp{}{\xi_5}, \fracp{}{\xi_5}}{\omega} = -\frac{1}{y_0}
\quad,\quad
\contrf{\fracp{}{\xi_6}, \fracp{}{\xi_6}}{\omega} = \frac{1}{y_0}
\mapob,
\end{moneq}
all other evaluations on basis tangent vectors being either zero or determined by graded skew-symmetry.
From these identities we} deduce that the (even) symplectic form is given as
\begin{moneq}
\omega = (y_0)\mo \cdot \bigl(\ \extder x_2 \wedge \extder x_1 - \tfrac12 \,\extder\xi_5 \wedge \extder\xi_5 +
\tfrac12 \,\extder\xi_6 \wedge \extder\xi_6\ \bigr)
\mapob.
\end{moneq}
Since $x_1$ and $x_2$ are coordinates on such an orbit, the different orbits are labeled by $\xnul_3$, $\xnulb_4$, $\xnulb_5$, $\xnulb_6$, $y_0$, thus giving a $5$-dimensional family of $2\vert2$-dimensional orbits with an even symplectic form.

\subhead{Case (iii)}
For $y_0=0$ and $\yb_7\neq0$, the base point $\mu_o$ is again determined by $7$ real values, but now  $\xnul_1$, $\xnul_2$, $\xnul_3$, $\xnulb_4$, $\xnulb_5$, $\xnulb_6$, $\yb_7$, while all other coordinates are zero.  The orbit dimension is still $2|2$, but now with even
coordinates $\xb_4,\xb_5$ and odd coordinates $\xib_1,\xib_3$ (as before, the other coordinates on $\LD[1]{\Liealg g}$ remain unchanged on this orbit). 
Substituting $(v,w) = (e_1,e_4)$ and $(e_3,e_5)$ in \recalf{DefiningFormulaForKKSForm} \KorteOfLangeTekst{allow us to}{we find the following identities
\begin{moneq}
\contrf{\fracp{}{\xib_1}, \fracp{}{\xb_4}}{\omega} = \frac{1}{\yb_7}
\quad,\quad
\contrf{\fracp{}{\xib_3}, \fracp{}{\xb_5}}{\omega} = \frac{1}{\yb_7}
\mapob,
\end{moneq}
while all other evaluations on basis tangent vectors being either zero or determined by graded skew-symmetry.
From this we} deduce that the (odd) symplectic form is given as
\begin{moneq}
\omega = (\yb_7)\mo \cdot \bigl(\ \extder\xb_4 \wedge \extder \xib_1 + \extder\xb_5 \wedge \extder\xib_3\ \bigr)
\mapob.
\end{moneq}
Here the $\xb_4$ and $\xb_5$ are coordinates on the orbit, so the different orbits are labeled by $\xnul_1$, $\xnul_2$, $\xnul_3$, $\xnulb_6$, $\yb_7$, giving a $5$-dimensional family of $2\vert2$-dimensional orbits with an odd symplectic form.

\subhead{Case (iv)}
In the last case $y_0\yb_7\neq0$ we have to be slightly more careful. The base point is determined by the $8$ real values $\xnul_1$, $\xnul_2$, $\xnul_3$, $\xnulb_4$, $\xnulb_5$, $\xnulb_6$, $y_0$, $\yb_7$, while all other coordinates are zero. But the combinations $s = y_0\xb_4 - \yb_7x_2$ and $\sigma = y_0\xib_3+\yb_7\xi_5$ are invariant under the group action. Since the base point of the orbit has real coordinates, $\sigma$ is zero at this point, and thus zero on the whole orbit. We introduce the change of coordinates on $\LD[1]{\Liealg g}$ involving only $x_2, \xb_4, \xib_3, \xi_5$ which change into $\xh_2, s, \sigma, \xih_5$ with $\xh_2 = x_2$ and $\xih_5 = \xi_5$. We then can use the coordinates $x_1, \xh_2, \xb_5, \xib_1, \xih_5, \xi_6$ as independent coordinates on the orbit, which thus has dimension $3\vert3$. 
In terms of these new coordinates (on $\LD[1]{\Liealg g}$) the fundamental vector fields $e_i^{\LD[1]{\Liealg g}}$ are given by 
\begin{align*}
e_1^{\LD[1]{\Liealg g}} &= -y_0\,\partial_{\xh_2}
\ ,\quad&
e_2^{\LD[1]{\Liealg g}} &= y_0\,\partial_{x_1}
\ ,\quad&
e_3^{\LD[1]{\Liealg g}} &= -\yb_7\,\partial_{\xb_5}
\\
e_4^{\LD[1]{\Liealg g}} &= \yb_7\,\partial_{\xib_1} 
\ ,\quad&
e_5^{\LD[1]{\Liealg g}} &= -y_0\,\partial_{\xih_5}
\ ,\quad&
e_6^{\LD[1]{\Liealg g}} &= y_0\,\partial_{\xi_6}
\mapob.
\end{align*}
As before, substituting suitable basis vectors for $(v,w)$ in the
formula \recalf{DefiningFormulaForKKSForm} for the symplectic form \KorteOfLangeTekst{allows us to obtain}{gives us the following identities
\begin{gather*}
\contrf{\fracp{}{x_1}, \fracp{}{\hat x_2} }{\omega} = \frac{1}{y_0}
\quad,\quad
\contrf{\fracp{}{\xib_1}, \fracp{}{\hat x_2} }{\omega} = \frac{1}{y_0}
\quad,\quad
\contrf{\fracp{}{\xb_5}, \fracp{}{\xih_5} }{\omega} = \frac{1}{y_0}
\\
\noalign{\vskip2\jot}
\contrf{\fracp{}{\xih_5}, \fracp{}{\xih_5}}{\omega} =
\frac{-1}{y_0}
\quad,\quad
\contrf{\fracp{}{\xi_6}, \fracp{}{\xi_6}}{\omega} = \frac{1}{y_0}
\mapob.
\end{gather*}
This results in} the (non-homogeneous) symplectic form
\begin{moneq}
\omega = (y_0)\mo \bigl(\,\extder\xh_2 \wedge \extder x_1 -\tfrac12\, \extder\xih_5 \wedge \extder\xih_5 + \tfrac12 \,\extder\xi_6 \wedge \extder \xi_6 + \extder\xh_2\wedge \extder\xib_1  + \extder\xih_5 \wedge \extder\xb_5\ \bigr)
\mapob.
\end{moneq}
In this case the invariants of the orbit are the components $\xnul_3$, $\xnulb_6$, $y_0$, $\yb_7$ and the combination $s^o = y_0\xnul_2 -\yb_7\xnulb_4$, giving a $5$-dimensional family of $3\vert3$-dimensional orbits with a mixed symplectic form.

\medskip\noindent\stress{Quantization of point orbits.}
Since the orbit is a point, the prequantum bundle is $\{\mu_o\} \times (\CA_0/d\ZZ)\times \CA_1$ together with the connection $1$-form $\Gamma = \extder t + \extder \tau$. Since the canonical momentum map takes the (constant) value $\mu_o$, the lifted vector fields $b^{Y}$ (determined by the conditions \recalf{LiftingEquationsOfGroupActionOnPrequantumBundle}) are given by $k_0^Y = k_7^Y = 0$ and
\begin{moneq}
e_i^Y = \xnul_i\,\partial_t
\ ,\ i=1,2,3
\qquad\text{and}\qquad
e_i^Y = \xnulb_i\,\partial_\tau
\ ,\ i=4,5,6
\mapob.
\end{moneq}
It follows immediately that the lifted action $\phi$ is given by
\begin{moneq}[liftedactioncasei]
\phi_\grel: \begin{pmatrix} t\\ \tau\end{pmatrix} \mapsto 
\begin{pmatrix} t  -a_1\xnul_1-a_2\xnul_2-a_3\xnul_3 \\
 \tau  -\alpha_4 \xnulb_4 -\alpha_5 \xnulb_5 -\alpha_6 \xnulb_6\end{pmatrix}
\mapob.
\end{moneq}

For a point orbit $\mu_o$, the stabilizer subgroup $G_{\mu_o}$ is the whole group. A polarization subalgebra $\Liealg h$ thus necessarily equals $\Liealg g$. 
Smooth functions $\psi$ on $Y=(\CA_0/d\ZZ) \times \CA_1$ satisfying condition \recalf{firstGQcondition} are determined by a complex number $C\in \CC$ as
\begin{moneq}
\psi(t,\tau) = C\cdot \eexp^{i\tau\kappa + it/\hbar}
\mapob,
\end{moneq}
whereas condition \recalf{secondGQcondition} is void. 
These functions constitute a $1$-dimensional graded vector space of dimension $1\vert0$. 
The induced action $\Phi$ of $G$ \recalf{InducedActionOnPrequantumFunctions} on such a function is given (using \recalf{liftedactioncasei}) as
\begin{align}
(\Phi_\grel \psi)(t,\tau)
&
=
\psi\bigl(\phi_\grel\mo(t,\tau)\bigr)
\notag
\\&
=
\psi(t  +a_1\xnul_1+a_2\xnul_2+a_3\xnul_3 ,
 \tau  +\alpha_4 \xnulb_4 +\alpha_5 \xnulb_5 +\alpha_6 \xnulb_6)
 \notag
\\&
=
\eexp^{i(\alpha_4 \xnulb_4 +\alpha_5 \xnulb_5 +\alpha_6 \xnulb_6)\kappa +
i(a_1\xnul_1+a_2\xnul_2+a_3\xnul_3)/\hbar} \cdot \psi(t,\tau)
\mapob.
\label{00dimorbitaction}
\end{align}

\medskip\noindent\stress{Quantization of orbits of dimension \texorpdfstring{$2\vert2$}{22} with an even symplectic form.}
As said in \recals{GQgeneral}, the prequantum bundle is just the direct product $Y=\Orbit_{\mu_o} \times (\CA_0/d\ZZ) \times \CA_1$. Using the explicit expression of the symplectic form, it is not hard to see that a connection form $\Gamma$ is given by
\begin{moneq}
\Gamma_0 
= 
\extder t + {(y_0)}\mo ( -x_1\, \extder x_2 - \tfrac12 \xi_5\, \extder \xi_5 + \tfrac12 \xi_6 \,\extder \xi_6 ) 
\qquad\text{and}\qquad
\Gamma_1 
=  
\extder \tau 
\mapob.
\end{moneq}
With this choice, the lifted fundamental vector fields $b^{Y}$ (determined by the conditions \recalf{LiftingEquationsOfGroupActionOnPrequantumBundle}) are given by $k_0^Y = y_0\,\partial_t$, $k_7^Y = 0$ and
\begin{align*}
e_1^Y & = -y_0\,\partial_{x_2}
&,\qquad
e_3^Y &= \xnul_3\,\partial_t
&,\qquad
e_5^Y &= -y_0\,\partial_{\xi_5} - \tfrac12\xi_5\,\partial_t +\xnulb_5\,\partial_\tau
\\
e_2^Y &= y_0\,\partial_{x_1} +x_2\,\partial_t
&,\qquad
e_4^Y & = \xnulb_4\,\partial_\tau
&,\qquad
e_6^Y &= y_0\,\partial_{\xi_6} - \tfrac12\xi_6\,\partial_t + \xnul_6\,\partial_\tau
\mapob,
\end{align*}
which integrate to an action $\phi$ of $G$ on $Y$ by
\begin{moneq}[liftedactioncaseii]
\phi_\grel:
\begin{pmatrix} x_1 \\ x_2  \\ \xi_5 \\ \xi_6 \\ t\\  \tau \end{pmatrix}
\mapsto
\begin{pmatrix}
x_1 - y_0a_2
\\ 
x_2 + y_0 a_1
\\
\xi_5 + y_0 \alpha_5
\\
\xi_6 - y_0 \alpha_6
\\
t- a_2x_2 - a_3 \xnul_3 + \tfrac12 \alpha_5 \xi_5 + \tfrac12 \alpha_6 \xi_6 -  c y_0
\\
\tau   - \alpha_4\xnulb_4 - \alpha_5\xnulb_5 -
\alpha_6 \xnulb_6 
\end{pmatrix}
\mapob.
\end{moneq}

For all these orbits the stabilizer subgroup $G_{\mu_o}$ is given by 
\begin{moneq}
G_{\mu_o} = \{\,(a, c, \alpha, \gamma )\mid a_1=a_2=\alpha_5=\alpha_6=0\,\}
\mapob.
\end{moneq}
The stabilizer subalgebra $\Liealg g_{\mu_o}$ of dimension $2\vert2$ thus is generated by the vectors $e_3, e_4, k_0, k_7$. If a subalgebra $\Liealg h\supset \Liealg g_{\mu_o}$ has even dimension $4$, then necessarily it contains all $4$ even vectors $e_1, e_2, e_3, k_0$. But then 
\begin{moneq}
\contrs{[e_1,e_2]}{\mu_o} = \contrs{-k_0}{\mu_o} = -y_0\neq0
\mapob.
\end{moneq}
Hence for a polarization subalgebra $\Liealg h$ the even dimension is at most $3$. 
Similarly, if $\Liealg h$ has odd dimension $4$, then it necessarily contains all odd vectors $e_4,e_5, e_6, k_7$ and in particular $e_5$. And then we have 
\begin{moneq}
\contrs{[e_5,e_5]}{\mu_o} = y_0\neq0
\mapob.
\end{moneq}
Since the same argument applies to $e_6$, it follows that the odd dimension of a polarization subalgebra $\Liealg h$ is at most $3$ and that $\Liealg h$ cannot contain the (isolated) vectors $e_5$ or $e_6$. 
On the other hand, it is easy to check that the subalgebra $\Liealg h$ of dimension $3\vert 3$ generated by the vectors $ae_1+be_2, e_3, e_4, e_5-\epsilon e_6, k_0, k_7$ with $\vert a\vert^2+ \vert b\vert^2\neq0$ and $\epsilon=\pm1$ satisfies the condition $\contrs{[\Liealg h, \Liealg h]}{\mu_o} = 0$, \ie, is a polarization subalgebra; with the given degrees of freedom in $a,b,\epsilon$ these are the only ones. 

We now choose the polarization subalgebra generated by $e_2, e_3, e_4, e_5-\epsilon e_6, k_0, k_7$.\footnote{Even though I did not check explicitly, I am convinced that other choices for $a,b\in\RR$ will yield equivalent representations. And even though I keep the $\epsilon$ dependence, we will see when comparing with the representations obtained in the decomposition of the left-regular representation that both give equivalent representations.} 
The horizontally lifted vector fields $b^h$, $b=e_2, e_3, e_4$, $e_5-\epsilon e_6, k_0, k_7$ determined by the conditions \recalf{DefinitionOfHorizontalLiftFundamentalVectorFields} are given by $k_0^h = k_7^h = e_3^h = e_4^h=0$ and
\begin{moneq}
e_2^h =  y_0\,\partial_{x_1}
\qquad\text{and}\qquad
e_5^h-\epsilon e_6^h =  -y_0\,\partial_{\xi_5} - \epsilon y_0\,\partial_{\xi_6} + \tfrac12(\xi_5-\epsilon\xi_6)\,\partial_t
\mapob.
\end{moneq}
Functions $\psi$ on the prequantum bundle $Y$ satisfying 
\KorteOfLangeTekst{\recalf{secondGQcondition} and \recalf{firstGQcondition}}{condition \recalf{firstGQcondition} are of the form
\begin{moneq}
\psi(x_1, x_2, \xi_5,\xi_6, t, \tau) = \eexp^{i\tau\kappa + it/\hbar} \cdot 
\psi(x_1, x_2, \xi_5,\xi_6,0,0)
\mapob. 
\end{moneq}
and condition \recalf{secondGQcondition} applied to functions of this form then gives the equations (at all points $(x_1, x_2, \xi_5,\xi_6,0,0)$)
\begin{moneq}
\fracp{\psi}{x_1} = 0
\qquad\text{and}\qquad
\fracp\psi{\xi_5} + \epsilon \fracp{\psi}{\xi_6} - i \cdot \Bigl( \frac{\xi_5-\epsilon\xi_6}{2y_0\hbar}  \Bigr) \cdot \psi = 0
\mapob.
\end{moneq}
The solutions to these equations} are determined by functions $\chi$ of one even and one odd variable according to
\begin{moneq}
\psi(x_1, x_2, \xi_5, \xi_6,t,\tau) = \eexp^{i\tau\kappa + it/\hbar} \cdot\eexp^{\frac{-i\epsilon}{2y_0\hbar}\xi_5\xi_6}\cdot \chi(x_2, \xi_5 - \epsilon \xi_6) 
\mapob.
\end{moneq}
It follows that the induced action $\Phi$ \recalf{InducedActionOnPrequantumFunctions} of $G$ on \KorteOfLangeTekst{such a function $\psi$}{a function $\psi$ on $Y$ satisfying \recalf{secondGQcondition} and \recalf{firstGQcondition}} is given (using \recalf{liftedactioncaseii}) by
\begin{align*}
\shifttag{5em}
(\Phi_\grel \psi)(x_1, x_2, 
\xi_5, \xi_6,t,\tau)
=
\psi\bigl(\phi_\grel\mo(x_1, x_2, 
\xi_5, \xi_6,t,\tau)\bigr)
\\&
=
\dontprintcomp{ 
\psi(x_1 + y_0a_2 ,
x_2 - y_0 a_1, 
\xi_5 - y_0 \alpha_5,
\xi_6 + y_0 \alpha_6,
\\&
\kern2em
t + a_2x_2 + a_3 \xnul_3 - \tfrac12 \alpha_5 \xi_5 - \tfrac12 \alpha_6 \xi_6 + 
c y_0 - y_0 a_1 a_2
,
\\&
\kern4em
\tau   + \alpha_4\xnulb_4 + \alpha_5\xnulb_5 +
\alpha_6 \xnulb_6 
)
\\&
=
\chi(x_2 - y_0 a_1, (\xi_5 - \epsilon \xi_6) - y_0 ( \alpha_5 + \epsilon \alpha_6)
)
\\&
\kern2em
\cdot
\eexp^{i(t + a_2x_2 + a_3 \xnul_3 - \frac12 \alpha_5 \xi_5 - \frac12 \alpha_6 \xi_6 + 
c y_0 - y_0 a_1 a_2
)/\hbar}
\\&
\kern4em
\cdot
\eexp^{i(\tau   + \alpha_4\xnulb_4 + \alpha_5\xnulb_5 +
\alpha_6 \xnulb_6 )\kappa}
\cdot
\eexp^{\frac{-i\epsilon}{2y_0\hbar}(\xi_5 - y_0 \alpha_5)(\xi_6 + y_0 \alpha_6)}
\\&
=
} 
\chi(x_2 - y_0 a_1, (\xi_5 - \epsilon \xi_6) - y_0 ( \alpha_5 + \epsilon \alpha_6)
)
\\&
\kern2em
\cdot
\eexp^{\frac{-i\epsilon}{2y_0\hbar}\xi_5\xi_6}
\cdot
\eexp^{i\tau\kappa + it/\hbar}
\cdot
\eexp^{i(\alpha_4\xnulb_4 + \alpha_5\xnulb_5 +
\alpha_6 \xnulb_6)\kappa}
\\&
\kern4em
\cdot
\eexp^{i(a_2x_2 + a_3 \xnul_3 -\frac12 (\alpha_5 - \epsilon \alpha_6)(\xi_5 - \epsilon\xi_6))/\hbar} 
\cdot
\eexp^{iy_0( 
c-a_1 a_2
+ \frac12\epsilon \alpha_5\alpha_6)/\hbar}
\mapob.
\end{align*}
\KorteOfLangeTekst{}{%
In terms of the functions $\chi$ the action thus becomes
\begin{align*}
(\Phi_\grel \chi)( x_2, 
\xi)
&
=
\chi\bigl(x_2 - y_0 a_1, \xi - y_0 ( \alpha_5 + \epsilon \alpha_6)
\bigr)
\\&
\kern2em
\cdot
\eexp^{i(a_2x_2 + a_3 \xnul_3 -\frac12 (\alpha_5 - \epsilon \alpha_6)\xi + 
c y_0 - y_0 a_1 a_2
+ \frac12\epsilon y_0 \alpha_5\alpha_6)/\hbar}
\\&
\kern4em
\cdot
\eexp^{i(\alpha_4\xnulb_4 + \alpha_5\xnulb_5 +
\alpha_6 \xnulb_6)\kappa}
\mapob.
\end{align*}
}%
If we introduce the change of variables $z_1 = (x_2 - \xnul_2)/y_0$ and $\zeta = \xi/y_0$, then, writing $\hat \chi(z_1,\zeta) = \chi(x_2,\xi)$, the action of $\grel\in G$ on the functions $\hat \chi$ becomes
\KorteOfLangeTekst{
\begin{align}
(\Phi_\grel \hat \chi)( z_1, 
\zeta)
&=
\hat \chi\bigl(z_1 -  a_1, \zeta -  ( \alpha_5 + \epsilon \alpha_6)
\bigr)
\cdot
\eexp^{i(a_2  z_1  -\frac12 (\alpha_5 - \epsilon \alpha_6)\zeta) y_0/\hbar}
\notag
\\&
\kern2em
\cdot
\eexp^{i(a_2 \xnul_2 + a_3 \xnul_3 )/\hbar}
\cdot
\eexp^{i(
c - a_1 a_2
+ \frac12\epsilon  \alpha_5\alpha_6)y_0/\hbar}
\cdot
\eexp^{i(\alpha_4\xnulb_4 + \alpha_5\xnulb_5 +
\alpha_6 \xnulb_6)\kappa}
\mapob.
\label{22evenaction}
\end{align}
}{
\begin{align}
(\Phi_\grel \hat \chi)( z_1, 
\zeta)
&=
\chi\bigl(x_2 - y_0 a_1, \xi - y_0 ( \alpha_5 + \epsilon \alpha_6)
\bigr)
\notag
\\&
\kern2em
\cdot
\eexp^{i(a_2x_2 + a_3 \xnul_3 -\frac12 (\alpha_5 - \epsilon \alpha_6)\xi + 
c y_0 - y_0 a_1 a_2
+ \frac12\epsilon y_0 \alpha_5\alpha_6)/\hbar}
\notag
\\&
\kern4em
\cdot
\eexp^{i(\alpha_4\xnulb_4 + \alpha_5\xnulb_5 +
\alpha_6 \xnulb_6)\kappa}
\notag
\\
&=
\hat \chi\bigl(z_1 -  a_1, \zeta -  ( \alpha_5 + \epsilon \alpha_6)
\bigr)
\cdot
\eexp^{i(a_2  z_1  -\frac12 (\alpha_5 - \epsilon \alpha_6)\zeta) y_0/\hbar}
\label{22evenaction}
\\&
\kern2em
\cdot
\eexp^{i(a_2 \xnul_2 + a_3 \xnul_3 )/\hbar}
\cdot
\eexp^{i(
c - a_1 a_2
+ \frac12\epsilon  \alpha_5\alpha_6)y_0/\hbar}
\cdot
\eexp^{i(\alpha_4\xnulb_4 + \alpha_5\xnulb_5 +
\alpha_6 \xnulb_6)\kappa}
\mapob.
\notag
\end{align}}%
This action depends upon the coordinates $\xnul_2$, $\xnul_3$, $\xnulb_4$, $\xnulb_5$, $\xnulb_6$, $y_0$ of the base point $\mu_o$. The last five are the parameters of the orbit, but the first is not an invariant of the orbit. 
To see that $\xnul_2$ is a spurious parameter in this action, it suffices to change the variable $z_1$ to $z_1 = x_2/y_0$, instead of $(x_2 - \xnul_2)/y_0$ in which case the $a_2 \xnul_2$ term in the exponential disappears. 
The \myquote{introduction} of the parameter $\xnul_2$ thus just gives an equivalent description of the same representation. We introduced it to facilitate future comparisons.

\medskip\noindent\stress{Quantization of orbits of dimension \texorpdfstring{$2\vert2$}{22} with an odd symplectic form.}
As before, the prequantum bundle is the direct product $Y=\Orbit_{\mu_o} \times (\CA_0/d\ZZ) \times \CA_1$. Here a connection form $\Gamma$ is given by
\begin{moneq}
\Gamma_0 
= \extder t 
\qquad\text{and}\qquad
\Gamma_1 
=  \extder \tau + {(\yb_7)}\mo ( \xb_5 \,\extder \xib_3 - \xib_1\, \extder \xb_4 ) 
\mapob.
\end{moneq}
The lifted fundamental vector fields $b^{Y}$ (determined by the conditions \recalf{LiftingEquationsOfGroupActionOnPrequantumBundle}) are given by $k_0^Y = 0$, $k_7^Y = \yb_7\,\partial_\tau$ and
\begin{align*}
e_1^Y & = -\yb_7\,\partial_{\xb_4} + \xnul_1\,\partial_t
&
e_3^Y &= -\yb_7\,\partial_{\xb_5} + \xnul_3\,\partial_t + \xib_3\,\partial_\tau
&
e_5^Y &= \yb_7\,\partial_{\xib_3}
\\
e_2^Y &= \xnul_2\,\partial_t
&
e_4^Y & = \yb_7\,\partial_{\xib_1} +\xb_4\,\partial_\tau
&
e_6^Y &=  \xnulb_6\,\partial_\tau
\mapob,
\end{align*}
which integrate to an action $\phi$ of $G$ on $Y$ by
\begin{moneq}[liftedactioncaseiii]
\phi_\grel:
\begin{pmatrix} \xib_1 \\ \xib_3 \\ \xb_4 \\ \xb_5 \\ t\\  \tau \end{pmatrix}
\mapsto
\begin{pmatrix}
\xib_1 - \yb_7 \alpha_4
\\
\xib_3 - \yb_7 \alpha_5
\\
\xb_4 + \yb_7a_1
\\ 
\xb_5+ \yb_7 a_3
\\
t- a_1\xnul_1 - a_2\xnul_2 - a_3 \xnul_3 
\\
\tau  
- a_3 \xib_3 - \alpha_4\xb_4 -
\alpha_6 \xnulb_6 - \gamma \yb_7 
\end{pmatrix}
\end{moneq}

For all orbits of this type the stabilizer subgroup $G_{\mu_o}$ is given by 
\begin{moneq}
G_{\mu_o} = \{(a, c, \alpha, \gamma)\mid a_1=a_3=\alpha_4=\alpha_5=0\}
\mapob.
\end{moneq}
The stabilizer subalgebra $\Liealg g_{\mu_o}$ of dimension $2\vert 2$ is generated by the vectors $e_2$, $e_6$, $k_0$, $k_7$. Since $y_0=0$, a polarization subalgebra $\Liealg h$ can have an even dimension $4$, \ie, containing all $4$ even vectors $e_1, e_2, e_3, k_0$. But then it cannot contain any combination of $e_4$ and $e_5$ as we have 
\begin{moneq}
\contrs{[ae_4+be_5, e_1]}{\mu_o} = a\yb_7
\qquad\text{and}\qquad
\contrs{[ae_4+be_5, e_3]}{\mu_o} = b\yb_7
\mapob.
\end{moneq}
Since $\yb_7\neq0$, the condition that these values are zero implies that we have $a=b=0$. The subalgebra $\Liealg h$ of dimension $4\vert 2$ generated by the vectors $e_1, e_2, e_3, e_6, k_0, k_7$ thus is a polarization subalgebra.

Similarly, a polarization subalgebra $\Liealg h$ can have an odd dimension $4$, \ie, containing all $4$ odd vectors $e_4, e_5,e_6, k_7$. But then it cannot contain any combination of $e_1$ and $e_3$ as we have 
\begin{moneq}
\contrs{[ae_1+be_3, e_4]}{\mu_o} = -a\yb_7
\qquad\text{and}\qquad
\contrs{[ae_1+be_3, e_5]}{\mu_o} = -b\yb_7
\mapob.
\end{moneq}
The subalgebra $\Liealg h$ of dimension $2\vert 4$ generated by the vectors $e_2, e_4, e_5, e_6, k_0, k_7$ thus is a polarization subalgebra.

If we now look at a subalgebra $\Liealg h\supset \Liealg g_{\mu_o}$ of dimension $3\vert 3$, it must be generated by $ae_1+be_3, e_2, ce_4+de_5, e_6, k_0, k_7$ for some constants $a,b,c,d$. The condition $\contrs{[\Liealg h, \Liealg h]}{\mu_o} = 0$ then implies in particular that we must have 
\begin{moneq}
\contrs{[ae_1+be_3, ce_4+de_5]}{\mu_o} = -ac-bd = 0
\mapob.
\end{moneq}%
Since the couples $(a,b)$ and $(c,d)$ can be multiplied by a constant without changing $\Liealg h$, it follows easily that the ``only'' solution is given by $c=b$ and $d=-a$ and that we have to require $\vert a\vert^2+\vert b\vert^2\neq0$. 
Since $\contrs{[\Liealg h, \Liealg h]}{\mu_o} = 0$ does not give any other conditions, it follows that the subalgebra $\Liealg h$ of dimension $3\vert 3$ generated by the vectors $ae_1+be_3, e_2, be_4-ae_5, e_6, k_0, k_7$ is a polarization subalgebra.
\KorteOfLangeTekst{}{\par As before, functions $\psi$ on the prequantum bundle $Y$ satisfying condition \recalf{firstGQcondition} are of the form
\begin{moneq}[caseiiifunctions]
\psi(\xb_4,\xb_5, \xib_1,\xib_3, t, \tau) = \eexp^{i\tau\kappa + it/\hbar} \cdot 
\psi(\xb_4,\xb_5, \xib_1,\xib_3,0,0)
\mapob. 
\end{moneq}
However, in order to proceed with quantization, we have to choose an invariant polarization. }

\medskip
\noindent$\bullet$
We start with the polarization subalgebra of dimension $4\vert2$ generated by $e_1, e_2, e_3$, $e_6, k_0, k_7$. 
The corresponding horizontally lifted vector fields determined by the conditions \recalf{DefinitionOfHorizontalLiftFundamentalVectorFields} are given by $k_0^h = k_7^h = e_2^h = e_6^h=0$ and
\begin{moneq}
e_1^h =  -\yb_7\,\partial_{\xb_4} - \xib_1\,\partial_\tau
\qquad\text{and}\qquad
e_3^h =  -\yb_7\,\partial_{\xb_5}
\mapob.
\end{moneq}
\KorteOfLangeTekst{Functions $\psi$ on the prequantum bundle $Y$ satisfying \recalf{secondGQcondition} and \recalf{firstGQcondition}}{Condition \recalf{secondGQcondition} applied to functions of the form \recalf{caseiiifunctions} then gives the equations (at all points $(\xb_4,\xb_5, \xib_1,\xib_3,0,0)$)
\begin{moneq}
\fracp{\psi}{\xb_4} + \frac{i\xib_1\kappa}{\yb_7} \psi = 0
\qquad\text{and}\qquad
\fracp{\psi}{\xb_5} = 0
\mapob.
\end{moneq}
The solutions} are determined by functions $\chi$ of two odd variables according to
\begin{moneq}
\psi(\xb_4,\xb_5, \xib_1,\xib_3,t,\tau)
=
\eexp^{i\tau\kappa + it/\hbar} \cdot\eexp^{-i\xb_4\xib_1\kappa/\yb_7}
\cdot \chi(\xib_1,\xib_3)
\mapob.
\end{moneq}
It follows that the induced action $\Phi$ \recalf{InducedActionOnPrequantumFunctions} of $G$ on \KorteOfLangeTekst{such a function $\psi$}{a function $\psi$ on $Y$ satisfying \recalf{secondGQcondition} and \recalf{firstGQcondition}} is given (using \recalf{liftedactioncaseiii}) by
\begin{align*}
\shifttag{8em}
(\Phi_\grel \psi)(\xb_4,\xb_5, 
\xib_1,\xib_3, t, \tau)
=
\psi\bigl(\phi_\grel\mo(\xb_4,\xb_5, 
\xib_1,\xib_3, t, \tau)\bigr)
\\&
\dontprintcomp{ 
=
\psi(\xb_4 - \yb_7a_1 ,
\xb_5- \yb_7 a_3, 
\xib_1 + \yb_7 \alpha_4 ,
\xib_3 + \yb_7 \alpha_5 , 
\\&
\qquad
t+ a_1\xnul_1 + a_2\xnul_2 + a_3 \xnul_3 ,
\\&
\qquad\qquad
\tau  
+ a_3 \xib_3 + \alpha_4\xb_4 +
\alpha_6 \xnulb_6 +\gamma \yb_7 - \yb_7 a_1 \alpha_4 + \yb_7 a_3 \alpha_5
)
\\
=\ 
&
\chi(\xib_1 + \yb_7 \alpha_4,\xib_3 + \yb_7 \alpha_5)
\cdot \eexp^{-i(\xb_4 - \yb_7a_1)(\xib_1 + \yb_7 \alpha_4)\kappa/\yb_7}
\\&
\qquad
\cdot \eexp^{i(t+ a_1\xnul_1 + a_2\xnul_2 + a_3 \xnul_3 )/\hbar}
\\&
\qquad\qquad
\cdot \eexp^{i(\tau  
+ a_3 \xib_3 + \alpha_4\xb_4 +
\alpha_6 \xnulb_6 + \gamma \yb_7 - \yb_7 a_1 \alpha_4 + \yb_7 a_3 \alpha_5)\kappa}
\\&
} 
=
\chi(\xib_1 + \yb_7 \alpha_4,\xib_3 + \yb_7 \alpha_5)
\cdot
\eexp^{i\tau\kappa + it/\hbar}
\cdot
\eexp^{-i\xb_4\xib_1\kappa/\yb_7}
\\&
\qquad
\cdot \eexp^{i( a_1\xnul_1 + a_2\xnul_2 + a_3 \xnul_3  )/\hbar}
\\&
\qquad\qquad
\cdot \eexp^{i(a_1\xib_1 + a_3 \xib_3  +
\alpha_6 \xnulb_6 + \gamma \yb_7 + \yb_7 a_3 \alpha_5)\kappa}
\mapob.
\end{align*}
\KorteOfLangeTekst{}{%
In terms of the functions $\chi$ the action thus becomes
\begin{align*}
(\Phi_\grel \chi)( 
&
\xib_1,\xib_3)
=
\chi(\xib_1 + \yb_7 \alpha_4,\xib_3 + \yb_7 \alpha_5)
\cdot \eexp^{i( a_1\xnul_1 + a_2\xnul_2 + a_3 \xnul_3 )/\hbar}
\\&
\qquad\qquad
\cdot \eexp^{i(a_1\xib_1 + a_3 \xib_3  +
\alpha_6 \xnulb_6 + \gamma \yb_7 + \yb_7 a_3 \alpha_5)\kappa}
\mapob.
\end{align*}
}%
If we introduce the change of variables $\zeta_4 = -\xib_1/\yb_7$ and $\zeta_5 = -\xib_3/\yb_7$, then, writing $\hat \chi(\zeta_4,\zeta_5) = \chi(\xib_1,\xib_3)$, the action of $\grel\in G$ on the functions $\hat \chi$ becomes
\begin{align*}
(\Phi_\grel \hat\chi)(\zeta_4,\zeta_5)
&=
\dontprintcomp{%
\chi(\xib_1 + \yb_7 \alpha_4,\xib_3 + \yb_7 \alpha_5)
\cdot \eexp^{i( a_1\xnul_1 + a_2\xnul_2 + a_3 \xnul_3 )/\hbar}
\\&
\qquad\qquad
\cdot \eexp^{i(a_1\xib_1 + a_3 \xib_3  +
\alpha_6 \xnulb_6 + \gamma \yb_7 + \yb_7 a_3 \alpha_5)\kappa}
\\&
=
\hat \chi(\zeta_4 - \alpha_4,\zeta_5 - \alpha_5)
\cdot \eexp^{i( a_1\xnul_1 + a_2\xnul_2 + a_3 \xnul_3 )/\hbar}
\\&
\qquad\qquad
\cdot \eexp^{i(-a_1\yb_7\zeta_4 - a_3\yb_7\zeta_5  +
\alpha_6 \xnulb_6 + \gamma \yb_7 + \yb_7 a_3 \alpha_5)\kappa}
\\&=}
\hat\chi(\zeta_4 - \alpha_4,\zeta_5 - \alpha_5)
\cdot \eexp^{-i(a_1\zeta_4 + a_3\zeta_5)\yb_7\kappa}
\\&
\qquad\qquad
\cdot \eexp^{i( a_1\xnul_1 + a_2\xnul_2 + a_3 \xnul_3 )/\hbar}
\cdot \eexp^{i\alpha_6 \xnulb_6\kappa}
\cdot \eexp^{i(\gamma + a_3 \alpha_5)\yb_7\kappa}
\mapob.
\end{align*}

\medskip\noindent$\bullet$
Next we choose the polarization subalgebra of dimension $2\vert 4$ generated by $e_2, e_4, e_5$, $e_6, k_0, k_7$.
The corresponding horizontally lifted vector fields determined by the conditions \recalf{DefinitionOfHorizontalLiftFundamentalVectorFields} are given by $k_0^h = k_7^h = e_2^h = e_6^h=0$ and
\begin{moneq}
e_4^h =  \yb_7\,\partial_{\xib_1}
\qquad\text{and}\qquad
e_5^h =  \yb_7\,\partial_{\xib_3} - \xb_5\,\partial_\tau
\mapob.
\end{moneq}
\KorteOfLangeTekst{Functions $\psi$ on the prequantum bundle $Y$ satisfying \recalf{secondGQcondition} and \recalf{firstGQcondition}}{Condition \recalf{secondGQcondition} applied to functions of the form \recalf{caseiiifunctions} then gives the equations (at all points $(\xb_4,\xb_5, \xib_1,\xib_3,0,0)$)
\begin{moneq}
\fracp{\psi}{\xib_1} = 0
\qquad\text{and}\qquad
\fracp{\psi}{\xib_3} - \frac{i\xb_5\kappa}{\yb_7} \cdot \psi = 0
\mapob,
\end{moneq}
whose solutions} are determined by functions $\chi$ of two even variables according to
\begin{moneq}
\psi(\xb_4,\xb_5, \xib_1,\xib_3,t,\tau)
=
\eexp^{i\tau\kappa+it/\hbar}\,\eexp^{i\xib_3\xb_5\kappa/\yb_7}
\cdot \chi(\xb_4, \xb_5)
\mapob.
\end{moneq}
It follows that the induced action $\Phi$ \recalf{InducedActionOnPrequantumFunctions} of $G$ on 
\KorteOfLangeTekst{such a function $\psi$}{a function $\psi$ on $Y$ satisfying \recalf{secondGQcondition} and \recalf{firstGQcondition}} is given (using \recalf{liftedactioncaseiii}) by
\KorteOfLangeTekst{
\begin{align*}
\shifttag{8em}
(\Phi_\grel \psi)(\xb_4,\xb_5, 
\xib_1,\xib_3, t, \tau)
=
\chi(\xb_4 - \yb_7a_1 ,
\xb_5- \yb_7 a_3)
\cdot
\eexp^{i\tau\kappa + it/\hbar}
\cdot
\eexp^{i\xib_3\xb_5\kappa/\yb_7}
\\&
\qquad
\cdot \eexp^{i(a_1\xnul_1 + a_2\xnul_2 + a_3 \xnul_3  )/\hbar}
\cdot \eexp^{i( \alpha_4\xb_4 + \alpha_5 \xb_5 + 
\alpha_6 \xnulb_6 + \gamma \yb_7 - \yb_7 a_1 \alpha_4)\kappa}
\mapob.
\end{align*}
}{%
\begin{align*}
\shifttag{8em}
(\Phi_\grel \psi)(\xb_4,\xb_5, 
\xib_1,\xib_3, t, \tau)
=
\psi\bigl(\phi_\grel\mo(\xb_4,\xb_5, 
\xib_1,\xib_3, t, \tau)\bigr)
\\&
=
\psi(\xb_4 - \yb_7a_1 ,
\xb_5- \yb_7 a_3, 
\xib_1 + \yb_7 \alpha_4 ,
\xib_3 + \yb_7 \alpha_5 , 
\\&
\qquad
t+ a_1\xnul_1 + a_2\xnul_2 + a_3 \xnul_3  ,
\\&
\kern2.5em
\tau  
+ a_3 \xib_3 + \alpha_4\xb_4 +
\alpha_6 \xnulb_6 +\gamma \yb_7 - \yb_7 a_1 \alpha_4 + \yb_7 a_3 \alpha_5)
\\&
=
\chi(\xb_4 - \yb_7a_1 ,
\xb_5- \yb_7 a_3)
\cdot \eexp^{i(\xib_3 + \yb_7 \alpha_5)(\xb_5- \yb_7 a_3)\kappa/\yb_7}
\\&
\qquad
\cdot \eexp^{i(t+ a_1\xnul_1 + a_2\xnul_2 + a_3 \xnul_3  )/\hbar}
\\&
\qquad\qquad
\cdot \eexp^{i(\tau  
+ a_3 \xib_3 + \alpha_4\xb_4 +
\alpha_6 \xnulb_6 + \gamma \yb_7 - \yb_7 a_1 \alpha_4 + \yb_7 a_3 \alpha_5)\kappa}
\\&
=
\chi(\xb_4 - \yb_7a_1 ,
\xb_5- \yb_7 a_3)
\cdot
\eexp^{i\tau\kappa + it/\hbar}
\cdot
\eexp^{i\xib_3\xb_5\kappa/\yb_7}
\\&
\qquad
\cdot \eexp^{i(a_1\xnul_1 + a_2\xnul_2 + a_3 \xnul_3  )/\hbar}
\\&
\qquad\qquad
\cdot \eexp^{i( \alpha_4\xb_4 + \alpha_5 \xb_5 + 
\alpha_6 \xnulb_6 + \gamma \yb_7 - \yb_7 a_1 \alpha_4)\kappa}
\mapob.
\end{align*}
}%
\KorteOfLangeTekst{}{%
In terms of the functions $\chi$ the action thus becomes
\begin{align*}
(\Phi_\grel h)(\xb_4,\xb_5 )
&=
\chi(\xb_4 - \yb_7a_1 ,
\xb_5- \yb_7 a_3)
\cdot \eexp^{i(a_1\xnul_1 + a_2\xnul_2 + a_3 \xnul_3 )/\hbar}
\\&
\qquad\qquad
\cdot \eexp^{i( \alpha_4\xb_4 + \alpha_5 \xb_5 + 
\alpha_6 \xnulb_6 + \gamma \yb_7 - \yb_7 a_1 \alpha_4)\kappa}
\mapob.
\end{align*}
}%
If we introduce the change of variables $\zb_1 = (\xb_4 - \xnulb_4)/\yb_7$ and $\zb_3 = (\xb_5 - \xnulb_5)/\yb_7$, then, writing $\hat \chi(\zb_1,\zb_3) = \chi(\xb_4,\xb_5)$, the action of $\grel\in G$ on the functions $\hat \chi$ becomes
\begin{align*}
(\Phi_\grel \hat\chi)(\zb_1,\zb_3 )
&=
\dontprintcomp{
\chi(\xb_4 - \yb_7a_1 ,
\xb_5- \yb_7 a_3)
\cdot \eexp^{i(a_1\xnul_1 + a_2\xnul_2 + a_3 \xnul_3 )/\hbar}
\\&
\qquad\qquad
\cdot \eexp^{i( \alpha_4\xb_4 + \alpha_5 \xb_5 + 
\alpha_6 \xnulb_6 + \gamma \yb_7 - \yb_7 a_1 \alpha_4)\kappa}
\\&
=
}
\hat\chi(\zb_1 - a_1 ,
\zb_3-  a_3)
\cdot \eexp^{i(a_1\xnul_1 + a_2\xnul_2 + a_3 \xnul_3  )/\hbar}
\\&
\qquad
\cdot \eexp^{i( \alpha_4 \zb_1 + \alpha_5 \zb_3)\yb_7\kappa}
\cdot \eexp^{i( \alpha_4\xnulb_4  + \alpha_5 \xnulb_5 + 
\alpha_6 \xnulb_6 )\kappa}
\cdot \eexp^{i(\gamma - a_1 \alpha_4)\yb_7\kappa}
\mapob.
\end{align*}
As for the case of the representation associated to an orbit of dimension $2\vert2$ with an even symplectic form, the appearance of the parameters $\xnulb_4$ and $\xnulb_5$ is artificial; they disappear when we use the coordinates $\zb_1 = \xb_4 /\yb_7$ and $\zb_3 = \xb_5 /\yb_7$.

\medskip\noindent$\bullet$
Our last choice is the polarization subalgebra of dimension $3\vert3$ generated by $e_2, e_3$, $e_4$, $e_6$, $k_0, k_7$.\footnote{As for the polarization subalgebras for the $2\vert2$-dimensional orbits with an even symplectic form, I am convinced (without having checked explicitly) that other choices for $a,b\in\RR$ will yield equivalent representations.} The corresponding horizontally lifted vector fields determined by the conditions \recalf{DefinitionOfHorizontalLiftFundamentalVectorFields} are given by $k_0^h = k_7^h = e_2^h = e_6^h=0$ and
\begin{moneq}
e_3^h =  -\yb_7\,\partial_{\xb_5}
\qquad\text{and}\qquad
e_4^h =  \yb_7\,\partial_{\xib_1} 
\mapob.
\end{moneq}
It follows immediately that functions $\psi$ on $Y$ satisfying \recalf{secondGQcondition} and \recalf{firstGQcondition} are determined by functions $\chi$ of one even and one odd variable according to
\begin{moneq}
\psi(\xb_4,\xb_5, \xib_1,\xib_3,t,\tau)
=
\eexp^{i\tau\kappa + it/\hbar}\,\eexp^{i\xb_4\xnul_1/\yb_7}\cdot \chi(\xb_4, \xib_3)
\mapob,
\end{moneq}
where we have separated a factor $\eexp^{i\xb_4\xnul_1/\yb_7}$ in order to get rid of the $\xnul_1$ dependence. 
The induced action $\Phi$ \recalf{InducedActionOnPrequantumFunctions} of $G$ on such a function $\psi$ thus is given (using \recalf{liftedactioncaseiii}) by
\begin{align*}
\shifttag{8em}
(\Phi_\grel \psi)(\xb_4,\xb_5, 
\xib_1,\xib_3, t, \tau)
\dontprintcomp{ 
=
\psi\bigl(\phi_\grel\mo(\xb_4,\xb_5, 
\xib_1,\xib_3, t, \tau)\bigr)
\\&
=
\psi(\xb_4 - \yb_7a_1 ,
\xb_5- \yb_7 a_3, 
\xib_1 + \yb_7 \alpha_4 ,
\xib_3 + \yb_7 \alpha_5 , 
\\&
\qquad
t+ a_1\xnul_1 + a_2\xnul_2 + a_3 \xnul_3  ,
\\&
\qquad\qquad
\tau  
+ a_3 \xib_3 + \alpha_4\xb_4 +
\alpha_6 \xnulb_6 + \gamma \yb_7 - \yb_7 a_1 \alpha_4 + \yb_7 a_3 \alpha_5
)
} 
\\&
=
\eexp^{i\tau\kappa + it/\hbar}\,\eexp^{i\xb_4\xnul_1/\yb_7}
\cdot \chi(\xb_4 - \yb_7a_1,\xib_3 + \yb_7 \alpha_5)
\cdot \eexp^{i( a_2\xnul_2 + a_3 \xnul_3 )/\hbar}
\\&
\qquad
\cdot \eexp^{i( a_3 \xib_3 + \alpha_4\xb_4 +
\alpha_6 \xnulb_6 + \gamma \yb_7 - \yb_7 a_1 \alpha_4 + \yb_7 a_3 \alpha_5)\kappa}
\mapob.
\end{align*}
\KorteOfLangeTekst{}{%
In terms of the functions $\chi$ the action $\Phi$ thus becomes
\begin{align*}
(\Phi_\grel \chi)(\xb_4,
\xib_3)
&=
\chi(\xb_4 - \yb_7a_1,\xib_3 + \yb_7 \alpha_5)
\cdot \eexp^{i( a_2\xnul_2 + a_3 \xnul_3  )/\hbar}
\\&
\qquad
\cdot \eexp^{i( a_3 \xib_3 + \alpha_4\xb_4 +
\alpha_6 \xnulb_6)\kappa}\,\eexp^{i (\gamma - a_1 \alpha_4 + a_3 \alpha_5)\yb_7\kappa}
\mapob.
\end{align*}
}%
If we introduce the change of variables $\zb_1 = (\xb_4 - \xnulb_4) /\yb_7$ and $\zeta_5 = -\xib_3/\yb_7$, then, writing $\hat \chi(\zb_1,\zeta_5) = \chi(\xb_4,\xib_3)$, the action of $\grel\in G$ on the functions $\hat \chi$ becomes
\begin{align}
(\Phi_\grel \hat\chi)(\zb_1,\zeta_5)
&=
\dontprintcomp{
\chi(\xb_4 - \yb_7a_1,\xib_3 + \yb_7 \alpha_5)
\cdot \eexp^{i( a_2\xnul_2 + a_3 \xnul_3  )/\hbar}
\notag
\\&
\qquad
\cdot \eexp^{i( a_3 \xib_3 + \alpha_4\xb_4 +
\alpha_6 \xnulb_6)\kappa}\,\eexp^{i (\gamma - a_1 \alpha_4 + a_3 \alpha_5)\yb_7\kappa}
\notag
\\&
=}
\hat\chi(\zb_1 - a_1,\zeta_5 - \alpha_5)
\cdot \eexp^{i(  \alpha_4\zb_1-a_3 \zeta_5) \yb_7\kappa}
\label{22odd33polaction}
\\&
\qquad
\cdot \eexp^{i( a_2\xnul_2 + a_3 \xnul_3  )/\hbar}
\cdot \eexp^{i( \alpha_4 \xnulb_4  +
\alpha_6 \xnulb_6 )\kappa}
\cdot \eexp^{i( \gamma - a_1 \alpha_4 + a_3 \alpha_5) \yb_7\kappa}
\mapob.
\notag
\end{align}
And as in previous cases, the appearance of the parameter $\xnulb_4$ is artificial and disappears when we use the coordinate $\zb_1 = \xb_4 /\yb_7$ instead of the shifted version $\zb_1 = (\xb_4 - \xnulb_4) /\yb_7$.

\medskip\noindent\stress{Quantization of orbits of dimension \texorpdfstring{$3\vert3$}{33}.}
As in the other cases, the prequantum bundle is the direct product $Y=\Orbit_{\mu_o} \times (\CA_0/d\ZZ) \times \CA_1$. Here a connection form $\Gamma$ is given by
\begin{align*}
\Gamma_0 
&
= \extder t + {(y_0)}\mo ( -x_1\, \extder \xh_2 - \tfrac12 \xih_5\, \extder \xih_5 + \tfrac12 \xi_6 \,\extder \xi_6 ) 
\\
\Gamma_1 
&
=  \extder \tau + {(y_0)}\mo ( -\xb_5 \,\extder \xih_5 - \xib_1\, \extder \xh_2 ) 
\mapob.
\end{align*}
The lifted fundamental vector fields $(v,z)^{Y}$ (as always determined by the conditions \recalf{LiftingEquationsOfGroupActionOnPrequantumBundle}) are given by $k_0^Y = y_0\,\partial_t$, $k_7^Y = \yb_7\,\partial_\tau$ and
\begin{align*}
e_1^Y & = -y_0\,\partial_{\xh_2} 
&
e_2^Y &= y_0\,\partial_{x_1} + \xh_2\,\partial_t
&
e_3^Y &= -\yb_7\,\partial_{\xb_5} + \xnul_3\,\partial_t + \xib_3\,\partial_\tau
\\
e_4^Y & = \yb_7\,\partial_{\xib_1} +\xb_4\,\partial_\tau
&
e_5^Y &= -y_0\,\partial_{\xih_5} - \tfrac12\xih_5\,\partial_t 
&
e_6^Y &= y_0\,\partial_{\xi_6} -\tfrac12\xi_6\,\partial_t + \xnulb_6\,\partial_\tau
\end{align*}
where the functions $\xb_4$ and $\xib_3$ must be seen as the functions $y_0\mo (s^o + \yb_7 \xh_2) = y_0\mo (y_0\xnulb_4- \yb_7\xnul_2 + \yb_7 \xh_2)$ and $-y_0\mo \yb_7 \xih_5$ respectively. These vector fields integrate to an action $\phi$ of $G$ on $Y$ by
\begin{moneq}[liftedactioncaseiv]
\phi_\grel:
\begin{pmatrix} x_1 \\ \xh_2 \\ \xb_5 \\ \xib_1 \\ \xih_5 \\ \xi_6 \\ t\\  \tau \end{pmatrix}
\mapsto
\begin{pmatrix}
x_1 - y_0a_2
\\ 
\xh_2 + y_0 a_1
\\ 
\xb_5+ \yb_7 a_3
\\
\xib_1 - \yb_7 \alpha_4
\\
\xih_5 + y_0 \alpha_5
\\
\xi_6 - y_0 \alpha_6
\\
t- a_2\xh_2 - a_3 \xnul_3 - c y_0 + \tfrac12 \alpha_5 \xih_5 + \tfrac12 \alpha_6 \xi_6  
\\
\tau 
- a_3 \xib_3 - \alpha_4\xb_4 -
\alpha_6 \xnulb_6 - \gamma \yb_7 
\end{pmatrix}
\end{moneq}

For all these orbits the stabilizer subgroup is given by 
\begin{moneq}
G_{\mu_o} = \{(a, c, \alpha, \gamma)\mid a_i=\alpha_j=0\}
\mapob.
\end{moneq}
The stabilizer subalgebra $\Liealg g_{\mu_o}$ of dimension $1\vert1$ thus is generated by $k_0,k_7$. 
Now suppose that $\Liealg h$ is a polarization subalgebra. 
If $\Liealg h$ has even dimension $4$, it contains all $4$ even basis vectors, but then 
\begin{moneq}
\contrs{[e_1, e_2]}{\mu_o} = -y_0\neq0
\mapob.
\end{moneq}
The conclusion is that the even dimension of $\Liealg h$ is at most $3$. Similarly, if the odd dimension is $4$, it contains in particular the vector $e_5$, but then 
\begin{moneq}
\contrs{[e_5,e_5]}{\mu_o} = y_0\neq0
\mapob.
\end{moneq}
Hence the odd dimension of $\Liealg h$ is at most $3$. A detailed analysis when the even dimension of $\Liealg h$ is $3$ shows that this is possible if and only if the even basis vectors are $ae_1+be_2, e_3, k_0$ (with $\vert a\vert^2+\vert b\vert^2=1$). Similarly, the only way to have three odd basis vectors in $\Liealg h$ is when they are $e_4, e_5+\epsilon e_6, k_7$ with $\epsilon=\pm1$. If $\Liealg h$ has dimension $3\vert 3$, it thus contains the vectors $ae_1+be_2, e_3, e_4, e_5+\epsilon e_6, k_0, k_7$. But then 
\begin{moneq}
\contrs{[e_3, e_5+\epsilon e_6]}{\mu_o} = -\yb_7\neq0
\mapob,
\end{moneq}
which shows that our polarization subalgebra can not have dimension $3\vert 3$ either. 
A detailed analysis of dimensions $3\vert 2$ and $2\vert 3$ for our polarization subalgebra $\Liealg h$ shows that the only possibilities are the following: for dimension $3\vert 2$ there is a single $\Liealg h$ possible, generated by $e_2, e_3, e_4, k_0, k_7$ and for dimension $2\vert 3$ there are exactly two possibilities, generated by $e_2, e_4, e_5+\epsilon e_6, k_0, k_7$, $\epsilon=\pm1$.
\KorteOfLangeTekst{}{As in the previous cases, functions $\psi$ on $Y$ satisfying condition \recalf{firstGQcondition} are of the form
\begin{moneq}[caseivfunctions]
\psi(x_1, \xh_2, \xb_5, \xib_1, \xih_5, \xi_6, t,\tau)
=
\eexp^{i\tau\kappa + it/\hbar} \cdot 
\psi(x_1, \xh_2, \xb_5, \xib_1, \xih_5, \xi_6,0,0)
\mapob. 
\end{moneq}
}

\medskip\noindent$\bullet$
We start with the invariant polarization subalgebra of dimension $3\vert2$ generated by $e_2, e_3, e_4, k_0, k_7$. 
The horizontally lifted vector fields $b^h$, $b=e_2, e_3, e_4, k_0, k_7$ determined by the conditions \recalf{DefinitionOfHorizontalLiftFundamentalVectorFields} are given by $k_0^h = k_7^h =0$ and 
\begin{moneq}
e_2^h = y_0\,\partial_{x_1}
\ ,\qquad
e_3^h = -\yb_7\,\partial_{\xb_5}
\qquad\text{and}\qquad
e_4^h = \yb_7\,\partial_{\xib_1}
\mapob.
\end{moneq}
It follows immediately that functions $\psi$ on $Y$ satisfying \recalf{secondGQcondition} and \recalf{firstGQcondition} are determined by functions $\chi$ of one even and two odd variables according to
\begin{moneq}
\psi(x_1, \xh_2, \xb_5, \xib_1, \xih_5, \xi_6,t,\tau)
=
\eexp^{i\tau\kappa+it/\hbar}\cdot \chi(\xh_2, \xih_5, \xi_6)
\mapob.
\end{moneq}
The induced action $\Phi$ \recalf{InducedActionOnPrequantumFunctions} of $G$ on such a function $\psi$ thus is given (using \recalf{liftedactioncaseiv}) by
\begin{align*}
\shifttag{5em}
(\Phi_\grel \psi)(x_1, \xh_2, \xb_5, 
\xib_1, \xih_5, \xi_6, t, \tau)
=
\psi\bigl(\phi_\grel\mo (x_1, \xh_2, \xb_5, 
\xib_1, \xih_5, \xi_6, t, \tau)\bigr)
\\
&=
\dontprintcomp{
\psi\bigl(x_1 + y_0a_2 ,
\xh_2 - y_0 a_1 ,
\xb_5 - \yb_7 a_3 ,
\\&
\kern2em
\xib_1 + \yb_7 \alpha_4 ,
\xih_5 - y_0 \alpha_5 ,
\xi_6 + y_0 \alpha_6 ,
\\&
\kern4em
t + a_2\xh_2 + a_3 \xnul_3 + (c-a_1a_2) y_0 - \tfrac12 \alpha_5 \xih_5 - \tfrac12 \alpha_6 \xi_6 ,
\\&
\kern6em
\tau 
+ a_3 \xib_3 + \alpha_4\xb_4 +
\alpha_6 \xnulb_6 + (\gamma-a_1\alpha_4+a_3\alpha_5) \yb_7 
\bigr)
\\&
=
\chi(\xh_2 - y_0 a_1 ,
\xih_5 - y_0 \alpha_5 ,
\xi_6 + y_0 \alpha_6 )
\\&
\kern2em
\cdot \eexp^{i(t + a_2\xh_2 + a_3 \xnul_3 + c y_0 - \frac12 \alpha_5 \xih_5 - \frac12 \alpha_6 \xi_6 - y_0 a_1 a_2)/\hbar}
\\&
\kern4em
\cdot \eexp^{i(\tau 
+ a_3 \xib_3 + \alpha_4\xb_4 +
\alpha_6 \xnulb_6 + \gamma \yb_7 - \yb_7 a_1 \alpha_4 + \yb_7 a_3 \alpha_5)\kappa}
\\&
=
}
\chi(\xh_2 - y_0 a_1 ,
\xih_5 - y_0 \alpha_5 ,
\xi_6 + y_0 \alpha_6 )
\cdot \eexp^{i\tau\kappa + it/\hbar}
\\&
\kern2em
\cdot \eexp^{i(a_2\xh_2 + a_3 \xnul_3 + c y_0 - \frac12 \alpha_5 \xih_5 - \frac12 \alpha_6 \xi_6 - y_0 a_1 a_2)/\hbar}
\\&
\kern4em
\cdot \eexp^{i(a_3 \xib_3 + \alpha_4\xb_4 +
\alpha_6 \xnulb_6 + \gamma \yb_7 -  \yb_7 a_1 \alpha_4 + \yb_7 a_3 \alpha_5)\kappa}
\mapob,
\end{align*}
where as before the dependent coordinates $\xib_3$ and $\xb_4$ represent the functions $-y_0\mo \yb_7 \xih_5$ and  $y_0\mo (s^o + \yb_7 \xh_2) = y_0\mo (y_0\xnulb_4- \yb_7\xnul_2 + \yb_7 \xh_2)$ respectively. 
\KorteOfLangeTekst{}{In terms of the functions $\chi$ the action then becomes
\begin{align*}
(\Phi_\grel \chi)( \xh_2, \xih_5, \xi_6)
&=
\chi(\xh_2 - y_0 a_1 ,
\xih_5 - y_0 \alpha_5 ,
\xi_6 + y_0 \alpha_6 )
\\&
\kern2em
\cdot \eexp^{i(a_2\xh_2 + a_3 \xnul_3 + c y_0 - \frac12 \alpha_5 \xih_5 - \frac12 \alpha_6 \xi_6 - y_0 a_1 a_2)/\hbar}
\\&
\kern4em
\cdot \eexp^{i(a_3 \xib_3 + \alpha_4\xb_4 +
\alpha_6 \xnulb_6 + \gamma \yb_7 -  \yb_7 a_1 \alpha_4 + \yb_7 a_3 \alpha_5)\kappa}
\mapob.
\end{align*}
}%
If we introduce the change of variables $z_1 = (\xh_2 - \xnul_2)/y_0$, $\zeta_5 = \xih_5/y_0$ and $\zeta_6 = -\xi_6/y_0$ and if we write $\hat \chi(z_1, \zeta_5, \zeta_6) = \chi(\xh_2, \xih_5, \xi_6)$, then the action of $\grel\in G$ on the functions $\hat \chi$  becomes
\begin{align}
(\Phi_\grel \hat \chi)( z_1, \zeta_5, \zeta_6)
&=
\dontprintcomp{ 
\chi(\xnul_2+y_0z_1 - y_0 a_1 ,
y_0\zeta_5 - y_0 \alpha_5 ,
-y_0\zeta_6 + y_0 \alpha_6 )
\notag
\\&
\kern2em
\cdot \eexp^{i(a_2\xh_2 + a_3 \xnul_3 + c y_0 - \frac12 \alpha_5 \xih_5 - \frac12 \alpha_6 \xi_6 - y_0 a_1 a_2)/\hbar}
\notag
\\&
\kern4em
\cdot \eexp^{i(a_3 \xib_3 + \alpha_4\xb_4 +
\alpha_6 \xnulb_6 + \gamma \yb_7 -  \yb_7 a_1 \alpha_4 + \yb_7 a_3 \alpha_5)\kappa}
\notag
\\&
=
\hat \chi(z_1 - a_1 , \zeta_5 - \alpha_5 , \zeta_6 - \alpha_6 )
\notag
\\&
\kern2em
\cdot \eexp^{i(a_2(\xnul_2+y_0z_1) + a_3 \xnul_3 + c y_0 - \frac12 \alpha_5 y_0\zeta_5 + \frac12 \alpha_6 y_0\zeta_6 - y_0 a_1 a_2)/\hbar}
\notag
\\&
\kern4em
\cdot \eexp^{i(a_3 (-y_0\mo \yb_7 \xih_5) + \alpha_4(y_0\mo (y_0\xnulb_4- \yb_7\xnul_2 + \yb_7 \xh_2)))\kappa}
\notag
\\&
\kern6em
\cdot \eexp^{i(\alpha_6 \xnulb_6 + \gamma \yb_7 -  \yb_7 a_1 \alpha_4 + \yb_7 a_3 \alpha_5)\kappa}
\notag
\\&
=
\hat \chi(z_1 - a_1 , \zeta_5 - \alpha_5 , \zeta_6 - \alpha_6 )
\notag
\\&
\kern2em
\cdot \eexp^{i(a_2(\xnul_2+y_0z_1) + a_3 \xnul_3 + c y_0 - \frac12 \alpha_5 y_0\zeta_5 + \frac12 \alpha_6 y_0\zeta_6 - y_0 a_1 a_2)/\hbar}
\notag
\\&
\kern4em
\cdot \eexp^{i(a_3 (-y_0\mo \yb_7 y_0\zeta_5) + \alpha_4(y_0\mo (y_0\xnulb_4- \yb_7\xnul_2 + \yb_7 (\xnul_2+y_0z_1))) 
)\kappa}
\notag
\\&
\kern6em
\cdot \eexp^{i(\alpha_6 \xnulb_6 + \gamma \yb_7 -  \yb_7 a_1 \alpha_4 + \yb_7 a_3 \alpha_5)\kappa}
\notag
\\&
=
\hat \chi(z_1 - a_1 , \zeta_5 - \alpha_5 , \zeta_6 - \alpha_6 )
\notag
\\&
\kern2em
\cdot \eexp^{i(a_2(\xnul_2+y_0z_1) + a_3 \xnul_3 + c y_0 - \frac12 \alpha_5 y_0\zeta_5 + \frac12 \alpha_6 y_0\zeta_6 - y_0 a_1 a_2)/\hbar}
\notag
\\&
\kern4em
\cdot \eexp^{i(-a_3 \yb_7 \zeta_5 + \alpha_4( \xnulb_4 + \yb_7 z_1) +
\alpha_6 \xnulb_6 + \gamma \yb_7 -  \yb_7 a_1 \alpha_4 + \yb_7 a_3 \alpha_5)\kappa}
\notag
\\&
=
}
\hat \chi(z_1 - a_1 , \zeta_5 - \alpha_5 , \zeta_6 - \alpha_6 )
\notag
\\&
\kern2em
\cdot \eexp^{i(a_2\xnul_2+a_3\xnul_3)/\hbar + i(\alpha_4\xnulb_4+\alpha_6\xnulb_6)\kappa} \cdot\eexp^{ i(a_2z_1  + c - \frac12 \alpha_5 \zeta_5 + \frac12 \alpha_6 \zeta_6 -  a_1 a_2)y_0/\hbar}
\notag
\\&
\kern4em
\cdot \eexp^{i(-a_3  \zeta_5 + \alpha_4 z_1 + \gamma  -   a_1 \alpha_4 +  a_3 \alpha_5)\yb_7\kappa}
\mapob.
\label{33mixed32polaction}
\end{align}
This action depends upon the base point parameters $\xnul_2$, $\xnul_3$, $\xnulb_4$, $\xnulb_6$, $y_0$, $\yb_7$. As in previous cases, the appearance of the parameters $\xnul_2$ and $\xnulb_4$ is artificial. 
Using the coordinate $z_1 = \xh_2/y_0$ instead of $z_1 = (\xh_2 - \xnul_2)/y_0$, the term $a_2\xnul_2$ in the exponent would have disappeared and the term $\alpha_4\xnulb_4$ would have become $\alpha_4 s^o/y_0$. 
And then it is obvious that the action only depends upon the orbit parameters $\xnul_3$, $s^o$, $\xnulb_6$, $y_0$, $\yb_7$. 
Another shift would have made the term $\alpha_4\xnulb_4$ disappear and changed the term $a_2\xnul_2$ to $-s^o a_2/\yb_7$.

\medskip\noindent$\bullet$
We finish with the polarization subalgebra of dimension $2\vert3$ generated by $e_2, e_4, e_5+\epsilon e_6, k_0, k_7$.\footnote{And as before, I did not check explicitly but I am convinced the two choices for $\epsilon$ give equivalent representations, as in the case of a $2\vert2$-dimensional orbit.} 
The horizontally lifted vector fields $b^h$ determined by the conditions \recalf{DefinitionOfHorizontalLiftFundamentalVectorFields} are given by $k_0^h = k_7^h =0$ and 
\begin{gather*}
e_2^h = y_0\,\partial_{x_1}
\ ,\qquad
e_4^h = \yb_7\,\partial_{\xib_1}
\\
e_5^h+\epsilon\, e_6^h 
=
-y_0\,\partial_{\xih_5} +\epsilon y_0\,\partial_{\xi_6} + \tfrac12(\xih_5+\epsilon\xi_6)\,\partial_t - \xb_5\,\partial_\tau
\mapob.
\end{gather*}
\KorteOfLangeTekst{Functions $\psi$ on the prequantum bundle $Y$ satisfying \recalf{secondGQcondition} and \recalf{firstGQcondition}}{Condition \recalf{secondGQcondition} applied to functions $\psi$ of the form \recalf{caseivfunctions} then gives the equations (at all points $(x_1, \xh_2, \xb_5, \xib_1, \xih_5, \xi_6,0,0)$)
\begin{moneq}
\fracp{\psi}{x_1} = 0
\quad,\quad
\fracp{\psi}{\xib_1} = 0
\qquad\text{and}\qquad
\fracp{\psi}{\xih_5} - \epsilon \fracp{\psi}{\xi_6} - i \cdot \frac{\xih_5+\epsilon \xi_6}{2 y_0\hbar} \cdot  \psi + i\cdot\frac{ \xb_5 \kappa}{y_0} \cdot \psi = 0
\mapob,
\end{moneq}
whose solutions} are determined by functions $\chi$ of two even and one odd variable according to
\begin{moneq}
\psi(x_1, \xh_2, \xb_5, \xib_1, \xih_5, \xi_6,t,\tau)
=
\eexp^{i\tau\kappa+it/\hbar} \cdot
\eexp^{\frac{i\epsilon}{2y_0\hbar}\xih_5\xi_6}
\cdot\eexp^{-\frac{i}{y_0}\xb_5\xih_5\kappa}\cdot \chi(\xh_2, \xb_5, \xih_5 + \epsilon \xi_6)
\mapob.
\end{moneq}
The induced action $\Phi$ \recalf{InducedActionOnPrequantumFunctions} of $G$ on \KorteOfLangeTekst{such a function $\psi$}{a function $\psi$ on $Y$ satisfying \recalf{secondGQcondition} and \recalf{firstGQcondition}} is given (using again \recalf{liftedactioncaseiv}) by
\begin{align*}
\shifttag{8em}
(\Phi_\grel \psi)(x_1, \xh_2, \xb_5, 
\xib_1, \xih_5, \xi_6, t, \tau)
\dontprintcomp{ 
\\&
=
\psi(x_1 + y_0a_2 ,
\xh_2 - y_0 a_1 ,
\xb_5 - \yb_7 a_3 ,
\\&
\kern2em
\xib_1 + \yb_7 \alpha_4 ,
\xih_5 - y_0 \alpha_5 ,
\xi_6 + y_0 \alpha_6 ,
\\&
\kern3em
t + a_2\xh_2 + a_3 \xnul_3 + c y_0 - \tfrac12 \alpha_5 \xih_5 - \tfrac12 \alpha_6 \xi_6 - y_0 a_1 a_2 ,
\\&
\kern4em
\tau 
+ a_3 \xib_3 + \alpha_4\xb_4 +
\alpha_6 \xnulb_6 + \gamma \yb_7 - \yb_7 a_1 \alpha_4 + \yb_7 a_3 \alpha_5
)
\\
&=
\chi(\xh_2 - y_0 a_1 ,
\xb_5 - \yb_7 a_3 ,
(\xih_5 +\epsilon\xi_6)  - y_0 (\alpha_5 - \epsilon \alpha_6 )
)
\\&
\kern2em
\cdot
\eexp^{\frac{i\epsilon}{2y_0\hbar}(\xih_5 - y_0 \alpha_5)(\xi_6 + y_0 \alpha_6)}
\cdot
\eexp^{-\frac i{y_0} (\xb_5 - \yb_7 a_3)(\xih_5-y_0\alpha_5)\kappa}
\\&
\kern3em
\cdot \eexp^{i(t + a_2\xh_2 + a_3 \xnul_3 + c y_0 - \frac12 \alpha_5 \xih_5 - \frac12 \alpha_6 \xi_6 - y_0 a_1 a_2)/\hbar}
\\&
\kern4em
\cdot \eexp^{i(\tau 
+ a_3 \xib_3 + \alpha_4\xb_4 +
\alpha_6 \xnulb_6 + \gamma \yb_7 - \yb_7 a_1 \alpha_4 + \yb_7 a_3 \alpha_5)\kappa}
} 
\\&
=
\chi(\xh_2 - y_0 a_1 ,
\xb_5 - \yb_7 a_3 ,
(\xih_5 +\epsilon\xi_6)  - y_0 (\alpha_5 -\epsilon \alpha_6 )
)
\\& 
\kern2em
\cdot \eexp^{i\tau\kappa + it/\hbar}
\cdot
\eexp^{\frac{i\epsilon}{2y_0\hbar}\xih_5\xi_6}
\cdot
\eexp^{-\frac i{y_0} \xb_5 \xih_5 \kappa}
\\&
\kern3em
\cdot \eexp^{i( a_2\xh_2 + a_3 \xnul_3 + c y_0 - \frac12 (\alpha_5 + \epsilon\alpha_6) ( \xih_5 + \epsilon \xi_6) - a_1 a_2 - \frac\epsilon2 y_0 \alpha_5\alpha_6)/\hbar}
\\&
\kern4em
\cdot \eexp^{i(  \alpha_4\xb_4 +  \alpha_5\xb_5 +
\alpha_6 \xnulb_6 + \gamma \yb_7 - \yb_7 a_1 \alpha_4)\kappa}
\mapob,
\end{align*}
where as before the dependent coordinates $\xib_3$ and $\xb_4$ represent the functions $-y_0\mo \yb_7 \xih_5$ and  $y_0\mo (s^o + \yb_7 \xh_2) = y_0\mo (y_0\xnulb_4- \yb_7\xnul_2 + \yb_7 \xh_2)$ respectively.
\dontprintcomp{ 
In terms of the functions $\chi$ the action thus becomes
\begin{align*}
(\Phi_\grel \chi)( \xh_2, \xb_5, \xi)
&=
\chi(\xh_2 - y_0 a_1 ,
\xb_5 - \yb_7 a_3 ,
\xi - y_0 (\alpha_5 - \epsilon \alpha_6 )
)
\\& 
\cdot \eexp^{i( a_2\xh_2 + a_3 \xnul_3 + c y_0 - \frac12 (\alpha_5 + \epsilon\alpha_6)  \xi - y_0 a_1 a_2 - \frac\epsilon2 y_0 \alpha_5\alpha_6)/\hbar}
\\&
\cdot \eexp^{i(  \alpha_4\xb_4 +  \alpha_5\xb_5 +
\alpha_6 \xnulb_6 + \gamma \yb_7 - \yb_7 a_1 \alpha_4)\kappa}
\mapob.
\end{align*}
}%
If we introduce the change of variables $z_1 = (\xh_2 - \xnul_2)/y_0$, $\zb_3 = (\xb_5 - \xnulb_5)/\yb_7$ and $\zeta = \xi/y_0$, then the action of $\grel\in G$ on the functions $\chi$ in terms of the new variables becomes
\begin{align*}
(\Phi_\grel \chi)( z_1, \zb_3, \zeta)
&=
\chi(z_1 - a_1 ,
\zb_3 -  a_3 ,
\zeta - (\alpha_5 -\epsilon \alpha_6 )
)
\cdot \eexp^{i\frac{y_0}{\hbar}( a_2 z_1 - \frac12 (\alpha_5 + \epsilon\alpha_6) \zeta)} 
\\& 
\kern2em
\cdot \eexp^{i(  \alpha_4 z_1 +  \alpha_5 \zb_3) \yb_7\kappa}
\cdot \eexp^{i( a_2\xnul_2 + a_3 \xnul_3 )/\hbar}
\cdot \eexp^{i(  \alpha_4\xnulb_4 +  \alpha_5\xnulb_5 +
\alpha_6 \xnulb_6 )\kappa}
\\&
\kern4em
\cdot \eexp^{i\frac{y_0}{\hbar}( c - a_1 a_2 - \frac\epsilon2 \alpha_5\alpha_6) }
\cdot \eexp^{i( \gamma - a_1 \alpha_4) \yb_7\kappa}
\mapob.
\end{align*}
And just as in previous cases, the dependence on the parameters $\xnul_2$, $\xnulb_4$ and $\xnulb_5$ is artificial. If we use the coordinates $z_1 = \xh_2 /y_0$ and $\zb_3 = \xb_5 /\yb_7$ instead of the shifted ones given above, the terms $a_2 \xnul_2$ and $\alpha_5 \xnulb_5$ in the exponents disappear and the term $\alpha_4 \xnulb_4$ becomes $\alpha_4 s^o/y_0$. And in that form this representation only depends upon the parameters $\xnul_3$, $\xnulb_6$, $y_0$, $\yb_7$, $s^o$ describing the orbit. Shifting $z_1$ differently would have made the term $\alpha_4 \xnulb_4$ disappear and changed the term $a_2 \xnul_2$ to $-a_2 s^o/\yb_7$.

\subsection{Recapitulation and comparison}

\vrule width0pt
\vskip0.5\baselineskip

\noindent
If we gather the results for all orbits, we get the following list: 
\begin{enumerate}
\item
For each $0\vert0$-dimensional orbit we have a ($1$-dimensional) representation.

\item
For each $2\vert2$-dimensional orbit with an even symplectic form we have a representation by functions on a $1\vert1$-dimensional manifold. 

\item
For each $2\vert2$-dimensional orbit with an odd symplectic form we have three different representations by functions on a manifold of dimension $2\vert0$, $0\vert2$ and $1\vert1$ respectively.

\item
For each $3\vert3$-dimensional orbit (with a mixed symplectic form) we have two different representations  by functions on a manifold of dimension $2\vert1$ and $1\vert2$ respectively.

\end{enumerate}
We now compare this with the four families of subrepresentations of the left-regular representation found in \recals{SUBSECTDecompLeftRegularRepHeisenbergLikeSuperGroup}.
\begin{enumerate}
\item
For the representation $\rhoh^{(0, \ell_2, \ell_3, \lambda_4, 0),(\ell_1, \lambda_5, \lambda_6)}$, comparison of the $G$ action given in \recalf{ell0zeroandlambda0zeroaction} with the action given in the representation associated to a $0\vert0$-dimensional orbit \recalf{00dimorbitaction} tells us that they are the same, provided we make the following identifications
\begin{moneq}
\begin{matrix}
\ell_1 &\kern-0.5em=\kern-0.5em& -\xnul_1/\hbar
&,\quad&
\ell_2 &\kern-0.5em=\kern-0.5em& -\xnul_2/\hbar 
&,\quad&
\ell_3 &\kern-0.5em=\kern-0.5em& {-\xnul_3/\hbar} 
\\ 
\lambda_4 &\kern-0.5em=\kern-0.5em& -\xnulb_4\kappa
&,\quad&
\lambda_5 &\kern-0.5em=\kern-0.5em& - \xnulb_5\kappa
&,\quad&
\lambda_6 &\kern-0.5em=\kern-0.5em& - \xnulb_6\kappa 
\mapob.
\end{matrix}
\end{moneq}

\item
For the representation $\rhoh^{(0, \ell_2, \ell_3, \lambda_4, \lambda_7),( \lambda_6)}$ with $\lambda_7\neq0$, comparing the action \recalf{ell0zerobutnotlambda0action} with the action \recalf{22odd33polaction} given in the representation associated to a $2\vert2$-dimensional orbit with an odd symplectic form and a polarization subalgebra of dimension $3\vert3$ tells us that they are the same, provided we make the identifications
\begin{moneq}
\begin{matrix}
0 &\kern-0.5em=\kern-0.5em& -\xnul_1/\hbar
&,\quad&
\ell_2 &\kern-0.5em=\kern-0.5em& -\xnul_2/\hbar 
&,\quad&
\ell_3 &\kern-0.5em=\kern-0.5em& -\xnul_3/\hbar 
\\ 
\lambda_4 &\kern-0.5em=\kern-0.5em& - \xnulb_4\kappa
&,\quad&
\lambda_6 &\kern-0.5em=\kern-0.5em& - \xnulb_6\kappa 
&,\quad&
\lambda_7 &\kern-0.5em=\kern-0.5em& -\yb_7\kappa
\mapob.
\end{matrix}
\end{moneq}

\item
For the representation $\rhoh^{(\ell_0, \ell_2, \ell_3, \lambda_4, 0),\pm}$ with $\ell_0\neq0$, comparing the action 
\recalf{RhohDecompositionPlusMoinsAtLambda7Zero} with the action \recalf{22evenaction} given in the representations associated to a $2\vert2$-dimen\-sio\-nal orbit with an even symplectic form tells us that they are the same, provided we make the identifications
\begin{moneq}
\begin{matrix}
\ell_0 &\kern-0.5em=\kern-0.5em& -y_0/\hbar
&,\quad&
\ell_2 &\kern-0.5em=\kern-0.5em& -\xnul_2/\hbar 
&,\quad&
\ell_3 &\kern-0.5em=\kern-0.5em& -\xnul_3/\hbar 
\\ 
\lambda_4 &\kern-0.5em=\kern-0.5em& - \xnulb_4\kappa
&,\quad&
0 &\kern-0.5em=\kern-0.5em& - \xnulb_5\kappa
&,\quad&
0 &\kern-0.5em=\kern-0.5em& - \xnulb_6\kappa 
\mapob.
\end{matrix}
\end{moneq}

\item
For the representation $\rhoh^{(\ell_0, \ell_2, \ell_3, \lambda_4, \lambda_7)}$ with $\ell_0\neq0\neq\lambda_0$, comparing the action 
\recalf{FirstDecompositionOfLeftRegularRepOfHLikeSuperGroup} with the action \recalf{33mixed32polaction} given in the representation associated to a $3\vert3$-dimensional orbit with a polarization subalgebra of dimension $3\vert2$ tells us that they are the same, provided we make the identifications
\begin{moneq}
\begin{matrix}
\ell_0 &\kern-0.5em=\kern-0.5em& -y_0/\hbar
&,\quad&
\ell_2 &\kern-0.5em=\kern-0.5em& -\xnul_2/\hbar 
&,\quad&
\ell_3 &\kern-0.5em=\kern-0.5em& -\xnul_3/\hbar 
\\ 
\lambda_7 &\kern-0.5em=\kern-0.5em& - \yb_7\kappa
&,\quad&
\lambda_4 &\kern-0.5em=\kern-0.5em& - \xnulb_4\kappa
&,\quad&
0 &\kern-0.5em=\kern-0.5em& - \xnulb_6\kappa 
\mapob.
\end{matrix}
\end{moneq}

\end{enumerate}
When we look at the required identifications, we see that they all have the same form. 
The even parameters $\ell_0$, $\ell_1$, $\ell_2$, $\ell_3$ have to be the same as the coordinates $y_0$, $\xnul_1$, $\xnul_2$, $\xnul_3$ of the basepoint $\mu_o$ of the orbit up to a factor $-\hbar$. 
And the odd parameters $\lambda_4$, $\lambda_5$, $\lambda_6$, $\lambda_7$ have to be the same as the coordinates $\xnulb_4$, $\xnulb_5$, $\xnulb_6$, $\yb_7$ of $\mu_o$ up to a factor $-\kappa$. 
But here the situation is not as simple as for the even parameters. For the even parameters we performed a honest Fourier transform of (smooth) functions of even coordinates. 
As such the even parameters $\ell_0$, $\ell_1$, $\ell_2$, $\ell_3$ thus are real. 
And identifying them up to a factor $-\hbar$ with the real coordinates $y_0$, $\xnul_1$, $\xnul_2$, $\xnul_3$ poses no problem. 
On the other hand, the parameters $\lambda_4$, $\lambda_5$, $\lambda_6$, $\lambda_7$ are odd, so identifying them with the real coordinates $\xnulb_4$, $\xnulb_5$, $\xnulb_6$, $\yb_7$ is not obvious. 
But as said, it should be up to the odd Fourier parameter $\kappa$. 
But then all four odd parameters $\lambda_4$, $\lambda_5$, $\lambda_6$, $\lambda_7$ should be real multiples of the same odd factor $\kappa$. Which would suggest that we cannot choose these four parameters freely. 
This is the same problem we already encountered when linking the representations obtained in the decomposition of the super Heisenberg group $\SuperH$ to coadjoint orbits.

\SwitchHats
\renewcommand\grel{{(\ah,\ch,\alphah,\gammah)}}

\begin{remarque}{Remark}
Although I did not pursue this line of investigation, it is worthwhile to note that a similar (with respect to the super Heisenberg group) phenomenon appears in this super group: we can decompose the left-regular representation in different ways with different components, components that are related to the same orbits, but with a different polarization. 
For instance, in the left-regular representation, functions of the form
\begin{moneq}
\psi(\adh,\cdh, \alphadh,\gammadh)
=
\chi(\adh_1, \adh_3, \alphadh_5 + \epsilon \alphadh_6) \cdot
\eexp^{-\frac i2\ell_0\epsilon\alphadh_5\alphadh_6}\cdot
\eexp^{i\ell_2\adh_2+i\alphadh_4\lambda_4+i\ell_0\cdh+i(\gammadh+\adh_3\alphadh_5)\lambda_7}
\mapob,
\end{moneq}
with $\epsilon=\pm1$ and $\chi$ a function of two even and one odd coordinate also form an invariant subspace.
The associated family of representations is (as in the original decomposition) linked to orbits of dimension $3\vert3$, but now with a polarization subalgebra of dimension $2\vert 3$. 
\end{remarque}

\section{Concluding remarks}

In the process of classifying \myquote{all} representations of a Lie group one tries to decompose a representation into smaller parts. And then the idea is to start with the smallest possible parts and to build any representation from them. 
This implies that we need a notion of \myquote{building} and a notion of \myquote{smallest.} 
For representations one uses the notion of an irreducible representation for \myquote{smallest} and one uses the notion of direct sum for \myquote{building.} 
But already for the regular representation of the additive real line $\RR$ just taking direct sums is not enough, one has to add direct integrals to the notion of \myquote{building.} 
Now of course these ideas are not perfect, as there exist non-decomposable representations that are not irreducible, which is one of the reasons to restrict attention to unitary representations. 

When looking at representations of super Lie groups, we have to ask ourselves the same questions: what are the smallest building blocks and what are the rules of building? 
In order to get an answer to these questions I performed my \myquote{experiments} to compare the super situation with the non-super situation. 
And we have seen that the standard analysis of the Heisenberg group carries over completely to the super Heisenberg group. 
However, in order to make sense of that analysis, we had to introduce direct Berezin integrals to our building process and we had to define a notion of families of representations depending upon odd parameters (the $p\vert q$-parameter families). 
As for the classification part by coadjoint orbits, even there the non-super analysis carries over to the super setting, provided I introduce the odd Berezin-Fourier parameter $\kappa$. 
Without it we can not link all representations obtained in the decomposition of $\SuperH$ to coadjoint orbits (it amounts to taking $\kappa=0$, or to take the super prequantum bundle as a circle bundle over the orbit, just as in the non-super case). 
When we do introduce it, we get a nice fit with coadjoint orbits, even the analogy with the position and momentum representations carries over. 
But this Fourier parameter $\kappa$ also poses problems as how to identify $p\vert q$-parameter families of representations with (families of) coadjoint orbits. 

The perfect link between the representations obtained in the decomposition of the super Heisenberg group $\SuperH$ and the Heisenberg like super group $G$ with the representations obtained by geometric quantization of their coadjoint orbits suggests that there must be a more profound reason for this correspondence, more profound than just happenstance. 
But then we have to understand how the single Berezin-Fourier parameter $\kappa$ can be transformed into a multitude of odd parameters in a family of representations. 
And we should not forget that we also need a way to obtain the inner product and super scalar product on these representations. 
I can only hope that someday somebody will find a good explanation for this nice link.

\section{Appendix: an overview of \texorpdfstring{$\CA$}{A}-manifold theory}
\label{appendixonAmanifoldssection}

This appendix\footnote{It is a truncated version of the appendix I wrote in \cite{Tuynman..2015}} is a very short overview of $\CA$-manifold theory, intended for readers with some familiarity with supermanifold theory, but not with this approach. It will not always be complete, and sometimes it might be slightly besides the truth, but it intends to give the gist of the theory, not the most precise formulation (which can be found in \cite{Tu04}).

\firstofmysubsection
\mysubsection{appendixonAmanifoldssection}
{The basic graded commutative ring}

The starting point of $\CA$-manifold theory is a graded-commutative ring $\CA$. We will fix it as being the exterior algebra of an infinite dimensional (real) vector space $V$: 
$$
\CA=\bigwedge V = \bigoplus_{k=0}^\infty \bigwedge{}{\kern-3pt\raise1.5ex\hbox{$\scriptstyle k$}\kern2pt} V =  \Bigl(\, \bigoplus_{k=0}^\infty \bigwedge {\kern-3pt\raise1.5ex\hbox{$\scriptstyle 2k$}\kern2pt} V \,\Bigr) \oplus \Bigl(\, \bigoplus_{k=0}^\infty \bigwedge{\kern-3pt\raise1.5ex\hbox{$\scriptstyle 2k+1$}\kern2pt}V \,\Bigr)
=
\CA_0\oplus \CA_1
\mapob.
$$
We denote by $\body:\CA\to \CA$ (sic) the canonical projection onto the direct summand $\RR\equiv\bigwedge^0V\subset \bigoplus_{k=0}^\infty\bigwedge^kV$; we will call $\body$ the \stresd{body map} and the image $\body(a)\in \RR$ the \stresd{body of $a\in \CA$}.
The set of nilpotent elements $\nilpotent=\bigoplus_{k=1}^\infty \bigwedge^kV$ (an ideal) is a supplement (in $\CA$) to $\RR=\body(\CA)=\bigwedge^0V$. 
We thus can identify the quotient $\CA/\nilpotent$ (which is a field) with $\body(\CA)\subset \CA$ and the canonical projection $\CA\mapsto \CA/\nilpotent$ with the body map $\body$.
Another feature of (this) $\CA$ (due to the fact that $V$ is infinite dimensional) is that for any $n\in \NN^*$ there exist elements $\xi_1, \dots, \xi_n\in \CA_1=\bigoplus_{k=0}^\infty \bigwedge^{2k+1}V$ such that the product $\xi_1\cdots\xi_n \equiv \xi_1\wedge\cdots\wedge \xi_n$ is non-zero.

\begin{definition}{Remark}
Actually, the particular form of $\CA$ is not that important, as long as it has the above mentioned features: a supplement to the nilpotent elements which is isomorphic to $\RR$ and for any $n\in \NN^*$ elements $\xi_1,\dots, \xi_n\in \CA_1$ whose product is non-zero. And even the last condition is slightly stronger than strictly needed: if an $\CA$-manifold has odd dimension $q$, then there should exist $q+1$ odd elements whose product is non-zero (and thus $V$ should have at least dimension $q+1$). However, as one does not wish to be restricted in the choice of the odd dimension, it is preferable to have this condition for any $n\in \NN^*$, not only for $q+1$. Hence the choice of an infinite dimensional vector space $V$.
\end{definition}

In the sequel we will never have any use for the particular form of $\CA$, so the symbol $V$ will no longer be reserved for the vector space whose exterior algebra is $\CA$.

\mysubsection{appendixonAmanifoldssection}
{\texorpdfstring{$\CA$}{A}-vector spaces}

As $\CA$ is not commutative, there is a difference between left- and right-modules over $\CA$. We will first concentrate on graded bi-modules, meaning an abelian group $M$ which is at the same time a left- and a right-module over $\CA$, and which splits as a direct sum $M=M_0\oplus M_1$ satisfying the conditions
\begin{moneq}
a\in \CA_\alpha \text{ and } m\in M_\beta
\quad\Rightarrow\quad
a\cdot m =(-1)^{\alpha\beta}m\cdot a \in M_{\alpha+\beta}
\mapob,
\end{moneq}
where the gradings $\alpha,\beta$ should be seen as belonging to $\ZZ/2\ZZ$ and thus $1+1=0$. 
(A slight abuse of naming, as by a \myquote{graded bi-module} one usually means a graded left- and right-module over $\CA$, without the condition that the left- and right-actions are related by $a\cdot m =(-1)^{\alpha\beta}m\cdot a$.) 
Associated to the direct sum decomposition $M=M_0\oplus M_1$ we define the \stresd{conjugation operator $\conjugate:M\to M$} by
\begin{moneq}[DefinitionOfConjugationOperator]
m_\alpha\in M_\alpha
\quad\Rightarrow\quad
\conjugate(m_0+m_1) = m_0-m_1
\quad\text{or}\quad
\conjugate (m_\alpha) = (-1)^\alpha m_\alpha
\mapob.
\end{moneq}

Again because $\CA$ is not commutative, there is a difference between left-linear homomorphisms and right-linear ones, even for graded bi-modules. 
A left-linear morphism that is even (meaning that it maps homogeneous elements of the source module to homogeneous elements of the target module of the same parity) is automatically also right-linear (and vice-versa), but a similar statement for non-even morphisms is not true. 
One thus has two morphism modules between two graded bi-modules, the left-linear ones and the right-linear ones. 
These morphism modules are in a natural way graded bi-modules over $\CA$ and they are in a natural way isomorphic.

\begin{definition}{Convention}
As we wish to adhere to the convention that interchanging two homogeneous elements (of whatever nature) induces a minus sign whenever both are odd, we are led to denote evaluation of left-linear morphisms on the left instead of on the right as is usual for maps. More precisely, if $f:M\to N$ is a left-linear morphism between the graded bi-modules $M$ and $N$, we will denote the image of the element $m\in M$ by the map $f$ as $\contrf{m}f$. In that way, left-linearity gets the form
$$
\contrf{a\cdot m}f = a\cdot(\contrf{m}f)
\mapob,
$$
instead of
$$
f(am) = af(m)
\mapob,
$$
which would violate the convention in case $f$ is not even.
This notational convention is already used in ordinary differential geometry when one evaluates\slash contracts a $k$-form with a tangent vector (to yield a $k-1$-form).

\end{definition}

The (for us) important constructions on vector spaces can be carried out also for graded bi-modules over $\CA$, such as direct sums, direct products, morphism modules (left or right), tensor products (over the graded-commutative ring $\CA$), exterior algebras and free bi-modules on a set of homogeneous generators. 

\medskip

We next turn our attention to a special case of graded bi-modules over $\CA$: those that are of the form $V\otimes_\RR \CA$ with $V=V_0\oplus V_1$ a graded vector space over $\RR$, which we call \stresd{$\CA$-vector spaces}. 
Inside such an $\CA$-vector space we have the (real) subspace $V\otimes_\RR \body(\CA)\cong V$ of elements that can be written in the form $v\otimes 1$ for some $v\in V$ and the projection (body map) $\body_V:V\otimes_\RR \CA\to V\otimes_\RR \body(\CA)\cong V$ defined by $\body_V(v\otimes a) = v\otimes \body(a)\cong \body(a)\cdot v$.

If $M$ is any (left) module over $\CA$, we can define the subset $\nilpotent_M\subset M$ of those elements $m$ for which there exists $a\in \CA$, $a\neq0$ such that $a\cdot m=0$.
As the reals are included in $\CA$, any module $M$ over $\CA$ is in particular a real vector space and one can show that $\nilpotent_M$ is a vector subspace of $M$ (as a real vector space).
For an $\CA$-vector space $M=V\otimes_\RR \CA$ we have
$$
\nilpotent_M = \{\,v\otimes n\mid v\in V\,,\, n\in\nilpotent\subset \CA\,\}
\mapob.
$$
Another way to define an $\CA$-vector space thus is as those graded bi-modules $M$ for which there exists a supplement $V\subset M$ for $\nilpotent_M$ in the category of graded vector spaces over $\RR$.
The body map then becomes the projection onto this summand $M=V\oplus \nilpotent_M$. 
In this way the definition of an $\CA$-vector space obtains the same flavor as the conditions imposed on our graded commutative ring $\CA$.

\begin{definition}{Remark}
The map $V\mapsto V\otimes_\RR \CA$ from the category of graded vector spaces over $\RR$ to the category of $\CA$-vector spaces is a functor (nearly an isomorphism of categories). In particular the constructions one can perform on these categories (direct sums, tensor products, exterior algebras) are preserved by this map. However, one has to make a choice what to do with morphisms, as in the category of graded vector spaces over $\RR$ there is no difference between left- and right-linear, whereas in the category of $\CA$-vector spaces there is.

\end{definition}

\mysubsection{appendixonAmanifoldssection}
{Smooth functions}
\label{subsecsmoothfunc}

In order to give a definition of smooth super functions from first principles, we first give an alternate description of smooth functions in the ordinary (non-super) case. We then just copy this alternate description to the super case to obtain our definition of super smooth functions.
The starting point is the formula
\begin{equation}\label{fminfequalsintoverder}
f(x) - f(y) = \sum_{i=1}^m (x_i-y_i)\cdot \left( \int_0^1 (\partial_if)\bigl(sx + (1-s)y\bigr)\,\extder s\right)
\mapob,
\end{equation}
valid for any function $f:O\subset \RR^m\to F$ of class $C^1$ defined on a convex set $O$ with values in a normed vector space $F$. If we define the functions $g_i:O\times O\to F$ by
$$
g_i(x,y) = \int_0^1 (\partial_if)\bigl(sx + (1-s)y\bigr)\,\extder s
\mapob,
$$
then we have:
\begin{equation}\label{smoothnesstrick}
\forall x,y\in O : 
f(x)-f(y) = \sum_{i=1}^m (x_i-y_i)\cdot g_i(x,y)
\mapob.
\end{equation}
It then it is easy to show that, for any $k\in \NN$, $f$ is of class $C^{k+1}$ if and only if there exist functions $g_i$ of class $C^k$ such that \recalf{smoothnesstrick} is valid.
Moreover, the partial derivatives of $f$ are given by $(\partial_if)(x) = g_i(x,x)$. As is well known, these partial derivatives are unique, while the functions $g_i$ themselves are not (apart from the case $n=1$).

\begin{definition}{\myquote{Classical} definitions}
A \stresd{smooth system} is an assignment $\smooths$ that associates, for any dimension $m\in\NN$, to any open set $O\subset \RR^m$ and any target (normed vector) space $F$ a collection of \stress{continuous} functions $\smooths(O,F)\subset C^0(O,F)$ verifying the property
\begin{multline*}
\forall f\in \smooths(O,F)
\quad
\exists g_i\in \smooths(O^2,F)
\quad
\forall x,y\in O\subset \RR^m
\quad:\quad
\\
f(x) - f(y) = \sum_{i=1}^m (x_i-y_i)\cdot g_i(x,y)
\mapob.
\end{multline*}
A smooth system $\smooths_1$ is said to be \stresd{smaller than} a smooth system $\smooths_2$, denoted as $\smooths_1\le \smooths_2$, if we have
$$
\forall O: \smooths_1(O,F)\subset \smooths_2(O,F)
\mapob.
$$

\end{definition}

\begin{proclaim}{Theorem}
$C^\infty$ is the (not a) maximal (with respect to the order $\le$) smooth system.
\end{proclaim}

\begin{definition}{Remark}
The above definitions and statements should not be taken literally, as they are wrong as stated. Formula \recalf{fminfequalsintoverder} is valid for convex sets, but not in general for arbitrary (open) sets. \stress{If} one disposes of smooth partitions of unity, one can re-establish such a result, but that would create, in the case of $\RR$, a circular definition, and in the case of $\CC$ it would be impossible. The solution is to use covers at every stage, but that complicates the notation (not the idea), so we left it out in this summary.

\end{definition}

In order to mimic these definitions in the super case, we have to start with a topology.

\begin{definition}[defoftheDeWittTopologyinAppendix]{Definition}
Let $V$ be any normed vector space over $\RR$ and let $E=V\otimes_\RR \CA$ be the associated $\CA$-vector space, \ie, $V=\body E$. 
The \stresd{DeWitt topology} on $E$ is the coarsest topology on $E$ for which the projection $\body:E\to \body E$ is continuous (when $\body E$ is equipped with the topology induced by the given norm; for finite dimensional $V$ this is the usual euclidean topology).
In particular $\CA$ itself is a (finite dimensional) $\CA$-vector space with $\body \CA=\RR\oplus\{0\}$ (i.e., $(\body\CA)_0=\RR$ and $(\body\CA)_1=\{0\}$), and thus is equipped with the DeWitt topology.
It should be noted that the DeWitt topology is not separated\slash Hausdorff and that for any open set $U\subset E$ we have $\body U\subset U$ and $U=\body\mo(\body U)$.

\end{definition}

Now let $E$ be a finite dimensional $\CA$-vector space and let $e_1, \dots, e_p,f_1, \dots, f_q\in \body E$ be a homogeneous basis, i.e., $e_1, \dots, e_p$ is a basis of the vector space (over $\RR$) $(\body E)_0\subset E_0$ and $f_1, \dots, f_q$ a basis of $(\body E)_1\subset E_1$. It follows that any element $e\in E$ is described uniquely by $p+q$ elements $x_1, \dots, x_p, \xi_1, \dots, \xi_q\in\CA$ according to
$$
e = \sum_{i=1}^p x_ie_i + \sum_{j=1}^q \xi_jf_j
\mapob.
$$
Now if $e$ belongs to the \stress{even} part of $E$, then the $x_i$ belong to $\CA_0$ and the $\xi_j$ to $\CA_1$. In other words, an element of $E_0$ is described by $p$ even \myquote{coordinates} and $q$ odd \myquote{coordinates.}
Moreover, a subset $U\subset E_0$ of the even part of $E$ is open (in the DeWitt topology) if and only if there exists an open set $O\subset \RR^p$ such that
\begin{equation}\label{defofopeninE0}
e = \sum_{i=1}^p x_ie_i + \sum_{j=1}^q \xi_jf_j \in U
\qquad\Longleftrightarrow\qquad
(\body x_1, \dots, \body x_p) \in O
\mapob.
\end{equation}
Another way to state this equivalence is the equality $U = (\body\caprestricted_{E_0})\mo(\body U)$ with $\body U = O\oplus \{0\}\subset (\body E)_0 \oplus (\body E)_1=\RR^p \oplus \RR^q$.

\begin{definition}[superdefinitions]{\myquote{Super} definitions}
$\bullet$
A \stresd{super smooth system} is an assignment $\smooths$ that associates to any open set $U\subset E_0$ in the even part of any finite dimensional $\CA$-vector space $E$ and to any $\CA$-vector space $F$ (the target space, equipped with the DeWitt topology induced by a norm on $\body F$) a collection of \stress{continuous} functions $\smooths(U,F)\subset C^0(U,F)$ verifying the two properties
$$
f(\body U) \subset \body F
$$
and
\begin{multline*}
\forall f\in \smooths(U,F)
\quad
\exists g_i\in \smooths(U^2,F)
\quad
\forall x,y\in U\subset E_0
\quad:\quad
\\
f(x) - f(y) = \sum_{i=1}^{p+q} (x_i-y_i)\cdot g_i(x,y)
\mapob.
\end{multline*}

$\bullet$
A super smooth system $\smooths_1$ is said to be \stresd{smaller than} a super smooth system $\smooths_2$, denoted as $\smooths_1\le \smooths_2$, if we have
$$
\forall U,F: \smooths_1(U,F)\subset \smooths_2(U,F)
\mapob.
$$

$\bullet$
$C^\infty$ is \stresd{the} maximal super smooth system with respect to the order $\le$.

\end{definition}

We now have a definition of super smooth functions that is a look-alike to a possible definition of smooth functions in the non-super case, but two questions remain: (1) why the additional condition $f(\body U)\subset \body F$ and (2) can we say anything interesting about these smooth functions?
One possible answer to (1) is that it allows us to give a positive answer to (2). 
Another possible answer to (1) is that we wish to stay as close as possible to non-super smooth functions and that this conditions is trivially satisfied in the non-super case.

\begin{proclaim}[Gextension]{Lemma}
Let $U\subset E_0$ be an open set in the even part of the $p\vert q$-dimensional $\CA$-vector space $E$, let $O=\body U\subset \RR^p$ be the corresponding open subset in $\RR^p$ (see \recalf{defofopeninE0}) and let $f:O\to \body F$ be an ordinary smooth function with values in the real vector space $\body F$ (where we thus ignore the grading).

With these data, we define the function $\mathbf{G}f:U\to F$ by the formula
$$
(\mathbf{G}f)(x_1, \dots, x_p, \xi_1, \dots, \xi_q)
=
\sum_{k=0}^\infty \frac1{k!} \bigl((D^kf)(r_1, \dots, r_p)\bigr)(\overbrace{n,\vrule width0pt height2.3ex\dots, n}^{k \text{ terms}})
\mapob,
$$
where $r_i=\body x_i\in \RR$, $n_i=x_i-r_i\in \nilpotent_0$, $n=(n_1, \dots, n_p)$, $(D^kf)(r_1, \dots, r_p)$ the $k$th order derivative of $f$ at $(r_1, \dots, r_p)\in O$ as $k$-linear symmetric map with values in $\body F$ and $\bigl((D^kf)(r_1, \dots, r_p)\bigr)({n,\dots, n})$ the formal evaluation of $(D^kf)(r_1, \dots, r_p)$ in the even nilpotent coefficients $n_i$.

This $\mathbf{G}f$ is a super smooth function.

\end{proclaim}

\begin{definition}[DefinitionRepeatedProducts]{Notation}
Let $\mathcal I$ be a finite set of integers, let $A_i$, $i\in \mathcal I$ be a family of objects that can be \myquote{multiplied\slash composed.} 
Then for any subset $I\subset \mathcal I$ we define $\repprod A^I$ as 
\begin{moneq}
\repprod A^I = A_{i_1}\cdot A_{i_{2}} \cdots A_{i_p}
\text{ where }
I = \{i_1, \dots, i_p\}
\text{ with }
i_1<\cdots<i_p
\mapob,
\end{moneq}
with the convention that $\repprod A^\emptyset = 1$ (or the identity map $\id$ in case it concerns maps with composition).  

\end{definition}

\begin{proclaim}[structureofsmoothfunctions]{Theorem}
Let $U\subset E_0$ be an open set in the even part of the $p\vert q$-dimensional $\CA$-vector space $E$, let $O=\body U\subset \RR^p$ be the corresponding open subset in $\RR^p$ and let $F$ be an  
$\CA$-vector space.
Then a function $f:U\to F$ belongs to $C^\infty(U,F)$ (i.e., is super smooth) if and only if there exist (ordinary) smooth functions $f_{I}:O\to \body F$ such that we have
\begin{moneq}
f(x_1, \dots, x_p, \xi_1, \dots, \xi_q)
=
\sum_{I\subset \{1, \dots, q\}} \repprod\xi^I \cdot (\mathbf{G}f_I)(x_1, \dots, x_p)
\mapob.
\end{moneq}

\end{proclaim}

\begin{proclaim}{Corollary}
Let $U\subset E_0$ be an open set in the even part of the $p\vert q$-dimensional $\CA$-vector space $E$, let $O=\body U\subset \RR^p$ be the corresponding open subset in $\RR^p$ and let $F$ be an  
$\CA$-vector space. Then we have the equality
\begin{moneq}
C^\infty(U,F)
=
C^\infty(O,\body F) \otimes \bigwedge \RR^q
\mapob.
\end{moneq}

\end{proclaim}

\begin{definition}{Remark}
Let $E$ and $F$ be two $\CA$-vector spaces (equipped with their DeWitt topologies induced by norms on $\body E$ and $\body F$ respectively). 
For any (left- or right-) linear map $f:E\to F$, we have to ask whether $f$ is continuous, and if so, whether it is smooth. 
Now if $E$ is finite dimensional (and we will always assume that the source space is!), $f$ will automatically be continuous, but it need not be smooth. 
It will be smooth if and only if we have the inclusion $f(\body E)\subset \body F$. 
This implies in particular that the space of smooth (left- or right-) linear maps from $E$ to $F$ is isomorphic to the space of continuous $\RR$-linear maps from $\body E$ to $\body F$. 

\end{definition}

Once we know the structure of super smooth functions, we can address the question of partial derivatives. For that we imitate the non-super case and want to define the partial derivative $\partial_if$ as
$$
(\partial_if)(x) = g_i\bigl(x,x)
\mapob,
$$
where the super smooth functions $g_i$ are given by the definition of super smoothness as
$$
f(x) - f(y) = \sum_{i=1}^{p+q} (x_i-y_i)\cdot g_i(x,y)
\mapob.
$$
But for this to be a coherent definition, one must show that the diagonal $g_i(x,x)$ is uniquely determined by $f$, even when the functions $g_i(x,y)$ are not unique. 
And it is here that one needs the fact that there exist $q+1$ odd elements in $\CA_1$ whose product is non-zero: if that condition is satisfied, the $\partial_if$ are well defined; if not, the diagonal functions $g_i(x,x)$ will not be unique and thus our super smooth functions will not have derivatives. 
(Of course one could define the partial derivatives by hand in terms of the functions $f_I$ of \recalt{structureofsmoothfunctions}, but then the whole purpose of giving an intrinsic definition of super smoothness would be superfluous.)

\medskip

With our choice of $\CA$, the partial derivatives are well defined and behave exactly as expected. If the index $i$ is associated to an odd coordinate, $\partial_i$ is an odd derivation of $C^\infty(U,F)$, whereas it is an even derivation when associated to an even coordinate. 
One should note that the $\partial_i$ are \stress{right}-derivations (as opposed to left-derivations), meaning that they are right-linear, which is a consequence of the choice in \recalt{superdefinitions} to write the \myquote{coefficients} $(x_i-y_i)$ to the left of $g_i(x,y)$ instead of to the right.

Once we know what smooth functions and their derivatives are, nearly all of the usual results of real analysis, and in particular the inverse function theorem, remain true (in most cases even the proofs are close copies).

\mysubsection{appendixonAmanifoldssection}
{\texorpdfstring{$\CA$}{A}-manifolds}
\label{subsecAmfds}

Once we know what open subset are and what (local) diffeomorphisms are, we can copy (nearly) all standard differential geometric constructions, such as manifolds (using an atlas and local diffeomorphisms for chart-changing), fiber bundles, vector bundles, the tangent bundle of a manifold, its cotangent bundle etc{\ae}tera. 
There is however one important difference that concerns vector bundles. 
An $\CA$-manifold is modelled by charts that are open subsets in the \stress{even} part of a finite dimensional $\CA$-vector space $E$; as such it has an even and an odd dimension (the dimensions of $(\body E)_0$ and $(\body E)_1$). 
On the other hand, the typical fiber of a vector bundle is \stress{all} of a finite dimensional $\CA$-vector space $F$. 
There are multiple reasons for doing so, one of them being that the standard constructions such as tensor products and exterior powers can be performed on $\CA$-vector spaces, but not on their even parts (at least not in a way that gives satisfactory results). 
Now a full $\CA$-vector space $F$ can be seen in a natural way as the even part of another $\CA$-vector space: the direct sum of $F$ with its parity dual $\prod F$ (the even part of $\prod F$ is the odd part of $F$ and the odd part of $\prod F$ is the even part of $F$), so vector bundles are still $\CA$-manifolds.

Concerning these constructions, one can show the following results. 
\begin{itemize}
\item
The body map can be extended to $\CA$-manifolds and their smooth maps, the result being an ordinary (non-super) manifold and an ordinary smooth map (taking the body map is essentially mapping all nilpotent elements in $\CA$ to zero).

\item
If $M$ is an $\CA$-manifold of graded dimension $p\vert q$, then $\body M$ is an ordinary manifold of dimension $p$ and for any (sufficiently small) open set $O\subset \body M$ with $U=\body\mo O\subset M$ we have the equality
$$
C^\infty(U, \CA) = C^\infty(O)\otimes \bigwedge \RR^q
\mapob,
$$
providing the link with the sheaf theoretic/ringed spaces approach (in this way we create a ringed space on the ordinary manifold $\body M$, and conversely, every appropriate ringed space appears this way).

\item
The set $C^\infty(M)$ of all smooth maps $f:M\to \CA$ on an $\CA$-manifold $M$ is a graded $\RR$-algebra and the set $\Gamma(B)$ of all smooth sections of a vector bundle $B\to M$ is a graded bi-module over $C^\infty(M)$.

\item
The tangent map $Tf:TM\to TN$ associated to a smooth map $f:M\to N$ between $\CA$-manifolds is naturally \stress{left}-linear and even. In particular we have, for local coordinate systems $x$ in $M$ and $y$ in $N$ with $y=f(x)$, the formula
$$
\iota\Bigl({\sum_{i=1}^{p+q}X_i \cdot \fracp{}{x_i}\bigrestricted_x}\Bigr){Tf} = \sum_{i=1}^{p+q}\sum_{j=1}^{p'+q'} X_i \cdot (\partial_if_j)(x)\cdot \fracp{}{y_j}\bigrestricted_y
\mapob,
$$
where $p\vert q$ is the graded dimension of $M$ and $p'\vert q'$ that of $N$.

\item
Any smooth and \stress{even} vector field $X$ on an $\CA$-manifold $M$ can be integrated to produce a flow, i.e., a $1$-parameter group of local diffeomorphisms $\Phi_t$ with $t\in \CA_0$, which has the usual properties of a flow.
According to \recalt{structureofsmoothfunctions} it has a local expression $\Phi_t(x,\xi)\equiv \Phi(t,x,\xi) = \sum_I (\mathbf{G}\Phi_I)(t,x)\cdot \xi^I$ in terms of ordinary smooth functions of real variables and products of odd coordinates. 
These functions can be computed by induction on the number of elements in $I$ (details can be found in \cite[Ch.V\S4]{Tu04}); the equation for $\Phi_\emptyset$ is a \myquote{standard} first order differential equation for the flow of a non-super vector field (it is the flow equation for the ordinary smooth vector field $\body X$ on the ordinary manifold $\body M$) and the other terms are determined by first order inhomogeneous linear differential equations (and as such they do not restrict the domain of definition of the flow). 

\item
A smooth \stress{odd} vector field $X$ can be integrated to a flow $\Phi_\tau$ with an odd time parameter $\tau\in \CA_1$ (but otherwise with the same properties of a flow) if and only if the auto-commutator $[X,X]$ is zero.
In local coordinates the expression for the flow essentially boils down to the formula
$$
\Phi_\tau(x) = x+\tau\cdot X(x)
\mapob.
$$

\item
Lie's third theorem and its converse are true: to any $\CA$-Lie group $G$ (i.e., an $\CA$-manifold equipped with a group structure such that multiplication and inversion are smooth maps) is associated an $\CA$-Lie algebra $\Liealg g$ (isomorphic to the tangent space at the identity and isomorphic to the space of all left-invariant vector fields). 
Conversely, for any finite dimensional $\CA$-Lie algebra $\Liealg g$ there exists a (unique up to isomorphisms) simply connected $\CA$-Lie group $G$ whose $\CA$-Lie algebra is (isomorphic to) $\Liealg g$.

\item
The notion of \aglgrp/ an $\CA$-Lie group is equivalent to the notion of a \myquote{Lie supergroup pair} or \myquote{super Harish-Chandra pair,} \ie, an \myquote{ordinary} Lie group $H$, an $\CA$-Lie algebra $\Liealg g$ and an action of $H$ on $\Liealg g$ that satisfy the two conditions that the Lie algebra $\Liealg h$ of $H$ is isomorphic to $\body\Liealg g_0$ and such that the restriction to $\Liealg h\cong\body\Liealg g_0$ of the action of $H$ on $\Liealg g$ is the adjoint action.
Starting from an $\CA$-Lie group $G$, one obtains the Lie supergroup pair directly by taking $H=\body G$ and as the action of $H$ on $\Liealg g$ the restriction to $H\subset G$ of the adjoint action of $G$ on $\Liealg g$.
Conversely, if we have a Lie supergroup pair $(H,\Liealg g)$, we first construct the simply connected $\CA$-Lie group $\widetilde G$ whose $\CA$-Lie algebra is $\Liealg g$. It then follows that $\body \widetilde G$ is a covering group of $H$.
If we denote by $Z\subset \body \widetilde G\subset \widetilde G$ the kernel of this covering map, then the condition that the action of $H$ on $\Liealg g$ extends the adjoint action of $H$ on $\Liealg h\cong \body\Liealg g_0$ implies that $Z$ is central (and discrete) in $\widetilde G$. 
One then shows easily that $G=\widetilde G/Z$ is an $\CA$-Lie group whose associated Lie supergroup pair (via the previous construction) is the given Lie supergroup pair.

\end{itemize}

\begin{definition}{A warning}
We have to discuss a negative side effect of the condition imposed on super smooth functions: that it maps the body of the source to the body of the target. 
On the positive side, this condition allows us to prove the structure of super smooth functions \recalt{structureofsmoothfunctions} and the fact that the body of an $\CA$-manifold is a subset of the $\CA$-manifold.
However, on the negative side we have the fact that fixing some coordinates in a super smooth function does not necessarily yield a super smooth function!

To explain the situation, we consider first a function $f:\RR^2\to \RR$, we take a fixed value $x\in \RR$ and we define the function $g:\RR\to \RR$ by
$$
g(y) = f(x,y)
\mapob.
$$
If $f$ is smooth, then $g$ is smooth. 
Let us next consider a function $f:(\CA_0)^2\to \CA$ and a fixed value $x\in \CA_0$, with which we define the function $g=\CA_0\to \CA$ again by $g(y) = f(x,y)$.
But now there is no guarantee that $g$ is smooth if $f$ is smooth!
The problem stems from the condition that for a smooth function we required the condition $g(\body U)\subset \body F$, which gives here the condition
\begin{equation}\label{restrictionbeingreal}
y\in \RR=\body \CA_0
\quad\Rightarrow\quad
g(y)=f(x,y)\in \RR=\body\CA
\mapob.
\end{equation}
But, since $f$ is smooth, we only know the property
$$
x,y\in \RR=\body\CA_0 
\quad\Rightarrow\quad
f(x,y)\in \RR=\body\CA
\mapob.
$$
It follows that if $x$ belongs to $\RR=\body\CA_0$, \stress{then} we can be sure that $g$ is super smooth. 
For any other $x$ it is highly unlikely that $g$ satisfies \recalf{restrictionbeingreal}, in which case it is not super smooth.
A very simple example (which can be interpreted as addition in the $\CA$-Lie group $\CA_1$ of dimension $0\vert1$) is the map $f:(\CA_1)^2\to \CA_1$ given by
$$
f(\xi,\eta) = \xi+\eta
\mapob.
$$
This map is smooth, but for any fixed non-zero element $\xi\in \CA_1$, the resulting map $g:\CA_1\to \CA_1$ given by
$$
g(\eta) = \xi+\eta
$$
is not smooth, as it maps $\{0\}=\body\CA_1$ to $\{\xi\}$, which is not the body of $\CA_1$.

More generally, fixing some coordinates to \stress{real} values in a super smooth function yields a super smooth function (of less variables), but when we restrict to non-real values (arbitrary values in $\CA_0$ or $\CA_1$), it is highly unlikely that the resulting function is super smooth.
For \stress{even} coordinates this is not a real problem, as a super smooth function of even coordinates is an extension of an ordinary smooth function of the same number of real coordinates (see \recalt{Gextension}).
On the other hand, for odd coordinates it is essential: if a super smooth function depends upon an odd coordinate, it is \stress{not} possible to fix this coordinate to an arbitrary value in $\CA_1$ (except $0$) and getting a super smooth function in the remaining coordinates.
This should not come as a surprise for readers familiar with the sheaf-theoretic\slash ringed spaces approach to supermanifolds. 
There a smooth (odd) coordinate function is a single element in the ring of local smooth superfunctions, it is not a function in the ordinary sense. 
Hence evaluating such a ``function'' in a point does not make sense.
On the other hand, we have shown that the local rings of smooth superfunctions are the same in the $\CA$-manifold approach and in the sheaf-theoretic\slash ringed spaces approach, so we should not be able to create more smooth superfunctions by restricting a superfunction to some (arbitrary) value.
\end{definition}

\section*{Acknowledgements}
I thank my colleague V. Thilliez for his patience, support and crucial help during the preparation of this paper. 
I also acknowledge the support of the CDP C2EMPI, as well as the French State under the France-2030 programme, the University of Lille, the Initiative of Excellence of the University of Lille, the European Metropolis of Lille for their funding and support of the R-CDP-24-004-C2EMPI project.

\providecommand{\bysame}{\leavevmode\hbox to3em{\hrulefill}\thinspace}
\providecommand{\MR}{\relax\ifhmode\unskip\space\fi MR }
\providecommand{\MRhref}[2]{%
  \href{http://www.ams.org/mathscinet-getitem?mr=#1}{#2}
}
\providecommand{\href}[2]{#2}

\end{document}